\documentclass[twoside,reqno,12pt, a4paper]{thesis}

\usepackage{amscd}
\usepackage[english,german]{babel}
\usepackage{a4}
\usepackage{amssymb}
\usepackage{amsmath}
\usepackage[a4paper]{geometry}
\usepackage{amsfonts}
\usepackage{amsthm}
\usepackage{graphicx}
\usepackage{epsfig}

\usepackage[dvipdf]{color}
\usepackage{multicol}
\definecolor{bs}{named}{Bittersweet}
\definecolor{og}{named}{OliveGreen}
\definecolor{BrickRed}{named}{BrickRed}
\definecolor{NavyBlue}{named}{NavyBlue}
\definecolor{Gray}{named}{Gray}

\usepackage[color,
 final]
 {showkeys}
\definecolor{refkeybis}{named}{Bittersweet}
\definecolor{labelkeybis}{named}{BrickRed}
{\makeatletter
 \def\SK@refcolor{\color{refkeybis}}
 \def\SK@labelcolor{\color{labelkeybis}}}

\input xy
\xyoption{all}
\def\bcat{\catcode`"=12}
\def\ecat{\catcode`"=13}

\numberwithin{equation}{chapter}

\theoremstyle{plain}   
\newtheorem{theorem}{\bf{Theorem}}[section]

\theoremstyle{definition}
\newtheorem{definition}[theorem]{\bf Definition}

\newcommand{\bea}{\begin{aligned}}
\newcommand{\eea}{\end{aligned}}
\newcommand{\bdm}{\begin{displaymath}}
\newcommand{\edm}{\end{displaymath}}
\newcommand{\barr}{\begin{array}}
\newcommand{\earr}{\end{array}}
\newcommand{\ben}{\begin{enumerate}}
\newcommand{\een}{\end{enumerate}}
\newcommand{\bde}{\begin{description}}
\newcommand{\ede}{\end{description}}
\newcommand{\nn}{\nonumber}

\newcommand{\fn}[1]{\footnote{#1}}

\numberwithin{equation}{chapter}
\def\noi{\noindent}

\newcommand{\hs}[1]{\hspace{#1cm}}
\newcommand{\vs}[1]{\vspace{#1cm}}

\def\Dots{,\,\dots\,,}
\def\eps{\varepsilon}
\let\Hat=\widehat

\def\spa{\hs{0.62}}
\def\spi{\vs{-0.1}\hs{0.62}}
\def\id{{\sf{id}}}
\def\ker{{\sf{ker}\,}}
\def\coker{{\sf{coker}\,}}
\def\Ker{{{\sf{ker}}}}
\def\im{{\sf{Im}}}
\def\pp{{\mathrm{p}}}
\def\pr{{\mathrm{pr}}}
\def\Pr{{\mathrm{Pr}}}
\def\rank{{\sf{rank}}}
\def\dim{{\sf{dim\,}}}
\def\rk{{\sf{rank}}}
\def\mon{{\sf{Mon}}}
\def\comp{
\circ
}
\def\wt{\widetilde}

\def\fii{\varphi}
\newcommand{\twice}[1]{{#1}^{\times 2}}
\newcommand{\threetimes}[1]{{#1}^{\times 3}}
\def\mod{/\hbox{\footnotesize{$\thicksim$}}}
\def\inc{\hookrightarrow}

\let\rinv=\overrightarrow
\let\linv=\overleftarrow

\def\inverse{^{\hbox{-\tiny{$1$}}}}
\def\pa{\partial}
\def\dd{\mathrm{d}}

\def\cif{\calC^\infty}
\def\frax{\mathfrak{X}}
\def\trans{\pitchfork}

\def\pish{\pi^\sharp}

\def\ohm{\omega}
\def\Ohm{\Omega}

\def\frag{\mathfrak{g}}

\def\gpd{\,\lower1.5pt\hbox{$\rightarrow$}\hskip-4.2mm\raise3.3pt
               \hbox{$\rightarrow$}\,}       
\def\dpg{\,\lower1pt\hbox{$\leftarrow$}\hskip-3.75mm\raise3pt
               \hbox{$\leftarrow$}\,}       

\def\acts{\ltimes}  
\newcommand{\ra}{\rangle}
\newcommand{\la}{\langle}

\newcommand{\poidd }[2]{#1\gpd #2}

\newcommand{\gammar}[1]{\rinv{\Gamma}{(#1)}}
\newcommand{\gammal}[1]{\linv{\Gamma}{(#1)}}

\def\rgamma{\rinv{\Gamma}{(\calG,\Ohm)}}
\def\lgamma{\linv{\Gamma}{(\calG,\Ohm)}}
\def\mgamma{\widehat{\Gamma}{(\calG,\Ohm)}}
\def\pgamma{\Gamma^{\downarrow}{(M,A)}}

\newcommand{\brasn}[2]{\pmb{[}{\,#1\,,\,#2\,}\pmb{]}}
\newcommand{\brak}[2]{[{\,#1\,,\,#2\,}]}

\newcommand{\brast}[2]{[{#1
{\overset{*}{,}}\,
#2}]}
\newcommand{\bracts}[2]{[{#1
{\overset{\acts}{,}}\,
#2}]}
\newcommand{\poib}[2]{\{{\,#1\,,\,#2\,}\}}
\newcommand{\pair}[2]{\langle{\,#1\,,\,#2\,}\rangle}
\newcommand{\ppair}[2]{\la\la{\,#1\,,\,#2\,}\ra\ra}

\newcommand{\RR}{\mathbb{R}}

\newcommand{\ZZ}{\mathbb{Z}}
\newcommand{\NN}{\mathbb{N}}
\newcommand{\KK}{\mathbb{K}}

\newcommand{\TT}{\mathbb{T}}
\newcommand{\FF}{\mathbb{F}}

\newcommand{\calA}{\mathcal{A}}
\newcommand{\calB}{\mathcal{B}}
\newcommand{\calC}{\mathcal{C}}
\newcommand{\calD}{\mathcal{D}}
\newcommand{\calF}{\mathcal{F}}
\newcommand{\calH}{\mathcal{H}}
\newcommand{\calV}{\mathcal{V}}
\newcommand{\calS}{\mathcal{S}}

\newcommand{\calN}{\mathcal{N}}
\newcommand{\calK}{\mathcal{K}}
\newcommand{\calG}{\mathcal{G}}
\newcommand{\calI}{\mathcal{I}}
\newcommand{\calJ}{\mathcal{J}}
\newcommand{\calL}{\mathcal{L}}

\newcommand{\calO}{\mathcal{O}}
\newcommand{\calR}{\mathcal{R}}

\newcommand{\calW}{\mathcal{W}}

\let\sf=\mathsf
\let\tsf=\textsf

\newcommand{\sfD}{\mathsf{D}}

\let\rm=\mathrm

\newcommand{\bfC}{\mathbf{C}}
\newcommand{\bfM}{\mathbf{M}}
\newcommand{\bfN}{\mathbf{N}}

\newcommand{\bfj}{\mathbf{j}}
\newcommand{\bfs}{\mathbf{s}}

\newcommand{\pmbgpdm}{$\poidd{\pmb\calG}{\bf M}\:\:$}

\def\pmbs{\pmb\sigma}
\def\pmbd{\pmb\delta}

\def\pmbe{\pmb\eps}
\def\pmbg{\hbox{$\pmb\calG$}}

\let\ol=\overline
\let\ul=\underline

\let\bgn=\begin

\def\it{\item[$i$)]}

\def\tt{\item[3)]}

\newcommand{\be}{\begin{eqnarray*}}
\newcommand{\ee}{\end{eqnarray*}}

\def\Lg{\mathcal{LG}}
\def\La{\mathcal{LA}}
\def\tlg{$\Lg$}
\def\tla{$\La$}
\def\tvb{$\mathcal{VB}$}

\def\pbk{^{\hbox{\tiny{$\pmb{+}$}}}\!}
\def\daga{^{\hbox{\tiny{$\pmb{++}$}}}\!}
\def\ppv{{\iota_{\hbox{\tiny{$\pa_V^+$}}}}^{\hs{-0.315}*}\hs{0.2}}
\def\pmv{{\iota_{\hbox{\tiny{$\pa_V^-$}}}}^{\hs{-0.315}*}\hs{0.2}}
\def\pph{{\iota_{\hbox{\tiny{$\pa_H^+$}}}}^{\hs{-0.315}*}\hs{0.2}}
\def\pmh{{\iota_{\hbox{\tiny{$\pa_H^-$}}}}^{\hs{-0.315}*}\hs{0.2}}
\def\ppmv{{\iota_{\hbox{\tiny{$\pa_V^\pm$}}}}^{\hs{-0.315}*}\hs{0.2}}
\def\ppmh{{\iota_{\hbox{\tiny{$\pa_H^\pm$}}}}^{\hs{-0.315}*}\hs{0.2}}

\def\ttar{\mathrm{\hbox{\small{$\mathrm{t}$}}}}
\def\tar{\mathrm{t}}
\def\ssor{\mathrm{\hbox{\small{$\mathrm{s}$}}}}
\def\sor{\mathrm{s}}
\def\vup{\hbox{\tiny{${}_{V}$}}}
\def\vip{{}^{
\mathrm{v}}}
\def\hup{\hbox{\tiny{${}_{H}$}}}

\def\ttv{\,{\ttar}\vup}
\def\tsv{\,{\ttar}_{v}}
\def\stv{\,{\ssor}\vup}
\def\ssv{\,{\ssor}_{v}}

\def\op{{}^{\rm{op}}}
\def\sc{{\sor}_c}

\def\tth{\,{\ttar}\hup}
\def\tsh{\,{\ttar}_{h}}
\def\sth{\,{\ssor}\hup}
\def\ssh{\,{\ssor}_{h}}

\def\nttv{{\mathrm{t}}\vup}
\def\ntsv{{\mathrm{t}}_{v}}
\def\nstv{{\mathrm{s}}\vup}
\def\nssv{{\mathrm{s}}_{v}}

\def\ntth{{\mathrm{t}}\hup}
\def\ntsh{{\mathrm{t}}_{h}}
\def\nsth{{\mathrm{s}}\hup}
\def\nssh{{\mathrm{s}}_{h}}

\def\etv{{\eps}\vup}
\def\esv{{\eps}_{v}}
\def\etv{{\eps}\vup}
\def\esv{{\eps}_{v}}

\def\eth{{\eps}\hup}
\def\esh{{\eps}_{h}}
\def\eth{{\eps}\hup}
\def\esh{{\eps}_{h}}

\def\itv{{\iota}\vup}
\def\isv{{\iota}_{v}}
\def\itv{{\iota}\vup}
\def\isv{{\iota}_{v}}

\def\ith{{\iota}\hup}
\def\ish{{\iota}_{h}}
\def\ith{{\iota}\hup}
\def\ish{{\iota}_{h}}

\def\mtv{{\mu}\vup}
\def\msv{{\mu}_{v}}
\def\mtv{{\mu}\vup}
\def\msv{{\mu}_{v}}

\def\mth{{\mu}\hup}
\def\msh{{\mu}_{h}}
\def\mth{{\mu}\hup}
\def\msh{{\mu}_{h}}

\def\epsh{\hat{\eps}}
\def\sh{\hat{\sor}}
\def\th{\hat{\tar}}
\def\iotah{\hat{\iota}}
\def\muh{\hat{\mu}}

\def\rhoh{\hat{\rho}}
\def\chih{\hat{\chi}}
\def\jh{\hat{j\,}}
\def\qh{\hat{\mathrm{q}}}
\def\sigmah{\hat{\sigma}}
\def\msigma{\pmb{\sigma}}

\def\hdot{\,\hat{\cdot}\,}

\def\tt{\tilde{t}}

\newcommand{\fib}[2]{\,{}_{#1}\hs{-0.05}\times\hs{-0.05}{}_{#2}\,}

\def\lagpd{\mathcal{LA}\hbox{-groupoid}}
\def\gpdm{\calG\gpd M}

\def\pgpd{(\calG,\Pi)\gpd M}

\def\jmod{\calJ\inverse(e_\star)/G}

\newcommand{\gr}[1]{\sf\Gamma({#1})}

\def\hhh{\,\hat{*}\,}
\def\lrtimes{\ltimes\hs{-0.325}\rtimes}
\begin{document}
\frontmatter
\title{On Morphic Actions and Integrability\\ of \tla-Groupoids}
%
\author{\textbf{ Dissertation}\vspace{+0.2cm}
\\zur
\\Erlangung der naturwissenschaftlichen Doktorw\"urde
\\(Dr. sc. nat.)\vspace{+0.2cm}\\
vorgelegt der
\\Mathematisch-naturwissenschaftlichen Fakult\"at 
\\der
\\Universit\"at Z\"urich\vspace{+0.2cm}
\\von\vspace{+0.3cm}
\\{\bf Luca Stefanini}\vspace{+0.1cm}
\\aus
\\Italien\vspace{1cm}
\\{\bf Promotionskomitee}\vspace{+0.2cm}
\\Prof. Dr. Alberto S. Cattaneo (Vorsitz)
\\Prof. Dr. Giovanni Felder
\\Prof. Dr. Kirill C. H. Mackenzie 
\vspace{+0.4cm}
\vfill Z\"urich 2008}
\maketitle
\newpage\renewcommand\thepage{}
\quad
%
\emph{Dedicated to Motata, Rita and Toto}
\newpage\quad\newpage
\renewcommand\thepage{}
\begin{center} \Large \bf Synopsis \end{center}
%
Lie theory for the integration of Lie algebroids to Lie groupoids, on the one hand, and of Poisson manifolds to symplectic
groupoids, on the other, has undergone tremendous developements in the last decade, thanks to the work of Mackenzie-Xu,
Moerdijk-Mr\v cun, Cattaneo-Felder and Crainic-Fernandes, among others.\\ 
In this thesis we study - part of - the categorified version of this story, namely the integrability of \tla-groupoids (groupoid
objects in the category of Lie algebroids), to double Lie groupoids (groupoid objects in the category of Lie groupoids)
providing a first set of sufficient conditions for the integration to be possible.  

Mackenzie's double Lie structures arise naturally from lifting processes,  such as the cotangent lift or the path
prolongation,  on ordinary Lie theoretic and Poisson geometric objects and we use them to study the integrability of
quotient Poisson bivector fields, the relation between ``local'' and ``global'' duality of Poisson groupoids and Lie
theory for Lie bialgebroids and Poisson groupoids.

In the first Chapter we prove suitable versions of Lie's 1-st and 2-nd theorem for Lie bialgebroids, that is, the
integrability of subobjects (coisotropic subalgebroids) and morphisms, extending earlier  results by Cattaneo and Xu,
obtained using different techniques.

We develop our functorial approach to the integration of \tla-groupoids \cite{07a} in the second Chapter, where we also
obtain partial results, within the program, proposed by Weinstein, for the integration of Poisson groupoids to symplectic
double groupoids. 

The task of integrating quotients of Poisson manifolds with respect to Poisson groupoid actions motivates the study we
undertake in third Chapter	 of what we refer to as morphic actions, i.e. groupoid actions in the categories of Lie
algebroids and Lie  groupoids, where we obtain general reduction and integrability results. 

In fact, applying suitable procedures \`a la Marsden-Weinstein zero level reduction to ``moment morphisms'', respectively
of Lie bialgebroids or Poisson groupoids, canonically associated to a Poisson $\calG$-space, we derive two approches to
the integration of the quotient Poisson bivector fields.\\
 The first, a kind of integration via symplectic double groupoids, is not always effective but reproduces the
``symplectization functor'' approch to Poisson actions of Lie groups, very recently developed by Fernandes-Ortega-Ratiu,
from quite a different perspective. We earlier implemented this approach successfully in the special case of complete
Poisson groups \cite{07b}.\\ 
The second approach, relying both on a cotangent lift of the Poisson $\calG$-space and on a prolongation of the original
action to an action on suitable spaces of Lie algebroid homotopies, produces necessary and sufficient integrability
conditions for the integration and gives a  positive answer to the integrability problem under the most natural assumptions.
\newpage
\quad
\newpage
\selectlanguage{german}
\begin{center}  {\Large \bf Zusammenfassung} \end{center}
%
%
%
Im letzten Jahrzehnt hat sich die Lie-Theorie insbesondere aufgrund ihrer Wichtigkeit
f\"ur die Integration von Lie Algebroiden nach Lie Gruppoiden, einerseits, und von Poisson Mannigfaltigkeiten
nach symplektischen Gruppoiden, andererseits, enorm weiterentwickelt. Arbeiten von Mackenzie-Xu, Moerdijk-Mr\v cun, 
Cattaneo-Felder und Crainic-Fernandes haben diese Entwicklung unter anderem entscheidend beeinflusst.

In der vorliegenden Dissertation interessieren wir uns diesbez\"uglich f\"ur den kategorientheoretischen Aspekt
der Integration von $\mathcal{L}\mathcal{A}$-Gruppoiden (n\"amlich von Gruppoid-Objekten in der Kategorie der Lie
Algebroide) nach doppelten Lie Gruppoiden (n\"amlich von Gruppoid-Objekten in der Kategorie der Lie Gruppoide). 
Hier erhalten wir hinreichende Bedingungen f\"ur die Integration.

Mackenzies doppelte Lie-Strukturen entstehen dabei ganz nat\"urlich aus Hochhebungsprozessen, wie der 
kotangentialen Hochebung oder der Pfad-Prolongation auf gew\"ohnlichen Lie-theoretischen oder Poisson-geometrischen
Strukturen. Mit ihrer Hilfe werden wir die Integrabilit\"at der Quotienten-Poisson-Bivektorfelder, 
die Beziehung zwischen lokaler und globaler Dualit\"at der Poisson-Gruppoide sowie die Lie-Theorie der
Lie-Bialgebroide und der Poisson-Gruppoide untersuchen.
Im ersten Kapitel beweisen wir gewisse Varianten des ersten und zweiten Satzes von Lie \"uber Lie-Bialgebroide, d.h. , 
\"uber die Integrabilit\"at von Unterobjekten (den ko-isotropischen Unteralgebroiden) und ihren Morphismen.
Hier verallgemeinern wir fr\"uhere Resultate von Cattaneo und Xu.

Im zweiten Kapitel entwickeln wir unseren funktorialen Ansatz zur Integration der $\mathcal{L}\mathcal{A}$-Gruppoide\cite{07a}.
Hier  erhalten wir  positive Teilergebnisse zur von Weinstein vorgeschlagenen Integration der
Poisson-Gruppoide nach doppelten symplektischen Gruppoiden.

Der Untersuchung der sogenannten morphischen Wirkungen, also Gruppoid-Wir-kungen in der Kategorie
der Lie Algebroide und der Lie Grupppoide, widmen wir uns im dritten Kapitel. Hier erhalten wir Reduktions- 
und Integrabilit\"atsresultate und wir gehen die  Aufgabe an, Quotienten von Poisson-Mannigfaltigkeiten in Bezug auf 
Wirkungen von Poisson-Gruppoiden zu integrieren.

Tats\"achlich erhalten wir durch die Anwendung geeigneter Varianten des Marsden-Weinstein-Reduktions-Verfahrens
zwei Ans\"atze zur Integration der Quotienten von Poisson-Bivektorfeldern.

Der erste Ansatz, eine Art von Integration durch doppelte symplektische Gruppoide, ist nicht immer erfolgreich. 
Dennoch gibt er im speziellen Fall der Wirkungen von Lie-Gruppen 
Fernendes-Ortega-Ratius Symplektisierungsfunktor-Ansatz wieder. Dieser Ansatz wurde schon zuvor von uns erfolgreich 
im speziellen Fall der Wirkungen von Poisson-Gruppen angewendet (siehe \cite{07b}).
Der zweite Ansatz basiert  schliesslich auf einer kotangentialen Hochhebung der Wirkung und ihrer Pfad-Prolongation nach einer
Wirkung auf geeigneten R\"aumen von Lie-Algebroid-Homotopien. \"Uber ihn erhalten wir notwendige sowie hinreichende
Integrabilit\"ats-bedingungen, womit wir unter kanonischen Voraussetzungen das Integrationsproblem vollst\"andig l\"osen.
\newpage
\quad
%
\selectlanguage{english}
\begin{center}  {\Large \bf Acknowledgements} \end{center}
For his guidance and constant support I am deeply indebted to Alberto Cattaneo; most
of all I wish to thank him for his persistently positive attitude, a decisive source
of encouragement throughout my Ph.D. I am very grateful to Kirill Mackenzie for
sharing with me many of his ideas through the years. I also wish to thank Ping Xu,
Marco Zambon, Rui Loja Fernandes, Serge Parmentier and Iakovos Androulidakis for
stimulating conversations, and Ping Xu, Kirill Mackenzie and Giovanni Felder for
interesting comments on final drafts of this dissertation.\\ 
Special thanks go to Mauro Carfora, Claudio Dappiaggi and Sergio Pini for helping me
find my way to Z\"urich. Nicola Kistler and his successor Felice Manganiello
deserve a big thank you for providing me in the office hours with good coffee, an
essential ingredient of my diet.
\newpage
\quad
%


%
%
%
%
\newpage\selectlanguage{german}
\begin{quotation}{\em\large K"onnte jeder brave Mann\\
solche Gl"ocken finden,\\
seine Feinde w"urden dann\\
ohne M"uhe schwinden}\footnote{Emanuel Schikaneder, Die Zauberfl"ote,  1791.} [\:?\,]\footnote{Pers\"onlische
Beischrift.} 
\end{quotation}
\newpage
\thispagestyle{empty}
\selectlanguage{english}
\tableofcontents
\thispagestyle{empty}
\newpage
\quad
\newpage
\renewcommand\thepage{\roman{page}}
\setcounter{page}{1}
\linespread{1.05}
\chapter*{Preliminary remark}
\chaptermark{ }
The first few pages of each chapter, ``section zero''s, can be read separately, as an introduction to each chapter,
or together as an introduction to the whole dissertation.

\chapter*{Notations and conventions}
\chaptermark{ }
Most of the notations and conventions are introduced in the main text, or should be clear from the context;
we list below a few remarks:
\bgn{itemize}
\item[$\cdot$] The symbol ``$\bullet$'' is used with various meanings throughout the text. It
can stand for the one point manifold, an unassigned variable, the
multiplication of the tangent groupoid,  etc \dots ; it shall be clear from the context which
meaning is intended.
\item[$\cdot$] For any function $f\in\cif(M)$ we denote with $X^f=\poib{f}{\cdot}$  the associated Hamiltonian
vector field,  when $M$ is symplectic or Poisson. The Poisson brackets are defined as 
$\{f,g\}:=\pi(\dd f,\dd g)$, respectively $\{f,g\}:=\ohm(X^g,X^f)$.
\item[$\cdot$] The symbol $\pair{\cdot}{\cdot}$ is used for the natural pairing of a vector bundle with its
dual.
\item[$\cdot$] $\pounds$ denotes the Lie derivative.
\item[$\cdot$] All manifolds are real smooth $\calC^\infty$ manifolds, unless otherwise stated (but not necessarily Hausdorff);
$\Gamma(M,E)$ or simply $\Gamma(E)$ is the space of $\calC^\infty$ sections of a vector bundle. We use the
 symbols $\gr{f}$, respectively $\gr{\thicksim}$, to denote the graph of a map $f$ (with respect to 
the standard ordering) or of an equivalence relation $\thicksim$. 
\item[$\cdot$] If $X\in\Gamma(E_-\otimes E_+)$, $X^\sharp:E_-^*\to E_+$ is the associated vector bundle map, for
all vector bundles $E_\pm$.
%
%
%
%
\item[$\cdot$] $\underset{i,j,k}\oint X_{ijk}$ is the sum over cyclic permutations. 
\item[$\cdot$] The symbol $\Delta$ is used with different meanings: 
to denote both singular and regular distributions, diagonal subsets and also the inclusion map of a diagonal
subset in a direct product.
\item[$\cdot$] The symbol $\fii$ is typically used to denote a morphism of Lie groupoids, $\phi$ instead to 
denote a morphism of Lie algebroids.
\end{itemize}


 
%
\mainmatter
\pagenumbering{arabic} 
\newpage
\thispagestyle{empty}
\linespread{1.05}
\chapter{Lie-Poisson duality}\label{chapi}
\spa Poisson algebras first arose in the mid 19-th century from Jacobi's algebraic study of mechanical  systems, embarked to
understand the  relation  of the brackets earlier discovered by Poisson with constrained dynamical systems and conservation 
laws (see e.g. \cite{ws98,xuonvaisman} and references therein for an historical review). The study of the geometry of these
brackets later led Lie to discover what nowadays we call Lie algebras; in fact, the dual of a Lie algebra is one of the
fundamental examples of a Poisson manifold. After a very long dormancy, Poisson manifolds were rediscovered across the 1960s
and 1970s in the work of Berezin, Kirillov, Kostant and Lichnerowicz, among others. After the introduction of Poisson
cohomology by Lichnerowicz in 1977 \cite{lnw77}, the unraveling of the local structure of Poisson manifolds by Weinstein in 1983, and the
discovery of symplectic groupoids towards the end of 1980s independently by Karas\" ev \cite{ks86}, 
Weinstein \cite{ws87} and Zakrzewski \cite{zr90},
a form of ``duality'' between Poisson brackets and Lie brackets begun to emerge in a more general setting, linking Poisson
geometry to the theory of Lie groupoids-Lie algebroids started earlier by Pradines in a series of papers
\cite{prd67,prd67b,prd68a,prd68b}.\\
Instances of this duality were further investigated in the late 1980s and through the
1990s to the present day, especially via Weinstein's coisotropic calculus, in the study of Poisson group-\emph{oids}, Poisson homogeneous
spaces, symplectic realizations, generalized moment(um) maps and the integrability of Lie algebroids.\\
\spi In the first two Sections of this mostly introductory Chapter we shall introduce, the preliminary material we will
need in the rest of this dissertation, by
reviewing well known facts about Poisson manifolds, Lie algebroids and Lie groupoids, respectively their morphisms, with
emphasis on the duality issues.\\
In Section \ref{algalg} we define actions of Lie groupoids and their infinitesimal counterpart.\\
In Section \ref{ilapm} we discuss the integrability of Lie algebroids to Lie groupoids and its dual counterpart, namely the
realization of a Poisson manifold via a symplectic groupoid. 
We choose to follow there a logical, rather then historical, order. A brief historical account of the tremendous
developments of Lie theory that took place in the last decade is nevertheless in order. The integrability of morphisms of 
Lie algebroids is due to Mackenzie and Xu \cite{mx00}. After that Moerdijk and Mr\v cun  provided in  \cite{mm02} every Lie
groupoid with a source 1-connected ``covering'' groupoid, given by the quotient of the monodromy groupoid of the associated
source foliation, and proved the integrability of Lie subgroupoids, also by foliation theoretic techniques. \eject 
About at the same
time, a reduction \`a la Marsden-Weinstein of the phase space of a topological field theory known as the Poisson sigma model,
performed by Cattaneo and Felder in \cite{cafe01}, produced a topological model for the symplectic groupoid of a Poisson
manifold; in the integrable case the construction yields the desired symplectic groupoid. Soon after Crainic and Fernandes
adapted Cattaneo and Felder's model, following an approach foreseen by \v Severa in \cite{sev05}, in order to characterize Moerdijk and
Mr\v cun's ``covering groupoid'' in terms of the Lie algebroid data only and obtained, by connection theoretic means, necessary
and sufficient conditions for the integrability of Lie algebroids \cite{crfs03}  and Poisson manifolds \cite{crfs04}. Crainic
and Fernandes' topological model, the so called \emph{Weinstein groupoid},  was in turn recovered via the reduction of the  Poisson
sigma model for the dual Poisson structure in \cite{cat04} by Cattaneo, where an integration of coisotropic submanifolds to
Lagrangian subgroupoids, namely the dual phenomenon to the  integrability of Lie subalgebroids and Lie algebroid morphisms,
was also given.\\
The original contributions presented in this Chapter can be found in the last Section. After introducing Lie
bialgebroids and Poisson groupoids, in a sense self-dual objects in the intersection of the Lie and Poisson worlds, we adapt
Mackenzie and Xu's approach for the integration of Lie bialgebroids to Poisson groupoids \cite{mx00} in order to prove the
integrability of coisotropic subalgebroids to coisotropic subgroupoids:
\bgn{thm}[\bf\ref{coisosub}] Let $\poidd{(\calG,\Pi)}{M}$ be a source 1-connected Poisson  groupoid
with Lie bialgebroid $(A,A^*)$ and $\poidd{\calC}{N}$ a source 1-connected Lie subgroupoid with Lie
algebroid $C$. Then $\calC\subset\calG$ is coisotropic if{f} so is $C\subset A$ for the dual Poisson
structure induced by $A^*$.  
\end{thm}
From this result it is then easy to derive the integrability of morphisms of Lie bialgebroids to morphisms of Poisson
groupoids and the equivalence of the category of Lie bialgebroids to that of source 1-connected Poisson groupoids.
\vs{0.5}
\section{Poisson manifolds, Lie algebroids and Lie groupoids}
\begin{quotation} 
In this Section we introduce the main definitions and basic facts about Poisson manifolds, Lie algebroids and Lie groupoids.
The geometric study of the Poisson brackets of classical mechanics leads naturally (and led historically) to Lie algebras
and Poisson manifolds. In fact, Lie algebras are dual to linear Poisson bivector fields and in turn  the dual object to a general
Poisson manifold can be understood as a Lie algebroid, the infinitesimal invariant of a Lie groupoid. We also discuss standard 
examples in some detail; a number of those  actually plays a role in the general theory of Poisson geometry and Lie theory
(see also, for instance, 
\cite{mackbook,vm94,cw99}).
\end{quotation}

\vs{0.1}
\spa A \tsf{Poisson bracket} on an associative algebra is given
by a Lie bracket for which the adjoint representation takes values in the derivations of the associative
product. The definition applies for any ground field, the algebra can be graded and the associative product need
not be commutative; still graded commutativity is a desirable property.
\bgn{definition} A \tsf{Poisson structure} on a manifold $M$ is a Poisson  algebra on $\cif(M)$, i.e. a Lie
bracket $\poib{}{}$ such that the Leibniz rule
$$
\poib{f}{g\cdot h}=g\cdot\poib{f}{h}+\poib{f}{g}\cdot h
$$
holds for all $f,g,h\in\cif(M)$.  
\end{definition}
The above Leibniz rule is equivalent to asking $\poib{}{}$ to be a skewsymmetric biderivation of the pointwise
product of $\cif(M)$, therefore setting
$$
\qquad
\pi(\dd f,\dd g):=\poib{f}{g}
\qquad,\qquad
f,g\in\cif(M)
\qquad,
$$
defines a bivector field $\pi\in\Gamma(\wedge^2TM)$, which we shall call a   \tsf{Poisson bivector}. A bivector field 
$\pi\in\frax^2(M)$ is the Poisson bivector of a Poisson structure if{f} the trivector field
$\brasn{\pi}{\pi}\in\Gamma(\wedge^3TM)$,
locally defined by the components
$$
\qquad
\brasn{\pi}{\pi}^{\alpha\beta\gamma}
=
\underset{\mu}\sum\oint_{\alpha,\beta,\gamma}\pi^{\alpha\mu}\pa_\mu\pi^{\beta\gamma}
\qquad
$$
in a coordinate patch, vanishes identically on $M$\fn{The notation is not incidental. The graded
commutative algebra $\frax^\bullet(M)=(\Gamma(\wedge^\bullet TM),\wedge)$ of multivector fields, carries a unique graded Lie bracket
of degree -1 extending the bracket of vector fields and making it a graded Poisson algebra
\cite{sch40,sch53,nij55}}
refer to the vector bundle map $\pish:T^*M\to TM$ induced by a Poisson structure as the  \tsf{Poisson anchor}.
The image $\im\,\pish$ of the Poisson anchor of a Poisson manifold M, regarded as a submodule of the local vector
fields over $M$, defines an integrable singular distribution in the sense of  Sussmann \cite{su73} and Stefan
\cite{st74}. That is, the points reached by the sequences of local flows of vector fields taking values in
$\im\,\pish$,  can be assembled  in a partition of $M$ in connected immersed submanifolds (not necessarily of the
same dimension), the \tsf{leaves} of $\pish$,  which can be pasted nicely (see \cite{vm94,balan} for details). The
vector fields  $X^f:=\pish\comp\dd f$, $f\in\cif(M)$, spanning the leaves of the distribution are called 
\tsf{Hamiltonian vector fields} in analogy with the classical Hamiltonian dynamics, since also in the Poisson case
$X^f=\poib{f}{\cdot}$, as a derivation of $\cif(M)$, and the subspace of  Hamiltonian vector fields is easily
seen to be a Lie subalgebra of $\frax(M)$. Moreover, the restriction of $\pish_\calL:T^*_\calL M\to T\calL$ to a
leaf $\calL$ has maximal rank equal to $\dim\calL$, by definition and one can check that setting
$$
\qquad
\ohm^\calL(\pish\alpha_+,\pish\alpha_-):=\pi(\alpha_-,\alpha_+)\qquad,\qquad \alpha_\pm\in\Ohm^1(M)\qquad,
$$
yields a well defined nondegenerate closed 2-form on $\calL$.
Probably, this is the most important property of a Poisson manifold, namely, being foliated in
\tsf{symplectic leaves} \cite{kv75} and might be regarded as
global a manifestation of the Weinstein-Darboux local structure theorem ~\cite{ws83}. Clearly, a symplectic manifold is a
Poisson manifold with only one leaf or equivalently a Poisson manifold for which the sharp map has maximal rank. %
%
Poisson manifolds are quite general objects interpolating between arbitrary manifolds and symplectic manifolds.
\bgn{example}\tsf{Poisson manifolds}  
\medskip\\
$i)$ Every manifold is a Poisson manifold for the zero bivector field, each point is a symplectic leaf.
\medskip\\
$ii)$ Every constant bivector field $B$ on $\RR^n$ is a Poisson bivector, the dimension of the symplectic leaves
is the rank of the matrix representing $B$.
\medskip\\
$iii)$ The sharp map of a Poisson bivector of maximal rank is the inverse of the sharp map of a symplectic
form, the whole manifold is the only symplectic leaf.
\medskip\\
$iv)$ Every smooth function $f$ on $\RR^2$, yields a Poisson structure by setting $\poib{x}{y}:=f(x,y)$ for the
coordinate functions $x$ and $y$ on $\RR^2$; the symplectic leaves are each point in the zero set $\sf Z(f)$ of
$f$ and the connected components of $\RR^2\backslash \sf Z (f)$.  For instance, 
the symplectic leaves of 
$$
f(x,y):=\exp\left(-\frac{1}{(1-(x^2 +y^2))^2}\right)
$$ 
are the interior of the unit disk, the complement of its closure and each point on the unit circle. 
\end{example}
Morally speaking, Poisson manifolds are nonlinear Lie algebras. In fact, the very idea of a Lie algebra emerged
from Sophus Lie's \cite{Lie} attempt to understand Poisson brackets in a geometric fashion,
by examining the simplest nontrivial examples. Consider a linear Poisson structure on $\RR^n$:
$$
\qquad
\poib{x^i}{x^j}=\sum_kf^{ij}{}_kx^k\qquad,\qquad i,j,k=1\Dots n
\qquad,
$$
for some constants $\{f^{ij}{}_k\}$, with $f^{ij}{}_k+f^{ji}{}_k=0$. The condition for
$\poib{}{}$ to be a Poisson bracket bracket writes
$$
\sum_{l=1}^{n}\underset{i,\,j,\,k}\oint f^{ij}{}_lf^{lk}{}_m=0\qquad,\qquad i,j,k,m=1\Dots n
\qquad,
$$
that is, of course, the requirement for the $f$'s to be the structure constants of a Lie algebra. Replacing
$\RR^n$ with a finite dimensional vector space $V$,  one can see that Lie algebras on $V$ are in bijective
correspondence with linear Poisson structures on $V^*$ by the formula
$$
\qquad
\pair{\xi}{\brak{v}{w}}=\poib{F_v}{F_w}(\xi)\qquad,\qquad
v,w\in V\simeq\cif_{\sf{lin}}(V^*)\quad,\quad \xi\in V^*
\qquad,
$$
where we identify a vector $v\in V$ with the associated linear functional $F_v$ on $V^*$.
\bgn{example} For any Lie group $G$ with Lie algebra $\frag$, the symplectic  leaves of the dual Poisson structure
on $\frag^*$ are the connected  components of the coadjoint orbits of $G$. Recall that the coadjoint action
$\sf{Ad}^*: G\to\sf{Diff}(\frag^*)$ is defined by
$$
\quad
\pair{\sf{Ad}^*_g\xi}{x}:=\pair{\xi}{\sf{Ad}_{g\inverse}x}
=
\pair{\xi}{\dd\calC_{g\inverse}x}\qquad,\qquad
g\in G \hbox{ , }\xi\in\frag^*\hbox{ and } x\in\frag
\quad,
$$
for the conjugation map $\calC_g:G\to G$. It is a basic exercise in classical Lie theory to see that the induced
infinitesimal action $\sf{ad}^*:\frag\to \frax(\frag^*)$ is
$$
\quad
\pair{\sf{ad}^*_x\xi}{y}=\pair{\xi}{\brak{y}{x}}=\poib{F_y}{F_x}(\xi)
\qquad,\qquad x,y\in\frag\hbox{ and } \xi\in\frag^*\quad,
$$
therefore the fundamental vector fields of the action are Hamiltonian, i.e. $\sf{ad}^*(x)_\xi=X^{F_x}_\xi$, and span the
symplectic foliation of $\frag^*$.
\end{example}
The Poisson bivector of a Poisson manifold $(M,\pi)$ may be regarded as a $\cif(M)$-bilinear operation
$\wedge^2\Ohm^1(M)\to\cif(M)$; according to the above discussion, it appears natural to ask, whether it is
possible to linearize $\pi$ in such a way to obtain a Lie bracket on $\Ohm^1(M)$. The answer is positive: set
$$
\qquad
\pair{\brak{\theta_+}{\theta_-}}{X}
=
(\pounds_X\pi) (\theta_+,\theta_-)
+
\pish\theta_+(\pair{\theta_-}{X})
-
\pish\theta_-(\pair{\theta_+}{X})
\qquad,
$$
for all $\theta_\pm\in\Ohm^1(M)$, $X\in\frax(M)$ or equivalently,
\bgn{eqnarray}\label{Koszulint}
\qquad
\brak{\theta_+}{\theta_-}
&=&
\dd\pi(\theta_+,\theta_-) +\iota_{\pish\theta_+}\dd\theta_- -
\iota_{\pish\theta_-}\dd\theta_+\\
&=&
\nn\pounds_{\pish\theta_+}\theta_- - \pounds_{\pish\theta_-}\theta_+
-
\dd\pi(\theta_+,\theta_-)
\qquad;
\end{eqnarray}
note that the $\cif(M)$-bilinearity is replaced by the Leibniz rule
$$
\qquad
\brak{\theta_+}{f\cdot\theta_-}=f\cdot\pi
\brak{\theta_+}{\theta_-}+\pish{\theta_+}(f)\cdot\theta_-
\qquad,\qquad f\in\cif(M)
\qquad,
$$
and $\brak{ }{}$ makes $\Ohm^1(M)$ a Lie algebra over $\RR$. This bracket we shall call the \tsf{Koszul bracket},
is very natural, for it is the unique extending 
$$
\qquad
\brak{\dd f_+}{\dd f_-}=\dd\poib{f_+}{f_-}
\qquad,
$$
$f_\pm\in\cif(M)$, according to the Leibniz rule above and therefore generalizing the Lie bracket dual to a linear
Poisson bracket on a vector space to the nonlinear case.%
\\
The Koszul bracket, together with the Poisson anchor endows $T^*M$, with a Lie algebroid structure.
\bgn{definition} A \tsf{Lie algebroid} structure on a vector bundle $A\to M$ is given by a $\RR$-bilinear Lie 
bracket $\brak{}{}$  on the space of sections $\Gamma(A)$ and a vector bundle map $\rho:A\to TM$ such that
the Leibniz rule
$$
\brak{a_+}{f\cdot a_-}=f\cdot\brak{a_+}{a_-}+\rho(a_+)(f)\cdot a_-
$$
holds for all $a_\pm\in\Gamma(A)$ and $f\in\cif(M)$. 
\end{definition} 
Before turning to examples and introducing the main features of Lie algebroids, we will show how  the duality
between Lie algebras and linear Poisson structures generalizes to a duality between Lie algebroids and fibrewise
linear Poisson structures. The space of  functions on $A^*$ is the completion of the subspace of  fiberwise
polynomial functions, which in turn is the associative algebra generated by the space 
$\cif_{\sf{lin}}(A^*)\oplus\cif_{\sf{cst}}(A^*)$ of fiberwise linear functions
$\cif_{\sf{lin}}(A^*)\simeq\Gamma(A)$ and fibrewise constant functions $\cif_{\sf{cst}}(A^*)\simeq\cif(M)$. In order to 
endow $A^*$ with a Poisson bracket $\poib{}{}_A$ it is then sufficient to define it on 
$\cif_{\sf{lin}}(A^*)\oplus\cif_{\sf{cst}}(A^*)$ and extend it according to  the Leibniz rule: setting, for all
$\xi\in\Gamma(A^*)$,
$$
\qquad
\poib{F_{a_+}}{F_{a_-}}_A(\xi):=\pair{\xi}{\brak{a_+}{a_-}}\qquad,\qquad
a_\pm\in\Gamma(A)
\qquad,
$$
in analogy with the case $M=\bullet$, where $F_{a_\pm}(\xi)=\pair{a_\pm}{\xi}$ are the associated fibrewise 
linear function on $A^*$, $\poib{}{}_A$ extends to  $\cif_{\sf{lin}}(A^*)\oplus\cif_{\sf{cst}}(A^*)$
compatibly with  the Leibniz rule is given by
\bgn{equation}\label{dfgh}
\poib{F_a}{\pr_{A^*}^*f}_A=\rho(a)(f)\comp \pr_{A^*}
\qquad\hbox{ and }\qquad
\poib{\pr_{A^*}^*f}{\pr_{A^*}^*g}=0\qquad,
\end{equation}
for all $a\in\Gamma(A)$, $f,g\in\cif(M)$.  
To check that the bivector $\pi_A$ on $A^*$ defined by equations (\ref{dfgh}) is Poisson, for example in the local
coordinates induced a choice of dual frames for $A$ and $A^*$, is then straightforward. Poisson structures of this
kind are called \emph{fibrewise linear}, in the sense that so is the restriction of the corresponding Poisson bivector to
the vertical subbundle of $\wedge^2A^*$.  Clearly, inverting the construction, a fibrewise linear Poisson structure
on a vector bundle  induces a Lie algebroid on the dual bundle.%
\\
Lie algebroids generalize various classes of  differential geometric  structures.
\bgn{example}\label{exala} \tsf{Lie algebroids}
\medskip\\
$i$) Every vector bundle is a Lie algebroid for the zero bracket and anchor. A \tsf{regular distribution} on a
manifold $M$,  i.e. a smooth subbundle $\Delta\to M$ of $TM$ or, equivalently, a singular distribution of constant
rank,  is a Lie algebroid if{f} it is integrable, i.e. if{f} the space of sections of $\Delta$ is closed under the Lie
bracket of vector fields. Notably, the tangent bundle itself is a Lie algebroid.
\medskip\\
$ii$) Lie algebras are Lie algebroids over the one point manifold: the anchor cannot be anything but zero. If a
Lie algebroid has zero anchor, setting
$$
\qquad
[x_+\,\overset{q},\:x_-]:=\brak{\wt{x}_+}{\wt{x}_-}_q \qquad,\qquad x_\pm\in A_q
\qquad,
$$
for any sections $\wt{x}_\pm\in\Gamma(A)$ with  $\wt{x}_\pm(q)=x_\pm$, yields a \tsf{bundle of Lie algebras}, i.e.
a smooth family of Lie algebras  $\{A_q,[\:\overset{q},\,\:]\}_{q\in M}$ on the fibres of $A$, parameterized by
$M$. 
\medskip\\
$iii$) Let $\sigma:\frag\to\frax(M)$ be an infinitesimal action on a manifold $M$ (i.e. a morphism of Lie
algebras). Then the flat bundle $\frag\times M$ carries a canonical Lie algebroid structure, the \tsf{action Lie
algebroid}  $\frag\ltimes M$. Note that  $\Gamma(\frag\times M)\simeq\cif(M,\frag)$ and setting
$$
\pair{\xi}{\wt{X}(F)}:=X(\pair{\xi}{F})\quad\qquad,\quad\qquad \xi\in\frag^*\quad,\quad
F\in\cif(M,\frag)\qquad,
$$
allows extending a vector field $X\in\frax(M)$ to a linear endomorphism of  $\cif(M,\frag)$. The bracket of
$\frag\ltimes M$ is defined by
$$
\bracts{F}{G}(q):=\brak{F(q)}{G(q)} 
+
(\wt{\sigma\comp F})_q(G)
-
(\wt{\sigma\comp G})_q(F)
\quad,\quad
F,G\in\cif(M,\frag)\quad,
$$
$q\in M$, inducing the evaluation $\rho_\ltimes$ of $\sigma$ as an anchor map: $\rho_\ltimes(x,q)=\sigma(x)_q$,
$(x,q)\in\frag\times M$.
\medskip\\
$iv$) A Lie algebroid on a flat line bundle $L\simeq\RR\times M$  is fully encoded by a vector field: one can
see, by fibrewise linearity of the Lie algebroid anchor, that $\poib{f}{g}=f\cdot X(g)-g\cdot X(f)$, where
$X\in\frax(M)$ is defined by $X=\rho_L(1)$.  
\end{example}
%
%
Analogously to the case of the Lie algebroid of a Poisson manifold, it is true in general that the image $\im\rho$
of the anchor of a Lie algebroid  $A\to M$, induces an integrable singular distribution.   On each leaf $\calL$ of
the anchor, $\rho$ has constant rank, thus there is a short exact sequence 
$$
\xymatrix{
0
\:\ar@{->}[r]&\:
\frag_\calL\:\ar@{->}[r]&\:
\left.A\right|_\calL\:\ar@{->}[r]\:&\:
T\calL
\:\ar@{->}[r]&\:
0
}
$$
of vector bundles over $\calL$, for $\frag_\calL:=\Ker_\calL\rho$. One can check that the Lie algebroid bracket on
$\Gamma(M,A)$ restricts to  $\Gamma(\calL,\left.A\right|_\calL)$ and induces a bundle of Lie algebras on 
$\frag_\calL$; the fibres of $\frag_\calL$ are called \tsf{isotropy Lie algebras}. Generalizing the Poisson case,
the anchor foliation on the base manifold, together with the associated short exact sequences, is to be regarded
as a global picture of Fernandes' local structure  theorem for Lie algebroids \cite{fs02}.\\
We describe below the leaves of the Lie algebroids of example \ref{exala}. 
\bgn{example} \tsf{Lie algebroid foliations}. 
\medskip\\
$i$) The leaves of the distribution $\Delta$: each point of $M$, respectively the whole of $M$, in the extreme cases;
\medskip\\
$ii$) All the points of the base manifold;
\medskip\\
$iii$) The orbits of the infinitesimal action;
\medskip\\
$iv$) The flows of $X$.
\end{example}
\eject

\spa Historically, Lie algebroids were discovered in the late sixties of last century by Pradines 
\cite{prd66,prd67}. Lie algebroids are the infinitesimal invariant of suitably smooth groupoids, just as Lie
algebras are for Lie groups, and an understanding of such structures was necessary to develop a geometric study of
groupoids endowed with a suitable smooth structure.
\bgn{definition} A \tsf{Lie groupoid} $\gpdm$ is a groupoid, i.e. a category all of whose arrows are invertible,
such that both the total space $\calG$ and the base space $M$ are smooth manifolds. Moreover the   the unit
section $\eps:M\to\calG$, the map assigning to each object the identity arrow, the inversion map
$\iota:\calG\to\calG$, the composition of arrows, regarded as a map  $\mu:\calG\fib{\sor}{\tar}\calG\to\calG$
and    the source and target maps $\sor,\tar:\calG\to M$ are required to be smooth. For
$\calG\fib{\sor}{\tar}\calG$ to be a manifold and the regularity condition on $\mu$ to make sense it is further 
assumed that the source map be submersive. 
\end{definition}
It follows directly from the definition that the inversion map $\iota:\calG\to\calG$ is a diffeomorphism and the
target map  $\tar:\calG\rightarrow M$ can be written as $\tar=\sor\comp\iota$, therefore it is also submersive; by
applying the inverse function theorem one can show that the unit section is a closed embedding.%
\\
The Lie algebroid of a Lie groupoid $\gpdm$ is obtained through the following construction. Since the groupoid
multiplication is only partially defined, so is the right translation $r_g:\sor\inverse(\tar(g))\to
\sor\inverse(\sor(g))$ by $g\in\calG$ and one cannot hope to reproduce the differentiation of Lie groups to Lie
algebras directly; however, mutatis mutandis, the procedure extends naturally. Being the source map submersive,
its fibres are smooth submanifolds of $\calG$, actually $(\calG,\sor)$ is a simple foliation, and the vectors
tangent to the source fibres $T^\sor\calG:=\ker\dd\sor\to \calG$ form a smooth vector subbundle. The total space
of the Lie algebroid $A$ of $\calG$ is the restriction $T^\sor_M\calG$ to $M$. Right translation is well defined
in $T^\sor\calG$ and the vector fields which are tangent to the source fibres are self-related under right
translation, called $\tsf{right invariant}$ vector fields,  form a $\sor^*\cif(M)$-submodule $\rinv{\frax}(\calG)$ of
$\frax(\calG)$,  isomorphic to the $\cif(M)$-module on $\Gamma(A)$. Analogously to the case of Lie groups and Lie
algebras, one can easily show that  $\rinv{\frax}(\calG)\subset\frax(\calG)$ is a Lie subalgebra,  endowing $A$
with a Lie algebroid bracket; the anchor is the restriction of tangent target $\dd\tar$ to $A$, compatibility with the Lie
bracket follows from the Leibniz rule for $T\calG$.
\\ 
The \tsf{orbit} $\calO_q$ of a Lie groupoid $\gpdm$ through $q\in M$ is by definition
$\calO_q=\tar(\sor\inverse(q))$ and one can show that each source fibre is a smooth principal
$\mathbb{G}_q$-bundle $\tar:\sor\inverse(q)\to\calO_q$ for the \tsf{isotropy group}
$\mathbb{G}_q:=\sor\inverse(q)\cap\tar\inverse(q)$. By construction of the Lie algebroid $A$ of $\calG$, for all
$p\in\calO_q$, $\im_p\rho=T_q\calO_p$, therefore the leaves of $A$ are connected components of the orbits of
$\calG$. We shall denote the orbit space of the foliation induced by $\calG$ with $M/\calG$.%
\\
We describe below the classic examples of Lie groupoids and their Lie algebroids in some detail.
\bgn{example} \tsf{Lie groupoids} 
\medskip\\
$i)$ Any manifold $M$ is trivially a Lie groupoid $\poidd{M}{M}$ over itself. Since necessarily $\sor=\id_M=\tar$,
the only possible multiplication is the identification $\Delta_M\to M$. On the other hand, the \tsf{pair groupoid}
$\poidd{M\times M}{M}$, given by $\sor=\pr_2$, $\tar=\pr_1$, the diagonal $M\to\Delta_M$ as unit section and 
$$
\qquad (x,y)\cdot(y,z)=(x,z) \qquad,\qquad (x,y)\inverse=(y,x) \qquad,\qquad x,y,z\in
M\qquad,
$$
for the multiplication and inversion, is a nontrivial Lie groupoid associated with every manifold. More
generally, to the graph $\gr{\thicksim}$ of any   equivalence relation $\thicksim$ on a smooth manifold $M$ one
can associate a subcategory of the pair groupoid $M\times M$ in the obvious way and conversely every subgroupoid
of $\poidd{M\times M}{M}$ with base $M$ is the graph of an equivalence relation; $\gr{\thicksim}$ is a Lie
groupoid if{f} $\thicksim$ is a \emph{regular equivalence relation}, i.e. if $\gr{\thicksim}$ is smooth and the
restriction of the first and second projection are submersive. The trivial groupoid over a manifold is the 
graph of the trivial equivalence relation, for the pair groupoid all points on the base  are equivalent; the
induced Lie algebroids are the trivial Lie algebroid and the tangent bundle, respectively. 
\medskip\\
$ii)$ A Lie group is a Lie groupoid over the one point manifold.  The \tsf{anchor} of a groupoid $\gpdm$ is the
map  $\chi=(\tar,\sor):\calG\to M\times M$. The image of $\chi$ always contains the diagonal submanifold; if
$\im\,\chi=\Delta_M$, $\calG$ is a  \tsf{smooth bundle of Lie groups}, in the sense that there are no arrows 
connecting different points and $\calG=\coprod_{q\in M}\mathbb{G}_q$. Smooth bundles of Lie groups differentiate
to bundles of Lie algebras.
\medskip\\
$iii)$ Let $\sigma:G\times M\to M$ be the action map of a  Lie group action. Then the product  $G\times M$ carries
a (unique) Lie groupoid structure over $M$, the \tsf{action groupoid} $G\acts M$,  whose  orbits coincide with the
orbits of the action. Source and target maps are $\sor_\acts(g,m)=m$ and $\tar_\acts(g,m)=g*m=\sigma(g,m)$,
$(g,m)\in G\times M$, thus composable pairs are those of the form $[(h,g*m);(g,m)]$, $h\in G$; the multiplication
$\mu_\acts$ is given by $(h,g*m)\comp (g,m)=(hg,m)$, hence the unit section $\eps_\acts$ and inversion
$\iota_\acts$ must be $ \eps_\acts(m)=(e,m)$ and  $\iota_\acts(g,m)=(g\inverse,g*m)$ for the unit element $e$ of
$G$. The Lie algebroid of $G\acts M$ is the action  algebroid of the induced infinitesimal action. The isotropies
of  ($\frag\ltimes M$) $G\acts M$, are precisely the isotropies of the (infinitesimal) action.
\medskip\\
$iv)$ For any left principal $G$-bundle $\pr:P\to M$ the quotient $(P\times P)/G$ for the diagonal action is a Lie
groupoid over $M$, known as the \tsf{gauge groupoid}. Target and source maps are $[p_+,p_-]\mapsto \pr(p_\pm)$,
the multiplication is 
$$
\qquad
[p_+,p_-]\cdot [q_+,q_-]=[g*p_+,q_-]
\qquad,
$$
for the unique $g\in G$ such that $p_-=g*q_+$, unit section and inversion are defined accordingly.  The Lie
algebroid of $(P\times P)/G$ can be identified with the Atiyah Lie algebroid $(TP)/G$ filling the short exact
sequence
\bgn{equation}\label{Atiyah}
\xymatrix{
0\:\ar@{->}[r]&
(P\times \frag)/G\:\ar@{->}[r]&\:
(TP)/G\ar@{->}[r]\:&\:
TM\ar@{->}[r]\:&\: 0
}
\end{equation}
for the adjoint bundle $(P\times \frag)/G$ (see \cite{mackbook} for details). If $M$ is a connected manifold with
universal cover $\wt{M}$, one can apply the above construction  to the covering projection $\pr:\wt{M}\to M$  for
the monodromy action of the fundamental group $\pi_1(M)$ and obtain the \tsf{fundamental groupoid}
$\poidd{\Pi{(M)}}{M}$. Since  $\pi_1(M)$ is discrete the Atiyah sequence (\ref{Atiyah}) yields an identification of
Lie algebroids  $(T\wt{M})/\pi_1(M)\equiv TM$. An alternative characterization of $\Pi(M)$ can be given in terms of
continuous paths in $M$ up to homotopy relative to the endpoints; in this description the groupoid multiplication
is induced by path-concatenation. The construction easily extends to a disconnected manifold $M$.  
\end{example}
\spa Let $\gpdm$ be a Lie groupoid. Since the source map is submersive, For a small enough open neighbourhood  $U'\subset
M$ there always exists a smooth section $\Sigma': U'\to\Sigma'(U')$ of the source map; by linear algebra, $\Sigma'$ can be
perturbed to a new smooth section  $\Sigma:U\to \Sigma(U)$, for some smaller neighbourhood $U\subset U'$, such that
$\tar\comp\sigma'$ is a diffeomorphism. In the above construction source and target can be exchanged and we shall call a
\tsf{local bisection} of $\calG$ a $\rank\,\sor=\rank\,\tar$-codimensional submanifold, on which both source and target
restrict to diffeomorphism to their images; when the domains of a local bisection $\Sigma$, regarded as a section either
way, coincide with $M$,  hence so do the images of $\sor\comp\Sigma$ and $\tar\comp\Sigma$, $\Sigma$ is called a
(\tsf{global})
\tsf{bisection}.
\bgn{example} \tsf{Global bisections}
\medskip\\
$i)$ The base manifold of any Lie groupoid is a bisection. Every section of a vector bundle $\pr: E\to M$ is a 
bisection for the \tsf{abelian groupoid} ($\sor=\!\pr\!=\tar$ and $\mu=+$).
\medskip\\
$ii)$ Bisections of a pair groupoid are in one to one correspondence to diffeomorphisms of the base. 
\medskip\\
$iii)$ The set of a bisections of a group is the group itself.
\medskip\\
$iv)$ See example \ref{patologo} for a groupoid with  points not admitting any global bisection.
\end{example}
Introducing local bisections allows one to deal with left and right translation of  arbitrary tangent vectors. For any
local bisection $\Sigma$ there are  well defined diffeomorphisms
$$
\qquad
\calL^\Sigma:\tar\inverse(\sor(\Sigma))\to \tar\inverse(\tar(\Sigma))
\qquad,\qquad
h\mapsto \Sigma(\tar(h))\cdot h
\qquad,
$$
where $\Sigma$ is regarded as a section of the  source map, and 
$$
\qquad
\calR^\Sigma:\sor\inverse(\tar(\Sigma))\to \sor\inverse(\sor(\Sigma))
\qquad,\qquad
h\mapsto h\cdot\Sigma(\sor(h))
\qquad,
$$
where $\Sigma$ is regarded as a section of the target map. For any  Lie groupoid $\gpdm$, there is a canonical
splitting
$$
\qquad
T_M\calG\simeq T_M^\sor\calG\oplus TM\qquad,\qquad \delta g=\delta g
-\dd\sor\delta g\oplus \dd s\delta g
\qquad,
$$
where we identify $M\simeq\eps(M)$; such a splitting extends to each point $g\in\calG$ off the base manifold, by
picking a local bisection $\Sigma$ through $g$:
$$
\qquad
T_g\calG\simeq T_g^\sor\calG\oplus T_g\Sigma\qquad,\qquad \delta g=\delta^\sor g
\oplus\delta^\sigma g
\qquad,
$$
where $\delta^\sigma g={\dd\calL^\Sigma}_g\dd\sor \delta g$ and $\delta^\sor g:=\delta g-\delta^\sigma g$.
Similarly, splittings subordinated to the target map can be obtained via right translation. One can check that the
composition  
$$
\Sigma_+\cdot\Sigma_-:=\mu(\Sigma_+\times\Sigma_-\cap \calG^{(2)})
$$
of local bisections $\Sigma_\pm$ and the inverse bisection $\Sigma\inverse:=\iota(\Sigma)$ of a local bisection
$\Sigma$ are also  local bisections; as a consequence left and right translations by local bisections form a pseudo
group of transformations: 
$$
\qquad
\bgn{array}{rclcrcl}
\calL^M&=&\id_\calG &\qquad &\calR^M&=&\id_\calG\\
& &\\
\calL^{\Sigma_+\cdot\Sigma_-}&=&\calL^{\Sigma_+}\comp\calL^{\Sigma_-}&\qquad &
\calR^{\Sigma_+\cdot\Sigma_-}&=&\calR^{\Sigma_-}\comp\calR^{\Sigma_+}\\
& &\\
\calL^{\Sigma\inverse}&=&(\calL^\Sigma)\inverse&\qquad &
\calR^{\Sigma\inverse}&=&(\calR^\Sigma)\inverse\\
\end{array}
\qquad, 
$$
the equalities in the second row hold, provided both sides are defined. In particular, restricting to global
bisections defines the  \tsf{group of bisections} $\sf{Bis}(\calG)$. In the following, we shall often  implicitly
make use of local bisections, especially in order tocount dimensions.
\spa A \tsf{morphism of Lie groupoids} is simply a smooth functor; in other words, a morphism of Lie groupoids is
given by a pair of smooth maps  $\fii:\calG^-\to\calG^+$, $f:M^-\to M^+$, which are equivariant with respect to
the groupoid structural maps in all possible ways, in particular so are the operations of left and right
translations: on the domains were the expressions below are defined,
$$
\qquad
\phi\comp l_g= l_{\phi(g)}\comp\phi
\qquad\hbox{ and }\qquad
\phi\comp r_g= r_{\phi(g)}
\qquad
$$ 
hold for all $g\in\calG^-$. Note that, in general, morphisms of Lie groupoids do not map bisections to
bisections. 
\bgn{remark}  For any Lie groupoids $\poidd{\calG^\pm}{M^\pm}$, a smooth map  $\fii:\calG^-\to\calG^+$ is a
morphism of Lie groupoids over the uniquely determined base map $f:M^-\to M^+$, 
$f=\sor_+\fii\comp \eps_-=\tar_+\fii\comp \eps_-$
if{f}  
$$\qquad
(\fii\times\fii\times\fii ) \,\gr{\mu_-}\subset\gr{\mu_+}
\qquad,$$
i.e. if $\fii$ preserves the graph of the partial multiplications.
\end{remark}
A smooth submanifold of a Lie groupoid which is also a subcategory is a \tsf{Lie subgroupoid} if{f} the restriction of the source map remains
submersive; a Lie subgroupoid $\calH$ of a Lie groupoid $\gpdm$ is called a \tsf{wide subgroupoid} if the source
map (hence the target map) of $\calG$ is surjective onto $M$, i.e. if it is a Lie groupoid over $M$. The
\tsf{direct product} of Lie groupoids is defined in the obvious way.
\bgn{example}\label{exaa}\tsf{Subgroupoids and morphisms of Lie groupoids}
\medskip\\
$i$) Let $N$ be a manifold and $\gpdm$ be a (Lie) groupoid. By picking an arbitrary point $n_o\in N$, every
morphism of (Lie) groupoids $\fii:N\times N\to\calG$ over $f:N\to M$ is of the form
$\fii(n_+,n_-)=\psi(n_+)\cdot\psi(n_-)\inverse$,  $n_\pm\in N$,
where $\psi(n)=\fii(n, n_o)$, $n\in N$. Conversely, given (smooth) maps $\psi:N\to \calG$ and $f:N\to M$, such
that $\sor\comp\psi(n)=n_o$, $n\in N$, inverting the procedure yields a morphism of (Lie) groupoids over
$\tar\comp\psi$.
\medskip\\
$ii$) The anchor  $\chi:\calG\to M\times M$ of any Lie groupoid is a morphism to the pair groupoid over the
identity. When it is surjective (thus, as one can check by diagram chasing, submersive), $\calG$ is  said
\tsf{transitive}, in the sense that the base foliation is the base manifold itself. On the other hand, when
$\im\chi=\Delta_M$, the base foliation consists of all the points of $M$ and $\calG$ is said \tsf{totally
intransitive}.
\medskip\\
$iii$) For any Lie groupoid $\gpdm$, the manifold $\calG\fib{\sor}{\sor}\calG$ defines a wide subgroupoid
$\calG\odot\calG$  of the  pair groupoid $\poidd{\calG\times\calG}{\calG}$. It is a Lie subgroupoid with source
fibres $\sor_\odot\inverse(g)=\sor\inverse(\sor(g))\times\{g\}$, $g\in\calG$; to see this, consider that, for all
$\delta g\in T_g\calG$ and $h\in\sor\inverse(\sor(g))$, $\delta
h:=(\dd\calL_\Sigma'\comp\dd\calL_\Sigma\inverse)\delta g$, is a tangent vector at $h$ such that $\dd\sor\delta
h=\dd\sor\delta g$, for any choice of local bisections $\Sigma$ through $g$ and $\Sigma'$ through $h$, therefore
$(\delta h,\delta g)\in T_{(g,h)}\calG\odot\calG$ and  $\dd\sor_\odot(\delta h,\delta g)=\dd\pr_2(\delta h,\delta
g)=\delta g$, i.e. $\sor_\odot$ is submersive. The \tsf{division map} $\delta:\calG\odot\calG\to\calG$, $(g,h)\mapsto
g\cdot h\inverse$, is a morphism of Lie groupoids over the target map of $\calG$.
\medskip\\
$iv$) For any morphism of Lie groupoids $\fii:\calG^-\to\calG^+$ over $f:M^-\to M^+$, the \tsf{kernel groupoid} 
$\poidd{\Ker\fii}{M^-}$ is $\Ker\fii:=\fii\inverse(\eps_+(M^+))$.  It is always a wide subgroupoid of
$\calG^-$, though, in general not a \emph{Lie} subgroupoid.  The kernel of the anchor $\chi$ of a groupoid
$\calG$, is a (topological) bundle of Lie groups $\mathbb{G}=\chi\inverse(\Delta_M)=\coprod_{q\in
M}\mathbb{G}_q$.  A kernel groupoid $\calN$ of a Lie groupoid $\poidd{\calG}{M}$ is a \tsf{normal subgroupoid}, in
the sense that it is wide and, for all $g\in\calG$
$$
g\cdot n\cdot g\inverse\in \NN_{\tar(g)}\quad,\quad\hbox{ for all } n\in\NN_{\sor(g)}\qquad;
$$
in fact, one can show that, taking the quotient of $\calG$ by the equivalence relation
$g\thicksim h$,  if $g=h\cdot n$  for some $n\in\calN$, $g,h\in\calG$,
induced by such a subgroupoid, yields a quotient groupoid  $\poidd{\calG/\calN}{M/\calN}$ over the orbit space of
$\calN$. A subgroupoid of a Lie groupoid is normal if{f} it is the kernel of a morphism of groupoids, since the
quotient projection $\calG\to\calG/\calN$ is a morphism of  groupoids by construction. 
\end{example}
Some constructions in the category of smooth manifolds extend straightforwardly to that of Lie groupoids; for an
instance, consider the next example.
\bgn{example} \tsf{The tangent prolongation groupoid}. If $\gpdm$ is a Lie groupoid, then so is the tangent
prolongation $\poidd{T\calG}{TM}$ for the tangent structural maps of $\calG$. Since the construction is functorial
and a groupoid structure is fully described in terms of diagrams,  $\poidd{T\calG}{TM}$ is clearly a smooth
groupoid; to further check that the tangent source map is submersive is a basic exercise in differential geometry
(see also remark \ref{bundlemaps1}). Computing the Lie algebroid of $\poidd{T\calG}{TM}$ yields the tangent
prolongation Lie algebroid $TA\to TM$ (see \cite{mkz00a,mkz92} for the Lie algebroid anchor and bracket).
\end{example}
Many other interesting constructions require an understanding of
morphisms in the categories of Lie algebroids and Poisson manifolds and 
duality therein. 
\vs{1}
\section{Morphisms and coisotropic calculus}
\begin{quotation} 
In the first part of this Section we first review (following \cite{hm90a,mackbook}) the notion of morphism and the
basic constructions in the category of Lie algebroids, such as direct products and pullbacks. 
The second part is devoted to recall nowadays standard facts about
Poisson maps and coisotropic submanifolds of Poisson manifolds; finally we discuss the dual Poisson-geometric
characterization of morphisms of Lie algebroids and Lie subalgebroids in terms of coisotropic submanifolds.
\end{quotation}

\vs{0.5}
\subsection{Morphisms, pullbacks and direct products of Lie
algebroids}\hfill

\vs{0.1}
\spa  Consider Lie algebroids $A^\pm\to M^\pm$ over different bases and a vector bundle map $\phi:A^-\to A^+$ over
$f:M^-\to M^+$. For such a map to be a morphism of Lie algebroids, the compatibility with the anchors is expressed by
the natural condition $\rho_+\comp\phi=\dd f\comp\rho_-$; this way leaves of $A^-$ are mapped to leaves of $A^+$.
Whenever $M^+=M^-$ and $f$ is the identity map, the bracket compatibility condition is also clear: the induced map of
sections $\Gamma(A^-)\to\Gamma(A^+)$ has to be a morphism of Lie algebras. The main complication encountered in
defining morphisms of Lie algebroids over different bases is met in not having at disposal  a natural way of mapping
sections of $A^-$ to sections of $A^+$.  Let us first consider a special case: suppose the vector bundle fibred
product $TM^-\fib{\dd f}{\rho_+}A^+$ exists; upon identifying the base manifold $M^-\fib{f}{} M^+\simeq\gr{f}$ with
the graph of $f$, yields a vector bundle  $f\daga A^+: TM^-\fib{\dd f}{\rho_+}A^+\to M^-$. 
\bgn{proposition}\label{morphismsla} Let $A^\pm\to M^\pm$ be Lie algebroids and  $\phi:A^-\to A^+$  a vector bundle
map over $f:M^-\to M^+$. Then there exists a unique Lie algebroid  on the vector bundle $f\daga A^+$, whose anchor
$\rho_{{}\daga}$ is the restriction of the first projection $TM^-\times A^+\to TM^-$.
\end{proposition}
Last proposition is a restatement of some of the results contained in  \cite{hm90a}; we shall only present the idea
of the proof. 
\bgn{proof}[Sketch of proof of proposition \ref{morphismsla}] A section $a\daga$  of  $f\daga A$ is given by a pair
$(X, \wt{a})\in\frax(M^-)\oplus\Gamma(f^{\sf +}A)$, satisfying $\dd f \comp X =\rho_+\comp\wt{a}$; on the other hand, the
space sections of the  pullback bundle\fn{Recall that, for any smooth map $f:N\to M$ and smooth vector bundle 
$\pr: A\to M$,  the space $\{(a,n)\:|\:\pr(a)=f(n)\}$ is canonically a smooth vector bundle over $M^-$.}
$f^{\sf +}A^+\to M^-$ is isomorphic to
$\cif(M^-)\otimes_{\cif(M^+)}\Gamma(A^+)$ as a left $\cif(M^-)$-module, where the tensor product is taken for the
right $\cif(M^+)$-module on $\cif(M^-)$ induced by precomposition with $f$. Then for any
$a_{1,2}\daga\in\Gamma(f\daga A)$ there exist decompositions of the form 
$$
\qquad
a_{1,2}\daga
= 
X_{1,2}\oplus\underset{k_{1,2}}\sum u_{1,2}^{k_{1,2}}(a_{1,2}^{{k_{1,2}}}\comp\phi)
\qquad,
$$
in terms of $X_{1,2}\in\frax(M^-)$, $\{u_{1,2}^{k_1,k_2}\}_{k_1,k_2}\subset\cif(M^-)$ and 
$\{a_{1,2}^{{k_{1,2}}}\}_{k_1,k_2}\subset\Gamma(A^+)$, which can be used to define
\bgn{eqnarray*}
\brak{a_1\daga}{a_2\daga}\daga
&:=\:&
\brak{X_1}{X_2}\\
&\:\oplus\:&\left\{
\underset{k_1,k_2}\sum u_{1}^{k_{1}}\cdot u_{2}^{k_{2}}(\brak{a_1^{k_{1}}}{a_2^{k_{2}}}\comp\phi)
\quad+\quad
\underset{k_2}\sum X_1(u_2^{k_2})(a_2^{k_2}\comp\phi)\right.\\
&\:-\:&\left.
\underset{k_1}\sum X_2(u_1^{k_1})(a_1^{k_1}\comp\phi)\right\}\hs{1.9}.
\end{eqnarray*}
It turns out that the above bracket does not depend on the choice of the decompositions  and yields a well defined bilinear 
skewsymmetric
operation on $\Gamma(f\daga A)$; the definition is tailored to have the Leibniz rule satisfied and one can check
that also the Jacobi identity holds.
\end{proof}
The Lie algebroid structure on $f\daga A^+\to M^-$ is that of a \tsf{pullback Lie algebroid} along $f$. When such a
Lie algebroid exists, for instance when $f$ is submersive or $\rho$ has maximal rank, it results then natural to ask the 
induced map $\Gamma(A^-)\to \Gamma(f\daga A^+)$, 
$a\mapsto \rho_-\comp a\oplus \phi\comp a$, to preserve the Lie brackets in order for $\phi$ to be a morphism of Lie
algebroids over $f$; provided the anchor compatibility holds, this requirement amounts to
\bgn{equation}\label{morphlaanc}
\qquad
\rho_-\comp\brak{a_1}{a_2}=\brak{\rho_-\comp a_1}{\rho_-\comp a_1}
\qquad,
\end{equation}
  for all 
$a_{1,2}\in\Gamma(A^-)$, on the $\frax(M^-)$ component  and
\bgn{eqnarray}\label{morphladec}
\nn\phi\comp\brak{a_1}{a_2}
&=&	
\underset{k_1,k_2}\sum 
u_{1}^{k_{1}}\cdot u_{2}^{k_{2}}
(\brak{a_1^{k_{1}}}{a_2^{k_{2}}}\comp\phi)
\quad+\quad
\underset{k_2}\sum \rho_-(a_1)(u_2^{k_2}) ( a_2^{k_2}\comp\phi)\\
&-&
\underset{k_1}\sum \rho_-(a_2)(u_1^{k_1}) ( a_1^{k_1}\comp\phi)
\end{eqnarray}
on the $\Gamma(f^{\pmb{+}}A^+)$ component, for any choice of decompositions
\bgn{equation}\label{dec}
\phi\comp a_{1,2}=
\underset{k_{1,2}}\sum u_{1,2}^{k_{1,2}} (a_{1,2}^{{k_{1,2}}}\comp\phi)
\qquad,
\end{equation}
in terms of $\{u_{1,2}^{k_1,k_2}\}_{k_1,k_2}\subset\cif(M^-)$ and 
$\{a_{1,2}^{{k_{1,2}}}\}_{k_1,k_2}\subset\Gamma(A^+)$. Even if the pullback Lie  algebroid along $f$ does
not exist, condition (\ref{morphladec}) still makes sense and does not depend on the choice of
decompositions, provided the anchor compatibility condition holds. This leads to the general definition of a
morphism of Lie algebroids.
\bgn{definition}\label{ordillo}  Let $A^\pm\to M^\pm$ be Lie algebroids and $\phi:A^-\to A^+$  a vector bundle map
over $f:M^-\to M^+$ which is compatible with the anchor maps in the sense that $\rho_+\comp\phi=\dd f\comp
\rho_-$.  Then $(\phi,f)$ is a \tsf{morphism of Lie algebroids} if the bracket compatibility (\ref{morphladec})
holds for all $a_{1,2}\in\Gamma(A^-)$ and any choice of decompositions  (\ref{dec}) of $\phi\comp a_{1,2}$.
\end{definition}
This notion of morphism is consistent with that of a category of Lie algebroids; even though it might seem quite
intractable, this definition is sufficiently effective to deal with general constructs such as subobjects and
direct products, as it was shown in \cite{hm90a}.
\bgn{example}\label{mphalgebroids} \tsf{Morphisms of Lie algebroids}. We list below some fundamental examples; the details, 
which can be found in \cite{mackbook}, are left to the reader as an exercise.  
\medskip\\
$i)$ If $\phi:A^-\to A^+$ is morphism of Lie algebroids over the identity map for $M^-=M^+$, the induced map of
sections  $\Gamma(A^-)\to\Gamma(A^+)$ is required to preserve the Lie brackets; in particular for, morphisms of Lie
algebras are morphisms of Lie algebroids ($M^-=M^+=\bullet$).
\medskip\\
$ii)$ For any smooth map $f\colon M\to N$, $\dd f\colon TM\to TN$ is a morphism of Lie algebroids; this can be checked using
decompositions and equating the sides bracket compatibility condition, regarding vector fields as derivations.
\medskip\\
$iii)$ The \tsf{inductor} $f\daga: f\daga A\to A$, $(X,a)\mapsto a$, for the pullback of a Lie algebroid of $A\to M$ along $f:N\to
M$, is a morphism of Lie algebroids.
\medskip\\
$iv)$ For any morphism of Lie groupoids $\fii:\calG^-\to\calG^+$ over $f:M^-\to M^+$, setting 
$$
\phi:=\dd\fii|_{T^{\sor_-}_{M^-}\calG^-}: A^-\to A^+
$$
yields a well defined vector bundle map over $f$. The anchor compatibility holds, since
$$
\quad
(\rho_+\comp\phi)(a)=(\dd\tar_-\comp\dd\fii)(a)=(\dd f\comp\dd\tar_-)(a)=
(\dd f\comp\rho_-)(a)
\quad,\quad a\in\Gamma(A^-)
\quad,
$$
and the bracket compatibility follows from the properties of right invariant vector fields.
\end{example}
Last example shows that there exist a \tsf{Lie functor} from the category of Lie groupoids to that of Lie
algebroids extending the classical Lie functor. 
\bgn{remark}\label{torsolo} Let $A^\pm\to M^\pm$ be Lie algebroids and $\phi:A^-\to A^+$  a vector bundle map over $f:M^-\to M^+$
which is compatible with the anchor maps in the sense of definition \ref{ordillo}. Then one can pick any
connection $\nabla$ for $A^+\to M^+$ to express the bracket compatibility more intrinsically. Denote with
$\ol{\nabla}$ the pullback connection induced by $\nabla$ on $f^{\pmb{+}}A^+\to M^-$, with
$\tau_\nabla\in\Gamma(A\otimes \wedge^2 A^*)$ the torsion tensor
$$
\tau_\nabla(a,b):=\nabla_{\rho^+(a)}b-\nabla_{\rho^+(b)}a - \brak{a}{b}\qquad,\qquad a,b\in\Gamma(A)
$$
and with $f^{\pmb{+}}\tau_\nabla$ its pullback to $f^{\pmb{+}}A^+$.  Then $(\phi,f)$ 
satisfies the bracket compatibility condition (\ref{morphladec}) if{f}
$$
\phi\comp\brak{a_1}{a_2}
=
\ol{\nabla}_{\rho^-(a_1)}\phi\comp a_2-\ol{\nabla}_{\rho^-(a_2)}\phi\comp a_1 -f^{\pmb{+}}\tau_\nabla({a_1},{a_2})
$$
holds for all $a_{1,2}\in\Gamma(A^-)$ \cite{hm90a}.
\end{remark}

\spa Let us now turn to subobjects (see \cite{mackbook} for details).
\bgn{definition} Let $A\to M$ be a Lie algebroid and $N\subset M$ an embedded submanifold. A Lie algebroid $B\to N$
on a vector subbundle of  $A\to M$  is a \tsf{Lie subalgebroid} if the inclusion $B\inc A$ is a morphism of Lie
algebroids over the inclusion $N\inc M$.
\end{definition}
From the last definition it is clear that restrictions of morphisms of Lie algebroids to Lie subalgebroids yield
morphisms of Lie algebroids.%
\\
As it is to be expected, Lie subalgebroids are precisely those vector subbundles for which the algebraic
operations of the ambient Lie algebroid suitably  restrict.
\bgn{lemma}\label{subalgebroid} A vector subbundle $B\to N$ of a Lie algebroid $A\to M$ is a Lie subalgebroid if{f}
the following conditions hold
\medskip\\
$1$. The anchor $\rho:A\to TM$ restricts to a bundle map $B\to TN$;
\medskip\\
$2$. For all $a_\pm\in\Gamma(A)$ such that $a_\pm|_N\in \Gamma(B)$, $\brak{a_+}{a_-}|_N\in \Gamma(B)$; 
%
%
\end{lemma}
\bgn{remark}\label{checksubal} It is always possible to restrict a Lie algebroid to a open submanifold of the base. Since
the conditions of lemma (\ref{subalgebroid}) are local, they can always be checked on some neighbourhood  $U$, open in
$M$, of an arbitrary point of $N$. Note, in particular, that for all $a_\pm$ such as in condition $2.$ above, it follows
that $\brak{a_+}{a_-}|_N=0$, whenever $a_-|_N=0$; that is, the bracket of sections of $\Gamma(B)$ can be computed using
arbitrary extensions.  Since for any frames $\{e^\alpha\}$, $\{e_\alpha\}$ in duality for $A$ and $A^*$ on (a restriction
of) $U$, the Leibniz rule implies
\be
\brak{a^1}{a^2}
&=&
\sum_\alpha a^2_\alpha\brak{a^1}{e^\alpha}+\rho(a^1)(a^2_\alpha)e^\alpha\\
&=&
\sum_\alpha a^2_\alpha\brak{a^1}{e^\alpha}+\pair{\dd a^2_\alpha}{\rho(a^1)}\hs{2},
\ee
where $a^2_\alpha:=\pair{e_\alpha}{a^2}\in\cif(U)$ vanishes on $U\cap N$ for all $\alpha$'s,  thus 
$\pair{\rho(a^1)}{\dd a^2_\alpha}$ also vanishes on $U\cap N$,  since the anchor restricts to a bundle map
$B|_{U\cap N}\to T(U\cap N)$.
\end{remark}
Bearing condition (\ref{morphladec})  in mind the proof of last lemma is straightforward by working locally and
picking extensions. As a consequence of example (\ref{mphalgebroids}, $iv$) Lie subgroupoids differentiate to Lie
subalgebroids.  There is a useful corollary to lemma \ref{subalgebroid}.
\bgn{corollary}\label{subalsubal} Given Lie algebroids $A$, $B$, $C$ and a sequence of vector bundle inclusions
$A\subset B\subset C$, if $A\subset C$ and $B\subset C$ are Lie subalgebroids, then so is $A\subset B$.
\end{corollary}
Lie subalgebras are obviously Lie subalgebroids over the one point manifold; we list a few other examples.
\bgn{example} \tsf{Lie subalgebroids}.
\medskip\\
$i$) An integrable regular distribution is a Lie subalgebroid of  the tangent bundle. More generally,  the tangent
bundle of a submanifold is a Lie subalgebroid of the ambient  manifold.
\medskip\\
$ii$) Let $P$ be a Poisson manifold and $C\subset P$ a submanifold; then the conormal bundle $N^*C\subset T^*P$ is
a Lie subalgebroid if{f} $C$ is coisotropic. We shall discuss this example in the next Section.
\medskip\\
$iii$) For any Lie algebroid $A$ and leaf $\calL$ of the anchor distribution the restriction $A|_\calL\to\calL$ is
a Lie subalgebroid. 
\medskip\\
$iv$) Let $A\to M$ and $B\to N$ be Lie algebroids and $\phi:A\to B$ a vector bundle map of constant rank, so that
its kernel is a vector subbundle of $A$; then, $\ker \phi\to M$ is a Lie subalgebroid whenever $\phi$ is a
morphism of Lie algebroids. In particular the bundle of isotropy Lie algebras $\frag_\calL$ of a leaf $\calL$ of a
Lie algebroid $A$ is a Lie subalgebroid of $A|_\calL$.
\end{example}
\spa Pullback Lie algebroids are not only useful to understand morphisms of Lie algebroids, but also play an
essential role to produce  general constructs, such as direct products, described below,  in the category of Lie
algebroids,  thanks to the following universal property.  
\bgn{proposition}\label{indio} Let $A\to M$ be a Lie algebroid and $f:M'\to M$ be a smooth map. Assume that the
pullback Lie algebroid $f\daga A\to N$ exists. Then for any morphism of Lie algebroids $\phi: B\to A$ over 
$g:N\to M$ factoring through $h:N\to M'$ there exists a unique morphism of Lie algebroids  $\psi:B\to
f^{\pmb{++}}A$, such that $f^{\pmb{++}}\comp\psi=\phi$.
\end{proposition}
Next, consider Lie algebroids $A^1$, $A^2$ and $B$ \emph{over the same base} $M$; given morphisms of Lie algebroids
$\phi_{1,2}: A^{1,2}\rightarrow B$ over the identity, such that the fibred product  
$A^1\fib{\phi_1}{\phi_2} A^2\to \Delta_M\simeq M$ is
a vector bundle, it is possible to introduce the  \tsf{product Lie algebroid}, over $B$ in this case,  for the vector
bundle  structure over $M$. The anchor is given by
$$
\qquad
\rho(a_1\oplus a_2)=\rho_B\comp\phi_1(a_1)
=\rho_B\comp\phi_2(a_2)\qquad,\qquad
a_1\oplus a_2\in A^1\fib{\phi_1}{\phi_2} A^2\qquad;
$$
the bracket is defined componentwise:
$$
\brak{a_1\oplus a_2}{b_1\oplus b_2}=\brak{a_1}{b_1}\oplus
\brak{ a_2}{b_2}\quad,\quad a_1\oplus a_2,b_1\oplus b_2\in
 \Gamma(M,A^1\fib{\phi_1}{\phi_2} A^2)
\quad.$$
Last construction is a straightforward generalization of the  \emph{fibred product of Lie
algebroids over the same base} in \cite{hm90a}, which is recovered  replacing $B$ with $TM$ and
$\phi_{1,2}$ with $\rho_{1,2}$. General fibred products shall be studied later on in Section
\ref{fplgla}.\\
Given Lie algebroids $A^{1,2}\rightarrow M^{1,2}$, denote with $M^{12}$ the direct product
$M^1\times M^2$ and with $\rm{pr}_{1,2}$ the projections onto $M^{1,2}$. Since the pullback
algebroids $\rm{pr}_{1,2}\daga A^{1,2}$ always exist and the  fibred product of manifolds 
$\rm{pr}_1\daga A^1\fib{\rho^{1}_{{\pmb{+}}{\pmb{+}}}}{\rho^{2}_{{\pmb{+}}{\pmb{+}}}}
\rm{pr}_2\daga A^2$ is to be identified with the vector bundle $A^1\times A^2\rightarrow
M^1\times M^2$, there is always a fibred product Lie algebroid over $TM^{12}$, the  \tsf{direct
product Lie algebroid} (denoted simply as $A^1\times A^2$) of $A^1$ and $A^2$; applying
proposition \ref{indio}, it is straightforward to check that it is indeed a direct 
product in the category of Lie algebroids, satisfying the relevant universal property.\\
The Lie bracket on $A^1\times A^2$ can be described explicitly as follows. Note that
$\Gamma(A^1\times A^2)=\Gamma(\pr_1^+A^1)\oplus\Gamma(\pr_2^+A^2)$;
for any choice of decompositions
$$
\alpha_{1,2}=\sum_{k_{1,2}}u^{1,2}_{k_{1,2}}(a^{1,2}_{k_{1,2}}\comp\pr_{1,2})
\qquad\hbox{and}\qquad
\beta_{1,2}=\sum_{l_{1,2}}v^{1,2}_{l_{1,2}}(b^{1,2}_{l_{1,2}}\comp\pr_{1,2})
$$
of $\alpha_{1,2},\beta_{1,2}\in\Gamma(\pr_{1,2}^+A^{1,2})$
with $\{u^{1,2}_{k_{1,2}}\},\{v^{1,2}_{l_{1,2}}\}\in\cif(M^1\times M^2)$ and 
$\{a^{1,2}_{k_{1,2}}\},\{b^{1,2}_{l_{1,2}}\}\in\Gamma(A^{1,2})$, the components of 
$[\alpha_1\oplus\alpha_2,\beta_1\oplus\beta_2]$ are given by ($i=1,2$)
\bgn{eqnarray}\label{directo}
[\alpha_1\oplus\alpha_2,\beta_1\oplus\beta_2]_i
\nn
&=&
\sum_{k_i,l_i}u^i_{k_i}v^i_{l_i}(\brak{a^i_{k_i}}{b^i_{l_i}}\comp\pr_i)\\
\nn&+&
\sum_{l_i}(\rho_1(\alpha_1)\times\rho_2(\alpha_2))(v^i_{l_i})({b^i_{l_i}}\comp\pr_i)\\
&-&
\sum_{k_i}(\rho_1(\beta_1)\times\rho_2(\beta_2))(u^i_{k_i})({a^i_{k_i}}\comp\pr_i)\hs{1}.
\end{eqnarray}
%
%
%
%
%
%
%
%
%
\subsection{Poisson maps and coisotropic calculus}\hfill

\vs{0.1}
\spa Morphisms in the category of Poisson manifolds are smooth maps inducing morphisms of Poisson
algebras.
\bgn{definition} Let $P_\pm$ be Poisson manifolds; a smooth map  $P_-\to P_+$ is called a
\tsf{Poisson map} if its pullback  $\cif(P_+)\to\cif(P_-)$ is a morphism of Lie algebras.
\end{definition} 
\bgn{remark} Not every notion in the category of Poisson manifolds extends
some symplectic notion: a symplectomorphism is not a Poisson map, and vice versa,
unless it is a local diffeomorphism.
\end{remark}
The notion of \emph{Poisson submanifold} is also fully established and natural, namely a
submanifold $Q$ of a Poisson manifold $P$, such that the Poisson bivector field of $P$ restricts to a
Poisson bivector field on $P'$; equivalently, a Poisson manifold $Q$, $Q\subset P$ is Poisson if the inclusion is a
Poisson map. However, in relation with duality issues, the role played by
subobjects is often inherent to \emph{coisotropic submanifolds}, in a sense we are to specify.
\bgn{definition} A submanifold $C$ of a Poisson manifold $(P,\pi)$ is called a \tsf{coisotropic
submanifold} if $\pi^\sharp N^*C\subset TC$.
\end{definition}
A coisotropic submanifold $C$ of a Poisson manifold $(P,\pi)$ comes equipped with a
\tsf{characteristic distribution} $\Delta^C$ spanned by the Hamiltonian vector fields of functions
in the \tsf{vanishing} (associative) \tsf{ideal} $\calI_C=\{f\in\cif(P)\:|\:f|_C=0\}$, i.e.
$$
\qquad
\Delta^C_c=\sf{span}\{\pish_c\dd f\:|\: f\in\calI_C\}
\qquad.
$$ The same symbol for $\Delta^C$ and the
corresponding subset of $TC$ shall occasionally be used. It is well known that $\Delta^C$ is always
integrable and we shall denote with $\ul{C}$ the leaf space of the associated singular foliation; if
$\ul{C}$ is a smooth manifold, one can show that $\cif(\ul{C})\simeq N{\calI_C}/\calI_C$, where
$$
N\calI_C=\{\:f\in\cif (P)\:\:|\:\: X(f)=0\:,\: X\in \Delta^C\:\}
$$ 
is the normalizer of $\calI_C$, i.e. the
smallest Poisson subalgebra of $\cif(P)$ making $\calI_C$ an ideal for the Poisson bracket.  
\bgn{proposition}\label{coisotropy} Let $(P,\pi)$ be a Poisson manifold and 
$C\subset P$ be a submanifold. Then, the following are equivalent:
\medskip\\$i)$ The conormal bundle $N^*C$ is a Lie subalgebroid of the
Koszul algebroid $T^*P$;
\medskip\\$ii)$ The characteristic ideal $\calI_C$ is a Poisson subalgebra of $\cif(P)$;
\medskip\\$iii)$ The characteristic distribution $\Delta^C$ is tangent to $C$;  
\medskip\\$iv)$ $C$ is coisotropic.
\end{proposition}
\bgn{proof}  ($ii$) $\Rightarrow$ ($iii$) $\Rightarrow$ ($iv$): obvious. 
($i$) $\Rightarrow$ ($ii$):  $f_\pm\in\calI_C$ if{f} $\dd
f_\pm|_C\in\Gamma(N^*C)$, then $ \poib{f_+}{f_-}=\pair{\dd f_-}{\pish\dd f_+} $ vanishes on $C$,
since $\pish$ is the anchor of $T^*P$.
($iv$) $\Rightarrow$ ($i$): for all $\theta\in\Ohm^1(P)$ and $X\in\frax(P)$, 
such that $\theta|_C\in\Gamma(N^*C)$  $X|_C\in\frax(C)$, $\theta(X)\in\calI_C$ and, by Cartan's
magic formula $\iota_X\dd\theta\in\Gamma(N^*C)$; use formula (\ref{Koszulint}), to check condition 2. of lemma
\ref{subalgebroid} (see also remark (\ref{checksubal})).
\end{proof}
Note that the direct product $P_1\times P_2$ of Poisson manifolds  $(P_{1,2},\pi_{1,2})$ is
canonically endowed with a product Poisson tensor $\pi_1\times\pi_2$; the Poisson bracket on
$\cif(P_1\times P_2)$ can be expressed as
$$
\qquad
\poib{F}{G}_{P_1\times P_2}(p_1,p_2)
=
\poib{F^1_{p_2}}{G^1_{p_2}}_{P_1}(p_1)
+
\poib{F^2_{p_1}}{G^2_{p_1}}_{P_1}(p_2)
\qquad,
$$
for all $F,H \in\cif(P_1\times P_2)$ and $p_{1,2}\in P_{1,2}$. 
Here $H^{1,2}_x\in\cif(P_{1,2})$ denotes the restriction to the
$P_{1,2}$-direction of any $H\in\cif(P_1\times P_2)$, $x\in P_{1,2}$.\\ 
For any Poisson manifold $(P, \pi)$, denote with $\ol{P}$ the Poisson manifold $(P, -\pi)$. 
\bgn{corollary} Let $P_\pm$ be Poisson manifolds and $\lambda:P_-\to P_+$ a smooth map. Then
$\lambda$ is Poisson if{f}  $\gr{\lambda}\subset P_-\times \ol{P}_+$ is coisotropic.
\end{corollary}
\bgn{example} For any Poisson manifold $P$, the diagonal $\Delta_P\subset P\times\ol{P}$ is the 
graph of the identity, therefore it is coisotropic, according to last corollary; on the other hand it
is straightforward to check coisotropicity directly.
\end{example}
In fact, coisotropic submanifolds of direct products 
are to be thought of as generalized morphisms of Poisson
manifolds.
\bgn{definition}~\cite{ws88} Let $P_\pm$ be Poisson manifolds. A \tsf{coisotropic relation} $R:P_-\to P_+$ is a
coisotropic submanifold $R\subset P_-\times \ol{P}_+$. For any Poisson manifolds $P_{1,2,3}$, the
\tsf{composition of coisotropic relations} $R_{12}:P_1\to P_2$ and $R_{23}:P_2\to P_3$ is 
$$
R_{12}\comp R_{23}=\pr_1\times\pr_4(R_{12}\times R_{23}\cap P_1\times \Delta_{P_2}\times
P_3)\subset P_1\times P_3\qquad.
$$
Coisotropic relations $R_{12}:P_1\to P_2$ and $R_{23}:P_2\to P_3$ are in \tsf{very clean position} if
\bgn{itemize}
\item[1.] $R_{123}:=(R_{12}\times R_{23}$ and $\Delta_{123}:=P_1\times \Delta_{P_2}\times P_3$ intersect cleanly;
\item[2.] The restriction of the projection $\Delta_{123}\to P_1\times P_3$ has constant rank;
\item[3.] $R_{12}\comp R_{23}\subset P_1\times P_3$ is a smooth submanifold; 
\item[4.] The projection $R_{123}\to R_{12}\comp R_{23}$ is a submersion.
\end{itemize}
\end{definition}
Note that the composition of coisotropic relations is not, in general a smooth submanifold. The
fundamental theorem of coisotropic calculus, which we give without proof, is due to Weinstein.
\bgn{theorem}\cite{ws88}\label{coisorel} For any coisotropic relation $R_{12}:P_1\to P_2$ and  $R_{23}:P_2\to P_3$ in
very clean position the composition $R_{23}:P_2\to P_3$ is a coisotropic relation $P_1\to P_3$.
\end{theorem}
\bgn{corollary} If $j:P_-\to P_+$ is a Poisson or anti-Poisson submersion and $C\subset P_+$ a
coisotropic submanifold, then so is $j\inverse (C)\subset P_-$.
\end{corollary}
\bgn{proof} $j\inverse (C)\times\bullet=\gr{j}\comp (C\times\bullet)\subset P_+\times\bullet$,  where $\gr{j}$ is a
coisotropic relation either $P_-\to P_+$ or $P_-\to\ol{P}_+$ (note that $C$ is coisotropic also in $\ol{P}_+$) and
submersivity of $j$ implies the cleanliness conditions.
\end{proof}
Lie subalgebroids and morphisms of Lie algebroids can be characterized  in the dual picture in terms of coisotropic
submanifolds.
\bgn{proposition}\label{subalandcoiso}\cite{xu95}
Let $A\rightarrow M$ be a Lie algebroid and $B\rightarrow N$ a smooth vector
subbundle. Then $B$ is a Lie subalgebroid if{f} the annihilator $B^o$ of $B$ is a
coisotropic submanifold of $A^*$.
\end{proposition}
\begin{proof} Suppose $B$ is a Lie subalgebroid; then, under the identification 
of 
$\Gamma(A)$ with the space of fibrewise linear functions $\cif_{\sf{lin}}(A^*)$,
condition (2.) of lemma (\ref{subalgebroid}) is equivalent to%
\medskip\\
$2^*$. For all 
$
F, G\in\cif_{\sf{lin}}(A^*)\cap\calI_{B^o}=:\calI^{\sf{lin}}_{B^o}$,
$\{F,G\}\in\calI^{\sf{lin}}_{B^o}$,
\medskip\\
for the induced Poisson bracket on $A^*$ and the vanishing ideal $\calI_{B^o}$.
Since $\calI_{B^o}$ is generated by 
$\calI^{\sf{lin}}_{B^o}$ as a 
$\cif(A^*)$-module, $2^*$ implies 
$$
\{\calI_{B^o},\calI_{B^o}\}\subset\calI_{B^o}
$$
by the Leibniz rule, i.e. $B^o$ is coisotropic. 
Conversely, to each $f\in\cif(M)$, associate the fibrewise constant function
$\pr_{A^*}^*f$, on $A^*$;
$\pr_{A^*}^*f\in\calI_{B^o}$ if{f} $f\in\calI_{N}$.
For any $a\in\Gamma(A)$, such that
$a|_N\in\Gamma(B)$,
$$
\{F_a,\pr_{A^*}^*f\}:=\pr_{A^*}^*({\rho(a)(f)})
=
\langle\dd\,f,\rho(a)\rangle\comp\pr_{A^*}
$$
if $B^o\subset A^*$ is coisotropic,   $\rho(a)(f)\in\calI_{N}$, for all $f\in\calI_{N}$, equivalently, $\rho(a)$
is a section of $TN=(N^*N)^o$, since the conormal bundle $N^*N$ is spanned by $\dd\,\calI_{N}$. That is, the
anchor of $A$ restricts to $B$.  Conditions (2.) and (1.) follow straightforwardly from the fact that 
$\calI_{B^o}$ is closed of the Poisson bracket and $\calI^{\sf{lin}}_{B^o}\subset\calI_{B^o}$ a Poisson
subalgebra, since the dual bracket on $A^*$ is fibrewise linear.
\end{proof}
\bgn{corollary}\label{caraci} Let $A^\pm\rightarrow M^\pm$ be Lie algebroids and $\phi:A^-\rightarrow A^+$ a
morphism of vector bundles over $f:M^-\rightarrow M^+$.  Then the following are equivalent:
\medskip\\
$i$) $\phi$ is a morphism of Lie algebroids;
\medskip\\
$ii$) $\gr{\phi}\subset  A^-\times A^+$ is a Lie subalgebroid;
\medskip\\
$iii$) ${\gr{\phi}}^o\subset A^{-*}\times A^{+*}$ is coisotropic for the 
induced linear Poisson structure.
\medskip\\
Moreover, if $\phi$ is base bijective, so that the transpose map $\phi^{\,\sf t}$ is well defined the three
statements above are equivalent to%
\medskip\\
$iv$) $\phi^{\sf{t}}: A^{+*}\rightarrow A^{-*}$ is Poisson.
\end{corollary}
\bgn{proof} For the equivalence of the first three statements, it  suffices to show the equivalence of ($i$) and
($ii$); this can be checked by picking
decompositions. 
When the transpose map is well defined $\gr{\phi}^o$ coincides with the vector bundle
$(-\id_{A^{-*}}\times\id_{A^{+*}})\,\gr{\phi^{\sf{+}}}\to \gr{f}$,  up to the exchange of $A^-$ with $A^+$; then
$\gr{\phi}^o$ is coisotropic in $A^{-*}\times A^{+*}$ if{f} $\gr{\phi^{\sf{+}}}$ is coisotropic in  $A^{+*}\times
\ol{A^{-*}}$, equivalently, if{f} $\phi^{\sf{+}}$ is Poisson.
\end{proof}
\spa We conclude this Subsection with a few examples of coisotropic submanifolds.
\bgn{example} \tsf{Coisotropic submanifolds}
\medskip\\
$i)$ Consider any Poisson bivector on $\RR^3$; it is a straightforward exercise to show that all  submanifolds of
the form $\{z=c(x,y)\}$, for some smooth function $c\in\cif(\RR^2)$ are coisotropic.
\medskip\\
$ii$) Recall that a submanifold $C$ of a symplectic manifold $(M,\ohm)$ is coisotropic if
  $T^\ohm C\subset TC$, for the symplectic orthogonal bundle 
$$
\qquad
T^{\ohm}_c C
=
\{\delta c\in T_c M
\:|\:
\ohm_c(\delta q ,\delta c)=0
\:,\:
\delta c\in T_c C\}\quad,\quad 
c\in C\qquad.
$$ 
Since $T^\ohm C=(\ohm^\sharp TC)^o$, the condition is equivalent to  $N^*C=T^oC\subset \ohm^\sharp TC$, i.e. a
submanifold of a symplectic manifold is coisotropic if{f} it is coisotropic for the inverse Poisson structure
$\pi^\sharp=\ohm^{\sharp-1}$.
\medskip\\
$iii)$  The symplectic leaves of a Poisson manifold are, essentially by definition, coisotropic submanifolds.
\medskip\\
$iv)$ Let $M$ be a symplectic manifold and $C$ a smooth submanifold defined by constraints $f_i=0$, 
$i=1,\dots,n$, $f_i\in\cif(M)$. The constraints are first class in the sense of Dirac if
$$
\poib{f_i}{f_j}=c_{ij}{^k} f_k
$$	
for some smooth functions $\{c_{ij}{^k}\}$ on $M$. Since $N^*C$ is spanned by the differentials $\{\dd f_i\}$,
the constraints are first class if{f} $C$ is coisotropic.
\end{example}
\vs{1}
\section{Actions of Lie groupoids and action Lie groupoids}\label{algalg}
\begin{quotation} 
We describe in this Section how Lie group actions and infinitesimal actions can be generalized replacing  Lie
groups--algebras with Lie groupoids--algebroids; we also present the natural conditions for the existence in the
smooth category of quotients with respect to such actions. We conclude by introducing principal bundles with
structure groupoid, which prepresent a handy tool in the study of the integrability and reduction of Lie groupoids and Lie
algebroids.
\end{quotation}

\vs{0.1}
\spa The main difference between the actions of Lie groups and
those of Lie group-oids consists in the fact that a Lie groupoid $\gpdm$ acts 
on maps to $M$ rather than on manifolds.
\bgn{definition} Let $\gpdm$ be a groupoid and consider a map $j:N\to M$. A
\tsf{left groupoid action} of $\calG$ on $j$ is given by a map
$\sigma:\calG\fib{\sor}{j}N\to N$, $(g,n)\mapsto g*n$, such that, for all
$n\in N$,
\medskip\\1. $j(g*n)=\tar(g)$ , whenever $g\in\sor\inverse(j(n))$;  
\medskip\\2. $\eps(j(n))*n=n$;
\medskip\\3. $(gh)*n=g*(h*n)$,
whenever $(g,h)\in\calG^{(2)}$ and $h\in\sor\inverse(j(n))$.
\medskip\\
The map $j$ is called a \tsf{moment map}. 
\end{definition} 
The action of a group on a set $M$ is a groupoid action with trivial  moment map $M\to\bullet$. Just like
the case of group actions,  to any left groupoid action one can associate an \tsf{action groupoid}
$\poidd{\calG\acts N}{N}$ with total space  $\calG\fib{\sor}{j}N$, the source map $\sor_\acts$ is the
restriction of the second projection, the target $\tar_\acts$ is the action map $\sigma$ itself, and the
multiplication $\mu_\acts$ maps any composable pair $[(g,h*n);(h,n)]$,   $(g,h)\in\calG^{(2)}$ and
$h\in\sor\inverse(j(n))$, to 
$$
\qquad
(g,h*n)\cdot(h,n)=(gh,n)
\qquad,
$$ thus unit bisection $\eps_\acts$ and inversion
$\iota_\acts$  must be $\eps_\acts(n)=(\eps(j(n)),n)$ and $(g,n)^{-\acts}=(g\inverse,g*n)$.\\ 
In what
follows we shall consider mostly actions of Lie groupoids on smooth maps; note that, in this case, the
domain of the action map is always a smooth manifold and  such an action is called a \tsf{Lie groupoid
action} if the action map is also smooth. Note that the source fibres
$\sor_\acts\inverse(n)=\sor\inverse(j(n))\times\{n\}$, $n\in N$, of an action groupoid $\calG\acts N$
have the same homotopy type as suitable source fibres  of $\calG$. The proof of the following lemma is
straightforward.
\bgn{lemma} The action groupoid of a Lie groupoid action is a Lie groupoid.
\end{lemma}
\spa Next, we shall explain the denomination ``moment map''\fn{ The ``controversy'' on the denomination \emph{moment}
vs \emph{momentum} is well known. From now on we shall use mostly the term moment map.}.   Consider an
\tsf{Hamiltonian action} of a Lie group $G$ (which we assume to be connected here) on a Poisson manifold $(M,\pi)$,
i.e. a Lie group action for which there exists a \tsf{moment map} $j:M\to\frag^*$, characterizing the induced 
infinitesimal action as $\sigma(x)=X^{j^*F_x}$, $x\in\frag$. Such a map is equivariant for the coadjoint
representation on the nose if{f} it is Poisson or, equivalently, the diagram
$$
\bcat
\xy
*+{}="0",    <-1cm,0cm>
*+{\frag}="1", <1cm,0cm>
*+{\frax(M)}="2", <0cm,1.4cm>
*+{\cif(M)}="3", 
\ar     @ {->} "1";"2"_{\sigma} 
\ar     @ {->} "1";"3"^{j^*}  
\ar     @ {->} "3";"2"^{X^\bullet}
\endxy
\ecat
$$
commutes in the category of Lie algebras, where $X^\bullet$ is the map  that assigns to each function on $M$ the
corresponding Hamiltonian vector field. When $M$ is symplectic this is equivalent to the classical notion of an 
equivariant momentum map \cite{moment} for an Hamiltonian $G$-space (see, for example,  \cite{cds01}
for more details).  
\\
Lie groupoid actions and moment maps are indeed a natural generalization of Hamiltonian actions and the momentum
maps of classical mechanics. Consider the tangent prolongation groupoid $\poidd{T^*G}{\frag^*}$ defined by the
structure maps
$$
\sh(\theta_g)=l_g{}^*\theta_g\qquad\th(\theta_h)=r_h{}^*\theta_h\qquad
\muh(\theta_g,\theta_h)=r_{h\inverse}{}^*\theta_g=l_{g\inverse}\theta_h
$$
(unit section and inversion are defined accordingly)\footnote{For any Lie groupoid $\calG$ with Lie algebroid $A$, it
is possible to define a cotangent prolongation groupoid $\poidd{T^*\calG}{A^*}$; we shall discuss the
construction in section \ref{ipm}.}.
An Hamiltonian action of $G$ on $M$ (symplectic or Poisson) can be lifted to a groupoid action of the 
$\poidd{T^*G}{\frag^*}$ along $j$ by a trivial fibrewise extension:  set $\theta_g\,\hat{*}\,m:=g*m$, for all $m\in
M$ and $\theta_g\in T^*G$, with $r_{g\inverse}{}^*\theta_g=j(m)$. By equivariance,  
$$
\pair{j(\theta_g\,\hat{*}\,m)}{x}=\pair{j(m)}{\sf{Ad}_{g\inverse}x} =\pair{\sh(\theta_g)}{\dd r_{g}\comp\dd
l_{g\inverse} x}=\pair{\th(\theta_g)}{x} \:\:,\:\: x\in\frag\:\:, 
$$ 
the remaining property of a groupoid action are manifest. 
\bgn{example} \tsf{Groupoid actions and action groupoids}.
\medskip\\
$i$) For any Lie groupoid $\gpdm$ there is a left action on the target map $\tar:\calG\to M$ and a right action on
the source map $\sor:\calG\to M$, the action map being the groupoid multiplication in both cases. In the latter
case the corresponding action groupoid is easily identified with the groupoid $\calG\odot\calG$ of example \ref{exaa},
in the former, with the same groupoid associated with $\calG^{\sf{op}}$ (i.e. the Lie groupoid corresponding to the
opposite category).
\medskip\\
$ii$) The action of a Lie group $G$  on a manifold $M$ always tangent lifts to an action of the tangent group
$TG$ on $TM$ and $TG\acts TM$ is canonically isomorphic to $T(G\acts M)$.
\medskip\\
$iii$) \cite{mackbook} For any closed Lie subgroup $H\subset G$ the action groupoid $G\acts G/H$ is canonically
isomorphic to the gauge groupoid $(G\times G)/H$.
\end{example}

\bigskip

\spa Consider a general Lie groupoid action such as above and let $A$ be the Lie algebroid of $\calG$. The action
map induces a morphism of Lie algebras $\Gamma(A)\to\frax(N)$, $a\mapsto X^{a}$
$$
\qquad
X^{a}_{n}:=\dd\sigma(a_{j(n)}, 0_n)\qquad,\qquad n\in N
\qquad.
$$
Note that the vector bundle underlying the Lie algebroid of $\calG\acts N$ is the pullback  $j^{\sf +}A\to N$ and
the space of sections $\Gamma(j^{\sf +}A)\simeq\cif(N)\otimes_{\cif(M)}\Gamma(A)$ can be endowed with a Lie bracket
by picking decompositions 
$\Gamma(j^{\sf +}A)\ni\alpha_{1,2}=\sum_{k_{1,2}} u_{1,2}^{k_{1,2}}(a_{1,2}^{{k_{1,2}}}\comp\phi)$ and setting
\bgn{eqnarray}\label{brackacts}
\nn\brak{\alpha_1}{\alpha_2}
&=&	
\underset{k_1,k_2}\sum 
u_{1}^{k_{1}}\cdot u_{2}^{k_{2}}
(\brak{a_1^{k_{1}}}{a_2^{k_{2}}}\comp j)
\quad+\quad
\underset{k_2}\sum \dd\sigma\comp(\alpha_1\times 0)(u_2^{k_2})(a_2^{k_2}\comp j)\\
&-&
\underset{k_1}\sum \dd\sigma\comp(\alpha_2\times 0)(u_1^{k_1})(a_1^{k_1}\comp j)\qquad.
\end{eqnarray}
One can check that the above expression is well defined and makes $j^{\sf +}A$ a Lie algebroid for the obvious
anchor. Moreover such a Lie algebroid is canonically isomorphic to the Lie algebroid of $\calG\acts N$.\\
We shall remark that it is possible to make sense of actions of Lie algebroids without making reference to actions
of Lie groupoids.
\bgn{definition} Let $A\to M$ be a Lie algebroid and consider a map $j:N\to M$. A \tsf{Lie algebroid action} of $A$
on $j$ is given by a map $\sigma:\Gamma(A)\to\frax(N)$, such that
\medskip\\
1. $\sigma$ is a morphism $\cif(M)$-modules, for the natural module on $\frax(N)$ induced by $j$ via pullback;
\medskip\\
2. $\sigma$ is a morphism of Lie algebras;
\medskip\\
Note that, thanks to condition 1., $\sigma$ induces a well defined bundle map $j^{\pmb +}\sigma:\Gamma(j^{\pmb
+}A)\to\frax(N)$, which is required to be compatible with $j$ in the sense that   
\medskip\\
3. $\dd j\comp j^{\pmb +}\sigma=\rho\comp j^!$
\end{definition}
Using the data in the above definition in the same way as those produced by a Lie groupoid action one can obtain a
Lie algebroid $A\acts N$, we shall refer to as the \tsf{action Lie algebroid}, with anchor 
$\rho^\acts:=j^{\pmb{+}}\sigma$ on the pullback $j^{\pmb +}A$. Essentially by construction the natural map
$j^!:j^{\pmb{+}}A\to A$ is a morphism of Lie algebroid.
\bgn{example} \tsf{Action Lie algebroids}
\medskip\\
$i$) An infinitesimal action $\frag\to\frax(M)$ is a Lie algebroid action on $M\to\bullet$ and  the corresponding
action Lie algebroid is the flat bundle $\frag\times M$ with the usual action Lie algebroid structure.
\medskip\\
$ii$) The action of a Lie groupoid $\poidd{\calG}{M}$ on $\tar:\calG\to M$ by {\em left}  translation induces an
action Lie algebroid on $\tar^{\pmb +}A$; the Lie algebroid anchor $\tar^{\pmb +}A\to T\calG$ is easily computed as
$(a_{\tar(g)},g)\mapsto \dd r_g a_{\tar(g)}$ and it is easy to identify such Lie algebroid with the action Lie
algebroid $A\acts\calG$ for the infinitesimal action $\Gamma(A)\to\rinv{\frax}(\calG)$ by {\em right} invariant
vector fields.
\end{example}

\bigskip

\spa By definition, an equivalence relation $\thicksim$ on a smooth manifold $M$ is regular if    the quotient space
$M\mod$ carries a smooth topology (the \emph{quotient topology})  making the quotient projection  $M\to M\mod$ a
submersion; the \tsf{graph of an equivalence relation} is  $\gr{\thicksim}:=\{(x,y)\:|x\thicksim y\:\}\subset
M\times M$. We shall repeatedly  make use of  the classical
\bgn{theorem}[Godement's criterion]\label{gode} An equivalence relation $\thicksim$ on a manifold $M$ is  regular
if{f}  
\medskip\\
$1.$ $\gr{\thicksim}$ is a submanifold of $M\times M$;
\medskip\\
$2.$ The restriction of the first or, equivalently, of the second projection  $\gr{\thicksim}\to M$ is submersive.
\end{theorem}
As a consequence the natural conditions for the orbit space $M/G$ of a Lie group action to be smooth consists in
asking the groupoid anchor $G\acts M\to M\times M$ to be an embedding, namely the action by right or left translation
of $G\acts M$ on $M$ to be free and proper. This motivates the following
\bgn{definition} A Lie groupoid $\poidd{\calG}{M}$ is called \tsf{$($isotropy$)$ free}, respectively \tsf{proper},
if one of the following equivalent conditions is fulfilled:
\medskip\\
$i)$ The action of $\calG$ on $M$ by left
translation is free, respectively proper, i.e.
$$
\tar\times\id_M:\calG\fib{\sor}{}M\to M\times M
$$
is an injective immersion, respectively a proper map;
\medskip\\
$ii)$ The action of $\calG$ on $M$ by right translation is free, respectively proper, i.e. 
$$
\id_M\times\sor:M\fib{}{\tar}\calG\to M\times M
$$
is an injective immersion, respectively a proper map;
\medskip\\
$iii)$ The groupoid anchor $\chi:\calG\to M\times M$ is an injective immersion, respectively a proper map.
\end{definition}
Accordingly, a Lie groupoid action of shall be called free,  respectively proper, if so is the associated action
Lie groupoid. Note that the orbit space $N/\calG$  and $N/(\calG\acts N)$ coincide.
\bgn{lemma} For any free and proper Lie groupoid action of $\poidd{\calG}{M}$ on $j:N\to M$ the orbit space
$N/\calG$  carries a unique smooth topology making the quotient projection  $N\to N/\calG$ a submersion.
\end{lemma}
\bgn{proof} Condition 1. for of theorem \ref{gode} holds by hypothesis; to check condition 2. pick a bisection of
$\calG\acts N$.
\end{proof}
Since action Lie groupoids differentiate to action Lie algebroids, the connected components of  a  $\calG$-orbit
$\calO^\calG$ on $N$ are the leaves of the associated action Lie algebroid, thus
$$
\qquad
T_n\calO^\calG=\dd\sigma (A_{j(n)}\times 0_n)\qquad, \qquad n\in\calO^\calG
\qquad,
$$
and it is straightforward to see that the quotient $N/\calG$, when smooth, has dimension
\bgn{eqnarray}\label{dime}
\qquad
\dim N/\calG =
\dim N +\dim M -\dim\calG
\qquad.
\end{eqnarray}
\bgn{example} For any isotropy free and proper Lie groupoid $\poidd{\calG}{M}$, the anchor $\chi$ is a closed
embedding, and so is the tangent anchor $\dd\chi:T\calG\to TM\times TM$.  That is, the tangent prolongation
groupoid is also isotropy free and proper.
\end{example}

\bigskip

\spa Apart from Lie group actions, many other differential geometric constructions
extend to Lie groupoids; in the following we shall need principal bundles with
groupoid structure.
\bgn{definition} A \tsf{left principal $\calG$-bundle} is given by a fibre
bundle $\tau\colon P\to N$, a groupoid $\poidd{\calG}{M}$ and a groupoid action
on a moment map $j\colon P\to M$, such that\nolinebreak
\medskip\\1. $\calG$ acts along the fibres of $\tau$;
\medskip\\2. The action is fibrewise free and transitive.
\medskip\\
A \tsf{right principal $\calG$-bundle} is defined similarly, replacing the left
action with a right action.
\end{definition}
We shall also say that a principal $\calG$-bundle is a  principal bundle with structure groupoid $\calG$.
Naturally, one can adapt this set theoretic definition to the category of topological spaces and smooth manifolds,
in the case of a Lie groupoid action. In the latter case, last two conditions can be expressed by requiring that
the anchor of the associated action groupoid be an isomorphism of Lie groupoids $\calG\acts P\to
P\fib{\tau}{\tau}P$, for the subgroupoid of the pair groupoid $\poidd{P\fib{\tau}{\tau}P}{P}$ defined by the
equivalence relation induced by $\tau$.
\bgn{example} \tsf{Principal groupoid bundles}
\medskip\\
$i$) Ordinary principal bundles are principal bundles in the sense of the definition above.
\medskip\\
$ii$) For any Lie groupoid $\gpdm$, $\sor:\calG\to M$ is a right principal $\calG$-bundle for the action by right
translation.
\end{example}
There are natural conditions for a groupoid action on a Lie groupoid to be compatible with the partial
multiplication.
\bgn{definition}\label{compaction} Let $\gpdm$ and $\poidd{\calH}{N}$ be Lie groupoids. If $\calG$ acts on 
moment maps $J:\calH\to M$ and $j:N\to M$ in such a way that the conditions
\medskip\\
1. $\sor_{\calH}(g*h)=g*\sor_{\calH}(h)$
\medskip\\
2. $\tar_{\calH}(g*h)=g*\tar_{\calH}(h)$
\medskip\\
3. $\iota_{\calH}(g*h)=g*\iota_{\calH}(h)$
\medskip\\
4. $\mu_{\calH}(g*h_+,g*h_-)=g*\mu_\calH(h_+,h_-)$
\medskip\\
hold whenever they make sense, we shall say that the action is \tsf{compatible with the groupoid structure}.
\end{definition}
Quotients with respect to compatible groupoid actions are well behaved.
\bgn{lemma}\label{pgbred}\cite{mm02} Let a Lie groupoid $\poidd{\calG}{M}$ act on a Lie groupoid 
$\poidd{\calH}{N}$ compatibly with the groupoid structure; if $N$ is a principal $\calG$-bundle, the quotient
$\calH/\calG$ carries a Lie groupoid structure over the quotient $N/\calH$.
\end{lemma}
In other words, whenever the action of $\calG$ on the base of $\calH$ is free and proper, so that the quotient
$N/\calG$ exists, the quotient $\calH/\calG$ is also smooth and it is easy to see that the   the groupoid structure
descends, thanks to the compatibility conditions of definition \ref{compaction}.\\
An example of last lemma plays a central role in the next Section.
\bgn{example} \tsf{Monodromy groupoids}
\medskip\\
$i$) For any regular foliation $\calF$ (i.e. corresponding to an integrable distribution of constant rank) on a
manifold $M$, the monodromy groupoid  $\poidd{\sf{Mon}(M,\calF)}{M}$ is given by the classes of foliated paths
(i.e. paths within the leaves of $\calF$) with respect to foliated homotopy  relative to the endpoints, for the
multiplication induced by concatenation of paths with matching endpoints. The source fibre over one point is easily
identified with the universal cover of that leaf, thus $\poidd{\sf{Mon}(M,\calF)}{M}$ is always Lie groupoid with
1-connected source fibres.
\medskip\\
$ii$) Consider the monodromy groupoid $\sf{Mon}{(\calG,\sor)}$ for the foliation induced by the source map of a Lie
groupoid with connected source fibres. On the one hand $\sf{Mon}{(\calG,\sor)}$ is  naturally acted on from the
right by $\calG$. On the other hand $\calG$ acts on itself also by right translation and makes  $\tar:\calG\to
M\equiv\calG/\calG$ a principal $\calG$-bundle. The quotient groupoid  $\poidd{\sf{Mon}{(\calG,\sor)}/\calG}{M}$ is then a 
\emph{Lie} groupoid with 1-connected source fibres.
\end{example}
\vs{1}
\section{Integrability of Lie algebroids and Poisson manifolds}\label{ilapm}
\begin{quotation} 
We devote this Section to a review of the fundamental results in Lie theory for Lie algebroids and Lie groupoids and
their dual counterparts in Poisson geometry; differently from the case of classical Lie theory, Lie algebroids
(Poisson manifolds) are not always integrable to Lie groupoids (symplectic groupoids). In the first part we present
the ``optimal'' generalizations of Lie's theorems \cite{Lie} about the integrability of subobjects and morphisms;
moreover we describe the construction of the Weinstein groupoid, a topological groupoid ``integrating 
nonintegrable Lie algebroids''. In the second part we recall the main properties of symplectic groupoids and their relation with
Poisson manifolds; here we also introduce cotangent prolongation groupoids, the symplectic groupoids integrating
Poisson structures which are dual to integrable Lie algebroids.	  
\end{quotation}
\vs{0.5}
\subsection{Integrability of Lie algebroids}\hfill

\vs{0.1}
\spa Let $\gpdm$ be a Lie groupoid with connected source fibres\fn{As it is customary, we shall call, for
short, source $n$-connected a groupoid with $n$-connected source fibres. Similarly a (morphism of ) groupoids are
said ``source-$x$'' if the property of being ``x'' holds sourcewise.} and Lie algebroid $A$; the quotient
$\wt{\calG}:=\mon(\calG,\sor)/\calG$ of the monodromy groupoid is a  cover of $\calG$ in a sense we are
about to make precise.\\
Since $\sor$ is a submersion,
$\calG\fib{\sor}{\sor}\calG\subset\calG\times\calG$ is a smooth submanifold and a wide Lie subgroupoid for
the pair groupoid on $\calG\times\calG$ (see example \ref{exaa} ($iii$)). Being $\calG$ source connected, the groupoid anchor
$\chi^{\mon}:\mon(\calG,\sor)\to\calG\times\calG$, selecting the endpoints,  is a source surjective and
submersive morphism of Lie groupoids onto $\calG\fib{\sor}{\sor}\calG$; thus, there is a short exact
sequence
$$
\xymatrix{
\calL(\calG,\sor)\:\ar@{^(->}[r]&\:
\mon(\calG,\sor)\ar@{->>}[r]\:&\:
\calG\fib{\sor}{\sor}\calG
}
$$
of Lie groupoids over $\calG$, where $\calL(\calG,\sor)=\ker\chi^{\mon}$ is the  wide normal Lie
subgroupoid of homotopy classes of foliated loops in $(\calG,\sor)$.\\ 
Note that $\calL(\calG,\sor)$ is stable
under the natural action of $\calG$  on $\mon(\calG,\sor)$ by right translation and
$\calG\fib{\sor}{\sor}\calG$ is acted on by $\calG$ from the right diagonally. The sequence above is
equivariant for these compatible Lie groupoid actions and yields a new short exact sequence
$$
\xymatrix{
\calL(\calG,\sor)/\calG\:\ar@{^(->}[r]&\:
\widetilde{\calG}\ar@{->>}[r]^{\kappa_\calG}\:&\:
\calG
}
$$
of Lie groupoids over $M$ by taking quotients. Computing the total space of the Lie algebroids of 
$\wt{\calG}$ yields the same vector bundle as $A$ and the morphism of the Lie algebroids induced by
$\kappa_\calG$ is the identity map; moreover, $\wt{\calG}$ is indeed the ``covering groupoid'' of $\calG$, in
the sense that all restrictions $\wt{\sor}\inverse(q)\to \sor\inverse(q)$, $q\in M$, of $\kappa_\calG$ 
to the source fibres are covering maps. Since the source connected component of any Lie groupoid $\calG$ 
is a Lie groupoid with the same Lie algebroid, we have  
\bgn{theorem}\label{0}{\em (0th Lie theorem)}\cite{mm02} Let $\poidd{\calG}{M}$ be a Lie groupoid with Lie algebroid
$A$.  Then, there exists a source 1-connected  Lie groupoid $\wt{\calG}$ over $M$, with Lie algebroid $\wt{A}$
and a morphism of Lie groupoids  $\kappa_\calG$, inducing an isomorphism $\wt{A}\to A$.
\end{theorem}
As it shall be clear after theorem \ref{2} and remark \ref{popeye} below, the covering groupoid is unique up to isomorphism.

\spa A Lie algebroid $A$ is called \tsf{integrable} if it is in the image of the Lie functor; in this case, we shall say 
that the covering groupoid of any groupoid differentiating to $A$, is \emph{the} (source 1-connected) Lie groupoid of
$A$.\\

\spa Let us consider now the integrability of Lie subalgebroids. The action Lie groupoid $\calG\acts\calG$ for the
action by left translation is isomorphic to $\calG\fib{\sor}{\sor}\calG$ for the mapping $(g,h)\mapsto (gh,h)$, thus
the Lie algebroid of $\calG\fib{\sor}{\sor}\calG$, with total space $T^\sor\calG$, is isomorphic to $A\acts\calG$ for
the bundle isomorphism $\calR':T^\sor\calG\to A\acts\calG$ given by right translation  and the right translation map
$\calR:=\tar^!\comp\calR':T^\sor\calG\to A$ a morphism of Lie algebroids over $\tar$, in fact it is the differential
of the division map  $\delta:\calG\fib{\sor}{\sor}\calG\to\calG$. Since $\calR$ is fibrewise a diffeomorphism any
subalgebroid $B\to N$ of $A$ can be pulled back to a vector subbundle $\calR\inverse(B)\to\tar\inverse(N)$ of 
$T^\sor\calG$ and the graph of  $\dd\tar|_{\calR\inverse(B)}:\calR\inverse(B)\to TN$ defines a \emph{regular}
distribution on $N\fib{}{\tar}\calG$. It is easy to identify the latter with the distribution spanned by the image of
the infinitesimal action 
$$
\sigma_B:\Gamma(B)\to\frax(N\fib{}{\tar}\calG)
\qquad,\qquad 
\sigma_B(b)_{(n,g)}=(\rho_n(b),\rinv{b}_g)
\qquad;
$$
$\dd \tar|_{\calR\inverse(B)}$ is thus an integrable distribution, since it coincides with the anchor distribution of
the action Lie algebroid $B\acts (N\times_{\tar}\calG)$. Moreover, since $\dd \tar|_{\calR\inverse(B)}$ has constant
rank, the associated foliation $\calF_B$ is regular and one can form the monodromy groupoid
$\sf{Mon}(N\times_{\tar}\calG,\calF_B)$; the latter is naturally acted on by  $\calG$ from the right, as is the
projection $N\times_{\tar}\calG\to N$, which is a principal $\calG$-bundle. It follows from lemma \ref{pgbred} that
$\sf{Mon}(N\times_{\tar}\calG,\calF_B)/\calG$ is a  source 1-connected Lie groupoid over $N$, whose Lie algebroid
turns out to be isomorphic to $B$.
\bgn{theorem}\label{1}{\em (1st Lie's theorem)}\cite{mm02} Let $\poidd{\calG}{M}$ be a Lie groupoid with Lie
algebroid $A$ and $B\to N$ a Lie subalgebroid of $A$. Then $B$ is integrable.
\end{theorem}
Note that the integration of $B$ in general is not obtained as a Lie subgroupoid of $\calG$; we shall clarify this
point after discussing the integrability of morphism of Lie algebroids (see remark \ref{popeye}).

\spa Consider a morphism of integrable Lie algebroids $\phi:A_-\to A_+$ over $f:M_-\to M_+$ and let 
$\calG_\pm$ be Lie groupoids inducing $A_\pm$, respectively. The vector fields on $\calG_-\times\calG_+$ 
of the form
$$
X_{(g_-,g_+)}=(\dd r_{g_-} a_{m_-}, \dd r_{g_+}\phi(a_{m_-}))
\qquad,\qquad \tar_-(g_-)=m_- \quad,\quad a\in\Gamma(A_-)
$$
span an integrable distribution $\Delta$ on $(\sor_-\times\sor_+)\inverse\gr{f}\subset\calG_-\times\calG_+$. 
Denote with $\calL_{m_-}$ the leaf of $\Delta$ through $(\eps_-(m_-),\eps_+(f(m_-))$, $m_-\in M_-$. 
One can check that the projection $\calL_{m_-}\to\sor\inverse_-(m_-)$ is a covering map, hence a diffeomorphism, if
$\calG_-$ is source 1-connected, therefore $\Delta$ induces a family of smooth maps 
$\phi_{m_-}:\sor\inverse_-(m_-)\to \sor\inverse_+(f(m_-))$. One can further check that the maps of this family
paste together nicely and in fact define a morphism of Lie groupoids $\fii:\calG_-\to\calG_+$ over 
$f:M_-\to M_+$, which, by construction, induces $\phi:A_-\to A_+$ (see \cite{mx00} for details).
\bgn{theorem}\label{2}{\em (2nd Lie's theorem)}\cite{mx00} Let $\poidd{\calG_\pm}{M_\pm}$ be Lie 
groupoids with Lie algebroids $A_\pm$ and $\phi\colon A_-\!\to\! A_+$ a morphism of Lie algebroids. If $\calG_-$ is
source 1-con-\linebreak nected, there exists a unique morphism of Lie groupoids $\calG_-\!\to\!\calG_+\!$ inducing
$\phi$.\fn{All the maps here are of class $\cif$. The proof can be adapted to morphisms of finite regularity.}
\end{theorem}
Uniqueness follows from the fact that the vector fields spanning the distribution are right invariant in
$\calG_-\times\calG_+$, thus $\Delta$, hence $\fii$, is determined by its span at $M_-\times M_+$, which in turn is
prescribed by $\phi$.
\bgn{remark}\label{popeye} 
For any Lie subalgebroid $A_-$ of an integrable Lie algebroid $A_+$, the inclusion $A_-\inc A_+$ can 
be integrated to a morphism of Lie groupoids $\calG_-\to\calG_+$, provided $\calG_-$ is source 1-connected, and
one can check by a right translation argument that this map is immersive; however it might fail to be
injective. This is not the case when $A_-\simeq A_+$, since the integration $\calG_-\to\calG_+$ is then a
bijective immersion; thus, the covering groupoid of theorem \ref{0} is unique up to isomorphism.
\end{remark}

\spa To any Lie algebroid $A$ one can associate a topological groupoid, the Weinstein groupoid 
$\calW(A)$, which fails to be smooth in general. 
If $A$ has an integrating groupoid $\calG$, thanks to theorem \ref{0} $\calG$ can be assumed to be 
source 1-connected and, in this case, $\calG\simeq\calW(A)$. We recall below the construction \cite{crfs03}
of the Weinstein groupoid.\\  
The source 1-connected cover $\wt{\calG}=\mon(\calG,\sor)/\calG$ can be equivalently
described as the quotient of the space  $P\calG$, of the so called (twice differentiable)
\tsf{$\calG$-paths}, namely source foliated class $\calC^2$ paths in $\calG$ starting at the unit section, by
\tsf{$\calG$-homotopy}, i.e. $\calC^2$ foliated homotopy relative to the endpoints
restricted to  $P\calG$. This is clear, since for each class in 
$\mon(\calG,\sor)/\calG$ there is a unique representative in $P\calG$. There
exists a natural map $P\calG\to PA$ taking values in the space of
class $\calC^1$ morphisms of Lie algebroids $TI\to A$, given by the right derivative 
$$
g\mapsto\delta_rg\qquad,\qquad \delta_rg(u)=\dd
r\inverse_{g(u)}\,\dot{g}(u)\quad,\quad 
u\in I\qquad.
$$
Here we identify paths $\alpha:I\to A$ with 1-forms $\alpha\,	\dd u:TI\to A$
taking values in the pullback bundle $\gamma^{\sf{+}}\!A\to I$, 
$\gamma:=\pr_A\comp \alpha$;
such paths are morphism of Lie algebroids if{f} the anchor compatibility condition
$\dd \gamma=\rho_\gamma\comp\alpha$ holds, for the tangent map 
$\dd \gamma\in\Ohm^1(I,\gamma^\sf{+}\!A)$, being the bracket compatibility a trivial
condition. The right derivative actually takes values in
$PA$, since, for any $\calG$-path $g$, $\delta_rg$ is a
vector bundle map $TI\rightarrow A$ with base map $\gamma=\tar\comp g$ and
$$
\rho_{\gamma(u)}\delta_rg(u)=\dd\tar^{\quad}_{g(u)}
\dd r\inverse_{g(u)}\,\dot{g}(u)=\left.\frac{\dd}{\dd v}\right|_{v=u}{\tar(g(v))}
$$
for all $u\in I$.\\
Since any  $\calC^1$ morphism of Lie algebroids $\alpha:TI\to A$ integrates to a
unique
$\calC^2$ morphism of Lie groupoids $a:I\times I\rightarrow\calG$, necessarily 
of the
form 
$$
a(u,v)=a(u,0)\cdot a(0,v)=g(u)\cdot g(v)\inverse
$$ 
for some $\calG$-path 
$g$,
and applying the right derivative to such a $\calG$-path, returns the
original $A$-path $\alpha$, there is a bijective correspondence between the
space of $A$-paths and that of $\calG$-paths.\\
The notion of $\calG$-homotopy translates to the following equivalence relation on 
the space of $A$-paths, the
so called \tsf{$A$-homotopy}: two $A$-paths $\alpha_\pm\in PA$ are 
$A$-homotopic if there exists a class $\calC^1$ morphism of Lie algebroids
$h:T\twice{I}\to A$, such that the boundary conditions 
$$
\ppmv h=\alpha_\pm\qquad\qquad\hbox{ and }\qquad\qquad\ppmh h=0
$$
hold, where $h\in\Ohm^1(\twice{I},X^+A)$ is regarded as a 1-form taking 
values in the pullback bundle for the base map $X$. According to remark \ref{torsolo}, the requirement to be a 
Lie algebroid
morphism can be expressed by the anchor compatibility condition
$\dd X=\rho_X\comp h$ and the Maurer-Cartan equation
$$\qquad
\mathrm{D}\,h+\frac{1}{2}[h
\,\overset{\wedge},\,
h]=0\qquad.
$$
Here 
$
\mathrm{D}:\Ohm^\bullet(\twice{I},X^{\sf{+}}A)\to 
\Ohm^{\bullet +1}(\twice{I},X^{\sf{+}}A)
$ is the covariant derivative of the pullback of an arbitrary linear connection 
$\nabla:\frax(M)\otimes\Gamma(A)\to\Gamma(A)$ for the vector bundle 
$A\rightarrow M$ and 
$
[\theta
\,\overset{\wedge},\,
\theta']=-X^{\sf{+}}\tau^\nabla(\theta\wedge\theta')
$ is the contraction with the pullback of the torsion tensor
$
\tau^\nabla\in\Gamma(A^*\otimes\wedge^2A)
$
. 
\bgn{proposition}\label{la-homotopy}\cite{crfs03} Let $\alpha_\pm\in PA$ be two $A$-paths and $g_\pm\in P\calG$ the
corresponding $\calG$-path. Then $\alpha_\pm$ are $A$-homotopic if{f} $g_\pm$ are $\calG$-homotopic. 
\end{proposition}
We give a considerably simpler proof of this fact, than that appeared in \cite{crfs03}.
\bgn{proof} A map $H:I\times I\to\calG$ is a $\calG$-homotopy from $g_-:=H\comp\pmv$ to $g_+:=H\comp\ppv$ if{f}
its tangent map takes values in $T^\sor\calG$ and 
$$
\qquad
\ppmh\dd H=0
\qquad.
$$
For the Lie algebroid induced by $\calG\fib{\sor}{\sor}\calG$, $T^\sor\calG\subset T\calG$ is a Lie subalgebroid
of $T\calG$, thus $\calR\comp\dd H:TI^{\times 2}\to A$ is a morphism of Lie algebroids, such that
\bgn{equation}\label{boundarui}
\qquad
\iota_\pa^*h=\calR\comp\dd (H\comp\iota_\pa)=\calR\comp\iota_\pa^*\dd H
\qquad,
\end{equation}
for any boundary component $\pa$, hence an $A$-homotopy from $\delta_r g_-$ to 
$\delta_r g_+$. Conversely, the integration $\wt{h}$ of any $A$-homotopy $h:TI\times I\to A$ from 
$\delta_r g_-$ to  $\delta_r g_+$ induces a map $H:I\times I\to\calG$, $H(u,v)=\wt{h}(u,v;0,0)$, whose image is
contained in a source fibre and a morphism of Lie algebroids $\dd H:TI^{\times 2}\to T^\sor\calG$, 
such that $\calR\comp\dd H=h$; $\dd H$ satisfies
the desired boundary condition because of (\ref{boundarui}), since $\calR$ is fibrewise a diffeomorphism.
\end{proof}
\spa Proposition (\ref{la-homotopy}) allows characterizing the source 1-connected  cover of a source connected Lie
groupoid $\calG$, in terms of the infinitesimal  data only, as the quotient 
$$
\wt{\calG}=PA\mod 
$$ 
of $A$-paths by $A$-homotopy. Since $A$-homotopy classes contain smooth representatives and reparameterizing 
such does not change
the class, the composition of $\calG$-paths in $P\calG\mod=\sf{Mon}(\calG,\sor)/\calG$ (induced, roughly, by
right translation and concatenation) can be translated to a groupoid multiplication on $PA\mod$, also given by
concatenation of good representatives (see \cite{crfs03} for details).  Most importantly, even when $A$ is not
integrable, $PA\mod$ is always a topological groupoid, the \tsf{Weinstein groupoid}
$\calW(A)$.
\spa The integrability conditions for a Lie algebroid, which we state below for completeness, are then to be 
understood as a characterization of the obstruction to the existence of a smooth structure on the 
associated Weinstein groupoid and depend on certain monodromy groups $\calN_\bullet(A)$: for all $q\in M$, 
$\calN_q(A)$ is the abelian subgroup for the isotropy algebra $\frag_q$ defined by the elements $\zeta$ in the
center of the isotropy algebra $\frag_q$, such that the constant
$A$-path $\zeta(u)\equiv\zeta$ is $A$-homotopic to the zero $A$-path $\zeta_0(u)\equiv 0_{q}$.
\bgn{theorem}\label{3}{\em (3rd Lie's theorem)} \cite{crfs03} A Lie algebroid $A$ is integrable to a Lie groupoid 
if{f} the monodromy groups $\calN_\bullet(A)$ are
\medskip\\
$1$. Discrete, i.e.  ${\sf d}_\calN(q))>0$ for all $q\in M$,
\medskip\\
$2$. Locally uniformly discrete, in the sense that 
$$
\qquad
{\underset{p{\underset{p\ni\calL_q}\longrightarrow} q }{\sf{liminf}}}\,{\sf d}_\calN(p)>0
\qquad,
$$ 
where ${\sf d}_{\calN}(q)=D(\calN_q,0)$ for the distance function $D$ of an arbitrary
norm on $A$.
\end{theorem}
See \cite{crfs03,crfs04,cafe01} and references therein for the discussion of nonintegrable examples.

\vs{0.5}
\subsection{Integrability of Poisson manifolds and coisotropic submanifolds}\label{ipm}\hfill

\vs{0.1}
\spa Symplectic groupoids are Lie groupoids with a compatible symplectic form and provide a desingularization, namely
a \emph{symplectic realization} \cite{cdw87},  of the symplectic foliation of certain Poisson manifolds
\bgn{definition}\label{sympgpd}\cite{ks86,ws87,zr90} Let $\gpdm$ be a Lie groupoid endowed with a symplectic form; we shall say
that the symplectic form is compatible with the groupoid structure and $\poidd{(\calG,\ohm)}{M}$ is  a \tsf{symplectic
groupoid} if the graph of the groupoid multiplication is Lagrangian in $T^*\calG\times
T^*\calG\times\ol{T^*\calG}$.
\end{definition}
The compatibility between the groupoid structure and the symplectic form of a symplectic groupoid can be
expressed equivalently in terms of the multiplicativity condition 
$$
\mu^*\ohm=\pr_1^*\ohm +\pr_2^*\ohm
$$
on the manifold  $\calG^{(2)}$ of composable elements, explicitly
\bgn{equation}\label{multisymp}
\ohm_{gh}(\delta g_+\bullet \delta h_+, \delta g_-\bullet \delta h_-)
=
\ohm_{g}(\delta g_+, \delta g_-)
+
\ohm_{h}( \delta h_+, \delta h_-)
\end{equation}
for the tangent multiplication and composable $\delta g_\pm\in T_g\calG$, $\delta h_\pm\in T_h\calG$. By
repeated application of (\ref{multisymp}) on suitable tangent vectors and by counting dimensions, one can show
\bgn{theorem}\label{sympgpdp}\cite{cdw87} Let $\gpdm$ be a symplectic groupoid. The following hold
\medskip\\
$1$. The unit section $M\to \calG$ is a Lagrangian embedding;
\medskip\\
$2$. The inversion $\calG\to\calG$ is an anti-symplectic
diffeomorphism;
\medskip\\
$3$. Source and target fibres are symplectic orthogonal to each other;
\medskip\\
$4$. The source invariant functions and the target invariant
functions form mutually commuting Poisson subalgebras of $\cif(M)$;
\medskip\\
as a consequence,
\medskip\\
$5$. There exists a unique Poisson structure on $M$ making the source
map Poisson and the target map anti-Poisson.
\end{theorem}
The statements above follow by counting dimensions from the analogous 
statements regarding Poisson groupoids (replace symplectic with Poisson and
Lagrangian with coisotropic), the proof of which shall
be outlined in the next Section. 
\eject
Dual Poisson structures of integrable Lie algebroids are always induced by symplectic
groupoids. The following results were obtained by Coste, Dazord and Weinstein.
\bgn{proposition}\cite{cdw87} Let $\gpdm$ be a Lie groupoid. Then
\medskip\\$i$) The cotangent bundle $T^*\calG$ carries a Lie groupoid structure
over the conormal bundle $N^*M$ making it a symplectic groupoid for the
canonical symplectic form;
\medskip\\ and, identifying $N^*M$ with the dual bundle $A^*$ to the Lie
algebroid $A$ of $\calG$,
\medskip\\$ii$) The symplectic groupoid $\poidd{T^*\calG}{A^*}$ induces the
dual Poisson structure on $A^*$ for the anticanonical symplectic form\fn{We need to change the sign of the
symplectic form since we use ``Poisson-friendly''
conventions: the Poisson bracket induced by the canonical symplectic has therefore opposite sign with
respect to the usual Poisson bracket of symplectic geometry.}.
\end{proposition}
The Lie groupoid of the above proposition is usually referred to as  the \tsf{cotangent prolongation groupoid}
and plays a central role throughout this dissertation.  Let us define the structure maps of the groupoid structure and
set some notations. Fix  $A=T^\sor_M\calG\simeq NM$ as a choice of a normal bundle. Source $\sh$  and target
$\th$ are the only natural maps $T^*\calG\to N^*M$ yielding vector bundle maps over $\sor$ and $\tar$; they are
defined by setting, for all $a\in A_{\sor(g)}$, respectively $b\in A_{\tar(g)}$,
\bgn{eqnarray*}
\pair{\sh(\theta_g)}{a}_{\sor(g)}&:=&\pair{\theta_g}{\dd l_g(a-\dd\eps(\rho(a)))}_g
\qquad,\qquad\:\hbox{respectively}\qquad,\\
\pair{\th(\theta_g)}{b}_{\tar(g)}&:=&\pair{\theta_g}{\dd r_g\, b}_g
\:\:\:\qquad\quad\qquad\qquad,\qquad
\:\:\theta_g\in T^*_g\calG\quad\qquad,
\end{eqnarray*}
for the natural pairings. The cotangent multiplication
$\muh:T^*\calG\fib{\sh}{\th}T^*\calG\to T^*\calG$ can be defined as a
vector bundle map over the groupoid multiplication
$\mu:\calG\fib{\sor}{\tar}\calG\to\calG$, by setting 
$$
\gr{\muh}:=(\id_{T^*\calG}\times
\id_{T^*\calG}\times-\id_{T^*\calG})N^*\Gamma(\mu)
$$
and checking that $(\id_{T^*\calG}\times \id_{T^*\calG}\times-\id_{T^*\calG})N^*\Gamma(\mu)$ is 
indeed the graph of a map. It follows that the cotangent unit section $\epsh:A^*\to
T^*\calG$ and inversion $\iotah:T^*\calG\to T^*\calG$ have to be the 
identification of $A^*$ with $N^*M$ and, respectively the antitranspose
$-\dd\iota^{\sf{t}}$ of the tangent inversion. Note that the cotangent
multiplication $\theta_g\hdot\theta_h$ of composable cotangent vectors
$(\theta_g,\theta_h)\in T^*\calG\fib{\sh}{\th}T^*\calG$
is completely defined by the formula
\bgn{equation*}\label{ctgliftcdw}
\langle\,\theta_g\hdot\theta_h,\delta g\bullet\delta h\,\rangle
=
\langle\,\theta_g\,,\,\delta g\,\rangle
+
\langle\,\theta_h\,,\,\delta h\,\rangle\qquad,\qquad (\delta
g,\delta h)\in T\calG\fib{\dd\sor}{\dd\tar}T\calG\qquad,
\end{equation*}
by means of which the groupoid compatibility conditions can easily  be checked 
(see the proof of theorem \ref{memento} for analogous computations); 
moreover the graph of $\muh$ is Lagrangian in $T^*\calG\times T^*\calG \times \ol{T^*\calG}$
(and $ \ol{T^*\calG}\times\ol{T^*\calG}\times{T^*\calG}$) 
by construction.\\
It is not hard to inspect when $T^*\calG$ is \emph{the} symplectic groupoid of $A^*$.
\bgn{remark}\label{homotopytype} The source fibres of the cotangent prolongation groupoid have the same homotopy type as
those of $\calG$. To see this, note that the kernel of  the cotangent source map is 
the annihilator $(T^{\tar}\calG)^o$, thus $\sh$ has maximal rank as a bundle map and the short exact
sequence
$$
\xymatrix{
0\:\ar@{->}[r]&\:
\ker\sh\:\ar@{->}[r]&\:
T^*\calG\ar@{->}[r]\:&\:
\sor^{\pmb{+}} A^*\ar@{->}[r]\:&\:
0\\
}
$$
of vector bundles over $\calG$ admits a splitting $\sigma:\sor^{\pmb{+}} A^*\to T^*\calG$, which allows
to describe each $\sh$-fibre as 
$$
\qquad
\sh\inverse(\xi)=\{\sigma_g(\xi)+\Ker_g\sh\}_{g\in\sh\inverse(\pr(\xi))}\qquad,\qquad \xi\in
A^*\qquad,
$$ 
i.e. as an affine subbundle of $T^*_{\sor\inverse(\pr(\xi))}\calG$. As a consequence all $\sh$-fibres homotopy
retract to suitable $\sor$-fibres and, in particular $T^*\calG$  is source 1-connected if{f} so is $\calG$.
\end{remark}
\bgn{example} \tsf{Cotangent prolongation groupoids}
\medskip\\
$i$) Regard a manifold $M$ as the trivial Lie groupoid; 
then the induced Lie algebroid is also trivial and the dual Poisson structure is
the zero bivector field over $M$. The construction above produces the abelian
Lie groupoid on $T^*M\to M$.
\medskip\\
$ii$) Consider the cotangent prolongation of the pair groupoid $M\times M$: a direct computation
shows that source and target are respectively $-\pr_2,\pr_1:T^*M\times T^*M\to T^*M$ and
\be
\pair{(\theta_x,-\theta_y)\hdot(\theta_y,\theta_z)}{(\delta x,\delta z)}
&=&
\pair{(\theta_x,-\theta_y)}{(\delta x,\delta y)}
+
\pair{(\theta_y,\theta_z)}{(\delta y,\delta z)}\\
&=&
\pair{(\theta_x,\theta_z)}{(\delta x,\delta z)}\hs{4.5}.
\ee
\end{example} 

\bigskip

\spa Let $\poidd{(\Lambda,\ohm)}{P}$ be a symplectic groupoid; denote with $A$ the Lie algebroid on
$T^\sor_P\Lambda$ induced by the bracket of right invariant vector fields. 
It is well known that $A$ is canonically isomorphic to the Koszul algebroid on $T^*P$. The
isomorphism is provided by the sharp map of the symplectic form: thanks to theorem \ref{sympgpdp} (3.), setting
$$
\qquad
\pair{\phi\inverse_\ohm(a)}{\delta p}:=-\ohm(a, \dd\eps \delta p)
\qquad,\qquad (a,\delta p)\in A\oplus TP\qquad 
\qquad,
$$
yields an isomorphism of vector bundles $\phi_\ohm: T^*P\to A$ over the identity map; the induced map on the
subalgebra of exact 1-forms $\Ohm^1_{\sf{exact}}(P)\to\Gamma(A)$ can be regarded as taking values in the 
algebra of right invariant vector fields. Define 
$\phi^{r}_\ohm(\dd f):=-X^{\tar^*f}$, thus
$$
\ohm(\phi^{r}_\ohm(\dd f)(\eps(p)),\dd\eps\dd p)
=
-\pair{(\tar^*\dd f)_{\eps{p}}}{\dd\eps \delta p}
=
-\pair{\dd f}{\delta p}=:\ohm({\phi_\ohm(\dd f_p)},{\dd\eps\delta p})
$$ 
for all $\delta p\in TP$; by symplectic orthogonality of source and target fibres, nondegeneracy of $\ohm$
implies that $\phi^{r}_\ohm(\dd f)$ coincides with the right invariant vector field 
$\rinv{\phi_\ohm\comp\dd f}$. Using this fact and the Leibniz rule
one can easily check that $\phi_\ohm$ preserves the canonical Lie brackets on $T^*P$ and $A$.\eject
It follows that the connected components of the orbits of a symplectic groupoid are the symplectic leaves of the
induced Poisson structure; in this sense, a symplectic groupoid is to be regarded as a desingularization of the
 symplectic foliation on its base manifold. The discussion above motivates the following
\bgn{definition} A Poisson bivector is called integrable if so is the induced Koszul algebroid. 
\end{definition}
%
%
As a corollary of integrability of Lie bialgebroids, to be discussed in the next Section, one can show that 
the source 1-connected Lie groupoid of an integrable Poisson manifold is always a symplectic groupoid.
\bgn{theorem}\label{inty}\cite{mx00} A Poisson manifold $(P,\pi)$ is integrable if{f} there exists a 
symplectic groupoid inducing $\pi$.
\end{theorem}
The correspondence between integrable Poisson manifolds and source 1-connected symplectic groupoids is however
not functorial. In fact, given source 1-connected symplectic groupoids $\Lambda_\pm$ over $P_\pm$
a Poisson map $f:P_-\to P_+$ does not induce a symplectomorphism 
$\Lambda_-\to\Lambda_+$, unless $f$ is a diffeomorphism. 
Coisotropic submanifolds are instead well behaved with respect to integration; roughly speaking, 
there exists a correspondence between coisotropic submanifolds and Lagrangian subgroupoids.
Rephrasing a result by Cattaneo \cite{cat04}, 
\bgn{theorem}\label{lagsub} Let $\poidd{\Lambda}{P}$ be a source 1-connected symplectic groupoid and
$\poidd{\calL}{C}$ a source 1-connected Lie subgroupoid with Lie algebroid $L$. Then $\calL\subset\calG$ is
Lagrangian if{f} $C\subset P$ is coisotropic and $L\simeq N^*C$.
\end{theorem}
The proof of last theorem given in \cite{cat04} is highly non trivial and involves symplectic reduction in a suitable
infinite dimensional weak symplectic manifold; an independent simple proof shall be obtained in the next Section 
We shall comment further on theorems \ref{inty}-\ref{lagsub}, after introducing Poisson groupoids.
\vs{0.75}
\section{Poisson groupoids and Lie bialgebroids}\label{pglb}
\begin{quotation} 
The first part of this Section consists in a quick introduction to Lie groupoids with a compatible Poisson
bivector and their infinitesimal invariant, namely Poisson groupoids and Lie bialgebroids. In the second part after
reviewing Mackenzie and Xu's integration of a Lie bialgebroid to a Poisson groupoid (Lie's third theorem holds
under the integrability conditions on the underlying Lie algebroid), we further develop Lie theory for Poisson
groupoids. In particular, we prove (theorem \ref{bipo}) the integrability of morphisms of Lie bialgebroids to morphisms of Poisson
groupoids via a version of Lie's first theorem (integrability of coisotropic subalgebroids, theorem
\ref{coisosub}).
\end{quotation}

\spa Poisson groupoids are simultaneous generalizations of Poisson groups and symplectic groupoids;  Weinstein's
coisotropic calculus provides a powerful technique to unify the theory of these two objects. 
\bgn{definition}\cite{ws88} A Poisson structure on a Lie groupoid $\gpdm$ is called compatible if
the graph of the groupoid multiplication
$$
\gr{\mu}\subset{\calG\times \calG\times\ol{\calG}}
$$
is coisotropic. A Lie groupoid endowed with a compatible Poisson structure is called a
\tsf{Poisson groupoid}.
\end{definition}
A Poisson groupoid $\poidd{G}{\bullet}$ over the one point manifold is a Poisson group in the
classical sense,  i.e. the compatibility condition is equivalent to asking the group multiplication
to be a Poisson map $G\times G\to G$; a Poisson groupoid with nondegenerate (hence symplectic)
Poisson bivector is a symplectic groupoid in the sense of definition \ref{sympgpd}, by counting dimensions.
\bgn{example} \tsf{Poisson groupoids}
\medskip\\$i$) Every Lie groupoid is a Poisson groupoid for the zero Poisson structure. A Poisson
manifold $P$ is \emph{not} a Poisson groupoid for the trivial Lie groupoid $\poidd{P}{P}$, since 
$\gr{\mu}=\{(p,p,p)\}$ is not coisotropic in $P\times P\times\ol{P}$.
The pair groupoid on  $\ol{P}\times P$ is a Poisson groupoid.
\medskip\\$ii$) Let $P$ be a \tsf{Poisson $G$-space}, i.e $G$ is a Poisson group acting on a Poisson manifold $P$ in
such a way that the action map $G\times P\to P$ is Poisson. If the action is free and proper $P\to P/G$ is a
principal $G$ bundle; it is an exercise in coisotropic calculus to check that the diagonal action of $G$ on
$\ol{P}\times P$ is Poisson and compatible with the pair groupoid. Then $(P\times\ol{P})/G$ carries a Lie groupoid
over $P/G$, namely the gauge groupoid. It is well known \cite{ws88} that smooth quotients of Poisson $G$ spaces always 
carry a unique Poisson structure making the quotient projection a Poisson submersion. For the quotient Poisson structure
$\poidd{(\ol{P}\times P)/G}{P/G}$ is a Poisson groupoid. 
\end{example}
The main properties of a Poisson groupoid were derived in \cite{ws88} and are listed in the following
theorem. The proof we sketch uses essentially the same techniques as in \cite{ws88}.
\bgn{theorem} Let $\pgpd$ be a Poisson groupoid. The following hold
\medskip\\
$1$. The unit section $M\to \calG$ is a coisotropic embedding;
\medskip\\
$2$. The inversion $\calG\to\calG$ is an anti-Poisson
diffeomorphism;
\medskip\\
$3$. The source invariant functions and the target invariant
functions form Poisson subalgebras of $\cif(M)$;
\medskip\\
$4$. There exists a unique Poisson structure on $M$ making the source
map Poisson and the target map anti-Poisson;
\medskip\\
$5$. The Hamiltonian vector fields of the source, respectively target, invariant functions 
are left, respectively right, invariant;  
\medskip\\
$6$. The subalgebras of source and target invariant functions commute.
\end{theorem}
\bgn{proof}[Sketch of proof] The first two properties follow easily from theorem \ref{coisorel}, the very cleanliness
assumptions can be shown to hold by direct inspection (1) $M\subset\calG$ is coisotropic since
$$
M\times\bullet=\gr{\mu}\comp(\Delta_\calG\times\bullet)\subset\calG\times\bullet
$$
is the composition of the coisotropic relations $\gr{\mu}:\calG\times\calG\to\calG$ and
$\Delta_\calG\times\bullet:\ol{\calG}\times\calG\to\bullet$. (2) $\gr{\iota}\subset\calG\times\calG$ is
coisotropic since
$$
\gr{\iota}\times\bullet=\gr{\mu}\comp(M\times\bullet)\subset\calG\times\calG\times\bullet
$$
is the composition of coisotropic the relations $\gr{\mu}:\calG\to\ol{\calG}\times\calG$ and
$M\times\bullet: \calG\to\bullet$. The proof of (3) is independent from that of (1) and (4) is a
straightforward consequence of (3). (3) By definition
$f\in\cif(\calG)$ is source invariant if $f=\sor^*\eps^*f$, equivalently if $f(gh)=f(h)$ for all composable
$g$, $h\in\calG$; that is, $f$ is source invariant if{f} the associated function 
$\wt{f}\in\cif(\threetimes{\calG)}$, $\wt{f}(g,h,k):=f(k)-f(h)$ vanishes on the characteristic ideal of
$\gr{\mu}$. Note that, for all $f_\pm\in\cif(\calG)$,
$$
\qquad
\wt{\poib{f_+}{f_-}_\calG}=\poib{\wt{f_+}}{\wt{f_-}}_{\calG\times\calG\times\ol{\calG}}
\qquad,
$$
thus source-invariant functions form a Poisson subalgebra, by coisotropicity of $\gr{\mu}$, and it
follows from (2) that also target-invariant functions do. (5) For any source invariant $f\in\cif(\calG)$ and
composable $g$, $h\in\calG$,
$$
(\Pi^\sharp\times\Pi^\sharp\times-\Pi^\sharp)_{(g,h,gh)}\dd\wt{f}=(0_g,-\Pi^\sharp_h\dd f,
-\Pi^\sharp_{gh}\dd f)\in T_{(g,h,gh)}\gr{\mu}=\gr{\dd\mu}_{(g,h,gh)}
\quad,
$$
thus $\dd\tar\Pi^\sharp_h\dd f=0_{\tar(h)}$ and  $\Pi^\sharp_{gh}\dd f=0_g\bullet\Pi^\sharp_h\dd f=\dd
l_g\Pi^\sharp_h\dd f$ for the tangent multiplication $\bullet$; the statement for right invariant vector fields is
proved analogously. (6) is a follows easily from (5).
\end{proof}

\spa Let $\poidd{(\calG,\Pi)}{M}$ be a Poisson groupoid and consider its Lie algebroid 
$A\equiv T^\sor\calG$. The dual bundle $A^*$ is to be canonically identified with the conormal 
bundle $N^*M\subset T^*\calG$ and, by coisotropicity of $M\subset\calG$, it carries a natural Lie
algebroid structure. The Lie algebroids on $A$ and $A^*$ are compatible in a way that makes sense without any
reference to the Poisson groupoid they are derived from.
\bgn{definition}\cite{mx94,ksb95} Let $A$ be a Lie algebroid and assume $A^*$ is also a Lie algebroid. The pair
$(A,A^*)$ is called a \tsf{Lie bialgebroid} if the Lie algebroid differential $\dd_{A^*}$ makes
$(\Gamma(\wedge^\bullet A), \wedge, \brasn{}{},\dd_{A^*}))$ a differential Gerstenhaber algebra\fn{Lie
bialgebroids were discovered by Mackenzie and Xu in \cite{mx94} as infinitesimal invariants of  Poisson groupoids. 
The definition \cite{ksb95} we use is equivalent to that of \cite{mx94}.}.
\end{definition}  
Le us explain the above definition. For any Lie algebroid $A\to M $, the Lie bracket on $\Gamma(A)$ can be
uniquely extended to a biderivation $\brasn{}{}$ of the graded commutative algebra on $\Gamma(\wedge^\bullet
A)=\oplus_{i\geq 0}\Gamma(\wedge^i A)$ by setting
\be\label{defu}
\brasn{f}{g}&:=&0\\
\brasn{a}{f}&:=&\rho(a)(f)\\
\brasn{a}{b}&:=&\brak{a}{b}\hs{3}f,g\in\cif(M)\:\:,\:\:a,b\in\Gamma(A)\quad,
\ee
and imposing the Leibniz rule\fn{We sometimes denote the degree of an homogeneous $a\in\Gamma(\wedge^k A)$
with the same symbol $a$
.}
\bgn{equation}\label{Leibrule}
\brasn{a}{b\wedge c}=\brasn{b}{a}\wedge c +(-)^{(a-1)b}b\wedge\brasn{a}{c}=0
\end{equation}
on all homogeneous $a,b,c\in\Gamma(\wedge^\bullet A)$. Remarkably $\brasn{}{}$ makes the suspension
$\Gamma(\wedge^\bullet A)[1]$ a graded Lie algebra, i.e.  $\brasn{}{}$ is graded skewsymmetric,
$$
\qquad
\brasn{a}{b}+(-)^{(a-1)(b-1)}\brasn{b}{a}=0
\qquad
$$
and the graded Jacobi identity 
$$
\brak{a}{\brak{b}{c}}=\brak{\brak{a}{b}}{c}+(-)^{(a-1)(b-1)}\brak{a}{\brak{b}{c}}  
$$
holds for all homogeneous elements $a,b,c\in\Gamma(\wedge^\bullet A)$; namely, 
$(\Gamma(\wedge^\bullet A), \wedge, \brasn{}{}))$ is a Gerstenhaber algebra. Dually, the Lie algebroid on $A$
induces a differential $\dd_{A}:\Gamma(\wedge^\bullet A^*)\to\Gamma(\wedge^{\bullet +1}A^*)$,
\be
\dd_A\eta(a_0,\dots, a_n)&:=&\sum_{i=0}^{n}(-)^i\rho(a_i)(\xi(a_1,\dots, \hat{a_i},\dots, a_n))\\
&\,+&
\sum_{i<j}^{}(-)^{i+j}\xi([a_i,a_j],a_1,\dots,a_n)\qquad,\qquad\eta\in\Gamma(\wedge^n A),
\ee
which generalizes the Chevalley-Eilenberg differential and the de Rham differential. Moreover it satisfies the 
usual rules of Cartan calculus; in particular $\dd_A$ is a derivation of the wedge product.
\bgn{remark}
If $(A,A^*)$ is a pair of Lie algebroids in duality the natural compatibility amounts therefore to asking the 
Lie algebroid differential $\dd_{A^*}$ induced by $A^*$ to be a derivation of the graded Lie bracket induced by $A$, i.e.
\bgn{equation}\label{lba}
\dd_{A^*}\brasn{a}{b}=\brasn{\dd_{A^*}a}{b}+(-)^{(a-1)}\brasn{a}{\dd_{A^*} b}
\end{equation}
for all homogeneous elements $a,b\in\Gamma(\wedge^\bullet A)$. The role played by $A$ and $A^*$ is symmetric, 
in the sense that $(A,A^*)$ is a Lie bialgebroid if{f} so is $(A,A^*)$.
\end{remark}
\bgn{example} \tsf{Lie bialgebroids}
\medskip\\
$i$) For any Lie algebra $\frag$, $\dd_\frag$ is the classical Chevalley-Eilenberg differential. If 
$(\frag,\frag^*)$ are Lie algebras denote with $\delta:\frag\to\frag\otimes\frag$ the map dual to the Lie
bracket on $\frag^*$; condition (\ref{lba}) is satisfied if{f} $\delta$ is a 1-cocycle in the Lie algebra
cohomology with coefficients in $\frag\otimes\frag$ for the adjoint representation, that is, if 
$(\frag,\frag^*)$ is a Lie bialgebra (see \cite{chpr} for details).
\medskip\\
$ii$)  If $(P,\pi)$ is a Poisson manifold the Jacobi identity for the Poisson bracket is equivalent
to $\brasn{\pi}{\pi}=0$ for the graded Lie bracket on $\frax^\bullet(P):=\Gamma(\wedge^\bullet TP)$ extending 
the Lie bracket of vector fields. By the graded Jacobi identity $\dd_\pi:=\brasn{\pi}{\cdot\,}$ makes 
$\frax^\bullet(P)$ a differential Gerstenhaber algebra and $(T^*P,TP)$ a Lie bialgebroid.\nolinebreak
\medskip\\
$iii$) More generally, for any Lie algebroid $A$ with a bisection $\Delta\in\Gamma(\wedge^2 A)$ such
that $\brasn{\Delta}{\Delta}=0$, the Lie bracket 
$$
\brak{\xi_+}{\xi_-}
:=
\dd_A\Delta(\xi_+,\xi_-) +\iota_{\Delta^\sharp \xi_+}\dd_A\xi_--
\iota_{\Delta^\sharp \xi_-}\dd_A\xi_+
\quad,\quad \xi_\pm\in\Gamma(A^*)
$$
makes $(A,A^*)$ a Lie bialgebroid, where the map $\Delta^\sharp:A^*\to A$ is 	induced by $\Delta$. 
The bisection $\Delta$ is called an $R$-matrix and $(A,A^*)$ a \tsf{triangular Lie bialgebroid}. If $A=\frag$ 
is a Lie algebra the condition on $\Delta\in\wedge^2\frag$ is the classical Yang Baxter equation. 
\medskip\\$iv$) \cite{ksb96} Let $(P,\pi)$ be a Poisson manifold and $\Phi:TP\to TP$ a (1,1) tensor with vanishing 
torsion
$$
\tau^\Phi(X,Y)=\brak{\Phi X}{\Phi Y}-\Phi\brak{\Phi X}{Y} -\Phi\brak{X}{\Phi Y} +\Phi^2\brak{X}{Y}\quad,\quad
X,Y\in\frax(P)
$$
Then $[X,Y]^\Phi=\brak{\Phi X}{Y} +\brak{X}{\Phi Y} -\brak{X}{Y}$ makes $TP$ a Lie algebroid $TP_\Phi$ with 
anchor $\Phi$ and $(T^*P, TP_\Phi)$ is a Lie bialgebroid if{f} 
$(P,\pi,N)$ is a \tsf{Poisson-Nijenhuis} manifold \cite{bv88,km90,bm94,ksb96}.  
\end{example}

\spa The main integrability result for Lie bialgebroids is the following theorem due to Mackenzie and Xu. We shall sketch 
below the idea of the proof.
\bgn{theorem}\label{mxiba}\cite{mx00} Let $(A,A^*)$ be a Lie bialgebroid over $M$ and $\gpdm$ a source 1-connected 
Lie groupoid with Lie algebroid $A$. Then there exists a unique Poisson structure on $\calG$ 
making it a Poisson groupoid and inducing the Lie algebroid $A^*$ on $N^*M$. 
\end{theorem}
First of all we shall need to introduce a certain canonical
antisymplectomorphism. 
\spa Note that, for any vector bundle $\pr: A\to M$, setting $\sf{r}_A:T^*A\to A^*$
$$
\qquad
{\sf{r}_A}(\theta)(a'):=\pair{\theta}{\iota\vup_a(a')}\qquad,\qquad a,a'\in A
\qquad,
$$
where $\iota\vup:\pr^{\sf +}A\simeq T\vup A\to TA$ is the inclusion of the vertical bundle, yields a
vector bundle map over $\pr:A\to M$. One can check that $\dd\pr: TA\to TM$ carries a vector bundle, whose
scalar multiplication and fibrewise addition are given by the tangent
of the corresponding maps of $A$ (see \cite{mkz92,mkz00b,mackbook} for the details); moreover the structure map of 
$\dd\pr$ can be suitably dualized to endow also ${\sf{r}_{A^*}}\colon T^*A^*\to A$ with a vector bundle structure. The diagram
\bgn{equation}\label{md}
\bcat
\xy
*+{}="0",    <-0.7cm,0.7cm>
*+{T^*A^*}="1", <0.7cm,0.7cm>
*+{A}="2", <-0.7cm,-0.7cm>
*+{A^*}="3", <0.7cm,-0.7cm>
*+{M}="4",
\ar  @ {->} "1";"2"^{\quad\sf{r}_{A^*}} 
\ar  @ {->} "1";"3"_{}
\ar  @ {->} "2";"4"^{}
\ar  @ {->} "3";"4"_{}
\endxy
\ecat
\end{equation}
is a double vector bundle in the sense of Ehresmann, i.e. a vector bundle in the category of vector
bundles. 
Remarkably, it is possible to define an isomorphism of double vector bundles
\bgn{equation}\label{mdm}
\bcat\xy
*+{}="0",    <-0.9cm,0.7cm>
*+{T^* A^*}="1", <0.7cm,0.7cm>
*+{A}="2", <-0.9cm,-0.7cm>
*+{A^*}="3", <0.7cm,-0.7cm>
*+{M}="4",     <2.6cm,-0.7cm>
*+{T^*A}="1'", <4.2cm,-0.7cm>
*+{A}="2'", <2.6cm,-2.1cm>
*+{A^*}="3'", <4.2cm,-2.1cm>
*+{M}="4'",  <5.2cm,-0.7cm>
*+{,}="2''"
\ar  @{->} "1";"2"^{} 
\ar  		   @{->} "1";"3"^{}
\ar                @{->} "2";"4"_{}
\ar     @{->} "3";"4"_{}
\ar    @{->} "1'";"2'"^{} 
\ar  		   @{->} "1'";"3'"^{}
\ar                @{->} "2'";"4'"_{}
\ar   @{->} "3'";"4'"^{}
\ar  		   @{->} "1";"1'"^{\tau}
\ar  		   @{->} "2";"2'"^{}
\ar  		   @{->} "3";"3'"^{}
\ar  		   @{->} "4";"4'"^{}
\endxy\ecat
\end{equation}
i.e. both $(\tau,\id_A)$ and $(\tau, \id_{A^*})$ form a morphism of vector bundles, which is also an antisymplectomorphism for
the canonical symplectic forms.  The definition of $\tau$ appeared in \cite{mx94} and generalizes Tulczyjew's canonical
antisymplectomorphism $T^*T^*M\to T^*TM$, encoding the Legendre transform geometrically; let us recall the construction of
$\tau$, since it will be needed in the following. Regard the canonical
pairing of $A^*$ with $A$ as a map  $\FF:A^*\oplus A\to\RR$; the submanifold  $L\subset T^*(A^*\times A)$ defined by those
elements  $(\ol{\theta},\theta)\in T^*_{A^*\oplus A}(A^*\times A)$ such that
\bgn{equation}\label{tmxrc}
\pair{(\ol{\theta},\theta)}{(\delta\xi,\delta a)}=\dd\FF_A(\delta\xi,\delta a)
\qquad,\qquad \delta\xi,\delta a\in T(A^*\oplus A)
\end{equation}
can easily be seen to be Lagrangian (choose Darboux coordinates adapted to $A^*\oplus A$). One can
further show that $L$ is the graph of a map, therefore (\ref{tmxrc}) completely defines an
antisymplectomorphism.\\
Whenever $(A,A^*)$ is a Lie algebroid over $M$, the composition  $\pi^\sharp_{A^*}\comp\tau:T^*A^*\to TA$ of the dual
Poisson anchor with the canonical antisymplectomorphism is a morphism of Lie algebroids
over the anchor $A^*\to TM$ for the Koszul and tangent prolongation Lie algebroids. In fact, this requirement on a
pair of Lie algebroids  $(A,A^*)$ is equivalent to $(A,A^*)$ being a Lie bialgebroid. If $A$ is integrable and
$\calG$ its source 1-connected Lie groupoid, $\pi_{A^*}^\sharp\comp\tau$ can then  be integrated to a morphism of  Lie
groupoids $\Pi^\sharp:T^*\calG\to T\calG$ for the tangent and  cotangent prolongation groupoids.  It turns out that
$\Pi^\sharp$, as the notation suggests, is indeed the sharp map of a Poisson bivector, which,  by construction
indeed, is compatible with the groupoid multiplication and
%
%
%
%
%
%
induces the given Lie algebroid on $A^*$.\eject
\bgn{example} \tsf{Integrability of Lie bialgebroids}
\medskip\\$i$) If $P$ is a Poisson manifold with connected components $\{P_i\}$ and 
$\widetilde{P_i}\to P_i$ are the associated covering principal bundles, the Poisson groupoid on 
$$
\Pi(P)\simeq{\underset{i\in\pi_{o}(P)}\coprod}(\ol{\wt{P_i}}\times{\wt{P_i}})/\pi_1(P_i)
$$
integrates the Lie bialgebroid on $(TP,T^*P)$ (the Poisson bivector on $\widetilde{P_i}$ is the
unique making the covering projection a Poisson submersion, see also lemma \ref{wer} below).
\medskip\\$ii$) \cite{lx96} For any Lie groupoid $\gpdm$ with Lie algebroid $A$ and $R$-matrix $\Delta$,
$$
\Pi=\rinv{\Delta}-\linv{\Delta}=\rinv{\Delta -\iota_*\Delta}
$$
is a Poisson bivector making $(\calG,\Pi)$ a Poisson groupoid and integrating the triangular Lie bialgebroid
associated with $(A,\Delta)$.
\end{example}

\spa The notion of morphism of Poisson groupoids is very natural; namely, for any Poisson groupoids
$\calG_\pm$, $\fii:\calG_-\to\calG_+$ is a \tsf{morphism of Poisson groupoids} if it is a Poisson map and a morphism of Lie
groupoids, equivalently if $\gr{\fii}\subset\calG_-\times\ol{\calG_+}$ is a coisotropic subgroupoid. Note that, if $(A_\pm,A^*_\pm)$ are Lie bialgebroids and $\phi:A_-\to A_+$ is a morphism of Lie
algebroids, there is no natural map $\Gamma(\wedge^\bullet A_-)\to\Gamma(\wedge^\bullet A_+)$. Nevertheless, by
transposition, there a map $\phi^{\wedge *}:\Gamma(\wedge^\bullet A_+^*)\to\Gamma(\wedge^\bullet A_-^*)$ and to
encode a compatibility condition for $\phi$ with the Lie algebroids on $A_\pm^*$ it is natural to ask 
$\phi^{\wedge *}$ to preserve the corresponding graded Lie brackets. Clearly,
$$
\phi^{\wedge *}(\xi\wedge\eta)=\phi^{\wedge *}\xi\wedge\phi^{\wedge *}\eta\qquad,\qquad
\xi,\eta\in\Gamma(\wedge^\bullet A_+^*)\qquad,
$$ 
thus, it is sufficient to check last condition on homogeneous elements of degree 0 and 1, equivalently  to
show that $\phi^*:\cif(A_+)\to\cif(A_-)$  preserves the Poisson brackets of pairs of functions which
are fibrewise constant or fibrewise linear for the dual Poisson structures. It makes therefore sense to call 
$\phi$ a \tsf{morphism of Lie bialgebroids} if it is a Poisson map and a morphism of Lie algebroids,
equivalently if $\gr{\phi}\subset A_-\times A_+$ is a coisotropic Lie subalgebroid; this definition clearly
yields a well defined category of Lie bialgebroids.\\
We devote the rest of this Section to prove that the categories of integrable Lie bialgebroids and source 1-connected Lie 
groupoids are isomorphic; the result is a corollary of the following
\bgn{theorem}\label{coisosub} Let $\poidd{(\calG,\Pi)}{M}$ be a source 1-connected Poisson  groupoid
with Lie bialgebroid $(A,A^*)$ and $\poidd{\calC}{N}$ a source 1-connected Lie subgroupoid with Lie
algebroid $C$. Then $\calC\subset\calG$ is coisotropic if{f} so is $C\subset A$ for the dual Poisson
structure induced by $A^*$.  
\end{theorem}
In other words, under the connectivity assumptions, coisotropic subalgebroids of integrable
Lie bialgebroids, integrate to coisotropic subgroupoids. Note, however, that the source
1-connected integration of a coisotropic subalgebroid of an integrable Lie bialgebroid $A$ is
in general only an immersed subgroupoid of the  Poisson groupoid integrating $A$. 
The only if{f} part of theorem \ref{coisosub} was proved in \cite{xu95} without any connectivity assumption, 
using a different approach; it is clear from our proof below that we also need no connectivity 
assumptions to obtain this implication.\\
Consider that for any morphism of Lie groupoids $\fii:\calG_-\to\calG_+$, the graph $\gr{\fii}$ is source 1-connected 
if{f} so is $\calG_-$; therefore, if $\calG_-\times\calG_+$ is also source 1-connected, we can apply theorem \ref{coisosub} and
conclude that, morphisms of integrable Lie bialgebroids, by the characterization in terms of graphs, integrate to morphisms of
Poisson groupoids.\\
The connectivity assumptions on the target groupoid may be removed:
\bgn{theorem}\label{bipo} Let $\poidd{\calG_\pm}{M_\pm}$ be Poisson groupoids and $(A_\pm,A^*_\pm)$ their Lie bialgebroids. 
If $\fii:\calG_-\to\calG_+$ is a morphism of Lie groupoids and $\calG_-$ is source $1$-connected, 
then $\fii$ is Poisson if{f} so is the induced morphism of Lie algebroids  $\phi:A_-\to A_+$ for the dual Poisson
bivector fields induced by $A^*_\pm$.
\end{theorem}
Let us start with the proof of theorem \ref{coisosub}, which we divide in three steps; we emphasize that the
 first two steps do not involve the Poisson bivector of $\calG$ or the Lie algebroid on $A^*$. 
\medskip\\
\emph{Step 1.} (Proposition \ref{dfg}) We characterize the conormal bundle of a coisotropic subgroupoid
$\calC\subset\calG$ as a Lagrangian subgroupoid of the cotangent prolongation $\poidd{T^*\calG}{A^*}$ and identify
its Lie algebroid with the conormal bundle $N^*C^o$ to the annihilator of the Lie algebroid of $\calC$;
\medskip\\
\emph{Step 2.} (Lemma \ref{wer}) We identify $N^*C^o$ with $N^*C$ using the canonical antisymplectomorphism;
\medskip\\
\emph{Step 3.} (End of proof) The result follows by diagram chasing.
\medskip\\
Note that for any Lie subgroupoid $\calC\subset\calG$ the inclusion of Lie algebroids $C\inc A$ induces an 
inclusion of Lie algebroids $N^*C^o\inc T^*A^*$. Next result depends on the underlying Lie groupoids--algebroids only.
\bgn{proposition}\label{dfg} Let $\poidd{\calG}{M}$ be a Lie groupoid with Lie algebroid $A$ and 
$\calC\subset\calG$ a Lie  subgroupoid over $N\subset M$ with Lie algebroid $C$. Then 
\medskip\\$i$) $N^*\calC$ is a Lagrangian subgroupoid of  $\poidd{T^*\calG}{A^*}$ with base $C^o$;
\medskip\\$ii$) The Lie algebroid of $\poidd{N^*\calC}{C^o}$ is $N^*C^o\to C^o$.
\end{proposition}\eject
\bgn{proof} For all $\nu\in N^*_c\calC$ and $x\in C_{\sor(c)}$,
$$
\qquad
\pair{\sh(\nu)}{x}=\pair{\nu}{\dd l_c(x-\rho(x))}=0
\qquad,
$$
since $x-\rho(x)$ is tangent to $\calC$ and so is its left translate by an element of $\calC$.
Denote with $\wt{\sor}$ the restriction of $\sh$ to $N^*\calG$. We have
$
\ker\wt{s}=(T^\tar\calG)^o\cap N^*\calC= (T^\tar\calG + T\calC)^o
$,
therefore $\wt{\sor}$ is a vector bundle map of rank
\be
\rank\,\wt{\sor}
&=&
(\dim\calG-\dim\calC)-\dim\calG +(\rank\tar +\dim N)\\
&=&
\rank A -\rank C\\
&=& 
\rank C^o\hs{3}\qquad,
\ee 
that is, $\wt{\sor}:N^*\calC\to C^o$ is fibrewise surjective over a surjective submersion, hence a
surjective submersion (see remark \ref{bundlemaps1}). ($i$) It remains to check that $\muh$ restricts to
$N^*\calC$, which is clear since the tangent multiplication restricts to a fibrewise surjective
multiplication on $T\calC$ and because of (\ref{ctgliftcdw}). ($ii$) Let $\ohm$ be the canonical
symplectic form on $T^*\calG$: we have to check that the isomorphism
$$
\phi_{-\ohm}: T^*A^*\to T^{\sh}_{A^*}T^*\calG\qquad,\qquad \theta\to
\ohm_{\epsh(\theta)}{}^\sharp{}\inverse\th^*\theta 
\qquad,
$$
restricts to give a map $N^*C^o\to T^{\wt{\sor}}N^*\calC$. For all $f\in\calI_{C^o}$ and 
$F\in\calI_{N^*\calC}$,
$$
\qquad
\rinv{\phi_{-\ohm}\comp\dd f}(F)=\poib{F}{\th^*f}=0
\qquad,
$$
since $\th^*f\in\calI_{N^*\calC}$ and being $N^*\calC\subset T^*\calG$ Lagrangian, that is, 
$\rinv{\phi_{-\ohm}\comp\dd f}|_{N^*\calC}\in\rinv{\frax}(N^*\calC)$.
\end{proof}
Next lemma does not depend on Lie algebroids or Poisson structures. 
\bgn{lemma}\label{wer} For any vector bundle $A\to M$ and smooth subbundle $C\to N$, the canonical
antisymplectomorphism $T^*A^*\to T^*A$ identifies $N^*C^o$ with $N^*C$.
\end{lemma}
\bgn{proof} For all $\nu\in N^*_{\ul{c}}C^o$ and $\delta c\in T_c C$ with 
$(\ul{c},c)\in C^o\oplus C$, we have
\bgn{equation}\label{qkx}
\pair{\tau(\nu_c)}{\delta c}=\dd \FF_A(\delta\xi,\delta c)-\pair{\nu_c}{\delta\xi}
\end{equation}
if $\dd\pr_{A^*}\delta\xi=\dd\pr_A\delta c$, $\delta\xi\in T_{\ol{c}}C^o$ . On the other hand, 
by picking a connection for $C^o$, we can always find one such $\delta\xi$; for such a choice 
 the right hand side of (\ref{qkx}) vanishes, since $\FF_A$ is constant on $C^o\oplus C$. 
\end{proof}
It is now easy to conclude the proof.
\bgn{proof}[Proof of theorem \ref{coisosub}] Assume $C\subset A$ is a coisotropic subalgebroid.
Because of lemma \ref{wer} and by coisotropicity of $C\subset A$, the diagram
$$
\bcat
\xy
*+{}="0",    <-0.9cm,0.7cm>
*+{T^*A^*}="1", <0.9cm,0.7cm>
*+{TA}="2", <-0.9cm,-0.5cm>
*+{N^*C^o}="3", <0.9cm,-0.5cm>
*+{TC}="4",
\ar  @ {->} "1";"2"^{\pi_{A^*}\comp\tau} 
\ar  @ {^(->} "3";"1"_{}
\ar  @ {^(->} "4";"2"^{}
\ar  @ {->} "3";"4"_{}
\endxy
\ecat
$$
commutes in the category of Lie algebroids. Moreover, 
the source fibres of both $\poidd{T\calC}{TM}$ and $\poidd{N^*\calC}{C^o}$  have the same homotopy type as those of 
$\calC$, this can be seen using a variation of the argument in remark \ref{homotopytype}.
By the connectivity assumptions the diagram above integrates to
$$
\qquad
\bcat
\xy
*+{}="0",    <-0.9cm,0.7cm>
*+{T^*\calG}="1", <0.9cm,0.7cm>
*+{T\calG}="2", <-0.9cm,-0.5cm>
*+{N^*\calC}="3", <0.9cm,-0.5cm>
*+{T\calC}="4",
\ar  @ {->} "1";"2"^{\Pi^\sharp} 
\ar  @ {^(->} "3";"1"_{}
\ar  @ {^(->} "4";"2"^{}
\ar  @ {->} "3";"4"_{}
\endxy
\ecat
\qquad,
$$
i.e. $\calC\subset\calG$ is coisotropic. Conversely, if $\calC\subset\calG$ is coisotropic invert the argument
and use the inverse antisymplectomorphism to show that $\pi_{A^*}^\sharp N^*C\subset TC$.
\end{proof}

\spa If $\poidd{(\Lambda, \ohm)}{P}$ is a symplectic groupoid,  one can check that the isomorphism of Lie algebroids
$\phi_\ohm:T^*P\to T^\sor_P\calG$ is also a Poisson map for the canonical symplectic form on $T^*P$ and the fibrewise
linear Poisson structure induced by $N^*P$ on $T^\sor_P\calG\simeq NP$ (which is then symplectic). That is, $(T^*P,
TP)$ and $(T^\sor_P\calG,N^*P)$ are isomorphic Lie bialgebroids and theorem \ref{inty} follows specializing the proof
of \ref{mxiba}.  
If $\calL\subset\calG$ is a Lagrangian subgroupoid, the Lie algebroid $L\subset T^\sor_P\calG$ of $\calL$, by direct
inspection, maps to $N^*C$ under  $\phi_\ohm\inverse$, thus making $N^*C\subset T^*P$ a Lagrangian subalgebroid and
$C$ coisotropic. Conversely if $C$ is coisotropic with $\phi_\ohm (N^*C)= L$, then  $L$ is Lagrangian and theorem
\ref{coisosub} implies that so is $\calL$, by counting dimensions; that is, we have recovered theorem
\ref{lagsub}.

\spa Before proving theorem \ref{bipo}, we shall need the following easy 
\bgn{lemma}\label{tech} Let $(P,\pi)$ be a Poisson manifold, $C\subset P$ a coisotropic submanifold and $\phi:M\to P$ a local
diffeomorphism. Then
\medskip\\
$1$. There exists a unique Poisson bivector $\phi^\star\pi$ on $M$ making $\phi$ a Poisson map;
\medskip\\
$2$. If $X\subset M$ is a submanifold, such that $\phi$ restricts to a local diffeomorphism $X\to C$, then $X$ is coisotropic
with respect to $\phi^\star\pi$.
\end{lemma}
\bgn{proof} It is easy to see and well known that the equation 
$$
\phi^\star\pi(\phi^*\dd f_+, \phi^*\dd f_-)=\poib{f_+}{f_-}\comp\phi
\qquad,\qquad
f_\pm\in\cif(P)
$$
defines the desired Poisson bivector. (2) By linear algebra we have
$
N^*_xX=\dd\phi_x{}^{\sf{t}}(\dd\phi_x T_xX)^o=\dd\phi_x{}^{\sf{t}}N^*\!\!{}_{\phi(x)}C 
$ for all $x\in X$, therefore
\be
(\phi^\star\pi)^\sharp N^*_xX
&=&
(\dd\phi_x{}^{\sf{t}})\inverse\pi^\sharp\!\!{}_{\phi(x)} N^*\!\!{}_{\phi(x)}C\\
&\subset&
(\dd\phi_x{}^{\sf{t}})\inverse T{}_{\phi(x)}C\\
&=&T_xX\qquad\qquad\qquad\qquad.
\ee
by coisotropicity of $C$.
\end{proof}
\bgn{proof}[Proof of theorem \ref{bipo}] If the target groupoid in the statement of theorem \ref{bipo} is actually
not source 1-connected, the integrating morphism $\fii:\calG_-\to\calG_+$ can always be factored through the covering
morphism $\kappa:\tilde{\calG}_+\to\calG_+$ and the integration $\widetilde{\fii}:\calG_-\to\tilde{\calG}_+$.
Let $\Pi_+$ be the Poisson bivector on $\calG_+$ and  $\phi:\tilde{A}_+\equiv T^{\tilde{\sor}_+}_M\tilde{\calG}_+\to
T^{\sor_+}_M\calG_+\equiv A_+ $ the  isomorphism of Lie algebroids induced by the covering morphism. On the one hand,
$\phi$ allows to induce a Lie bialgebroid on 
$
(\widetilde{A}_+, \widetilde{A}_{+}^*\equiv N^*M_{\sf{rel\:}\tilde{\calG}_+})
$
from that on $(A_+,A_+^*\equiv N^*M_{\sf{rel\:}\calG})$ via $\phi^{\sf{t}}$  ($N^*X_{\sf{rel\:}Y}$ denotes the
conormal bundle of $X$ as a submanifold of $Y$); denote with $\tilde{\Pi}$ the Poisson bivector on $\calG$ making it
a Poisson groupoid and inducing $(\widetilde{A}_+, \widetilde{A}_+^*)$.  Since
both $\calG_-$ and $\tilde{\calG}_+$ are source 1-connected, $\widetilde{\fii}$ is Poisson and it suffices to show
that the covering morphism is also Poisson.
On the other hand, because of lemma
\ref{tech},  $\kappa^\star\Pi_+$ is also a compatible Poisson bivector on $\tilde{\calG}_+$, since  $\kappa^{\times
3}$ restricts to a local diffeomorphism $\gr{\tilde{\mu}_+}\to\gr{\mu_+}$.  We claim that the Lie bialgebroid of
$(\widetilde{\calG}_+,\kappa^\star\Pi_+)$ coincides  (not only up to isomorphism!) with 
$(\widetilde{A}_+, \widetilde{A}_{+}^*)$; it follows
by uniqueness (theorem \ref{inty}) that $\tilde{\Pi}$ coincides with $\kappa^\star\Pi_+$ and $\kappa$ is therefore
Poisson.  We have to show that  $\phi^{\sf{t}}:N^*M_{\sf{rel\:}\calG}\to N^*M_{\sf{rel\:}\tilde{\calG}_+}$ is a
morphism of Lie algebroids for the bracket on $N^*M_{\sf{rel\:}\tilde{\calG}_+}$ induced by $\kappa^\star\Pi_+$. The
anchor compatibility holds since
$$
\kappa^\star\Pi_+^\sharp\phi^{\sf{t}}
=
\left.\kappa^\star\Pi_+^\sharp\dd\kappa^{\sf{t}}\right|_{N^*M_{\sf{rel\:}\calG}}
=
\left.(\dd\kappa^{\sf{t}})\inverse \Pi^\sharp\right|_{N^*M_{\sf{rel\:}\calG}}
=
\left.\Pi^\sharp\right|_{N^*M_{\sf{rel\:}\calG}}\qquad,
$$
where we have used coisotropicity of $M\subset(\calG,\Pi)$ for the last equality. For sections of the form
$\left.\dd F_{1,2}\right|_M\in\Gamma(A^*_+)$, $F_{1,2}\in\calI_{M}\subset\cif(\calG_+)$, we have
\be
\phi^{\sf{t}}\brak{\left.\dd F_1\right|_M}{\left.\dd F_{2}\right|_M}_{A^*_+}
&=&
\dd\kappa^{\sf{t}}\left.\dd\poib{F_1}{F_2}_{\tilde{\Pi}}\right|_M
=
\left.\dd\poib{\kappa^*F_1}{\kappa^*F_2}_{\kappa^\star\Pi_+}\right|_M\\
&=&
\brak
{\phi^{\sf{t}}\left.\dd F_1\right|_M}
{\phi^{\sf{t}}\left.\dd F_{2}\right|_M}_{\widetilde{A}_+^*}
\ee
where the bracket in the second line is that induced by $\kappa^\star\Pi_+$.
Since $\Gamma(A^*_+)$ are locally finitely generated over $\cif(M)$ by differentials of functions in the vanishing
ideal, it follows from the Leibniz rule(s) that $\phi^{\sf{t}}$ preserves the bracket of arbitrary sections.
\end{proof}

\chapter{Double structures in Lie theory and Poisson geometry}\label{chapii}
\spa Ehresmann's categorification of a groupoid is a groupoid object in the category of groupoids; this is a symmetric
notion and it makes sense to regard such a structure as a ``double groupoid''.  A double Lie groupoid is, essentially, a ``Lie
groupoid in the category of Lie groupoids''; one can apply the Lie functor to either object of a double Lie groupoid, to obtain
an \tla-groupoid, i.e. a  ``Lie groupoid in the category of Lie algebroids''. The application of the Lie functor can still be
iterated; the result, a double Lie algebroid, is the best approximation to what one would mean as a ``Lie algebroid in the
category of Lie algebroids''. Such double structures do arise in various areas of mathematics, especially from Poisson geometry and the theory of
Poisson actions.
Topological double groupoids arose naturally from homotopy theory \cite{bs76,bh78}, while double Lie algebroids have found
applications, for instance, in Lu's extension of Drinfel`d's work \cite{dfd93} on the classification of Poisson homogeneous spaces
\cite{lu97}. \tla-groupoids represent the intermediate object between double Lie groupoids and double Lie algebroids,
and provide an equivalent characterization of the compatibility for Poisson bivector and group multiplication in a Poisson
groupoid \cite{mkz92,mkz99}.\linebreak
\spa Double Lie groupoids, \tla-groupoids and double Lie algebroids were introduced by Mackenzie \cite{mkz92}, who also
developed Lie theory ``from double Lie groupoids to double Lie algebroids'' \cite{mkz92,mkz00a}. The relation between Poisson
groupoids and double structures was foreseen by Weinstein in \cite{ws88}, where a program for the integration of Poisson
groupoids to symplectic double groupoids was proposed, and further investigated by Mackenzie in \cite{mkz99}.\\ 
Instances of Lie theory for
the integration of double Lie algebroids to double Lie groupoids have appeared in literature. Karas{\"e}v noticed in \cite{ks86}
that to a 1-connected complete Poisson group one can always associate a symplectic double groupoid and Lu-Weinstein
obtained an alternative construction in \cite{lw89}, which applies, under the same
connectivity assumptions also in the noncomplete case. Recently, Li and Parmentier studied in \cite{lp05} a certain class of coboundary dynamical Poisson 
groupoids, introduced by Etingof and Varchenko \cite{ev98} in relation with the classical dynamical Yang-Baxter equation, and
 produced an integrating symplectic double groupoid in a special case. Mackenzie and Xu's integration of Lie bialgebroids to
 Poisson groupoids can be equivalently regarded as the integration of the cotangent double (Lie algebroid) of a Lie
 bialgebroid \cite{mkz98} to a  cotangent prolongation \tla-groupoid.

\spa In this Chapter, after clarifying the notions of fibred products in the category of Lie groupoids and Lie algebroids  in
Section \ref{fplgla} (our techniques are borrowed from Higgins and Mackenzie \cite{hm90a,hm90b}), and introducing  double Lie
groupoids and \tla-groupoids in Section \ref{dlglag}, we address the integration problem of an \tla-groupoid to a double Lie
groupoid. In Section \ref{ilag} the study of the integrability of fibred products of Lie algebroids, 
leads us to derive conditions for endowing the differentiable graph (always) integrating an \tla-groupoid with a further
compatible multiplication making it a double Lie groupoid. Our integrability conditions are to be understood as
Lie-algebroid-homotopy-lifting conditions along suitable morphisms of Lie algebroids. The main result of this Chapter is the
following
\bgn{thm}[\bf{\ref{intla}}] Let $\sf{\Ohm}$ be an \tla-groupoid with integrable top Lie algebroid. If source and
target of the top Lie groupoid of $\sf{\Ohm}$ are strongly transversal, there exists a unique vertically source
1-connected double Lie groupoid integrating $\sf\Ohm$.
\end{thm}
In the final Section we specialize our results to the case of the \tla-groupoid canonically associated with a Poisson groupoid
and show that, under our conditions, the integration always yields a symplectic double groupoid:
\bgn{thm}[\bf\ref{double}] Let $\poidd{(\calG,\Pi)}{M}$ be an integrable Poisson groupoid with weak dual Poisson groupoid
$(\calG^\star,\Pi^\star)$. If cotangent source and target map of $\calG$ are strongly transversal, the symplectic
groupoid $\calS$ of $\calG$ carries a further Lie groupoid making it a symplectic groupoid for $\calG^\star$ and
$$
\bcat
\xy
*+{}="0",    <-0.7cm,0.7cm>
*+{\calS}="1", <0.7cm,0.7cm>
*+{\calG^\star}="2", <-0.7cm,-0.7cm>
*+{\calG}="3", <0.7cm,-0.7cm>
*+{M}="4",
\ar  @ <0.07cm>   @{->} "1";"2"^{} 
\ar  @ <-0.07cm>  @{->} "1";"2"_{}  
\ar  @ <-0.07cm>  @{->} "1";"3"_{}
\ar  @ <0.07cm>   @{->} "1";"3"^{}
\ar  @ <0.07cm>   @{->} "2";"4"^{}
\ar  @ <-0.07cm>  @{->} "2";"4"_{}
\ar  @ <-0.07cm>  @{->} "3";"4"_{}
\ar  @ <0.07cm>   @{->} "3";"4"^{}
\endxy
\ecat
$$
a symplectic double groupoid.
\end{thm}
We conclude this Chapter checking the integrability conditions in the case of complete Poisson groups, 
which allows for explicit computations, to extend Lu and Weinstein's result, in the complete case, by dropping all the 
connectivity assumptions on the Poisson group.
\bgn{thm}[\bf\ref{mainddd}] For any complete Poisson group $G$, the source  1-connected
symplectic groupoid $\calS$ of $G$ carries a unique Lie groupoid structure over the
1-connected dual Poisson group $G^\star$ making it a symplectic groupoid for the dual
Poisson structure and a double of $G$.
\end{thm}
We shall obtain further examples of integrable \tla-groupoids, arising from Poisson actions, in the next Chapter.
\newpage
\subsection*{Notations and remarks on groupoid objects}\hfill

\vs{0.1}
\spa Let $\gpdm$ be a groupoid; the \tsf{nerves} of the underlying category, i.e. the
strings of composable elements
\be
\calG^{(n)}&:=&\{(g_0,\ldots g_n)\in\calG^{\times
n}\:|\:\tar({g_i})=\sor({g_{i-1})}\:,\:i=1,\ldots n\}
\qquad,\qquad n\geq 1,\\
\calG^{(0)}&:=& M 
\ee
can also be inductively defined as fibred products
$$
\calG^{(n)}=\calG\fib{\sor}{\tar\comp\pp_1}\calG^{(n-1)}\:\:,\:\:n>1\:\:,
$$
where $\pp_1:\calG^{(n-1)}\to\calG$ is the restriction of the first projection for all
$n$.
The groupoid compatibility conditions may be grouped in two sets of commuting
diagrams: those for the \tsf{graph compatibility}
$$
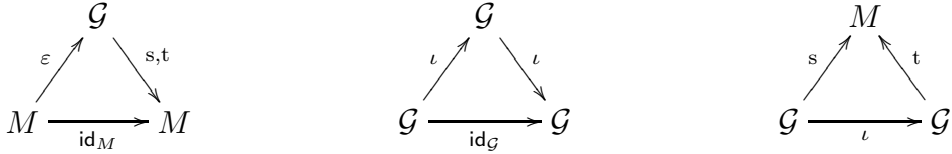
\bcat
\xy
*+{}="0",    <-1cm,0cm>
*+{M}="1", <1cm,0cm>
*+{M}="2", <0cm,1.4cm>
*+{\calG}="3", 
\ar     @ {->} "1";"2"_{\id_M} 
\ar     @ {->} "1";"3"^{\eps}  
\ar     @ {->} "3";"2"^{\sor,\tar}
\endxy
\ecat
\qquad\qquad\qquad
\bcat
\xy
*+{}="0",    <-1cm,0cm>
*+{\calG}="1", <1cm,0cm>
*+{\calG}="2", <0cm,1.4cm>
*+{\calG}="3", 
\ar     @ {->} "1";"2"_{\id_\calG} 
\ar     @ {->} "1";"3"^{\iota}  
\ar     @ {->} "3";"2"^{\iota}
\endxy
\ecat
\qquad\qquad\qquad
\bcat
\xy
*+{}="0",    <-1cm,0cm>
*+{\calG}="1", <1cm,0cm>
*+{\calG}="2", <0cm,1.4cm>
*+{M}="3", 
\ar     @ {->} "1";"2"_{\iota} 
\ar     @ {->} "1";"3"^{\sor}  
\ar     @ {->} "2";"3"_{\tar}
\endxy
\ecat
$$
and those for the multiplication,
$$
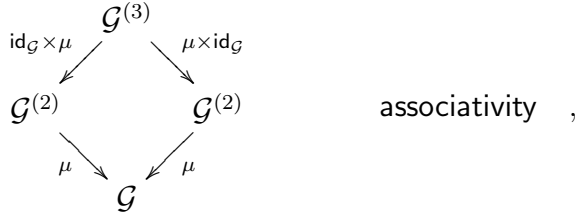
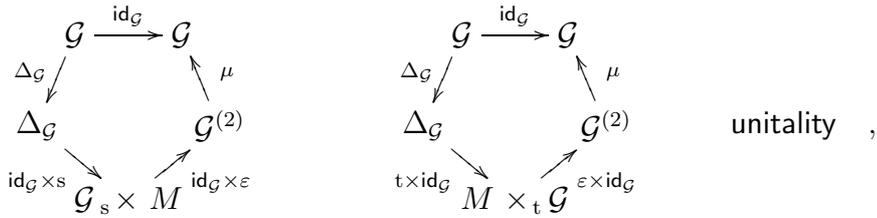
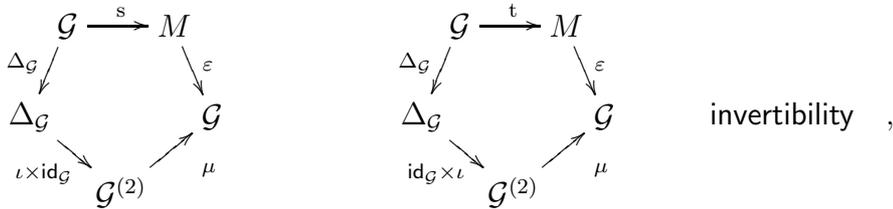
\qquad
\bcat
\xy
*+{}="0",    <0cm,1.2cm>
*+{\calG^{(3)}}="t", <0cm,-1.2cm>
*+{\calG}="b", <-1.2cm,0cm>
*+{\calG^{(2)}}="l", <+1.2cm,0cm>
*+{\calG^{(2)}}="r",  
\ar     @ {->} "t";"l"_{\id_\calG\times\mu} 
\ar     @ {->} "t";"r"^{\mu\times\id_\calG}  
\ar     @ {->} "l";"b"_{\mu}
\ar     @ {->} "r";"b"^{\mu}
\endxy
\ecat
\qquad\qquad
\hbox{\tsf{associativity}}\quad,
$$
$$
\bcat
\xy
*+{}="0",    <-0.7cm,1.7cm>
*+{\calG}="tl",  <0.7cm,1.7cm>
*+{\calG}="tr",    <-1.2cm,0.5cm>
*+{\Delta_\calG}="l", <1.2cm,0.5cm>
*+{\calG^{(2)}}="r", <0cm,-0.5cm>	
*+{\calG\fib{\sor}{}M}="b",
\ar     @ {->} "tl";"l"_{\Delta_\calG}
\ar     @ {->} "r";"tr"_{\:\:\mu}  
\ar     @ {->} "l";"b"_{\id_\calG\times\sor\:\:\:}  
\ar     @ {->} "b";"r"_{\:\:\:\id_\calG\times\eps}
\ar     @ {->} "tl";"tr"^{\id_\calG}
\endxy
\ecat
\qquad\qquad
\bcat
\xy
*+{}="0",    <-0.7cm,1.7cm>
*+{\calG}="tl",  <0.7cm,1.7cm>
*+{\calG}="tr",    <-1.2cm,0.5cm>
*+{\Delta_\calG}="l", <1.2cm,0.5cm>
*+{\calG^{(2)}}="r", <0cm,-0.5cm>	
*+{M\fib{}{\tar}\calG}="b", <3.8cm, 0.5cm>
*+{\hbox{\tsf{unitality}}\quad,}="z"
\ar     @ {->} "tl";"l"_{\Delta_\calG}
\ar     @ {->} "r";"tr"_{\:\:\mu}  
\ar     @ {->} "l";"b"_{\tar\times\id_\calG\:\:\:}  
\ar     @ {->} "b";"r"_{\:\:\:\eps\times\id_\calG}
\ar     @ {->} "tl";"tr"^{\id_\calG}
\endxy
\ecat
$$
and
$$
\bcat
\xy
*+{}="0",    <-0.7cm,1.7cm>
*+{\calG}="tl",  <0.7cm,1.7cm>
*+{M}="tr",    <-1.2cm,0.5cm>
*+{\Delta_\calG}="l", <1.2cm,0.5cm>
*+{\calG}="r", <0cm,-0.5cm>	
*+{\calG^{(2)}}="b",
\ar     @ {->} "tl";"l"_{\Delta_\calG}
\ar     @ {->} "tr";"r"^{\eps}  
\ar     @ {->} "l";"b"_{\iota\times\id_\calG}  
\ar     @ {->} "b";"r"_{\quad\mu}
\ar     @ {->} "tl";"tr"^{\sor}
\endxy
\ecat
\qquad\qquad\quad
\bcat
\xy
*+{}="0",    <-0.7cm,1.7cm>
*+{\calG}="tl",  <0.7cm,1.7cm>
*+{M}="tr",    <-1.2cm,0.5cm>
*+{\Delta_\calG}="l", <1.2cm,0.5cm>
*+{\calG}="r", <0cm,-0.5cm>	
*+{\calG^{(2)}}="b", <3.8cm, 0.5cm>
*+{\hbox{\tsf{invertibility}}\quad,}="z"
\ar     @ {->} "tl";"l"_{\Delta_\calG}
\ar     @ {->} "tr";"r"^{\eps}  
\ar     @ {->} "l";"b"_{\id_\calG\times\iota}  
\ar     @ {->} "b";"r"_{\quad\mu}
\ar     @ {->} "tl";"tr"^{\tar}
\endxy
\ecat
$$
where we have denoted with $\Delta_\calG$ both the diagonal in $\calG\times\calG$ and
the diagonal map $\calG\to\calG\times\calG$. 
It makes therefore sense to consider
groupoid objects in any small category with direct products and fibred products.
\bgn{definition}\label{gpdobj}
Let $\mathbf{C}$ be a small category with direct products and pullbacks.
Consider objects $\pmb{\calG},\mathbf{M}\in{\mathrm{Obj}(\mathbf{C})}$,
such that the set underlying $\pmb{\calG}$ carries a groupoid structure over the set
underlying $\mathbf{M}$, with structural maps 
$(\pmb{\sor},\pmb{\tar},\pmb{\eps},\pmb{\iota},\pmb{\mu})$; 
$\poidd{\pmb{\calG}}{\mathbf{M}}$ is a \tsf{groupoid object in}  $\mathbf{C}$ if all the
structural maps are arrows in $\mathbf{C}$. 
\end{definition}
Note that, for any groupoid $\gpdm$, given the source map, the inversion map and
the division map $\delta:\calG\fib{\sor}{\sor}\calG\to\calG$, the remaining
structural maps are recovered by considering the diagrams
\bgn{equation}\label{groupoidobjects}
\bcat
\xy
*+{}="0",    <-0.7cm,1.7cm>
*+{\calG}="tl",  <0.7cm,1.7cm>
*+{\calG}="tr",    <-1.2cm,0.5cm>
*+{\Delta_\calG}="l", <1.2cm,0.5cm>
*+{\calG\fib{\sor}{\sor}\calG}="r", <0cm,-0.5cm>	
*+{M\fib{}{\sor}\calG}="b",
\ar     @ {->} "tl";"l"_{\Delta_\calG}
\ar     @ {->} "r";"tr"_{\quad\delta}  
\ar     @ {->} "l";"b"_{\sor\times\id_\calG}  
\ar     @ {->} "b";"r"_{\:\:\eps\times\id_\calG}
\ar     @ {->} "tl";"tr"^{\iota}
\endxy
\ecat
\quad\qquad
\bcat
\xy
*+{}="0",    <-1.2cm,1.2cm>
*+{\calG}="l",  <1.2cm,1.2cm>
*+{M}="r",    <0cm,-0.2cm>
*+{\calG}="b", 
\ar     @ {->} "l";"r"^{\tar}
\ar     @ {->} "l";"b"_{\iota}  
\ar     @ {->} "b";"r"_{\sor}  
\endxy
\ecat
\quad\qquad
\bcat
\xy
*+{}="0",    <-1.2cm,1.2cm>
*+{\calG^{(2)}}="l",  <1.2cm,1.2cm>
*+{\calG}="r",    <0cm,-0.2cm>
*+{\calG\fib{\sor}{\sor}\calG}="b", 
\ar     @ {->} "l";"r"^{\mu}
\ar     @ {->} "l";"b"_{\id_\calG\times\iota}  
\ar     @ {->} "b";"r"_{\delta}  
\endxy
\ecat
\quad.
\end{equation}
Thus, a groupoid $\gpdm$ could be alternatively defined in terms of the data 
$(\calG,M,\sor,\eps,\delta)$, imposing the commutativity of the 
suitable diagrams; that is however lengthy and not illuminating.
Nevertheless this description is quite handy when dealing with groupoid objects.\\
According to the above discussion, we have
\bgn{proposition}\label{obietto} Let $\mathbf{C}$ a small category with direct products and  fibred
products. For any objects $\pmb{\calG},\mathbf{M}\in\mathrm{Obj}(\mathbf{C})$ a groupoid 
$\poidd{\pmb{\calG}}{\mathbf{M}}$ is a groupoid object in $\mathbf{C}$  if{f}  the source map, the
unit section and the division map of $\poidd{\pmb{\calG}}{\mathbf{M}}$ are arrows in
$\mathbf{C}$.
\end{proposition}
\bgn{proof} If $\poidd{\pmb{\calG}}{\mathbf{M}}$ is a groupoid object in $\mathbf{C}$, the
division map $\pmb{\delta}:=\pmb{\mu}\comp(\id_{\pmb{\calG}}\times \pmb{\iota})$ is also an
arrow. For the opposite implication use the diagrams (\ref{groupoidobjects}) in
sequence to recover $\pmb{\iota}$, $\pmb{\tar}$ and $\pmb{\mu}$ in terms of compositions
of arrows.
\end{proof}
If $\mathbf{C}$ does not possess general fibred products, the same idea applies, 
provided the relevant fibred products exist.
\bgn{example} A Poisson groupoid $\poidd{(\calG,\Pi)}{M}$ is \emph{not} a groupoid object in the category of Poisson
manifolds, for many reasons. First of all the target map is anti-Poisson for the induced Poisson structure on the
base manifold and the inversion map is always anti-Poisson. Enlarging the category of Poisson manifolds by adding the
anti-Poisson maps to the space of arrows does not solve the problem, since the domain of the groupoid multiplication
is not canonically endowed with a Poisson structure. 
\end{example} 
\bgn{example}\label{esempiuccio} For any group $G$ the total space of the fundamental groupoid  $\poidd{\Pi(G)}{G}$
also carries a group structure induced by pointwise multiplication of paths in $G$. For such structures it is a
groupoid object in the category of groups, since it is also possible to pointwise multiply homotopies, therefore
pointwise multiplication of homotopy classes is well defined and commutes with  concatenation up to
reparametrization. Note that the space of composable pairs  $\Pi(G)^{(2)}$ is described by all triples
$
([l], l(0)=g=r(1),[r])
$,
where $[l]$ and $[r]$ are homotopy classes of paths in $G$; setting  for all 
$([l_\pm], g_\pm,[r_\pm])\in \Pi(G)^{(2)}$
$$
\qquad
([l_+], g_+,[r_+])\cdot([l_-], g_-,[r_-])
:=
([l_+\cdot l_-], g_+g_-,[r_+\cdot r_-])
\qquad,
$$
where homotopy classes are multiplied pointwise for any choice of representatives, yields a well defined group
multiplication. Similarly are defined the groups on the relevant fibred products.
\end{example}
\vs{1}
\section{Fibred products of Lie groupoids and Lie algebroids}\label{fplgla}
\begin{quotation}
In this technical Section we clarify the notions of fibred products in the category of Lie groupoids and of Lie
algebroids, showing that they exist under the natural transversality conditions. In the case of Lie groupoids, the
proof of existence is essentially independent of the groupoid multiplication. For this reason, and for later
purposes, we introduce here differentiable graphs, namely ``Lie groupoids without a multiplication'' and study the
construction of fibred products in their category.
\end{quotation}

\vs{0.2}
\spa\quad
\bgn{definition} A $($\tsf{differentiable}$)$ \tsf{graph} on a pair of manifolds  $(\Gamma,M)$ is given by an
immersion $\eps:M\to\Gamma$ and a pair of submersions $\sor,\tar:\Gamma\to M$, which are both left inverses to
$\eps$: $\sor\comp\eps=\id_M$  and $\tar\comp\eps=\id_M$%
\fn{Pradines \cite{p84} uses the term \emph{differentiable graph}, when only a pair of submersions is given; 
similar objects are called bisubmersions in \cite{aspp}. We chose to use ``\emph{graph}'' for the sake of economy.}.
\end{definition}
Clearly Lie groupoids are differentiable graphs in the sense of the above definition,  hence so are vector bundles
and Lie algebroids for the underlying abelian groupoids; thus we shall use the same terminology for graphs as for Lie
groupoids: for any differentiable graph $(\Gamma,M)$ we shall call $\eps:M\to\Gamma$ the unit section and 
$\sor:\Gamma\to M$, respectively $\tar:\Gamma\to M$, the source map, respectively the target map. The natural notion
of a \tsf{morphism of graphs} $(\Gamma',M')\to(\Gamma, M)$, namely a pair of smooth maps $\fii:\Gamma'\to\Gamma$ and
$f:M'\to M$ such that
$$\qquad
\eps\comp f=\fii\comp\eps'
\qquad\qquad
\sor\comp\fii=f\comp\sor'
\qquad\qquad
\tar\comp\fii=f\comp\tar'
\qquad,
$$
makes the set of differentiable graphs a category with direct products. For any
differentiable graph $(\Gamma,M)$, a pair of submanifolds $\Gamma^o\subset\Gamma$ and
$M^o\subset M$ is called a $($\tsf{differentiable}$)$ \tsf{subgraph} if source and
target of $\Gamma$ restrict to submersions $\Gamma^o\to M^o$; a subgraph shall be
called wide if $M^o=M$.\\
Fibred products of differentiable graphs and  preimages of differentiable subgraphs
under morphisms of graphs (in particular kernels, to be defined  below) are not differentiable
graphs in general; however they are, under natural transversality conditions.\\   
Recall (for example, from \cite{amrbook,lang,conlon}) that a smooth map $f:N\to M$ is
\tsf{transversal to a submanifold}  $M^o\subset M$ if the tangent map induces a surjection $T_p
N\to T_{f(p)}M/T_{f(p)}M^o$ for all $p\in N$, equivalently if
$$
\qquad
\dd fT_p N + T_{f(p)}M^o=T_{f(q)}M\qquad,\qquad p\in N
\qquad.
$$
If $f:N\to M$ is transversal to $M^o\subset M$, we shall write $f\trans M^o$; in that
case $f\inverse(M^o)\subset N$ is a submanifold and 
$$
T_{p_o}f\inverse(M^o)=(\dd f)_{p_o}\inverse T_{f(p_o)}M^o
\qquad,\qquad
p_o\in f\inverse(M^o)
\qquad.
$$
Note that transversality is not a necessary condition for $f\inverse(M^o)$ to be a submanifold;
however, last formula holds true whenever $f\inverse(M^o)$ is a smooth submanifold.
\bgn{example} Consider $f:\RR^2\to\RR^2$, $f(x,y)=(xy,0)$ and $M^o=\{x\neq 0, y=0\}$:
$f\inverse(M^o)=\RR^2\backslash\{xy=0\}\subset\RR^2$ is a smooth submanifold even if $f$ is not
transversal to $M^o$, 
$$
\qquad
\dd f_{(x,y)}T_{(x,y)}\RR^2=T_{xy}\RR\oplus\{0_y\}=T_{xy}M^0\neq\RR^2
\qquad,
$$ 
and
$Tf\inverse(M^o)=T_{f\inverse(M^o)}\RR^2=(\dd f)\inverse TM^0$.
\end{example}
Two smooth maps $f_{1,2}:M^{1,2}\to M$ are said transversal, and we
shall write $f_1\pitchfork f_2$, if $f_1\times f_2:M^1\times M^2\to
M\times M$ is transversal to the diagonal $\Delta_M\subset M\times M$,
equivalently 
$$
\qquad
\dd f_1T_{q_1} M^1 +\dd f_2T_{q_2}M^2
=
T_{q}M
\qquad,
$$
for all 
$(q_1,q_2)\in M^{12}$, $q=f_{1,2}(q_{1,2})$.
In that case the fibred product 
$
M^{12}
:=
M^1\fib{f_1}{f_2}M^2\equiv(f_1\times f_2)\inverse(\Delta_{M})
\subset 
M^1\times M^2$ is a smooth submanifold and
$$
\qquad
T_{(q_1,q_2)}M^{12}=T_{q_1}M^{1}\fib{\dd f_1}{\dd f_2}T_{q_1}M^{1}
\qquad,\qquad 
(q_1,q_2)\in M^{12}
\qquad,
$$ 
where the fibred product on the right hand side is taken in the category of vector
spaces.
\bgn{proposition}\label{les} Let $\fii_{1,2}:\Gamma^{1,2}\to \Gamma$ be morphisms of graphs over
$f_{1,2}:M^{1,2}\to M$ such that $\fii_1\pitchfork \fii_2$ and $f_1\pitchfork f_2$ (so
that the fibred products 
$
\Gamma^{12}:=
\Gamma^1\fib{\fii_1}{\fii_2}\Gamma^2\subset\Gamma^1\times\Gamma^2
$
and
$
M^{12}:=
M^1\fib{\fii_1}{\fii_2}M^2\subset M^1\times M^2
$
are smooth submanifolds). Then $(\Gamma^{12},M^{12})$ is a subgraph of the direct
product $(\Gamma^1\times\Gamma^2, M^1\times M^2 )$ if{f}  the {\em \tsf{ source
transversality condition}}
\bgn{equation}\label{sourcetrans}
\dd\fii_1T_{x_1}^{\sor_1}\Gamma^1 +\dd\fii_2T_{x_2}^{\sor_2}\Gamma^2
=
T_{q}^{\sor}\Gamma
\end{equation}
holds for all $(x_1,x_2)\in \Gamma^{12}$ and $x=\fii_{1,2}(x_{1,2})$.
\end{proposition}
Note that condition (\ref{sourcetrans}) above is equivalent to asking
all the restrictions of $\phi_1$ and $\phi_2$ to the source fibres to be transversal,
whenever they map to the same source fibre of $\Gamma$.
\bgn{proof} Since $\sor_{12}:(\sor_1,\sor_2)|_{\Gamma^{12}}$ 
takes values in $M^{12}$ it suffices to check its submersivity, the same argument
applies to $\tar_{12}$. Applying the snake lemma (see, for instance, \cite{weibel}) to the exact commuting diagram
$$
\xymatrix{
0\:\ar@{->}[r]&\:
T_{x_1}^{\sor_1}\Gamma^1\oplus T_{x_2}^{\sor_2}\Gamma^2\ar@{->}[r]\ar@{->}_{}[d]&
T_{x_1}\Gamma^1\oplus
T_{x_2}\Gamma^2\ar@{->}^{\hbox{\tiny$\bgn{array}{c}\dd\sor_1\oplus\dd\sor_2\\\\\end{array}$}}[r]
\ar@{->}_{\dd\fii_1-\dd\fii_2}[d]&
T_{q_1}M^1\oplus T_{q_2}M^2\ar@{->}[r]\ar@{->}_{\dd f_1-\dd f_2}[d]&
0\\
0\:\ar@{->}[r]&\:
T_{x}^{\sor}M\:\ar@{->}[r]&\:
T_x\Gamma\ar@{->}^{\dd\sor}[r]&\:
T_qM\:\ar@{->}[r]&
0\\
}
$$ 
yields the long exact sequence
\bgn{equation}\label{snakefib}
\xymatrix{
0
\:\ar@{->}[r]&\:
\ker\,(\dd\fii_1-\dd\fii_2)|_{T_{x_1}^{\sor_1}\Gamma^1\oplus T_{x_2}^{\sor_2}\Gamma^2}\ar@{->}[r]\:&\:
T_{(g_1,g_2)}\Gamma^{12}\ar@{->}[r]&\:
T_{(q_1,q_2)}M^{12}
\:\ar@{->}[dll]\\
&
\coker\,(\dd\fii_1-\dd\fii_2)|_{T_{x_1}^{\sor_1}\Gamma^1\oplus T_{x_2}^{\sor_2}\Gamma^2}
\:\ar@{->}[r]&
C_\Gamma
\:\ar@{->}[r]&
C_M
\ar@{->}[r]&\:0\:,
}
\end{equation}
where the rightmost arrow on the top row is the restriction of $\dd\sor_1\times\dd\sor_2$ and
the cokernels $C_\Gamma=\coker\,(\dd\fii_1-\dd\fii_2)$ and $C_M=\coker\,(\dd f_1-\dd f_2)$
vanish.
\end{proof}
By specializing last result we obtain transversality conditions for preimages of
differentiable subgraphs under morphisms of graphs to be also differentiable subgraphs. 
\bgn{corollary} Let $\fii:\Gamma'\to\Gamma$ be a morphism of differentiable graphs over
$f:M'\to M$ and $(\Gamma^o,M^o)\subset(\Gamma, M)$ a subgraph such that $f\pitchfork
M^o$ and $\phi\pitchfork \Gamma^o$ (so that $f\inverse(M^o)\subset M'$ and
$\phi\inverse(\Gamma^o)\subset\Gamma'$ are smooth submanifolds). Then 
$(\phi\inverse(\Gamma^o),f\inverse(M^o))\subset(\Gamma',M')$ is a differentiable 
subgraph if{f} 
$$
\dd\fii T_{q'}^{\sor'}\Gamma' +T_{f(q')}^{\sor_o}\Gamma^o
=
T_{f(q')}^{\sor}\Gamma
$$
for all $q'\in f\inverse(M^o)$. 
\end{corollary}
\bgn{proof} Identify $\fii\inverse(\Gamma^o)$ with the fibred product
$\Gamma^1\fib{\fii}{\iota_o}\Gamma^o$ for the inclusion $\iota_o:\Gamma^o\to\Gamma$.
\end{proof}
For any morphism $\fii:\Gamma'\to\Gamma$ of differentiable graphs over
$f:M'\to M$, we shall say that the preimage $\ker\fii:=(\fii\inverse(\eps(M)), M)$ of the 
trivial wide subgraph of $(\Gamma,M)$ is the \tsf{kernel} of $(\fii,f)$. The
transversality condition of last corollary on $f$ is trivial and the source transversality
condition reduces to 
\bgn{equation}\label{ssub}
\qquad
\dd\fii T_{q'}^{\sor'}\Gamma'=T_{f(q')}^{\sor}\Gamma
\qquad,\qquad q'\in M'\qquad.
\end{equation}  
We shall say that a morphism of graphs satisfying (\ref{ssub}) is \tsf{source submersive}.
If $\fii$ is \emph{not} transversal to $\eps(M)$ but $\ker\fii$ is however smooth, it could still
happen to be a subgraph.
\bgn{lemma}\label{kergr} Let $\fii:\Gamma'\to\Gamma$ be a morphism of differentiable graphs over
$f:M'\to M$, such that $\fii\inverse(\eps(M))\subset \Gamma'$ is a smooth submanifold,
then $\ker\fii$ is a smooth graph if{f} 
\bgn{equation}\label{ga}
\dd\fii  T_{q'}^{\sor'}\Gamma\simeq \dd\fii T_{f(q)}\Gamma/\dd\eps\dd f(T_{q'} M')
\end{equation}
for all $q'\in \fii\inverse(\eps(M))$, $f(q')=q$.
\end{lemma}
\bgn{proof} The long exact sequence in the proof of proposition \ref{les} takes the form
$$
\qquad
\xymatrix{
0
\:\ar@{->}[r]&\:
\ker\,\dd\fii|_{ T_{q'}^{\sor'}\Gamma'}\ar@{->}[r]\:&\:
T_{g'}\ker\,\fii
	\ar@{->}[r]&\:
T_{q'}M'
\:\ar@{->}[dll]\\
&
\coker\,\dd\fii|_{ T_{q'}^{\sor'}\Gamma'}
\:\ar@{->}[r]&
C_\Gamma
\:\ar@{->}[r]&
C_M
\ar@{->}[r]&\:0\qquad\qquad,
}
$$
where the cokernels now are $C_M=0$, 
$C_\Gamma=T_{f(q)}\Gamma/(\dd\fii T_{q'}\Gamma' +\dd\eps(T_q M))\simeq 
T_{f(q)}^{\sor}\Gamma/(\dd\fii T_{q'}\Gamma'/\dd\eps\dd f(T_q' M'))$ 
and
$\coker\,\dd\fii|_{ T_{q'}^{\sor'}\Gamma'}=T_{f(q)}^{\sor}\Gamma/\dd\fii T_{q'}^{\sor'}\Gamma'$.
\end{proof}
A similar argument could be applied to the long exact sequence for fibred products of arbitrary
graphs, but the resulting condition is not illuminating, since the cokernels do not have, in
general, a nice description.
\bgn{example} Consider the action of $S^1$ on $\RR^2$ by rotation. The anchor of the action groupoid $S^1\acts\RR^2$
is a morphism of differentiable graphs to the pair groupoid  $\RR^2\times\RR^2$, which is not source submersive. Its
kernel is the isotropy groupoid, whose fibres are $S^1$ over $(x,y)=0$ and $\{1\}$ over $(x,y)\neq 0$, hence not a
differentiable subgraph.
\end{example}
\bgn{example} Consider the action of $\RR$ on the cylinder $S^1\times\RR$ by rotation.  The action groupoid 
$\RR\acts (S^1\times\RR)$ is a graph with differentiable kernel  $\cup_{k\in\ZZ}\{2\pi k\}\times S^1\times\RR$ for
the groupoid anchor  $\chi:\RR\acts (S^1\times\RR)\to S^1\times\RR\times  S^1\times\RR$, 
$(x;\theta,h)\mapsto(x+\theta,h;\theta,h)$. The anchor $\chi$ is not transversal to the diagonal, nor it
is source submersive, but 
\be
\dd\chi\, T_{(2\pi k,\theta,h)}\RR\acts (S^1\times\RR)
&\simeq&
\RR\oplus 0\oplus 0\oplus 0 +\Delta_{\RR\oplus\RR}\\
&=&
\dd\chi  T^\sor_{(2\pi k,\theta,h)}\RR\acts (S^1\times\RR) 
+
\Delta_{T_{(\theta,h)}S^1\times\RR}\qquad,
\ee
that is, condition (\ref{ga}) is satisfied.
\end{example}

\spa Fibred products of Lie groupoids are naturally endowed with a groupoid structure, therefore they
stay in the category, provided the transversality conditions for the underlying graphs hold. The
groupoid multiplication allows reducing the source transversality condition to a requirement on
the induced Lie algebroids.
\bgn{theorem}\label{fibgpd} Let $\fii_{1,2}:\calG^{1,2}\to \calG$ be morphisms of Lie groupoids
over $f_{1,2}:M^{1,2}\to M$ such that $\fii_1\pitchfork \fii_2$ and $f_1\pitchfork f_2$ (so that
the fibred products 
$
\calG^{12}:=
\calG^1\fib{\fii_1}{\fii_2}\calG^2\subset\calG^1\times\calG^2
$
and
$
M^{12}:=
M^1\fib{\fii_1}{\fii_2}M^2\subset M^1\times M^2
$
are smooth submanifolds). Then $(\calG^{12},M^{12})$ is a Lie subgroupoid of the direct
product $(\calG^1\times\calG^2, M^1\times M^2 )$ if{f}  the source
transversality condition
\bgn{equation}\label{sourcetransgpd}
\dd\fii_1T_{g_1}^{\sor_1}\calG^1 +\dd\fii_2T_{g_2}^{\sor_2}\calG^2
=
T_{g}^{\sor}\calG
\end{equation}
holds for all $(g_1,g_2)\in \calG^{12}$ and $g=\fii_{1,2}(g_{1,2})$, equivalently if{f} 
\bgn{equation}\label{sourcetransagd}
\phi_1A^1_{q_1} +\phi_2A^2_{q_2}
=
A_q
\end{equation}
for all $(q_1,q_2)\in M^{12}$ and $q=\fii_{1,2}(q_{1,2})$, where $A^{1,2}$, $A$ are
the Lie algebroids of $\calG^{1,2}$, $\calG$ and $\phi_{1,2}:A^{1,2}\to A$ the
morphisms of Lie algebroids induced by $\fii_{1,2}$. 
\end{theorem}
\bgn{proof} It remains to show that the two transversality conditions are equivalent.
Evaluating (\ref{sourcetransgpd}) on $g_{1,2}=\eps_{1,2}(q_{1,2})$ yields
(\ref{sourcetransagd}). Conversely, by equivariance under right translation, we have, for all $(g_1,g_2)\in \calG^{12}$,
\be
\dd\fii_1T_{g_1}^{\sor_1}\calG^1 +\dd\fii_2T_{g_2}^{\sor_2}\calG^2
&=&
\dd r_{\fii_1(g_1)}\inverse\phi_1A^1_{\tar_1(g_1)} 
+
\dd r_{\fii_2(g_2)}\inverse\phi_2A^2_{\tar_2(g_2)}\\
&=&
\dd r_{g}\inverse (\phi_1A^1_{\tar_1(g_1)} + \phi_2A^2_{\tar_2(g_2)})
\ee
and $T_{g}^{\sor}\calG=\dd r_{g}\inverse A_{\tar(g)}$, for $g=\fii_{1,2}(g_{1,2})$. 
Then (\ref{sourcetransgpd}) follows from (\ref{sourcetransagd}).
\end{proof}
Stronger transversality conditions for the existence of fibred products were given in \cite{mackbook}.
Note that the source transversality condition (\ref{sourcetransgpd}) and the infinitesimal
condition (\ref{sourcetransagd}) are equivalent, even dropping the transversality conditions on
$\fii_{1,2}$ and $f_{1,2}$.
We state a special case of last result for further reference.
\bgn{corollary}\label{pregpd} Let $\fii:\calG'\to\calG$ be a morphism of Lie groupoids over
$f:M'\to M$ and $(\calG^o,M^o)\subset(\calG, M)$ a Lie subgroupoid such that $f\pitchfork M^o$
and $\phi\pitchfork \calG^o$ (so that $f\inverse(M^o)\subset M'$ and
$\phi\inverse(\calG^o)\subset\calG'$ are smooth submanifolds). Then 
$(\phi\inverse(\calG^o),f\inverse(M^o))\subset(\calG',M')$ is a Lie subgroupoid if{f} 
\bgn{equation}\label{transpre}
\dd\fii T_{g'}^{\sor'}\calG' +T_{f(g')}^{\sor_o}\calG^o
=
T_{f(q')}^{\sor}\calG
\end{equation}
for all $g'\in \fii\inverse(\calG^o)$, equivalently
$$
\phi A'_{q'} +A^o_{f(q')}
=
A_{f(q')}
$$
for all $q'\in f\inverse(M^o)$, where $A$, $A'$ and $A^0$ are the Lie algebroids of $\calG'$,
$\calG$ and $\calG^o$.  
\end{corollary}
In view of theorem \ref{fibgpd} and corollary \ref{pregpd} we shall say that a morphism of Lie groupoids
$\fii:\calG'\to\calG$ over $f:M'\to M$ is \tsf{transversal to a Lie subgroupoid} $\poidd{\calG^o}{M^o}$ of
$\poidd{\calG}{M}$ if $\phi\pitchfork \calG^o$, $f\pitchfork M^o$ and the source transversality condition
(\ref{transpre}) holds. Similarly, we shall say that  $\fii_{1,2}:\calG^{1,2}\to \calG$ are \tsf{transversal
morphisms of Lie groupoids} over $f_{1,2}:M^{1,2}\to M$  if $\fii_1\times \fii_2:\calG^1\times\calG^2\to
\calG\times\calG$ is transversal to the diagonal subgroupoid $\poidd{\Delta_\calG}{\Delta_M}$ of the direct product
groupoid  $\calG\times\calG$, equivalently if $\fii_1\pitchfork \fii_2$, $f_1\pitchfork f_2$ and the source
transversality condition (\ref{sourcetransgpd}) holds.
\bgn{remark} 
If a morphism of Lie groupoids $\fii:\calG'\to\calG$ over $f:M'\to M$ is transversal to the trivial subgroupoid
$\eps(M)\subset\calG$, its kernel  $\ker\fii=\fii\inverse(\eps(M))\subset\calG$ is Lie subgroupoid if{f}  $\phi$ is
source submersive in the sense of (\ref{ssub}). If $\fii$ is {\em not} transversal to $\eps(M)$, $\ker\fii$ could be a
Lie subgroupoid even if $\fii$ is not source submersive: e.g. for any Lie group $G$ consider  $\fii:G\to G$,
$g\mapsto e$. In this case lemma \ref{kergr} gives a necessary and sufficient condition for a smooth kernel groupoid
to be a Lie subgroupoid. 
\end{remark} 

\spa Next we shall consider preimages and fibred products in the category of Lie algebroids. Differently from the
case of Lie groupoids, preimages of Lie subalgebroids under morphisms always carry a Lie algebroid structure,
provided they are smooth vector bundles. 
\bgn{proposition}\label{preagd} Let $A\to M$ and $A'\to M'$ be Lie algebroids and $B\to N$ a Lie subalgebroid of
$A\to M$. Consider a morphism of Lie algebroids $\phi:A'\to A$ over $f:M'\to M$ such that $f\trans N$ and
$\phi\inverse(B)\to f\inverse(N)$ is a smooth vector subbundle. Then  $\phi\inverse(B)\subset A'$ is a Lie
subalgebroid.  
\end{proposition}
\bgn{proof} 
First of all note that the anchor $\rho':A'\to TM'$ restricts to a bundle map $\phi\inverse(B)\to Tf\inverse(N)$,
since, for all $b'\in \phi\inverse(B)$, $\dd f\rho'(b')=\rho_B(\phi(b'))\in TN$, for the anchor $\rho_B$ on $B$
induced by the anchor of $A$; we can then check the condition 2. of lemma \ref{subalgebroid} (see also remark
\ref{checksubal}) to show that $\phi\inverse(B)\subset A'$ is a Lie subalgebroid. The condition is local  and we can
restrict to coordinate charts $U$ of $M'$ and $V$ of $M$, such that $f(U)\subset V$; setting $U_o:=U\cap
f\inverse(N)$, $V_o:=V\cap N$ and $f_o:=f|_{U_o}$ we can always assume that $U_o\subset U$ and $V_o\subset V$ are
smooth submanifolds, since $f\trans N$. Upon restriction, fix trivializing frames in duality  $\{e^\alpha\}$ for $A$
and $\{e_\alpha\}$ for $A^*$ over $U$. For any $\beta\in\Gamma(U,f^{\sf{+}} A)$ such that $\beta|_{U_o}=o$,
$$
\beta=\sum_\alpha\pair{e_\alpha\comp f}{\beta}(e^\alpha\comp f)
=
\sum_\alpha\beta_\alpha(e^\alpha\comp f)
$$ 
yields a decomposition $\{\beta_\alpha\otimes e^\alpha\}$  in $\cif(U)\otimes_{\cif(V)}\Gamma(V,A)$  with
$\beta_\alpha|_{U_o}=0$.  Let $a^{1,2}\in\Gamma(U,A)$ be local sections such that
$a^{1,2}|_{U_o}\in\Gamma(U,\phi\inverse(B))$, i.e. $\phi\comp a^{1,2}|_{U_o}\in\Gamma(U_o,f^{+}_o(B))$ and pick
decompositions
$$
\phi\comp a^{1,2}|_{U_o}=\sum_{i_{1,2}}u_{i_{1,2}}^{1,2}(b_{i_{1,2}}^{1,2}\comp f_o)
$$
with $\{u_{i_{1,2}}^{1,2}\}\subset\cif(U_o)$ and $\{b_{i_{1,2} }^{1,2}\}\subset\Gamma(V_o,B)$;  for any choice
of extensions 
$\{\wt{u}_{i_{1,2}}^{1,2}\}\subset\cif(U)$ of $\{u_{i_{1,2}}^{1,2}\}\subset\cif(U_o)$ and 
$\{\wt{b}_{i_{1,2}}^{1,2}\}\subset\Gamma(V, A)$ of $\{b_{i_{1,2} }^{1,2}\}\subset\Gamma(V_o,B)$,
$$
\left.
{\left(\phi\comp a^{1,2} 
-
\sum_{i_{1,2}}\wt{u}_{i_{1,2}}^{1,2}(\wt{b}_{i_{1,2}}^{1,2}\comp f)
\right)}
\right|_{U_o}=0
$$
therefore there exist decompositions
\bgn{equation}\label{belledeco}
\phi\comp a^{1,2} 
=
\sum_{i_{1,2}}\wt{u}_{i_{1,2}}^{1,2}(\wt{b}_{i_{1,2}}^{1,2}\comp f)
+
\sum_{l_{1,2}}v_{l_{1,2}}^{1,2}(d_{l_{1,2}}^{1,2}\comp f)
\end{equation}
with $\{d_{l_{1,2}}^{1,2}\}\subset\Gamma(V, A)$ and $\{v_{l_{1,2}}^{1,2}\}\subset\cif(U)$, with the property
that $v_{l_{1,2}}^{1,2}|_{U_o}=0$ for all $l_{1,2}$. Evaluating the bracket compatibility condition for
$\phi$ on the decompositions (\ref{belledeco}), yields
{\small
\be
\phi\comp\brak{a^1}{a^2}\hs{-0.3}
&=&\hs{-0.3}
\underset{i_1,i_2}\sum 
\wt{u}^{1}_{i_{1}}\wt{u}^{2}_{i_{2}}
(\brak{\wt{b}_{i_{1}}^{1}\!}{\wt{b}_{i_{2}}^{2}}\!\comp\!f)
\!+\!
\underset{i_2}\sum\, \rho'(a^1)(\wt{u}_2^{i_2})(\wt{b}_{i_{2}}^{2}\!\comp \!f)
\!-\!
\underset{i_1}\sum\, \rho'(a^2)(\wt{u}_1^{i_1})(\wt{b}_{i_{1}}^{1}\!\comp\! f)\\
\hs{-0.3}
&+&\hs{-0.3}
\underset{l_1,l_2}\sum 
v^{1}_{l_{1}}v^{2}_{l_{2}}
(\brak{d^1_{l_{1}}\!\!}{d^2_{l_{2}}}\!\comp\!f)
\!+\!
\underset{i_1,l_2}\sum 
\wt{u}^{1}_{i_{1}}v^{2}_{l_{2}}
(\brak{\wt{b}_{i_{1}}^{1}\!\!}{d^2_{l_{2}}}\!\comp\!f)
\!+\!
\underset{l_1,i_2}\sum 
v^{1}_{l_{1}}\wt{u}^{2}_{i_{2}}
(\brak{d^1_{l_{1}}\!\!}{\wt{b}_{i_{2}}^{2}}\!\comp\!f)\\
\hs{-0.3}
&+&\hs{-0.3}
\underset{l_2}\sum \rho'(a^1)(v^2_{l_2})(d^2_{l_2}\!\comp \!f)
\!-\!
\underset{l_1}\sum \rho'(a^2)(v^1_{l_1})(d^1_{l_1}\!\comp\! f)\hs{3}.\\
\ee
}%
\noindent Restricting last expression to $U_o$ the terms in the second line vanish, since so do the $v$'s,
those on the third line also vanish since $\rho'(a^{1,2})|_{U_o}$ is tangent to $U_o$ and $\dd
v_{1,2}^{l_{1,2}}\in\Gamma(N^*U_o)$; for the same reason the terms of the form $\rho(a)(\wt{u})$ do not depend on the
choice of extension. Then
\be
\phi\comp\brak{a^1}{a^2}|_{U_o}
&=&
\underset{i_1,i_2}\sum 
u^{1}_{i_{1}}u^{2}_{i_{2}}
(\brak{\wt{b}_{i_{1}}^{1}}{\wt{b}_{i_{2}}^{2}}|_{U_o}\!\comp\!f_o)
+
\underset{i_2}\sum \rho'(a_1)(u_2^{i_2})(b_2^{i_2}\comp f_o)\\
&-&
\underset{i_1}\sum \rho'(a_2)(u_1^{i_1})(b_1^{i_1}\comp f_o)
\ee
takes values in $B|_{U_o}$, since the local sections $\{\wt{b}_{i_{1}}^{1}\}$ extend sections of $B|_{U_o}$ and
$B\subset A$ is a Lie subalgebroid, thus 
$\brak{\wt{b}_{i_{1}}^{1}\!}{\wt{b}_{i_{2}}^{2}}|_{U_o}\in\Gamma(U_o,B)$. 
\end{proof}
\bgn{remark} Specializing last proposition to the trivial subalgebroid $B=M'$, we obtain that
the kernel of a morphism of Lie algebroids $\phi:A'\to A$ is always a Lie algebroid, 
provided $\phi$ has constant rank.
\end{remark}
The existence of fibred products of Lie algebroids under natural transversality conditions was stated without proof by Higgins
and Mackenzie in \cite{hm90a} and  can now be deduced from proposition
\ref{preagd}.
\bgn{theorem}\label{dracula} Let $\phi_{1,2}:A^{1,2}\to A$ be morphisms of Lie algebroids over $f_{1,2}:M^{1,2}\to M$
such that $f_1\pitchfork f_2$ and the fibred product
$
A^{12}
:=
A^1\fib{\fii_1}{\fii_2}A^2
\to
M^{12}
:=
M^1\fib{\fii_1}{\fii_2}M^2
$
is a smooth subbundle of the direct product $A^1\times A^2\to M^1\times M^2$.  Then $A^{12}\subset A^1\times A^2$ is
a Lie subalgebroid, hence a fibred product in the category of Lie algebroids. 
\end{theorem}
\bgn{proof} $(\phi_1\times\phi_2):A^1\times A^2\to B^{\times 2}$ is a morphism of Lie algebroids.
To see this, consider that $\psi_{1,2}:A^1\times A^2\to B$, $\psi_{1,2}:=\phi_{1,2}\comp\pr_{1,2}$ are morphisms
of Lie algebroids, for the projections $\pr_{1,2}:A^1\times A^2\to A^{1,2}$ and $\phi_1\times\phi_2$ is the unique
morphism of vector bundles whose composition with $\pr_{1,2}$ give $\psi_{1,2}$. Since $B^{\times 2}$ is a direct
product in the category of Lie algebroids there exists a unique morphism of Lie algebroids $\phi$ such that
$\pr_{1,2}\comp\phi=\psi_{1,2}$, since it is in the category of vector bundles,  $\phi=\phi_1\times\phi_2$. The
diagonal subbundle $\Delta_B\subset B^{\times 2}$ is the graph of the identity, it follows  from
proposition \ref{preagd} and corollary \ref{caraci} that the vector  bundle
$A^{12}=(\phi_1\times\phi_2)\inverse\Delta_B$ is a Lie subalgebroid of the direct product $A^1\times A^2$. For the
induced Lie algebroid structure, $A^{12}$ is indeed a fibred product: let $\chi_{1,2}:C\to A^{1,2}$ be morphisms
of Lie algebroids such that $\phi_1\comp \chi_1=\phi_2\comp \chi_2$, then there exist a unique morphism of Lie
algebroids $C\to A^1\times A^2$ such that $\pr_{1,2}\comp\chi=\chi_{1,2}$; since
$(\phi_1\times\phi_2)\comp\chi=(\phi_1\comp\chi,\phi_2\comp\chi)$ takes values in $\Delta_B$, $\chi$ takes values
in $A^{12}$. That is, $\chi:C\to A^{12}$ lifts $\chi_1$ and $\chi_2$. Regarding any other $\chi':C\to A^{12}$ with
the same property as a morphism  $C\to A^1\times A^2$, shows that $\chi=\chi'$.
\end{proof}
Typically, the fibred product Lie algebroid $A^{12}$ exists when $\phi_1\pitchfork\phi_2$, $f_1\pitchfork f_2$
and the infinitesimal linear transversality condition (\ref{sourcetransagd}) holds. In that case we shall say that
$\phi_1$ and $\phi_2$ are \tsf{transversal morphism of  Lie algebroids}.
\bgn{example} For any morphism of Lie algebroids $\phi:A'\to A$ over $f:M'\to M$, the (ordinary) graph
$\gr{\phi}\equiv A'\fib{\phi}{\id_A}A$ always carries a Lie algebroid structure over  $\gr{f}\equiv
M'\fib{f}{\id_M}M$, isomorphic to $A'\to M'$ by the projection $A'\fib{\phi}{\id_A}A\to A'$; the identity $A\to  A$
is transversal to any morphism of Lie algebroids $A'\to A$.
\end{example}
\bgn{example} For any map $f:N\to M$ the pullback algebroid $f\daga A\to N$ of a Lie algebroid $A\to M$, provided it
exists, is, by construction indeed, the fibred product $TN\fib{\dd f}{\rho}A$.
\end{example}
\bgn{remark} From the two equivalent transversality conditions of theorem \ref{fibgpd} it is clear that a fibred
product of Lie groupoids $\calG^1\fib{\fii_1}{\fii_2}\calG^2$, when it is  Lie, differentiates to a fibred 
product of Lie algebroids. In fact, the
source transversality condition (\ref{sourcetransgpd}), implies that the induced morphisms $\phi_{1,2}$ of 
Lie algebroids are transversal in the sense of the definition above; therefore 
the associated fibred product of Lie algebroids $A^1\fib{\phi_1}{\phi_2}A^2$ exists
and can be easily seen to coincide with the Lie algebroid of $\calG^1\fib{\fii_1}{\fii_2}\calG^2$ as a vector
bundle. It also does as a Lie algebroid by uniqueness.
\end{remark}
\vs{1}
\section{Double Lie groupoids and \tla-groupoids}\label{dlglag}
\begin{quotation} 
We define in the first part of this Section double Lie groupoids, their morphisms and derive conditions for kernels
of morphisms to be double Lie groupoids. Characterizing double Lie groupoids as suitably smooth groupoid objects in the
category of Lie groupoids, leads us directly in the second part to their infinitesimal invariant, namely
\tla-groupoids, i.e. suitably smooth groupoid objects in the category of Lie algebroids; we further obtain
conditions for kernels of morphisms of \tla-groupoids to stay in the category.
\end{quotation}
\vs{0.5}
\subsection{Double Lie groupoids}\hfill

\vs{0.1}
\spa A \tsf{double groupoid} in the sense of Ehresmann
\bgn{equation}\label{dbgpd}
\bcat
\xy
*+{}="0",    <-0.7cm,0.7cm>
*+{\calD}="1", <0.7cm,0.7cm>
*+{\calV}="2", <-0.7cm,-0.7cm>
*+{\calH}="3", <0.7cm,-0.7cm>
*+{M}="4",
\ar  @ <0.07cm>   @{->} "1";"2"^{} 
\ar  @ <-0.07cm>  @{->} "1";"2"_{}  
\ar  @ <-0.07cm>  @{->} "1";"3"_{}
\ar  @ <0.07cm>   @{->} "1";"3"^{}
\ar  @ <0.07cm>   @{->} "2";"4"^{}
\ar  @ <-0.07cm>  @{->} "2";"4"_{}
\ar  @ <-0.07cm>  @{->} "3";"4"_{}
\ar  @ <0.07cm>   @{->} "3";"4"^{}
\endxy
\ecat
\end{equation}
is a groupoid object in the category of groupoids. The definition is symmetric. Assume that $\poidd{\calD}{\calV}$, 
is a groupoid object in the category of groupoids with structure maps  ($\nsth$, $\ntth$, $\eth$, $\ith$, $\mth$). 
Then $\poidd{\calD}{\calH}$ and $\poidd{\calV}{M}$ are groupoids, with structure maps ($\nstv$, $\nttv$, $\etv$,
$\itv$, $\mtv$) and ($\nssv$, $\ntsv$, $\esv$, $\isv$, $\msv$), making ($\nsth$, $\ntth$, $\eth$, $\ith$, $\mth$)
morphisms of groupoids. The base diagrams for the groupoid structure of $\poidd{\calD}{\calV}$ define a  groupoid on
$\calH$ over $M$ given by the base maps ($\nssh$, $\ntsh$, $\esh$, $\ish$, $\msh$). The compatibility
conditions for the structure maps of $\poidd{\calD}{\calV}$ to be morphisms of groupoids over those of 
$\poidd{\calH}{M}$ with respect to the vertical groupoids are equivalent to the compatibility conditions for the
structure maps of $\poidd{\calD}{\calH}$ to be morphisms over those of $\poidd{\calV}{M}$ with respect to the
vertical groupoids. For example the compatibility of $\nsth:\calD\to\calV$ with the vertical multiplications, namely
$\nsth\comp\mtv=(\nsth\times\nsth)\comp\msv$, is the compatibility of $\mtv:\calD\fib{\stv}{\ttv}\calD\to\calD$ with
the horizontal source map; the other compatibility conditions can be read similarly in two ways. Then 
$\poidd{\calD}{\calH}$ is also  a groupoid object in the category of groupoids.
\\
We shall refer to the groupoid structures of a general double groupoid such as  (\ref{dbgpd}) as \emph{top
horizontal}, {\em top vertical},  {\em side horizontal} and {\em side vertical}, with the obvious meaning.
\begin{definition}\cite{mkz92}
A \tsf{\textsf{double Lie groupoid}} $\sf{D}:=(\calD,\calH,\calV,M)$ such as  (\ref{dbgpd})  is a double groupoid,
such that  $\poidd{\calD}{\calV}$, $\poidd{\calD}{\calH}$, $\poidd{\calH}{M}$, $\poidd{\calV}{M}$ are Lie groupoids
and the double source map  
$$
\mathbb{S}\doteq(\nstv,\nsth):\calD\rightarrow\calH\fib{\ssh}{\ssv}\calV 
$$
is  submersive\fn{%
In \cite{mkz92} the double source map is required to be also surjective; this condition does not really play a role
in the study of the internal structure of a double Lie groupoid and the descent to double  Lie algebroids. Moreover,
there are interesting examples, such as Lu and Weinstein's double of a Poisson group (\ref{lwd}) for instance, which
do not fulfill the double source surjectivity condition.
}.
\end{definition} 
We remark that the double source submersivity is required to make the domains of the
top multiplications,
$
\poidd{\calD\hbox{\tiny{${}_H^{\raise0.7mm\hbox{$\,(2)$}}$}}}
{\calV\hbox{\tiny{${}^{\raise0.7mm\hbox{$\,(2)$}}$}}}
$
and 
$
\poidd{\calD\hbox{\tiny{${}_V^{\raise0.7mm\hbox{$\,(2)$}}$}}}
{\calH\hbox{\tiny{${}^{\raise0.7mm\hbox{$\,(2)$}}$}}}
$,  Lie groupoids; thus, in particular, a double Lie groupoid is a groupoid object in the category of Lie groupoids.
 Let us introduce a special class of
morphisms of Lie groupoids.
\bgn{definition}%
A morphism of Lie groupoids $\fii:\calG'\to\calG$ over $f:M'\to M$ is called an
\tsf{\tlg-fibration} if%
\medskip\\
$1$. It is submersive and base submersive;%
\medskip\\
$2$. The characteristic map
$
(\fii,\sor'):\calG'\longrightarrow\calG\fib{\sor}{f}M'
$
is submersive, equivalently $\fii$ is source submersive.%
\medskip\\
A \tsf{strong \tlg-fibration} is an \tlg-fibration such that%
\medskip\\
$1$'. It is surjective and base surjective;%
\medskip\\
$2$'. The characteristic map is surjective, equivalently $\fii$ is source surjective.
\end{definition}
A strong \tlg-fibration is a fibration of Lie groupoids in the sense of Higgins and Mackenzie \cite{hm90b,mackbook}.\\
For a double groupoid (\ref{dbgpd}) all of whose side and top groupoids are Lie groupoids, the following
conditions are easily seen to be  are equivalent
\medskip\\
$i$) The double source map 
$
\mathbb{S}\doteq(\nstv,\nsth):\calD\rightarrow \calH\fib{\ssh}{\ssv}\calV
$
is  submersive,\medskip\\
$ii$) The top horizontal source map $\nsth:\calD\to\calV$ is an
\tlg-fibration,\medskip\\
$iii$) The top vertical source map $\nstv:\calD\to\calH$ is an
\tlg-fibration.\medskip\\
In particular, the top horizontal and vertical source maps of a double Lie groupoid are
\tlg-fibrations and the transversality conditions of theorem \ref{fibgpd} to make the fibred products
$\calD\fib{\sth}{\tth}\calD$ and $\calD\fib{\stv}{\ttv}\calD$ Lie groupoids are met.\\
We list below the typical examples; more interesting ones shall be studied throughout the rest of
this dissertation.
\bgn{example} \tsf{Double Lie groupoids}
\medskip\\$i$) Any Lie groupoid $\gpdm$ is the top horizontal groupoid of a double Lie
groupoid for the trivial groupoids on $\calG$ and $M$ as vertical groupoids.
\medskip\\$ii$) The \tsf{pair double Lie groupoid}
$$
\bcat
\xy
*+{}="0",    <-1cm,0.7cm>
*+{\calG\times\calG}="1", <1cm,0.7cm>
*+{M\times M}="2", <-1cm,-0.7cm>
*+{\calG}="3", <1cm,-0.7cm>
*+{M}="4",
\ar  @ <0.07cm>   @{->} "1";"2"^{} 
\ar  @ <-0.07cm>  @{->} "1";"2"_{}  
\ar  @ <-0.07cm>  @{->} "1";"3"_{}
\ar  @ <0.07cm>   @{->} "1";"3"^{}
\ar  @ <0.07cm>   @{->} "2";"4"^{}
\ar  @ <-0.07cm>  @{->} "2";"4"_{}
\ar  @ <-0.07cm>  @{->} "3";"4"_{}
\ar  @ <0.07cm>   @{->} "3";"4"^{}
\endxy
\ecat
$$
of any Lie groupoid $\gpdm$ is defined by the pair groupoid on the top vertical edge
and the direct product groupoid on the top horizontal edge. 
\medskip\\$iii$)  A double vector bundle in the sense of Ehresmann, i.e. a vector bundle
in the category of vector bundles, is a double Lie groupoid for the abelian structures
on the four sides. 
\medskip\\$iv$) For any Lie group $G$,
$$
\bcat
\xy
*+{}="0",    <-0.7cm,0.7cm>
*+{\Pi(G)}="1", <0.7cm,0.7cm>
*+{\bullet}="2", <-0.7cm,-0.7cm>
*+{G}="3", <0.7cm,-0.7cm>
*+{\bullet}="4",
\ar  @ <0.07cm>   @{->} "1";"2"^{} 
\ar  @ <-0.07cm>  @{->} "1";"2"_{}  
\ar  @ <-0.07cm>  @{->} "1";"3"_{}
\ar  @ <0.07cm>   @{->} "1";"3"^{}
\ar  @ <0.07cm>   @{->} "2";"4"^{}
\ar  @ <-0.07cm>  @{->} "2";"4"_{}
\ar  @ <-0.07cm>  @{->} "3";"4"_{}
\ar  @ <0.07cm>   @{->} "3";"4"^{}
\endxy
\ecat
$$
is a double Lie groupoid, for the pointwise 
multiplication of homotopy classes of paths described in  example \ref{esempiuccio}.
\end{example}
Morphisms and subobjects of double Lie groupoids are defined in the obvious way.
\bgn{definition}\label{mdlg} Let $\sfD^\pm=(\calD^\pm,\calH^\pm,\calV^\pm,M^\pm)$
be double Lie groupoids; a \tsf{morphism of double Lie groupoids} 
$\sf\Phi:\sfD^-\to\sfD^+$ 
\bgn{equation}\label{morphismdlg}
\bcat\xy
*+{}="0",    <-0.7cm,0.7cm>
*+{\calD^-}="1", <0.7cm,0.7cm>
*+{\calV^-}="2", <-0.7cm,-0.7cm>
*+{\calH^-}="3", <0.7cm,-0.7cm>
*+{M^-}="4",     <2.8cm,-0.7cm>
*+{\calD^+}="1'", <4.2cm,-0.7cm>
*+{\calV^+}="2'", <2.8cm,-2.1cm>
*+{\calH^+}="3'", <4.2cm,-2.1cm>
*+{M^+}="4'",
\ar  @ <0.07cm>   @{->} "1";"2"^{} 
\ar  @ <-0.07cm>  @{->} "1";"2"_{}  
\ar  @ <-0.07cm>  @{->} "1";"3"_{}
\ar  @ <0.07cm>   @{->} "1";"3"^{}
\ar  @ <0.07cm>   @{->} "2";"4"^{}
\ar  @ <-0.07cm>  @{->} "2";"4"_{}
\ar  @ <-0.07cm>  @{->} "3";"4"_{}
\ar  @ <0.07cm>   @{->} "3";"4"^{}
\ar  @ <0.07cm>   @{->} "1'";"2'"^{} 
\ar  @ <-0.07cm>  @{->} "1'";"2'"_{}  
\ar  @ <-0.07cm>  @{->} "1'";"3'"_{}
\ar  @ <0.07cm>   @{->} "1'";"3'"^{}
\ar  @ <0.07cm>   @{->} "2'";"4'"^{}
\ar  @ <-0.07cm>  @{->} "2'";"4'"_{}
\ar  @ <-0.07cm>  @{->} "3'";"4'"_{}
\ar  @ <0.07cm>   @{->} "3'";"4'"^{}
\ar               @{->} "1";"1'"^{\Phi}
\ar               @{->} "2";"2'"^{\fii_v}  
\ar               @{->} "3";"3'"^{\fii_h}  
\ar               @{->} "4";"4'"^{f}     
\endxy\ecat
\end{equation}
is given by four maps $(\Phi,\fii_h,\fii_v,f)$ such that all sides of the diagram above are morphisms of Lie
groupoids. $\sfD^-$ is called a \tsf{double Lie subgroupoid} if all of its side groupoids are Lie subgroupoids of the
corresponding sides of $\sfD^+$.
\end{definition}
It is straightforward to see that, with this notion of morphism, double Lie groupoids form a category. We shall
say that $(\Phi,\fii_v)$, respectively  $(\Phi,\fii_h)$, is the top horizontal, respectively vertical, component
of $\sf{\Phi}$ and, similarly, that $(\phi_h,f)$, respectively $(\phi_v,f)$, is the side horizontal, respectively
vertical, component. The following lemma gives a characterization of the kernel of a morphism of double Lie
groupoids; it is not, in general, a double Lie groupoid. However, when the top components are suitably well
behaved, the double Lie groupoid structure is inherited by the kernel.
\bgn{lemma}\label{kerdlg} For any morphism $(\ref{morphismdlg})$ of double Lie groupoids:
\medskip\\
$i)$ The top vertical kernel $\calK\vup:=\ker(\Phi,\fii_h)$ has a natural groupoid structure over the side
vertical kernel $\calK_v:=\ker(\fii_v,f)$, making
\bgn{equation}\label{biker}
\bcat\xy
*+{}="0",    <-0.7cm,0.7cm>
*+{\calK\vup}="1", <0.7cm,0.7cm>
*+{\calK_v}="2", <-0.7cm,-0.7cm>
*+{\calH}="3", <0.7cm,-0.7cm>
*+{M}="4",
\ar  @ <0.07cm>   @{->} "1";"2"^{} 
\ar  @ <-0.07cm>  @{->} "1";"2"_{}  
\ar  @ <-0.07cm>  @{->} "1";"3"_{}
\ar  @ <0.07cm>   @{->} "1";"3"^{}
\ar  @ <0.07cm>   @{->} "2";"4"^{}
\ar  @ <-0.07cm>  @{->} "2";"4"_{}
\ar  @ <-0.07cm>  @{->} "3";"4"_{}
\ar  @ <0.07cm>   @{->} "3";"4"^{}
\endxy\ecat
\end{equation}
a double subgroupoid;
\medskip\\
$ii)$ Assume that $\Phi$ and $\fii_v$ are \tlg-fibrations, respectively over $\fii_h$ and $f$ (so that  $\calK\vup$
and $\calK_v$ are Lie subgroupoids for the vertical  structures of $\sf{D}$); then, $\poidd{\calK\vup}{\calK_v}$ is a
Lie groupoid if{f} 
\bgn{equation}\label{regularbiker}
(\Phi,\sth^-):\calD^-\rightarrow\calD^+\fib{\sth^+}{\fii_v}\calV^-
\end{equation}
is submersive;
\medskip\\
$iii)$ Under the hypotheses of $($ii$)$, the double groupoid $(\ref{biker})$ is a  double Lie subgroupoid of $\sf{D}$ if{f} 
\bgn{equation}\label{regularbiker1}
(\Phi,\stv^-):\calD^-\rightarrow\calD^+\fib{\stv^+}{\fii_h}\calH^-
\end{equation}
is an \tlg-fibration over 
\bgn{equation}\label{regularbiker2}
(\fii_h,\ssv):\calV^-\rightarrow\calV^+\fib{\stv^+}{f}M^-
\end{equation}
\end{lemma}
Note that the fibred product Lie groupoid in ($iii$) is always well defined, since the top vertical source
$\nstv^+$is an \tlg-fibration, and $(\Phi,\stv^-)$ is always a morphism of Lie groupoids over $(\fii_h,\ssv^-)$.
\bgn{proof} ($i$) For any horizontally composable $k_{1,2}\in\calK\vup$
\be
\Phi(\mth^-(k_1,k_2))
&=&
\mth^+(\Phi(k_1),\Phi(k_2))\\
&=&
\mth^+(\etv^+(\fii_h(\stv^-(k_1))),\etv^+(\fii_h(\stv^-(k_1))))\\
&=&
\etv^+(\msh^+(\fii_h(\stv^-(k_1)),\fii_h(\stv^-(k_1))))\qquad,
\ee
thus the top horizontal multiplication restricts to $\calK\vup$; in the same way, it is
easy to see that all the top horizontal maps of $\calD^-$ restrict to
$\poidd{\calK\vup}{K_v}$. ($ii$) Taking the long exact sequence for the diagram
$$
\xymatrix{
0\:\ar@{->}[r]&\:
T_k\calK\vup\:\ar@{->}[r]\ar@{->}[d]&\:
T_{k}\calD^-\ar@{->}^{\dd\Phi}[r]\ar@{->}_{\dd\sth^-}[d]\:&\:
T_{\Phi(k)}\calD^+\:\ar@{->}[r]\ar@{->}_{\dd\sth^+}[d]&
0\\
0\:\ar@{->}[r]&\:
T_{\sth^-(k)}\calK_v\:\ar@{->}[r]&\:
T_{\sth^-(k)}\calV^-\ar@{->}^{\dd\fii_v}[r]\:&\:
T_{\fii_v(\sth^-(k))}\calV^+\:\ar@{->}[r]&
0\\
}
$$ 
yields $\coker\,\dd\sth^-|_{T\calK\vup}=0$ if{f}  $\coker\dd\Phi|_{T^{\sth^-}\calD^-}=0$, that is 
$\nsth:\calK\vup\to\calK_v$ is
submersive if{f}  so is  $(\Phi,\sth^-)$.
($iii$) Consider that, for
all $k\in\calK\vup$, 
$$
\qquad
T^{\stv^-}_k\calK\vup=T^{(\Phi,\stv^-)}_k\calD^-\qquad\hbox{ and }\qquad 
T^{\ssv^-}_{\sth^-(k)}\calK_v=T^{(\fii_v,\ssv^-)}_{\sth^-(k)}\calV^-\qquad,
$$ 
therefore the source submersivity condition on $(\Phi,\nsth^-)$ is equivalent to the
submersivity condition on the double source map of (\ref{biker}).
\end{proof}
%
%
%
%
%
%
%
%
%
%
%
%
%
%
%
%
%
%
%
%
%
%
%
%
%
%
\subsection{\tla-groupoids}\hfill

\vs{0.1}
\spa \tla-groupoids are groupoid objects in the category of Lie algebroids and 
the first order infinitesimal invariant of double Lie groupoids.
\begin{definition}\cite{mkz92}
An $\mathcal{LA}$-\tsf{groupoid} $\sf{\Ohm}:=(\Ohm,A,\calG,M)$
\bgn{equation}\label{typlagpd}
\bcat
\xy
*+{}="0",    <-0.7cm,0.7cm>
*+{\Omega}="1", <0.7cm,0.7cm>
*+{A}="3", <-0.7cm,-0.7cm>
*+{\calG}="2", <0.7cm,-0.7cm>
*+{M}="4",
\ar  @ <-0.07cm>   @{->} "1";"3"^{\Hat{}} 
\ar  @ <0.07cm>    @{->} "1";"3"^{}  
\ar  		   @{->} "1";"2"^{}
\ar                @{->} "3";"4"_{}
\ar  @ <0.07cm>    @{->} "2";"4"^{}
\ar  @ <-0.07cm>   @{->} "2";"4"_{}
\endxy
\ecat
\end{equation}
is a groupoid object in the category of Lie algebroids, such that  $\poidd{\Omega}{A}$ and
$\poidd{\calG}{M}$ are Lie groupoids and  the double source map 
$$
\$\doteq(\sh,\rm{Pr}):
\Omega\rightarrow A\fib{\rm{pr}}{\sor}\calG
$$
is surjective\fn{Then it is also submersive, as it shall be clear from example 
\ref{bundlemaps1} below.}.
\end{definition}
We shall refer to the Lie groupoid and Lie algebroid structures of an \tla-groupoid as top and
side, with the obvious meaning and  denote with $(\sh,\th,\epsh,\iotah,\muh)$ the top groupoid
structural maps of a typical \tla-groupoid such as (\ref{typlagpd}). Note that an \tla-groupoid
is a double groupoid for the abelian groupoids on the side Lie algebroids and the double source
map $\$$ should be understood with respect to this structure.\\
Before explaining the surjectivity condition on $\$$, which ensures that the domain of the top
multiplication be a Lie algebroid indeed, let us introduce the infinitesimal analog of
\tlg-fibrations. 
\bgn{definition} A morphism of Lie algebroids $\phi:A'\to A$ over
$f:M'\to M$ is called an \tsf{\tla-fibration} if
\medskip\\
$1$. It is submersive and base submersive;
\medskip\\
$2$. The characteristic map
$
(\phi,\pr'):A'\longrightarrow A\fib{\pr}{f}M'
$
is surjective, equivalently $\phi$ is fibrewise surjective.
\medskip\\
A \tsf{strong \tla-fibration} is an \tla-fibration which is also base surjective, hence
surjective.
\end{definition}
A strong \tla-fibration is a fibration of Lie algebroids in the sense of Higgins and Mackenzie
\cite{hm90a,mackbook}. For any vector bundles $E$ and $E'$, we shall say that a
vector bundle map $E'\to E$  satisfying the conditions above is a (\tsf{strong})
\tsf{\tvb-fibration};  such maps are  \tla-fibrations for the zero algebroid structures.\\
Let us describe an obvious
example for further reference.
\bgn{remark}\label{bundlemaps1} A base submersive and fibrewise surjective vector bundle map $\phi$
over $f$ is a \tvb-fibration,  since for any choice of trivializations the Jacobian $\sf J(\phi)$
has  the form 
\bgn{equation}\label{reppo}
\sf J(\phi)= \left(\barr{cc}
\mathbb{F}&*\\ 0&\sf J (f)\\ \earr\right) 
\end{equation}  
for the matrix $\mathbb{F}$ representing $\sf \Phi$  and therefore has maximal rank. From 
(\ref{reppo}) is clear that a \tvb-fibration is also fibrewise submersive. 
\end{remark}
The source surjectivity condition on the double source map of a \tla-groupoid is  equivalent to asking the top source
map $\sh$  to be fibrewise surjective; since it is also base submersive and base surjective,  the top source map is a
strong \tla-fibration and so is the top target map $\th=\sh\comp\iotah$. In fact the following are equivalent:
\medskip\\
$i$) The double source map 
$
\$\doteq:\Omega\rightarrow A\fib{\rm{pr}}{\sor}\calG
$
is  surjective,
\medskip\\
$ii$) The top source map $\sh:\Ohm\to A$ is a strong \tla-fibration,
\medskip\\
$iii$) The top vector bundle projection $\Pr:\Ohm\to \calG$ is a strong \tlg-fibration.
\medskip\\
It is then clear that condition ($ii$) above makes the domain of the top multiplication
$\Ohm\fib{\sh}{\th}\Ohm$ of an \tla-groupoid (\ref{typlagpd}) a Lie algebroid.
\bgn{example} \tsf{\tla-groupoids}%
\medskip\\
$i)$ Any Lie groupoid $\gpdm$ is an \tla-groupoid for the rank zero algebroids on $\calG$ and
$M$. Any Lie algebroid $A\to M$ is an \tla-groupoid for the trivial Lie groupoid on $A$ and $M$.%
\medskip\\
$ii)$ The prototypical example of an \tla-groupoid is the \tsf{tangent prolongation
\tla-groupoid} 
$$
\bcat
\xy
*+{}="0",    <-0.7cm,0.7cm>
*+{T\calG}="1", <0.7cm,0.7cm>
*+{TM}="2", <-0.7cm,-0.7cm>
*+{\calG}="3", <0.7cm,-0.7cm>
*+{M}="4",
\ar  @ <-0.07cm>   @{->} "1";"2"_{} 
\ar  @ <0.07cm>    @{->} "1";"2"^{}  
\ar  		   @{->} "1";"3"^{}
\ar                @{->} "2";"4"_{}
\ar  @ <0.07cm>    @{->} "3";"4"^{}
\ar  @ <-0.07cm>   @{->} "3";"4"_{}
\endxy
\ecat
$$
of a Lie groupoid $\gpdm$. Tangent maps are always morphisms of Lie algebroids, then the tangent
prolongation is a groupoid object. It is sufficient to check that the double source map is a
surjective submersion; surjectivity is clear (pick a bisection) and submersivity follows from
example (\ref{bundlemaps1}).
\medskip\\
$iii)$ Any double vector bundle is an \tla-groupoid for the horizontal
abelian group\-oids and the vertical algebroids with zero anchor and bracket. In,
general, an \tla-groupoid is a double Lie groupoid, replacing the Lie algebroid
structures with abelian groupoids.%
\medskip\\
$iv)$ Consider example ($ii$) in the special case of a Lie group $G$; the Lie algebroid $TG\to G$ of $\poidd{\Pi(G)}{G}$ is 
an \tla-\emph{group}, for the tangent group on the top vertical side and the tangent Lie algebroid on the top vertical side.
\end{example}

To see that applying the Lie functor to a double Lie groupoid yields an \tla-groupoid, consider
the following obvious lemma.
\bgn{lemma} Let $\fii:\calG'\to\calG$ over $f:M'\to M$ be an \tlg-fibration, then the induced
morphism of Lie algebroids $\phi:A'\to A$ is an \tla-fibration.
\end{lemma}
%
%
Given any double Lie groupoid (\ref{dbgpd}) the top horizontal source map $\nsth:\calD\to\calV$
differentiates to an \tla-fibration $\sh:A\vup(\calD)\to V$ from the Lie algebroid $A\vup(\calD)$
of $\poidd{\calD}{\calH}$ to the Lie algebroid $V$ of $\calV$. Then differentiating the diagrams
for the top horizontal groupoid $\poidd{\calD}{\calV}$, yields the Lie groupoid  
$\poidd{A\vup(\calD)}{V}$, which is also a groupoid object in the the category of Lie
algebroids with surjective double source map.

\spa Also morphisms of \tla-groupoids are defined in the obvious way.
\bgn{definition}\label{mlag} Let $\sf{\Ohm}^\pm=(\Ohm^\pm,\calG^\pm,A^\pm,M^\pm)$  be 
\tla-groupoids; a \tsf{morphism of \tla-groupoids}  $\sf\Phi:\sf{\Ohm}^-\to\sf{\Ohm}^+$ 
\bgn{equation}\label{morphismlag}
\bcat
\xy
*+{}="0",    <-0.7cm,0.7cm>
*+{\Omega^-}="1", <0.7cm,0.7cm>
*+{A^-}="2", <-0.7cm,-0.7cm>
*+{\calG^-}="3", <0.7cm,-0.7cm>
*+{M^-}="4",     <2.8cm,-0.7cm>
*+{\Omega^+}="1'", <4.2cm,-0.7cm>
*+{A^+}="2'", <2.8cm,-2.1cm>
*+{\calG^+}="3'", <4.2cm,-2.1cm>
*+{M^+}="4'",
\ar  @ <-0.07cm>   @{->} "1";"2"^{} 
\ar  @ <0.07cm>    @{->} "1";"2"^{}  
\ar  		   @{->} "1";"3"^{}
\ar                @{->} "2";"4"_{}
\ar  @ <0.07cm>    @{->} "3";"4"^{}
\ar  @ <-0.07cm>   @{->} "3";"4"_{}
\ar  @ <-0.07cm>   @{->} "1'";"2'"^{} 
\ar  @ <0.07cm>    @{->} "1'";"2'"^{}  
\ar  		   @{->} "1'";"3'"^{}
\ar                @{->} "2'";"4'"_{}
\ar  @ <0.07cm>    @{->} "3'";"4'"^{}
\ar  @ <-0.07cm>   @{->} "3'";"4'"_{}
\ar  		   @{->} "1";"1'"^{\hat{\phi\,}}
\ar  		   @{->} "2";"2'"^{\phi}
\ar  		   @{->} "3";"3'"^{\fii}
\ar  		   @{->} "4";"4'"^{f}
\endxy
\ecat
\end{equation}
is given by four maps  $(\hat{\phi\,},\fii,\phi,f)$ such that the vertical faces of the diagram
above are morphisms of Lie algebroids, while the horizontal faces are morphisms of Lie
algebroids. 
\end{definition}
The same remarks as after definition (\ref{mdlg}) apply, mutatis mutandis, and we shall use the
analogous nomenclature for the components of a morphism of \tla-groupoids, as for morphisms of
double Lie groupoids. 
\bgn{lemma}\label{kerla} For any morphism $(\ref{morphismlag})$ of \tla-groupoids, such that both
$\hat{\phi\:}$ and $\fii$ are \tla-fibrations, respectively over
$\fii$ and over $f$ (so that the kernel Lie algebroids
$\hat{K}:=\ker\hat{\phi}\rightarrow\calG$ and $K:=\ker\phi\rightarrow M$ exist),
\medskip\\
$i$) $\hat{K}$ has a natural  groupoid structure over $K$, making
\bgn{equation}\label{occoc}
\bcat
\xy
*+{}="0",    <-0.7cm,0.7cm>
*+{\hat{K}}="1", <0.7cm,0.7cm>
*+{K}="2", <-0.7cm,-0.7cm>
*+{\calG}="3", <0.7cm,-0.7cm>
*+{M}="4",
\ar  @ <-0.07cm>   @{->} "1";"2"^{} 
\ar  @ <0.07cm>    @{->} "1";"2"^{}  
\ar  		   @{->} "1";"3"^{}
\ar                @{->} "2";"4"_{}
\ar  @ <0.07cm>    @{->} "3";"4"^{}
\ar  @ <-0.07cm>   @{->} "3";"4"_{}
\endxy
\ecat
\end{equation}
a sub-groupoid-object of $\sf{\Ohm}^-$;
\medskip\\
$ii)$ $\poidd{\hat{K}}{K}$ is a Lie groupoid if{f}  
\bgn{equation}\label{ss}
(\hat{\phi\:},\sh^-):\Ohm^-\rightarrow\Ohm^+\fib{\sh^+}{\phi} A^-
\end{equation}
is submersive;
\medskip\\
$iii$) In addition, $(\ref{occoc})$ is an \tla-groupoid if{f}  $(\ref{ss})$ is surjective. 
\end{lemma}
\bgn{proof} The proof of ($i$) is a simple 
exercise in diagram chasing and the proof of ($ii$) is the same as the analogous
statement in lemma (\ref{kerdlg}). ($iii$) By applying the snake lemma fibrewise to the diagram
$$
\xymatrix{
0\:\ar@{->}[r]&\:
\hat{K}\:\ar@{->}[r]\ar@{->}[d]&\:
\Ohm^-\ar@{->}^{\hat{\phi}}[r]\ar@{->}_{\sh^-}[d]\:&\:
\Ohm^+\:\ar@{->}[r]\ar@{->}_{\sh^+}[d]&
0\\
0\:\ar@{->}[r]&\:
K\:\ar@{->}[r]&\:
A^-\ar@{->}^{\phi}[r]\:&\:
A^+\:\ar@{->}[r]&
0\\
}
$$ 
the restriction of the top source map of
$\Ohm^-$ to $\hat{K}$ is seen to be fibrewise surjective onto $K$ if{f}  the restriction of $\hat{\phi}$ to
$\ker\sh^-$ is fibrewise surjective onto $\ker\sh^+$; last condition is equivalent to surjectivity of 
\ref{occoc} by right translation in the abelian groupoids. 
\end{proof}
\vs{1}
\section{Integrability of \tla-groupoids}\label{ilag}
\begin{quotation}
In this Section we develop our functorial approach \cite{07a} to the integrability of \tla-groupoids. In the first
part we obtain the natural integrability result for fibred products of Lie algebroids to fibred products of Lie
groupoids  and moreover we produce necessary and sufficient conditions for an integrating fibred product to be
source  (1-)connected. This allows us to derive easily integrability conditions for \tla-groupoids and morphisms of
\tla-groupoids  by diagrammatics in the second part; a major part of the technical work required to this aim has
already been done in the first Section of this Chapter. The conditions presented here are slightly more general than
those obtained in \cite{07a}. 
\end{quotation}
\spa Unlike the case of Lie groupoids, the vertically source connected component of a double Lie groupoid need not be
a double Lie subgroupoid and, even if it is, it  might not have a vertically source 1-connected cover. This problem
makes the integration of \tla-groupoids no straightforward corollary to the integration of Lie algebroids; namely, 
the top groupoid of an \tla-groupoid does not, in general, induce a groupoid on  the Weinstein groupoids of its
vertical Lie algebroids. In fact there are examples of \tla-groupoids which are integrable to double Lie groupoids
and do not admit any vertically source 1-connected integration  (see example \ref{patologo}).\\ 
We provide below (theorem \ref{intla}) a criterion, which is often computable in the  examples, for the integration 
of an \tla-groupoid $\sf{\Ohm}$, with integrable top and side Lie algebroids, to a vertically source 1-connected double Lie groupoid,
depending on the infinitesimal data only. Our conditions can be understood as \tla-homotopy lifting properties for
the top source and target map of $\sf{\Ohm}$.\\  
The idea is the following. Source, target and unit section of the top
groupoid of an \tla-groupoid are integrable to morphisms of Lie groupoids making the graph compatibility diagrams
(\ref{graphc}) commute. Moreover, the regularity condition on the double source map of $\sf\Ohm$ can also be
``integrated'' to the submersivity condition on the double source of the resulting graph $\sf{\Gamma}$, 
whose (vertical) nerves are
then Lie groupoids integrating the top vertical nerves of $\sf{\Ohm}$. To obtain a double Lie groupoid, it is then
sufficient to integrate the top multiplication of $\sf{\Ohm}$ to a compatible multiplication for $\sf{\Gamma}$; 
this can be done provided suitable vertical nerves of $\sf{\Gamma}$ are source 1-connected.
%
%
%
\subsection{Integrability of fibred products}\hfill

\vs{0.1}
\spa Under the natural transversality conditions, fibred products of integrable Lie algebroids integrate to fibred
products of Lie groupoids. Next result is a corollary of theorem \ref{fibgpd}.
\bgn{theorem} Let $\phi_{1,2}:A^{1,2}\to A$ be transversal morphisms of  Lie algebroids over
$f_{1,2}:M^{1,2}\to M$ (so that  the fibred product Lie algebroid
$
A^{12}:=
A^1\fib{\fii_1}{\fii_2}A^2
\to
M^{12}:=
M^1\fib{\fii_1}{\fii_2}M^2
$
exists). Assume that $A^{1,2}$, $A$ are integrable; then any integrations 
$\fii_{1,2}:\calG^{1,2}\to\calG$ of $\phi_{1,2}$, are transversal morphisms of Lie groupoids and
the fibred product Lie groupoid $\poidd{\calG^1\fib{\fii_1}{\fii}\calG^2}{M^1\fib{f_1}{f_2}M^2}$
exists.  
\end{theorem}
\bgn{proof} In the long exact sequence (\ref{snakefib}) for the fibred product of the graphs
underlying $\calG^{12}$ the cokernels of $\dd f_1 -\dd f_2$ and
$(\dd\fii_1-\dd\fii_2)|_{T^{\sor_1}\calG^1\oplus T^{\sor_2}\calG^2}$ vanish everywhere, the first
one by hypothesis, the second by right translation and the linear transversality condition
(\ref{sourcetransagd}); then the cokernel of $\dd\fii_1-\dd\fii_2$ also vanishes everywhere.
\end{proof}
Note that a fibred product $\calG^1\fib{\fii_1}{\fii_2}\calG^2$ of source (1-)connected Lie groupoids 
might fail to be source (1-)connected; however, it is possible to encode the source (1-)connectivity 
of a fibred product in terms of the induced morphisms of Lie algebroids $\phi_{1,2}$.
\bgn{theorem}\label{fibbiaconnessa} Let  $\calG^{1,2}$ and $\calG$ be source 1-connected Lie
groupoids and  $\fii_{1,2}:\calG^{1,2}\to \calG$ be morphisms of Lie groupoids over
$f_{1,2}:M^{1,2}\to M$ such that the fibred product Lie groupoid 
$\poidd{\calG^1\fib{\fii_1}{\fii_2}\calG^2}{M^1\fib{f_1}{f_2}M^2}$ exists; denote with
$\phi_{1,2}:A^{1,2}\to A$ the induced morphisms of Lie algebroids. Then  
$\poidd{\calG^1\fib{\fii_1}{\fii_2}\calG^2}{M^1\fib{f_1}{f_2}M^2}$ is source
connected if{f} 
\medskip\\
$0$. For any $A^{1,2}$-paths $\alpha_{1,2}^-$ such that
$\phi_1\comp\alpha_1^-$ is $A$-homotopic to $\phi_2\comp\alpha_2^-$, there exist $A^{1,2}$-paths
$\alpha_{1,2}^+$, which are  $A^{1,2}$-homotopic to $\alpha_{1,2}^-$, such that
$\phi_1\comp\alpha_1^-=\phi_2\comp\alpha_2^-$;
\medskip\\
furthermore it is source 1-connected if{f} 
\medskip\\
$1$. For any $A^{1,2}$-paths $\alpha_{1,2}$, which are $A^{1,2}$-homotopic to the constant
$A^{1,2}$-paths 
$\alpha_{1,2}^o\equiv 0^{A^{1,2}}_{\pr_{1,2}(\alpha_{1,2}(0))}$
and such that
$\phi_1\comp\alpha_1=\phi_2\comp\alpha_2$, there exist $A^{1,2}$-homotopies $h_{1,2}$ from
$\alpha_{1,2}$ to the constant $A^{1,2}$-paths 
$\alpha_{1,2}^o$, 
such that
$\phi_1\comp h_1=\phi_2\comp h_2$.
\end{theorem}
\bgn{proof} First of all note that for any $\calG^{1,2}$-paths $\gamma^-_{1,2}$ and $u\in I$
\be
\delta_r(\fii_{1,2}\comp\gamma_{1,2})(u)
&=&
\dd r\inverse_{\fii_{1,2}(g_{1,2}(u))}\dd\fii_{1,2}\dot{\gamma}_{1,2}(u)\quad=\quad\dd\fii_{1,2}\dd
r\inverse_{\gamma_{1,2}(u)}\dot{\gamma}_{1,2}(u)\\
&=&(\phi_{1,2}\comp\delta_r\gamma_{1,2})(u)\qquad\qquad.
\ee
Assume $\calG^{12}$ is source connected, let $\alpha_{1,2}^-$ be $A^{1,2}$-paths such that $\phi_1\comp\alpha_1^-$ is
$A$-homotopic to $\phi_2\comp\alpha_2^-$ and denote with $\gamma^-_{1,2}$ the unique corresponding
$\calG^{1,2}$-paths; we have $\delta_r(\fii_{1,2}\comp\gamma^-_{1,2})=\phi_{1,2}\comp\alpha_{1,2}^-$, thus
$\fii_{1,2}\comp\gamma^-_{1,2}$ are $\calG$-homotopic $\calG$-paths and 
$\fii_1(\gamma^-_1(1))=\fii_2(\gamma^-_2(1))$. Since $\calG^{12}$ is source connected, one can find
$\calG^{1,2}$-paths $\gamma^+_{1,2}$ such that $\gamma^+_{1,2}(1)=\gamma^-_{1,2}(1)$, 
(i.e. $\calG^{1,2}$-homotopic to $\gamma^-_{1,2}$) with $\fii_1\comp\gamma_1^+=\fii_2\comp\gamma_2^+$; 
the unique corresponding $A^{1,2}$-paths
$\alpha_{1,2}^+:=\delta_r\gamma_{1,2}^+$ are then homotopic to $\alpha_{1,2}^-$ and satisfy
$\phi_1\comp\alpha_1^-=\phi_2\comp\alpha_2^-$.\\
Conversely, assume that $(0.)$ holds, let $(g_1,g_2)\in\calG^{12}$ and pick $A^{1,2}$-paths $\alpha_{1,2}^-$ with
$A^{1,2}$-homotopy classes $[\alpha_{1,2}^-]=g_{1,2}$. Then
$
[\phi_1\comp\alpha^-_1]=\fii_1(g_1)=\fii_2(g_2)=[\phi_2\comp\alpha^-_2]
$, i.e. $\phi_1\comp\alpha_1$ is  $A$-homotopic to $\phi_2\comp\alpha_2$ and one can find $A^{1,2}$-paths
$\alpha_{1,2}^+$, which are $A^{1,2}$-homotopic to $\alpha_{1,2}^-$ and satisfy
$\phi_1\comp\alpha_1^-=\phi_2\comp\alpha_2^-$. The unique corresponding $\calG^{1,2}$-paths $\gamma_{1,2}$ satisfy
$\gamma_{1,2}(1)=g_{1,2}$ and  $\delta_r(\fii_1\comp\gamma_1)=\delta_r(\fii_2\comp\gamma_2)$, therefore 
$\fii_1\comp\gamma_1=\fii_2\comp\gamma_2$, by uniqueness, and the pair $(\gamma_1,\gamma_2)$ form a $\calG^{12}$-path
starting from the unit section and reaching $(g_1,g_2)$; that is, $\calG^{12}$ is source connected.\\
Assume now that $\calG^{12}$ is source 1-connected and let $\alpha_{1,2}$ be $A^{1,2}$-paths, which are
$A^{1,2}$-homotopic to the constant $A^{1,2}$-paths $\alpha_{1,2}^o$, with 
$\phi_1\comp\alpha_1=\phi_2\comp\alpha_2$. The unique corresponding $\calG^{1,2}$-paths are in fact loops
$\lambda_{1,2}$ in the source fibres of $\calG^{1,2}$ starting form  $\eps_{1,2}(\pr_{1,2}(\alpha_{1,2}(0)))$, such
that $\fii_1\comp\lambda_1=\fii_2\comp\lambda_2$, i.e. the pair $(\lambda_1,\lambda_2)$ forms a $\calG^{12}$-loop.
Since $\calG^{12}$ is source 1-connected, one can find homotopies $H^{1,2}$ within the source fibres of 
$\pr_{1,2}(\alpha_{1,2}(0)$ from $\lambda_{1,2}$ to the constant $\calG^{1,2}$ loops
$\lambda_{1,2}^o\equiv\eps_{1,2}(\pr_{1,2}(\alpha_{1,2}(0)))$ such that $\fii_1\comp H^1=\fii_2\comp H^2$. The unique corresponding $A^{1,2}$-homotopies
$h^{1,2}$ satisfy the boundary conditions
\be
\pmh h^{1,2} (u)&=&\delta_rH^{1,2}(u,0)=\delta_r\lambda_{1,2}(u)=\alpha_{1,2}\\
\pph h^{1,2} (u)&=&\delta_rH^{1,2}(u,1)=\delta_r\lambda_{1,2}^o(u)=0\qquad,
\ee
thus are $A^{1,2}$-homotopies from  $\alpha_{1,2}$ to the constant $A^{1,2}$-paths
$\alpha_{1,2}^o$; moreover
\be
\phi_{1,2}\comp h^{1,2}
&=&
\phi_{1,2}\delta_r^u H^{1,2}(u, \eps)\cdot\dd u +\phi_{1,2}\delta_r^\eps H^{1,2}(u,\eps)\cdot\dd\eps\\
&=&
\dd r\inverse_{\fii_{1,2}(H^{1,2}(u,\eps))}\dd\fii_{1,2}\frac{\pa}{\pa u}H^{1,2}(u, \eps)\cdot\dd
u\\ 
&+&
\dd r\inverse_{\fii_{1,2}(H^{1,2}(u,\eps))}\dd\fii_{1,2}\frac{\pa}{\pa \eps}H^{1,2}(u, \eps)\dd\eps\\
&=&
\delta_r^u (\fii_{1,2}\comp H^{1,2})(u, \eps)\cdot\dd u 
+
\delta_r^\eps (\fii_{1,2}\comp H^{1,2})(u,\eps)\cdot\dd\eps\qquad,
\ee
for the partial right derivatives $\delta_r^{\eps,u}$, i.e. $\phi_1\comp h^1=\phi_2\comp h^2$,\\
Conversely, assume that $(1.)$ holds and $\calG^{1,2}$-loops $\lambda_{1,2}$ such that
$\fii_1\comp\lambda_1=\fii_2\comp\lambda_2$ are  assigned; we have to find homotopies $H^{1,2}$
from $\lambda_{1,2}$ to the constant loops $\lambda_{1,2}^o$ within the source fibres, such that
$\fii_1\comp H^1=\fii_2\comp H^2$.   The unique $A^{1,2}$-paths $\alpha_{1,2}$ corresponding to
$\lambda_{1,2}$ are $A^{1,2}$-homotopic to the constant $A^{1,2}$-paths $\alpha_{1,2}^o$, since
$\lambda_{1,2}$ are $\calG^{1,2}$-homotopic to the constant $\calG^{1,2}$-loops
$\lambda_{1,2}^o$; moreover $\phi_1\comp\alpha_1=\phi_2\comp\alpha_2$, thus we can find
$A^{1,2}$-homotopies $h^{1,2}$ from $\alpha_{1,2}$ to $\alpha_{1,2}^o$ with 
$\phi_1\comp h_1=\phi_2\comp h_2$, which we shall regard as morphisms of Lie algebroids 
$h^{1,2}:TI^{\times 2}\to A^{1,2}$. The unique corresponding $\calG^{1,2}$-homotopies $H^{1,2}$
are given by $H^{1,2}(u,\eps)=\wt{h}^{1,2}(u,\eps;0,0)$, where 
$\wt{h}^{1,2}:I^{\times 2}\times I^{\times 2}\to\calG^{1,2}$ are the unique morphisms of Lie 
groupoids integrating  $h^{1,2}$. Then $\fii_{1,2}\comp\wt{h}^{1,2}$ are the unique morphisms of
Lie groupoids integrating $\phi_{1,2}\comp h^{1,2}$ and $\fii_1 \comp H^1=\fii_2 \comp H^2$.
\end{proof}
Last result motivates the following
\bgn{definition} Two morphisms of Lie algebroids $\phi_{1,2}:A^{1,2}\to A$ are \tsf{strongly transversal}
if{f} they are transversal morphisms of Lie algebroids and the lifting conditions (0., 1.) of theorem \ref{fibbiaconnessa} hold.
\end{definition}
It straightforward to translate conditions 0. and 1. above to the case of the kernel groupoid of 
of a morphism of Lie groupoids; it turns out that condition 1. is equivalent to the vanishing of
suitable loop groups. Let $\phi:A\to A'$ be a morphism of Lie algebroids over $f:M\to M'$ and
assume $\ker\phi\subset A$ is a Lie subalgebroid. For any $q\in M$, consider the space of class
$\cif$ $A$-loops based in $q$, taking values in $\ker\phi$ and $A$-homotopic to the null constant path,
 modulo $\cif$ $A$-homotopy 
taking values in $\ker\phi$, namely $\cif$
$\ker\phi$-loops based in $q$ modulo $\cif$ $\ker\phi$-homotopy:
$$
\qquad
\mathbb{K}_q(\phi):=\{\alpha\in \ker\phi\hbox{\,-\,}paths\:|\:\alpha(0)=\alpha(1)=0_q,\,\alpha\thicksim_A\alpha^o\}/\ker\phi
\hbox{\,-\,}homotopy
\qquad
$$
Note that, since $\ker\phi\subset A$ is a Lie subalgebroid the composition of smooth
$\ker\phi$-paths is well defined up to smooth $\ker\phi$-homotopy and induces a
group multiplication on the loop spaces above (we shall not need this group multiplication 
here or in the following).
\bgn{corollary}\label{lupo} Let $\fii:\calG\to\calG'$ be a morphism of Lie groupoids over $f:M'\to M$ such
that the kernel groupoid $\ker\fii\subset\calG'$ is a Lie subgroupoid and denote with $\phi:A\to
A'$ the induced morphism of Lie algebroids. If $\calG$ and $\calG'$ are source 1-connected, 
then $\ker\fii$ is source connected if{f} 
\medskip\\
$0$. For any $A$-path $\alpha^-$ such that
$\phi\comp\alpha^-$ is $A'$-homotopic to  the constant
$A'$-path
$(\phi\comp\alpha^-)^o\equiv 0_{\pr'(\phi(\alpha^-(0))}$, 
there exists an  $A$-path
$\alpha^+$, which is  $A$-homotopic to $\alpha^-$, such that
$\phi\comp\alpha^+=(\phi\comp\alpha^-)^o$;
\medskip\\
furthermore it is source 1-connected if{f} 
\medskip\\ 
$1$. For any $A$-path $\alpha$, which is $A$-homotopic to the constant $A$-path 
$\alpha^o\equiv 0_{\pr(\alpha(0))}$ and such that $\phi\comp\alpha=(\phi\comp\alpha)^o$
is the constant $A'$-path, there exists an  $A$-homotopy $h$ form $\alpha$ to the constant
$A$-path  $\alpha^o\equiv 0_{\pr(\alpha(0))}$, such that $\phi\comp h=0$
\medskip\\
equivalently if{f} 
\medskip\\
$1'$. The loop groups $\mathbb{K}_q(\phi)$ are trivial for all $q\in M$.
\end{corollary}
\bgn{proof} Set $K:=\ker\phi$. Conditions (0.) and (1.) are just restatements of the corresponding  conditions in theorem
\ref{fibbiaconnessa} for the case $\fii_1=\fii:\calG\to\calG'$ and $\fii_2=\eps': M'\to\eps'(M')$. The $K$-loops
representing elements of  $\mathbb{K}_q(\phi)$ are in particular $A$-paths satisfying the hypothesis of (1.), thus
(1.) implies $(1'.)$. For the opposite implication, consider that any $K$-path is $K$-homotopic to a smooth $K$-paths
and reparametrization of smooth $K$-paths does not change the $K$-homotopy class. Thus condition $1.$ holds, 
for all smooth $\alpha$ with compact support, which in particular represent elements of $\mathbb{K}_q(\phi)$.
\end{proof}
We conclude this Subsection introducing a special class of morphisms of Lie algebroids.
\bgn{definition}\label{liftingprop}\cite{07a} A morphism $\phi:A\to B$ of Lie algebroids has the
\medskip\\
$l_0$) \tsf{0-\tla-homotopy lifting property} if, for any $A$-path $\alpha_-$ and $B$-path
$\beta_+$,  which is $B$-homotopic to $\beta_-:=\phi\comp\alpha_-$,  there exists an $A$-path
$\alpha_-$, which is $A$-homotopic to $\alpha_+$ and satisfies  $\phi\comp\alpha_+=\beta_+$;
\medskip\\ 
$l_1$) \tsf{1-\tla-homotopy lifting property} if,  for any $A$-path $\alpha$, which is
$A$-homotopic to the constant $A$-path $\alpha_o\equiv 0_{\pr_A(\alpha(0))}$,  and $B$-homotopy $h$
from  $\beta:=\phi\comp\alpha$ to the constant $B$-path  $\beta_o\equiv 0_{\pr_B(\beta(0))}$, there
exists an $A$-homotopy $\hat{h}$ from $\alpha$ to $\alpha_o$,  such that $\phi\comp\hat{h}=h$.
\end{definition}
The \tla-homotopy lifting conditions above on a morphism of integrable Lie
algebroids  $A\to B$, translate to the infinitesimal level path lifting conditions
in the source fibres of the source 1-connected Lie groupoids $\calA$ and $\calB$ of
$A$ and $B$ along the integration $\calA\to\calB$; therefore, with the same notations and
assumptions of theorem \ref{fibbiaconnessa} we have
\bgn{lemma}\label{intelift} The fibred product $\poidd{\calG^1\fib{\fii_1}{\fii_2}\calG^2}{M^1\fib{f_1}{f_2}M^2}$
is%
\medskip\\
$i$) source connected if $\phi_1:A^1\to A$ has the $0$-\tla-homotopy lifting property
\medskip\\
furthermore, it is%
\medskip\\
$ii$) source 1-connected if $\phi_1:A^1\to A$ has the $1$-\tla-homotopy lifting property
\end{lemma}
Clearly, he statement holds true replacing $\phi_1$ with $\phi_2$.
\bgn{proof} Condition 0. of theorem \ref{fibbiaconnessa} follows from ($i$) since one can set
$\alpha_2^+=\alpha_2^-$ and lift $\phi_2\comp\alpha_2^+$ along $\phi_1$. Similarly, condition
1. of   \ref{fibbiaconnessa} follows from ($ii$).
\end{proof}
\vs{0.5}
\subsection{Integrability of \tla-groupoids}\hfill

\vs{0.1}
\spa Here and in the following we shall repeatedly use of next easy
\bgn{lemma}\label{fibrel} Let $\fii:\calG'\to\calG$ be a morphism of Lie groupoids over $f:M'\to M$
integrating $\phi:A'\to A$. Then  $\fii$ is an \tlg-fibration if{f}  $\phi$ is an \tla-fibration.
\end{lemma} 
\bgn{proof} By right translation, fibrewise surjectivity of $\phi$ is equivalent to source
submersivity of $\fii$. Then the implication to the right is clear and the implication to the left
follows by submersivity of $f$ and the 5-lemma.
\end{proof}
Let us consider the integrability of the graph underlying an \tla-groupoid. We shall say that a differentiable graph 
$(\Gamma, M;\sor,\tar)$
is an \tsf{invertible graph} if its total space is endowed with an automorphism $\iota$ such that $\iota^2=1$ and 
$\tar\comp\iota=\sor$.
\bgn{proposition}\label{integraph} The top differentiable graph of any \tla-groupoid with integrable
top Lie algebroid integrates to an invertible differentiable graph in the category of source
1-connected  Lie groupoids.
\end{proposition}
\bgn{proof} Since $A$ is a Lie subalgebroid of $\Ohm$ and $\Ohm$ is integrable (theorem \ref{1}), $A$ is also
integrable. Denote with $\calA$ and $\Gamma$ the source 1-connected integrations of $A$ and $\Ohm$,
respectively. The graph compatibility diagrams
\bgn{equation}\label{graphc}
\bcat
\xy
*+{}="0",    <-1cm,0cm>
*+{A}="1", <1cm,0cm>
*+{A}="2", <0cm,1.4cm>
*+{\Ohm}="3", 
\ar     @ {->} "1";"2"_{\id_A} 
\ar     @ {->} "1";"3"^{\epsh}  
\ar     @ {->} "3";"2"^{\sh,\th}
\endxy
\ecat
\qquad\qquad\qquad
\bcat
\xy
*+{}="0",    <-1cm,0cm>
*+{\Ohm}="1", <1cm,0cm>
*+{\Ohm}="2", <0cm,1.4cm>
*+{\Ohm}="3", 
\ar     @ {->} "1";"2"_{\id_\Ohm} 
\ar     @ {->} "1";"3"^{\iotah}  
\ar     @ {->} "3";"2"^{\iotah}
\endxy
\ecat
\qquad\qquad\qquad
\bcat
\xy
*+{}="0",    <-1cm,0cm>
*+{\Ohm}="1", <1cm,0cm>
*+{\Ohm}="2", <0cm,1.4cm>
*+{\calA}="3", 
\ar     @ {->} "1";"2"_{\iotah} 
\ar     @ {->} "1";"3"^{\sh}  
\ar     @ {->} "2";"3"_{\th}
\endxy
\ecat
$$
integrate to commuting diagrams
$$
\bcat
\xy
*+{}="0",    <-1cm,0cm>
*+{\calA}="1", <1cm,0cm>
*+{\calA}="2", <0cm,1.4cm>
*+{\Gamma}="3", 
\ar     @ {->} "1";"2"_{\id_\calA} 
\ar     @ {->} "1";"3"^{\eth}  
\ar     @ {->} "3";"2"^{\sth,\tth}
\endxy
\ecat
\qquad\qquad\qquad
\bcat
\xy
*+{}="0",    <-1cm,0cm>
*+{\Gamma}="1", <1cm,0cm>
*+{\Gamma}="2", <0cm,1.4cm>
*+{\Gamma}="3", 
\ar     @ {->} "1";"2"_{\id_\Gamma} 
\ar     @ {->} "1";"3"^{\ith}  
\ar     @ {->} "3";"2"^{\ith}
\endxy
\ecat
\qquad\qquad\qquad
\bcat
\xy
*+{}="0",    <-1cm,0cm>
*+{\Gamma}="1", <1cm,0cm>
*+{\Gamma}="2", <0cm,1.4cm>
*+{\calA}="3", 
\ar     @ {->} "1";"2"_{\ith} 
\ar     @ {->} "1";"3"^{\sth}  
\ar     @ {->} "2";"3"_{\tth}
\endxy
\ecat
\end{equation} 
in the category of Lie groupoids; thus, $\ith$ is a diffeomorphism, being the inverse to itself,
$\eth$ is injective, having left inverses and $\nsth$, $\ntth$ are surjective, being left inverses
to $\eth$. Moreover the tangent diagrams
$$
\bcat
\xy
*+{}="0",    <-1cm,0cm>
*+{T_a\calA}="1", <1cm,0cm>
*+{T_a\calA}="2", <0cm,1.4cm>
*+{T_{\eth(a)}\Gamma}="3"
\ar     @ {->} "1";"2"_{\id_{T_a\calA}} 
\ar     @ {->} "1";"3"^{\dd\eth}  
\ar     @ {->} "3";"2"^{\dd\sth,\dd\tth}
\endxy
\ecat
$$
commute for all $a\in\calA$, therefore $\eth$ is immersive. It remains to show that $\sth$ and
$\tth$ are submersive, which follows from lemma \ref{fibrel}, since both $\sh$ and $\th$ are
\tla-fibrations.
\end{proof}
Consider the top nerves of an \tla-groupoid $\sf{\Ohm}=(\Ohm,\calG;A,M)$
$$
\barr{lcl}
\Ohm^{(0)}&:=&A\\
\Ohm^{(1)}&:=&\Ohm\\
\Ohm^{(n+1)}&:=&\Ohm\fib{\sh}{\th\comp\pp_1}\Ohm^{(n)}
\qquad,\qquad n\geq 1\qquad,
\earr
$$
where $\pp_1:\Ohm^{(n)}\to\Ohm$ denotes the restriction of the first projection $\Ohm^{\times n}\to
\Ohm$. Since the top source map $\sh$ is an \tla-fibration each top nerve $\Ohm^{(n)}$ carries a Lie
algebroid structure over the corresponding side nerve $\calG^{(n)}$ making it a Lie subalgebroid of
the direct product $\Ohm^{\times n}\to\calG^{\times n}$, $n\geq 1$ (theorem \ref{dracula}). Assume now that
the top Lie algebroid of $\sf{\Ohm}$ is integrable and consider the integrating graph 
$(\Gamma,\calA)$ of lemma \ref{integraph} given by the source 1-connected integrations 
$\poidd{\Gamma}{\calG}$ and
$\poidd{\calA}{M}$ of the top and side Lie algebroids $\Ohm\to \calG$ and $A\to M$. Each nerve, 
$$
\barr{lcl}
\Gamma^{(0)}&:=&\calA\\
\Gamma^{(1)}&:=&\Gamma\\
\Gamma^{(n+1)}&:=&\Gamma\fib{\sth}{\tth\comp\pp_1}\Gamma^{(n)}\subset\Gamma^{\times n}
\qquad,\qquad n\geq 1\qquad,
\earr
$$
of $\Gamma$ is a smooth submanifold of the corresponding direct product, where $\nsth,\ntth:\Gamma\to\calA$ are
the integrations of top source and target of $\Gamma$ and  $\pp_1:\Gamma^{(n)}\to\Gamma$ denotes the restriction
of the first projection    $\Gamma^{\times n}\to\Gamma$. Moreover, since \tla-fibrations integrate to 
\tlg-fibrations   (lemma \ref{fibrel}), the transversality conditions of theorem \ref{fibbiaconnessa} are met and 
$\Gamma^{(n)}$ carries a Lie
groupoid structure over the corresponding side nerve $\calG^{(n)}$ making it a Lie subgroupoid of the direct
product $\Gamma^{\times n}\to\calG^{\times n}$, $n\geq 1$. Note that if $\nsth$ has the 0- and 
1-homotopy lifting properties, it follows from lemma \ref{intelift} that \emph{all} the nerves
$\poidd{\Gamma^{(\bullet)}}{\calG^{(\bullet)}}$ are source  1-connected.
\bgn{theorem}\label{intla} Let $\sf{\Ohm}$ be an \tla-groupoid with integrable top Lie algebroid. If source and
target of the top Lie groupoid of $\sf{\Ohm}$ are strongly transversal, there exists a unique vertically source
1-connected double Lie groupoid integrating $\sf\Ohm$.
\end{theorem}
\bgn{remark} In particular, the transversality condition holds when the top source, equivalently target, map
of $\sf{\Ohm}$ enjoys the \tla-homotopy lifting conditions of definition \ref{liftingprop}.
\end{remark}
\bgn{proof}[proof of theorem \ref{intla}] With the same notations as above,  $\Gamma\fib{\sth}{}M\simeq
M\fib{}{\tth}\Gamma\simeq\Delta_\Gamma\simeq\Gamma$ and $\Gamma^{(2)}$ are source 1-connected Lie groupoids, then
the unitality diagrams
$$
\bcat
\xy
*+{}="0",    <-0.7cm,1.7cm>
*+{\Gamma}="tl",  <0.7cm,1.7cm>
*+{\Gamma}="tr",    <-1.2cm,0.5cm>
*+{\Delta_\Gamma}="l", <1.2cm,0.5cm>
*+{\Gamma^{(2)}}="r", <0cm,-0.5cm>	
*+{\Gamma\fib{\sth}{}M}="b",
\ar     @ {->} "tl";"l"_{\Delta_\Gamma}
\ar     @ {->} "r";"tr"_{\mth}  
\ar     @ {->} "l";"b"_{\id_\Gamma\times\sth\quad}  
\ar     @ {->} "b";"r"_{\quad\id_\Gamma\times\eth}
\ar     @ {->} "tl";"tr"^{\id_\Gamma}
\endxy
\ecat
\qquad\qquad
\bcat
\xy
*+{}="0",    <-0.7cm,1.7cm>
*+{\Gamma}="tl",  <0.7cm,1.7cm>
*+{\Gamma}="tr",    <-1.2cm,0.5cm>
*+{\Delta_\Gamma}="l", <1.2cm,0.5cm>
*+{\Gamma^{(2)}}="r", <0cm,-0.5cm>	
*+{M\fib{}{\tth}\Gamma}="b", 
\ar     @ {->} "tl";"l"_{\Delta_\Gamma}
\ar     @ {->} "r";"tr"_{\mth}  
\ar     @ {->} "l";"b"_{\tth\times\id_\Gamma\quad}  
\ar     @ {->} "b";"r"_{\quad\eth\times\id_\Gamma}
\ar     @ {->} "tl";"tr"^{\id_\Gamma}
\endxy	
\ecat
$$
and invertibility diagrams
$$
\bcat
\xy
*+{}="0",    <-0.7cm,1.7cm>
*+{\Gamma}="tl",  <0.7cm,1.7cm>
*+{\calA}="tr",    <-1.2cm,0.5cm>
*+{\Delta_\Gamma}="l", <1.2cm,0.5cm>
*+{\Gamma}="r", <0cm,-0.5cm>	
*+{\Gamma^{(2)}}="b",
\ar     @ {->} "tl";"l"_{\Delta_\Gamma}
\ar     @ {->} "tr";"r"^{\eth}  
\ar     @ {->} "l";"b"_{\ith\times\id_\Gamma}  
\ar     @ {->} "b";"r"_{\quad\mth}
\ar     @ {->} "tl";"tr"^{\sth}
\endxy
\ecat
\qquad\qquad\quad
\bcat
\xy
*+{}="0",    <-0.7cm,1.7cm>
*+{\Gamma}="tl",  <0.7cm,1.7cm>
*+{\calA}="tr",    <-1.2cm,0.5cm>
*+{\Delta_\Gamma}="l", <1.2cm,0.5cm>
*+{\Gamma}="r", <0cm,-0.5cm>	
*+{\Gamma^{(2)}}="b", 
\ar     @ {->} "tl";"l"_{\Delta_\Gamma}
\ar     @ {->} "tr";"r"^{\eth}  
\ar     @ {->} "l";"b"_{\id_\Gamma\times\ith}  
\ar     @ {->} "b";"r"_{\quad\mth}
\ar     @ {->} "tl";"tr"^{\tth}
\endxy
\ecat
$$
commute, since the integrate the corresponding diagrams for the top Lie groupoid of $\sf{\Ohm}$.
The third nerve $\Gamma^{(3)}$ is also source 1-connected, hence the associativity diagram
$$
\bcat
\xy
*+{}="0",    <0cm,1.2cm>
*+{\Gamma^{(3)}}="t", <0cm,-1.2cm>
*+{\Gamma}="b", <-1.2cm,0cm>
*+{\Gamma^{(2)}}="l", <+1.2cm,0cm>
*+{\Gamma^{(2)}}="r",  
\ar     @ {->} "t";"l"_{\id_\Gamma\times\mth} 
\ar     @ {->} "t";"r"^{\mth\times\id_\Gamma}  
\ar     @ {->} "l";"b"_{\mth}
\ar     @ {->} "r";"b"^{\mth}
\endxy
\ecat
$$
commutes. Then $(\Gamma,\calG;\calA,M)$ is a groupoid object in the category of Lie groupoids, thanks to
proposition \ref{integraph}; the source submersivity condition follows from lemma \ref{fibrel}.
\end{proof}
The transversality conditions for integrability imply that all top horizontal nerves of the integrating double Lie
groupoid are source 1-connected; however, only the role of the second and third nerve is essential.
\bgn{remark}\label{cvgt}
In the proof of last theorem, even if $\Gamma^{(3)}$ is only source connected, 
last diagram still commute since
$$
\bcat
\xy
*+{X}="0",    <0cm,1.2cm>
*+{\Gamma^{(3)}}="t", <0cm,-1.2cm>
*+{\Gamma}="b", <-1.2cm,0cm>
*+{\Gamma^{(2)}}="l", <+1.2cm,0cm>
*+{\Gamma^{(2)}}="r",  
\ar     @ {->} "t";"l"_{\id_\Gamma\times\mth} 
\ar     @ {->} "t";"r"^{\mth\times\id_\Gamma}  
\ar     @ {->} "l";"b"_{\mth}
\ar     @ {->} "r";"b"^{\mth}
\ar     @ {->} "0";"t"^{\kappa}
\ar     @ {->} "0";"l"^{\mu_l}
\ar     @ {->} "0";"r"_{\mu_r}
\endxy
\ecat
$$
commutes for the covering morphism $\kappa:X\to\Gamma^{(3)}$ and the morphisms  $\mu_{l,r}:X\to\Gamma^{(2)}$ from the
source 1-connected cover $X$ of $\Gamma^{(3)}$ integrating $\id_\Gamma\times\mth$ and $\mth\times\id_\Gamma$
respectively. Thus to in order to integrate an \tla-groupoid with integrable top algebroid to a double Lie groupoid,
it is sufficient to have $\Gamma^{(2)}$ source 1-connected and $\Gamma^{(3)}$ source connected.
\end{remark}
In principle, an integrable \tla-groupoid could admit a vertically source 1-connected integration which does not
possess source (1-)connected top horizontal nerves; we could not find examples of this kind.  
\bgn{example} Consider the tangent prolongation \tla-groupoid 
$$
\sf{T\TT^2}=(T\TT^2,\TT^2;TS^1,S^1)
$$
of the pair groupoid on the torus $\TT^2\colon\poidd{S^1\times S^1}{S^1}$; then the fundamental
groupoid $\Pi(\TT^2)=(\RR^2\times\RR^2)/\ZZ^2$ carries a further natural Lie groupoid (induced by the
direct  product of pair groupoids $\poidd{\RR^2=\RR\times\RR}{\RR}$) over
$\Pi(S^1)=(\RR\times\RR)/\ZZ$,  making
$
(\Pi(\TT^2),\TT^2;\Pi(S^1),S^1)
$
a vertically source 1-connected Lie groupoid integrating $\sf{T\TT^2}$. The top source map of 
$\sf{T\TT^2}$ is the second projection $TS^1\times TS^1\to TS^1$, which clearly satisfies the
\tla-homotopy lifting conditions $l_{0,1}$.
\end{example}
Let us consider an example of a tangent prolongation \tla-groupoid, which has no vertically 
source 1-connected integration.
\bgn{example}\label{patologo} Consider the Lie groupoid given by $\calG=\RR^2\times\ZZ$,
$M=\RR^2\equiv\{0\}\times\RR^2$, $\sor(k;x,y)=(0, x-k,y)$, $\tar(k;x,y)=(0, x+k,y)$,
$$
\quad
\mu((l;x+k+l,y),(k;x,y))=(l+k;x+l,y)
\quad\hbox{ and } \quad
\iota(k;x,y)=(-k;x,y)
\quad.
$$
Removing lattices from all the layers of $\calG$ but the base yields a Lie subgroupoid 
$\poidd{\calH}{M}$,
$$
\calH:=\RR^2\times\{0\} \hbox{{\small $\coprod$}}
(\underset{k\neq0}{\hbox{{\small $\coprod$}}}\RR^2\backslash\ZZ^2\times \{k\})
$$
such that the differentiable graph $(\Pi(\calH), \Pi(M))$ (source, target, unit section and
inversion being induced by composition) carries no compatible groupoid
multiplication, therefore $(\Pi(\calH),\calH;\Pi(M), M)$ is never a double Lie groupoid. 
To see this, consider the the equivalence classes in $\Pi(\calH)$ of the paths
$$
\qquad
\gamma_\pm(t)=(1;\pm 1 \pm\sin(\pi t)/2,1/2 + t)\qquad,\qquad t\in
I\qquad:
$$
$\sor\comp\gamma_-(t)=(0;-\sin(\pi t)/2,1/2 + t)$ and 
$\tar\comp\gamma_+(t)=(0;\sin(\pi t),1/2 + t)$ define homotopic paths in $M$, but there is no
representative in the homotopy class of $\gamma_+$ which is pointwise  composable with some
representative of the class of $\gamma_-$; the source map of $\Pi(\calH)$ does not satisfy the
\tla-homotopy lifting condition $l_0$.
\end{example}
Integrability of morphisms of \tla-groupoids follows by diagrammatics in the same fashion.
\bgn{corollary} Let $\sf{\Phi}:\sf{\Ohm}^-\to\sf{\Ohm}^+$ be a morphism of \tla-groupoids $(\ref{morphismlag})$.
 Assume that $\sf{\Ohm}^\pm$ are integrable to double Lie groupoids $\sf{\Gamma}^\pm$; moreover assume that top
horizontal source and target of $\Ohm^-$ are strongly transversal. Then there exist a unique morphism of double Lie
groupoids $\wt{\sf{\Phi}}$ integrating $\sf{\Phi}$.   
\end{corollary}
\bgn{proof} The integration $\wt{\Phi}:\Gamma^-\to\Gamma^+$ over $\fii:\calG^-\to\calG_+$ of the top vertical 
component $\Phi:\Ohm^-\to\Ohm^+$ makes the diagram
\bgn{equation}\label{uyt}
\bcat
\xy
*+{}="0",    <-1.7cm,0.7cm>
*+{(\Gamma^-)^{(2)}}="1", <0.7cm,0.7cm>
*+{\Gamma^-}="3", <-1.7cm,-0.7cm>
*+{(\Gamma^+)^{(2)}}="2", <0.7cm,-0.7cm>
*+{\Gamma^+)}="4",
\ar  @ {->} "1";"3"^{\mth^-}  
\ar  @ {->} "1";"2"_{\wt{\Phi}\times\wt{\Phi}}
\ar  @ {->} "3";"4"^{\wt{\Phi}}
\ar  @ {->} "2";"4"_{\mth^+}
\endxy
\ecat
\end{equation}
commute.
\end{proof}
In the same spirit as in remark \ref{cvgt}, we shall comment on the connectivity requirements on the top horizontal
nerves.
\bgn{remark}
In order to integrate a morphism of \tla-groupoids such as above it is sufficient to have the second top horizontal
nerve of $\sf{\Gamma}^-$ source connected. Commutativity of diagram (\ref{uyt}) is induced by that of the similar
diagram for the covering groupoid.
\end{remark}
\vs{1}
\section{Double \!structures, \!duality \!and \!integrability \!of \!Poisson \!groupoids}\label{dsdipg}
\begin{quotation} 
Here we specialize the integrability results of last Section to the \tla-groupoid canonically associated with a
Poisson groupoid. We show that when the Poisson bivector is integrable and the strong transversality conditions on
the associated \tla-groupoid are met, the integration produces a symplectic double groupoid realizing a strong
duality between the given Poisson groupoid and its unique source 1-connected weak dual. Finally we study a class of
examples where the transversality conditions can be computed explicitly and show that all complete Poisson groups are
integrable to symplectic double groupoids.
\end{quotation}

\vs{0.1}
\spa Recall from Section \ref{pglb} that $(A,A^*)$ is a Lie bialgebroid if{f}  so is $(A^*,A)$. Moreover changing the
signs of the anchor and bracket of a Lie algebroid $A$, yields another Lie algebroid $-A$; it is then easy to see
that $(A,A^*)$ is a Lie bialgebroid if{f}  so is the \tsf{flip} $(A^*,-A)$.
\bgn{definition}\label{weakduality} Two Poisson groupoids $\poidd{\calG_\pm}{M}$ with Lie algebroids 
$(A_\pm,A^*_\pm)$ are in \tsf{weak duality} if $(A_-, A^*_-)$ and the flip $(A^*_+,-A_+)$ are isomorphic Lie
bialgebroids\footnote{This notion was introduced by Weinstein in \cite{ws88} simply as ``duality'' for Poisson
groupoids}. 
\end{definition}
This notion of duality is a symmetric relation: if $\phi:A^*_+\to A_-$ induces an isomorphism of Lie bialgebroids 
$(A^*_+,-A_+)\to(A_-, A^*_-)$,  $-\phi^\sf{t}$ induces an isomorphism $(A^*_-,-A_-)\to(A_+, A^*_+)$. However, a
Poisson groupoid might not have any weak dual and weak duals are not unique (for instance $-\phi$ makes 
$(A^*_+,A_+)$ weakly dual to $(-A_-,-A^*_-)$
\bgn{example} For any Poisson manifold $(P,\pi)$, the Lie bialgebroid of $\poidd{\ol{P}\times P}{P}$ is (canonically) 
$(TP,-T^*P)$  and there is no Poisson groupoid integrating its flip $(-T^*P, -TP)$ when $\pi$ is not integrable. If $\pi$ is
integrable  to a symplectic groupoid $\poidd{\Lambda}{P}$, $\ol{\Lambda}$ is the canonical choice of a weak dual to
$\ol{P}\times P$. 
\end{example}

\bgn{remark} The convention of declaring Poisson groupoids $\calG_\pm$ in weak duality if 
$(A_+, A^*_+)\simeq(A^*_-,A_-)$ is also frequent in the literature, especially for Poisson groups. To recover this notion
of duality one has to change the sign of a Poisson bivector, that is $\calG_\pm$
are dual to each other in this sense if{f}  $\calG_+$ is dual to $\ol{\calG_-}$ in the sense of definition
\ref{weakduality}.
\end{remark}
%
%
A stronger notion of duality for Poisson groupoids arises from symplectic double groupoids.
\bgn{definition} A \tsf{symplectic double groupoid} is a double Lie groupoid
\bgn{equation}\label{sdggen}
\bcat
\xy
*+{}="0",    <-0.9cm,0.7cm>
*+{\calS}="1", <0.7cm,0.7cm>
*+{\calG_-}="2", <-0.9cm,-0.7cm>
*+{\calG_+}="3", <0.7cm,-0.7cm>
*+{M}="4",
\ar  @ <0.07cm>   @{->} "1";"2"^{} 
\ar  @ <-0.07cm>  @{->} "1";"2"_{}  
\ar  @ <-0.07cm>  @{->} "1";"3"_{}
\ar  @ <0.07cm>   @{->} "1";"3"^{}
\ar  @ <0.07cm>   @{->} "2";"4"^{}
\ar  @ <-0.07cm>  @{->} "2";"4"_{}
\ar  @ <-0.07cm>  @{->} "3";"4"_{}
\ar  @ <0.07cm>   @{->} "3";"4"^{}
\endxy
\ecat
\end{equation}
endowed with a symplectic form which is compatible with both top groupoid structures. 
\end{definition}
Next we shall compute the \tla-groupoid of a symplectic double groupoid. Before doing that we shall need the
following
\bgn{proposition}\label{ctglla} For any Poisson groupoid $\gpdm$, the cotangent prolongation groupoid
induces an \tla-groupoid 
\bgn{equation}\label{ctgla}
\qquad\qquad
\bcat
\xy
*+{}="0",    <-0.9cm,0.7cm>
*+{T^*\calG}="1", <0.7cm,0.7cm>
*+{A^*}="3", <-0.9cm,-0.7cm>
*+{\calG}="2", <0.7cm,-0.7cm>
*+{M}="4",
\ar  @ <-0.07cm>   @{->} "1";"3"^{\Hat{}} 
\ar  @ <0.07cm>    @{->} "1";"3"^{}  
\ar  		   @{->} "1";"2"^{}
\ar                @{->} "3";"4"_{}
\ar  @ <0.07cm>    @{->} "2";"4"^{}
\ar  @ <-0.07cm>   @{->} "2";"4"_{}
\endxy\qquad\qquad.
\ecat
\end{equation}
\end{proposition}
\bgn{proof} By definition, the graph of the cotangent multiplication is a Lie subalgebroid of
$\threetimes{T^*\calG}$; moreover it is easy to see that
$\gr{\iotah}=N^*\gr{\iota}\subset\twice{T^*\calG}$ is also a Lie subalgebroid, since $\iota:\calG\to\ol{\calG}$ is
Poisson. The graph of the
unit section $\gr{\epsh}\subset A^*\times T^*\calG$ is a smooth subbundle and 
$A^*\times T^*\calG\equiv N^*M\times T^*\calG\subset\twice{T^*\calG}$ a Lie subalgebroid; on
the other hand we can identify $\gr{\epsh}$ with the diagonal  
$\Delta_{N^*M}\subset\twice{T^*\calG}$ as a subbundle, hence with a Lie subalgebroid. It
follows then from lemma \ref{subalsubal} that $\gr{\epsh}\subset A^*\times T^*\calG$ is a Lie subalgebroid.
Note that one can identify $\gr{\th}$ with the subbundle 
$\{\epsh(\th(\theta_g)),\theta_g,\theta_g\}$ of $\gr{\muh}$, i.e.
$$
\qquad
\gr{\th}\simeq(\pr_2\times\pr_3)|_{\gr{\muh}}\inverse(\Delta_{T^*\calG})
\qquad,
$$
thus with a Lie subalgebroid, due to proposition \ref{preagd}. Then
$$
\gr{\th}\subset\gr{\muh}\subset T^*\calG\times T^*\calG\times T^*\calG
$$
is a sequence of Lie subalgebroids; since 
$$
T^*\calG\times A^*\simeq A^*\times T^*\calG\simeq
N^*M\times \Delta_{T^*\calG}\subset T^*\calG\times T^*\calG\times T^*\calG
$$
is a Lie subalgebroid, it also
follows from lemma \ref{subalsubal} that $\gr{\th}\subset T^*\calG\times A^*$ is a Lie subalgebroid.
Analogously one shows that $\gr{\sh}\subset T^*\calG\times A^*$ is a Lie subalgebroid. Finally
note that $\sf{ker}_g\,\sh=N^*_g\tar\inverse(g)$, for all $g\in\calG$, then $\sh$ is a bundle map of
rank
$$
\rank_g\,\sh=\dim\calG-(\dim\calG-\rk\,\tar)=\rank_{\tar(g)}A^*
$$
over $\sor$, thus an \tla-fibration.
\end{proof}
Applying the Lie functor vertically to a symplectic double groupoid (\ref{sdggen}) yields an \tla-groupoid, which
is to be canonically identified with the cotangent prolongation \tla-groupoid associated with its side horizontal
Poisson groupoid $\calG_+$. Let $A\vup(\calS):=T^{\stv}_{\calG_+}\calS$ and $A_-:=T^{\ssv}_M\calG_-$ denote the 
Lie algebroids of the vertical groupoids and $\Ohm$ be the symplectic form on $\calS$. Note that the isomorphism
of Lie algebroids $\phi_\Ohm:T^*\calG_+\to A\vup(\calS)$ identifies  $A^*_+\equiv N^*M\subset T^*\calG_+$ with
$A_-$: in fact, for all  $(\delta g_-,\delta q)\in T^{\ssv}_M\calG_-\oplus TM$, we have
\be
\pair{\phi_\Ohm\inverse(\dd\eth\delta g_-)}{\dd\esh\delta q}
&=&-\Ohm(\dd\eth\delta g_-,\dd\etv\dd\esh\delta q)\\
&=&-\Ohm(\dd\eth\delta g_-,\dd\eth\dd\esv\delta q)\\
&=& 0\hs{3}\qquad,
\ee
since $\calG_-\subset \calS$ is Lagrangian. Moreover $\phi_\Ohm$ is an isomorphism of Lie groupoids from 
$\poidd{A\vup(\calS)}{A_-}$ to $\poidd{T^*\calG_+}{A_+^*}$: for all composable pairs 
$(\delta x,\delta y)\in A\vup(\calS)$ and $(\delta g_-,\delta h_-)\in T\calG_-$,
\be
\pair{\phi_\Ohm\inverse(\delta x\cdot_{A\vup(\calS)}\delta y))}{\delta g_+\cdot_{T\calG_+}\delta h_+)}
\!\!&=&\!\!
-\Ohm(\delta x\cdot_{A\vup(\calS)}
\delta y,\dd\etv(\delta g_+\cdot_{T\calG_+}\delta h_+))\\
\!\!&=&\!\!
-\Ohm(\delta x\cdot_{A\vup(\calS)}
\delta y,\dd\etv(\delta g_+)\cdot_{A\vup(\calS)}\dd\etv(\delta h_+))\\
\!\!&=&\!\!
-\Ohm(\delta x,\dd\etv(\delta g_+))-\Ohm(\delta y,\dd\etv(\delta h_+))\\
\!\!&=&\!\!
\pair{\phi_\Ohm\inverse(\delta x)}{\delta g_+}
+
\pair{\phi_\Ohm\inverse(\delta y)}{\delta h_+}\\
\!\!&=&\!\!
\pair{\phi_\Ohm\inverse(\delta x)\cdot_{T^*\calG_+}\phi_\Ohm\inverse(\delta y)}{\delta g_+\cdot_{T\calG_+}\delta
h_+}\quad.
\ee
In other words $\phi_\Ohm: T^*\calG\to A\vup(\calS)$ induces an isomorphism of \tla-groupoids. Note that 
$\phi_\Ohm$ is a Poisson map for the fibrewise linear structure induced by $A\vup(\calS)^*\equiv N^*\calG_+$ and 
the \emph{canonical} symplectic form on $T^*\calG_+$; therefore the restriction $A_+^*\to A_-$ is Poisson for the
Poisson structure induced by $-A_+$ on $A_+$, i.e. an isomorphism of Lie bialgebroids $(A_+^*, -A_+)\to (A_-,
A^*_-)$.\\
Therefore we have
\bgn{proposition}\cite{mkz99}\label{999} The side groupoids of a symplectic double groupoid are weakly dual Poisson groupoids.
\end{proposition}

Even if a Poisson groupoid might not admit a weak dual, it is easy to see that any Poisson groupoid,  which is
integrable as a Poisson manifold, has a unique source 1-connected dual.
%
%
%
\bgn{lemma} Any integrable Poisson groupoid $\calG$ has a unique canonical weakly dual source 1-connected 
Poisson groupoid $\calG^\star$.
\end{lemma}
\bgn{proof} It follows from proposition \ref{ctglla} that the unit section $A^*\to T^*\calG$ of the 
cotangent prolongation groupoid is a closed embedding of Lie algebroids. Thus $A^*$ is integrable, since so is
$T^*\calG$ (theorem \ref{1}), and the source 1-connected integration $\calG^\star$ of $A^*$ carries a unique
compatible Poisson bivector inducing the Lie bialgebroid $(A^*, -A)$, thanks to theorem \ref{mxiba}.     
\end{proof}
In the following we shall refer to the Poisson groupoid $\calG^\star$ of last lemma as \emph{the} weak dual of $\calG$.
The above arguments motivate the following
\bgn{definition} Two Poisson groupoids $\poidd{(\calG_\pm,\Pi_\pm)}{M}$ are in \tsf{strong
duality} if there exists a symplectic double groupoid inducing the given Poisson structures $\Pi_\pm$; a
\tsf{double} of an integrable Poisson groupoid $\calG$ is a symplectic double groupoid realizing a strong duality 
between $\calG$ and its unique weak dual $\calG^\star$.\footnote{This notion was suggested to the author by K.
Mackenzie in a private conversation (2005).}
\end{definition}
\bgn{example} For any symplectic groupoid $\poidd{\Lambda}{P}$, $(\ol{\Lambda}\times\Lambda;\Lambda,\ol{P}\times P;
P)$ is a symplectic double groupoid. The weak duality is canonically realized by the isomorphism of Lie bialgebroids
$(T^\sor_P\Lambda, N^*P)\to (T^*P, TP)$, since the target is the flip of the Lie bialgebroid of  $\poidd{\ol{P}\times
P}{P}$. 
\end{example}

\spa Last example shows that symplectic groupoids always admit strongly dual Poisson groupoids (canonically). 
Specializing theorem \ref{intla} to the case of the cotangent prolongation \tla-group\-oid associated with a Poisson
groupoid we obtain a criterion for integrability of Poisson groupoids to symplectic double groupoids and for weak
duality to imply strong duality. 
\bgn{theorem}\label{double} Let $\poidd{(\calG,\Pi)}{M}$ be any integrable Poisson groupoid with weak dual Poisson groupoid
$(\calG^\star,\Pi^\star)$. If cotangent source and target map of $\calG$ are strongly transversal, the symplectic
groupoid $\calS$ of $\calG$ carries a further Lie groupoid making it a symplectic groupoid for $\calG^\star$ and
$$
\bcat
\xy
*+{}="0",    <-0.7cm,0.7cm>
*+{\calS}="1", <0.7cm,0.7cm>
*+{\calG^\star}="2", <-0.7cm,-0.7cm>
*+{\calG}="3", <0.7cm,-0.7cm>
*+{M}="4",
\ar  @ <0.07cm>   @{->} "1";"2"^{} 
\ar  @ <-0.07cm>  @{->} "1";"2"_{}  
\ar  @ <-0.07cm>  @{->} "1";"3"_{}
\ar  @ <0.07cm>   @{->} "1";"3"^{}
\ar  @ <0.07cm>   @{->} "2";"4"^{}
\ar  @ <-0.07cm>  @{->} "2";"4"_{}
\ar  @ <-0.07cm>  @{->} "3";"4"_{}
\ar  @ <0.07cm>   @{->} "3";"4"^{}
\endxy
\ecat
$$
a symplectic double groupoid.
\end{theorem}
Note that we make no source connectivity assumptions on $\calG$.
\bgn{proof} It remains to show ($1$) that the integrating double Lie groupoid is symplectic for the top horizontal
multiplication and ($2$) the side vertical groupoid is indeed the canonical weakly dual groupoid.
($1$) The symplectic form $\ohm$ of $\calS$ induces an isomorphism of Lie bialgebroids
$$
\Phi:=(\phi_\ohm\times\phi_\ohm\times\phi_{-\ohm})\inverse: 
A\vup(\calS)\times A\vup(\calS)\times A\vup(\calS)\overset{\thicksim}
\longrightarrow 
T^*\calG\times T^*\calG\times T^*\calG
$$
for the Lie bialgebroids $(A\vup(\calS),A\vup(\calS)^*)$, respectively $(T^*\calG,T\calG)$, on the first two
components and $(A\vup(\calS),-A\vup(\calS)^*)$,  respectively $(-T^*\calG,T\calG)$ on the last (note that the first
is the Lie bialgebroid of the symplectic groupoid $\poidd{\ol{\calS}}{\calG}$ and the second is the Lie bialgebroid
associated with the Poisson manifold $\ol{\calG}$). Let $\sf{L}(\mth)$ be the graph of the top multiplication of  the
\tla-groupoid of $\calS$, namely the Lie algebroid of the graph of the top horizontal multiplication $\mth$ of
$\calS$, and regard it as a Lie groupoid over $N^*M$.   For all composable $a_\pm\in\sf{L}(\mth)$
and $\delta g_\pm\in T\calG$, we have
\be
\pair{\phi\inverse_{-\ohm}(a_+\cdot_{A\vup(\calS)}a_-)}{\delta g_+\cdot_{T\calG}\delta g_-}
&=&
\ohm({a_+\cdot_{A\vup(\calS)}a_-},{\dd\etv(\delta g_+\cdot_{T\calG}\delta g_-)})\\
&=&
\ohm(a_+,\dd\etv\delta g_+)
+
\ohm(a_-,\dd\etv\delta g_-)\\
&=&
-\pair{\phi\inverse_{\ohm}(a_+)}{\delta g_+}
-\pair{\phi\inverse_{\ohm}(a_-)}{\delta g_-}
\qquad,
\ee
that is, $\Phi$ identifies  $\sf{L}(\mth)\subset A\vup (\calS\times\calS\times\ol{\calS})$ with  
$(T\gr{\mu})^o\equiv N^*\gr{\mu}\subset T^*(\calG\times\calG\times\ol{\calG})$. 
Therefore, by  coisotropicity of $\gr{\mu}$   (use theorem
\ref{lagsub}, the connectivity requirements are implied by the strong transversality assumption), 
$\gr{\mth}\subset\calS\times\calS\times\ol{\calS}$ is also Lagrangian and $\calS$ is a symplectic double groupoid.
The Poisson structure induced on (the groupoid underlying) $\calG^\star$ makes it weakly dual to $\calG$ by
proposition \ref{999}, thus it must coincide with $\Pi^\star$ by uniqueness (theorem \ref{mxiba}).
\end{proof}
The symplectic double groupoid of last theorem, when it exists, is a vertically source 1-connected double
of its side horizontal Poisson groupoid; we shall discuss below a class of examples which allow for checking the
transversality conditions of theorem \ref{double}.
\vs{0.5}
\subsection{The case of complete Poisson groups}\hfill

\vs{0.1}
\spa Let $(G,\Pi)$ be a Poisson group. The (\tsf{left}) \tsf{dressing action} of $\frag^*$
on $G$ is the infinitesimal action 
$$
\qquad
\Upsilon:\frag^*\to\frax(G)\qquad,\qquad \Upsilon(\xi):=\Pi^\sharp\linv{\xi}
\qquad,
$$
where $\linv{\xi}$ is the left invariant 1-form on $G$ associated with $\xi\in\frag^*$:
$\linv{\xi}_g=l_{g\inverse}^*\xi$, $g\in G$.  The (\tsf{left}) \tsf{dressing vector fields} on $G$ are those in
the image of $\frag$ under the left dressing action map. Left dressing vector fields do  not have, in general,
complete flows; when they have, e.g. when $G$ is compact, $G$ is called a \tsf{complete} Poisson group.\\
The left trivialization  $T^*G\to\frag^*\acts G$ is always an isomorphism of Lie algebroids  to the action Lie
algebroid.  When $G$ is complete, the infinitesimal dressing action integrates to a global action,
$\wt{\Upsilon}:G^\star\times G\to G$ of the dual 1-connected Poisson group $G^\star$ on $G$ and  the action
groupoid $\poidd{G^\star\acts G}{G}$ is the source 1-connected integration of $\frag\acts G$. Provided $G$ is
1-connected (in this case $G^\star$ is also complete \cite{luth}) the same argument applies to the integration of the
infinitesimal right dressing action  $\frag^*\to\frax(G)$, defined using $\Pi^\star$, to a global action 
$G^\star\times G\to G^\star$.  One can show that the 1-connected integration $D$ of  the Drinfel'd double
$\mathfrak{d}$ is isomorphic to bitwisted product  ${G^\star{\lrtimes} G}$ carrying two further
symplectic groupoid structures over $G$ and $G^\star$ \cite{ks86,luth} and making 
$(G^\star{\lrtimes} G, G ; G^\star ,\bullet)$ a double of the Poisson group on $G$. A construction
of a double in the  noncomplete case, for a 1-connected $G$ was given in \cite{lw89} by Lu and Weinstein. 
\bgn{example}\label{lwd} \tsf{Lu and Weinstein's double} \cite{lw89}. 
When $G$, $G^\star$ and $D$ are 1-connected there exist integrations $\lambda:G\hookrightarrow D$, respectively
$\rho:{G^\star}\hookrightarrow D $,  of $\frag\hookrightarrow\mathfrak{d}$, respectively
$\frag^*\hookrightarrow\mathfrak{d}$. Moreover, one can show that $D$ has a compatible
Poisson structure $\pi_D$, which happens to be nondegenerate on the submanifold of elements $d$, admitting a
decomposition $d=\lambda(g_+)\rho(g^\star_+)=\rho(g^\star_-)\lambda(g_-)$, $g_\pm\in G$ and
$g^\star_\pm\in{G^\star}$. Moreover there is also a natural double Lie groupoid
$$\qquad\qquad\qquad\qquad\qquad\qquad\qquad
\bcat
\xy
*+{}="0",    <-0.7cm,0.7cm>
*+{\calD}="1", <0.7cm,0.7cm>
*+{G^\star}="2", <-0.7cm,-0.7cm>
*+{G}="3", <0.7cm,-0.7cm>
*+{\bullet}="4",
\ar  @ <0.07cm>   @{->} "1";"2"^{}
\ar  @ <-0.07cm>  @{->} "1";"2"_{}
\ar  @ <-0.07cm>  @{->} "1";"3"_{}
\ar  @ <0.07cm>   @{->} "1";"3"^{}
\ar  @ <0.07cm>   @{->} "2";"4"^{}
\ar  @ <-0.07cm>  @{->} "2";"4"_{}
\ar  @ <-0.07cm>  @{->} "3";"4"_{}
\ar  @ <0.07cm>   @{->} "3";"4"^{}
\endxy\ecat
\qquad\qquad\calD=G\times{G^\star}\times{G^\star}\times G\quad.
$$
It turns out that the double subgroupoid, whose total space is
$$
\calS=\{(g_+,g^\star_+,g^\star_-,g_-)\:|\:\lambda(g_+)\rho(g^\star_+)
=
\rho(g^\star_-)\lambda(g_-)\}
$$
carries a compatible symplectic form, inducing the Poisson structures on $G$ and
${G^\star}$, which is the inverse of the pullback of $\pi_D$, under
the natural local diffeomorphism $\calS\rightarrow D$.
\end{example}

\spa
We  remark that the total space of Lu-Weinstein's double is  only locally diffeomorphic  to $D$, unless $G$ is
complete, therefore it is in general neither source (1-)connected over $G$, nor over $G^\star$.  Moreover, the
1-connectivity of $G$ is essential in both constructions.\\ 
In the complete case it is possible to drop the connectivity assumptions on $G$  and, nevertheless, to obtain
vertically source 1-connected double.
\bgn{theorem}\cite{07b}\label{mainddd} For any complete Poisson group $G$, the source  1-connected
symplectic groupoid $\calS$ of $G$ carries a unique Lie groupoid structure over the
1-connected dual Poisson group $G^\star$ making it a double of $G$.
%
\end{theorem}
The transversality conditions of theorem \ref{double} are met if the top source map of the cotangent prolongation
groupoid has the \tla-homotopy lifting properties \ref{liftingprop}, therefore theorem \ref{mainddd} is a
consequence of the following
\bgn{lemma}\label{lift} For any complete Poisson group $G$, the cotangent source map 
$\sh:T^*G\to \frag^*$ has the 0- and 1-\tla-homotopy lifting properties.
\end{lemma}
\bgn{proof} Note that the diagram
$$
\bcat
\xy
*+{T^*G}="0",    <2cm,0cm>
*+{\frag^*\acts G}="1", <1cm,-1cm>
*+{\frag^*}="2"
\ar  @ {->} "0";"1"^{\thicksim\quad}  
\ar  @ {->} "1";"2"^{\pr}
\ar  @ {->} "0";"2"_{\sh}
\endxy
\ecat
$$
commutes in the category of Lie algebroids for the left trivialization on the top edge: it suffices to prove the
statement for the projection   $\frag^*\acts G\to\frag^*$.  $\frag^*\acts G$-paths are pairs $(\xi,\gamma)$ of
$\frag^*$-paths and paths in $G$, such that $\Upsilon(\xi(u))_{\gamma(u)}=\dd\gamma_u$, $u\in I$. A Lie algebroid
homotopy  $h_\acts: TI^{\hbox{\tiny{$\times 2$}}}\to\frag^*\acts G$ is a pair $(h, X)$ for which $h$ is a
$\frag^*$-homotopy and $X:I\to G$  satisfies  $\dd X=\Upsilon\comp h$. To see this, note that $h_\acts$ takes
values in the pullback  bundle  $X^+(\frag^*\acts G)\simeq \frag^*\times I^{\hbox{\tiny{$\times 2$}}}$.  The
bracket compatibility for  $h_\acts$ can be written choosing the de Rham  differential on $G$  with coefficients
in $\frag^*$ as a linear connection $\nabla$ for   $\frag^*\acts G$; this way the covariant derivative for the
pullback connection is simply the de Rham differential on $I^{\hbox{\tiny{$\times 2$}}}$  with coefficients in
$\frag^*$. The torsion tensor of $\nabla$ is	
\be
\tau^\nabla(F_+,F_-)
&=&
{\dd F_-}_g(\rho^\acts(F_+))-{\dd F_+}_g(\rho^\acts(F_-))-\bracts{F_+}{F_-}_g\\
&=&-\brast{F_+(g)}{F_-(g)}
\quad,
\ee
$F_\pm\in\cif(G,\frag^*)$, thus the pullback of $\tau^\nabla$ induces the  canonical
graded Lie bracket, also denoted by $\brast{\:}{\,}$, on 
$\Ohm^\bullet(I^{\hbox{\tiny{$\times 2$}}},\frag^*)$
and the bracket compatibility condition for $h_\acts$ reduces to the classical 
Maurer-Cartan equation
$$
\dd h+\frac{1}{2}\brast{\,h\,}{h\,}=0
$$
for $h\in\Ohm^1(I^{\hbox{\tiny{$\times 2$}}},\frag^*)$. Suppose now $(\xi_-,\gamma_-)$
is a fixed $\frag^*\acts G$-path, let $h$ be a $\frag^*$-homotopy form $\xi_-$ to some
other $\frag^*$-path $\xi_+$ and $H:I^{\hbox{\tiny{$\times 2$}}}\to G^\star$ the unique
$G^\star$-homotopy integrating $h$, i.e. such that 
$$
h= \delta_r^uH \cdot du + \delta_r^\eps H\cdot d\eps
$$
for the partial right derivatives. We claim that $h_\acts:=(h,X)$, where
\be
X(u,\eps)&=&(H(u,\eps)\cdot H(u,0)\inverse )*\gamma_-(u)\\
         &=&H(u,\eps)*(H(u,0)\inverse *\gamma_-(u))
\qquad,\qquad 
\ee
$u$, $\eps\in I$, is a $\frag^*\acts G$-homotopy; here $*$ denotes the integrated 
dressing action map.
We postpone to the end of the proof this straightforward but lengthy check.
The lifting conditions follow: by construction
\be
\iota^*_{\pa^\pm_H}(h,X)&=&(\iota^*_{\pa^\pm_H}h,X\comp\iota_{\pa^\pm_H})=0\\
\iota^*_{\pa^\pm_V}(h,X)&=&(\xi_\pm, \gamma_\pm) 
\ee
for some path $\gamma_+$ in $G$, 
($l_0$) If $h$ is a fixed $\frag^*$-homotopy to some fixed  $\frag^*$-path $\xi_+$,
$(\xi_+,\gamma_+)$ is the desired lift. 
($l_1$) In particular if $\xi_-$ is $\frag^*\acts G$-homotopic to the constant 
$\frag^*\acts G$-path $\xi_o\equiv 0$ and $h$ a $\frag^*$-homotopy from $\xi_-$ to
the constant $\frag^*$-path $\xi_+\equiv 0$, we have
$$
\qquad
\dd\gamma_+|_u=\Upsilon(\xi_+(u))_{\gamma_+(u)}=0_{\gamma_+(u)}
\qquad, 
$$
hence $ \gamma_+(u)\equiv\gamma_+(0)=\gamma_-(0)$, since the base paths of homotopic Lie algebroid paths are homotopic relatively to the
endpoints; $(h,X)$ is then  the desired homotopy. 
 In order to prove our claim it remains to check the anchor compatibility condition for $(h,X)$; 
set $a=\delta_r^uH$ and  $b=\delta_r^\eps H$. The derivative of $X$ in the
$\eps$-direction is thus
\be
\pa_\eps X(u,\eps)&=&\Upsilon(\pa_\eps H(u,\eps))_{H(u,0)\inverse *\gamma_-(u)}
\:=\:
\Upsilon(\dd r_{H(u,\eps)}b(u,\eps))_{H(u,0)\inverse *\gamma_-(u)}\\
&=&\Upsilon(b(u,\eps))_{X(u,\eps)}\hs{3},
\ee
since for all $g\in G$, $h^\star\in G^\star$, $\xi=\dot{g}^\star(o)\in\frag^*$ and path
$g^\star$ in $G^\star$
$$
\Upsilon(\dd r_{h^\star}\xi)
=
\left.
\frac{d}{d\alpha}\right|_{\alpha=o}\wt{\Upsilon}(g^\star(\alpha)\cdot h^\star, g)
=
\left.\frac{d}{d\alpha}\right|_{\alpha=o}\wt{\Upsilon}(g^\star(\alpha),h^\star*g)
=
\Upsilon(\xi)_{h^\star*g}
\quad.
$$
The computation of the derivative in the $u$-direction is more involved. We have
$$
\qquad
\pa_u X(u,\eps)=\dd\wt{\Upsilon}_{(H(u,\eps)\cdot H(u,0)\inverse, \gamma_-(u))}
		(\delta	 H, \dot{\gamma}_-(u))
\qquad,
$$
where 
$\quad
\dot{\gamma}_-(u)
=
\Upsilon(\xi_-(u))_{\gamma_-(u)}
=
\Upsilon((\iota_{\pa\hup^{-}}^* h)(u))
=
\Upsilon(a(u,0))_{X(u,0)}
\quad$ and
\be
\delta H 
&=&
\dd\mu^\star_{(H(u,\eps),H(u,0)\inverse)}
(\pa_u H(u,\eps),\dd\iota_\star\pa_u H(u,0))\\
&=&
\dd\mu^\star_{(H(u,\eps),H(u,0)\inverse)}
(\dd r_{H(u,\eps)}a(u,\eps),\dd\iota_\star\dd r_{H(u,0)}a(u,0))\\
&=&
\dd\mu^\star_{(H(u,\eps),H(u,0)\inverse)}
(0_{H(u,\eps)},\dd\iota_\star\dd r_{H(u,0)}a(u,0))\\
&+&
\dd\mu^\star_{(H(u,\eps),H(u,0)\inverse)}
(\dd r_{H(u,\eps)}a(u,\eps),0_{H(u,0)\inverse})\\
&=& \delta H_+ +\delta H_-\hs{5},
\ee
since
$\dd\mu^\star:TG^\star\times TG^\star\to TG^\star$ is fibrewise linear, with 
\be
\delta H_+
&=&
\dd l_{H(u,\eps)}\dd\iota_\star\dd r_{H(u,0)}a(u,0)
\:\:=\:\:
\dd l_{H(u,\eps)}\dd l_{H(u,0)\inverse}\dd\iota_\star a(u,0)
\\
&=&
\dd l_{H(u,\eps)\cdot H(u,0)\inverse}\dd\iota_\star a(u,0)
\\
\delta H_-&=&\dd r_{H(u,0)}\inverse\dd r_{H(u,\eps)}a(u,\eps)=
\dd r_{H(u,\eps)\cdot H(u,0)\inverse}a(u,\eps)\hs{2}.
\ee
The tangent action map $\dd\wt{\Upsilon}:TG^\star\times TG\to TG$ is also fibrewise
linear, hence
\be
\pa_u X(u,\eps)
&=&
\dd\wt{\Upsilon}_{(H(u,\eps)\cdot H(u,0)\inverse,\gamma_-(u))}
(\delta H_+,\Upsilon(a(u,0))_{\gamma_-(u)})\\
&+&
\dd\wt{\Upsilon}_{((H(u,\eps)\cdot H(u,0),\gamma_-(u))}
(\delta H_-,0_{\gamma_-(u)});
\ee
the first term of last expression vanishes, since it can be rewritten as 
\be\!\!\!\!\!\ \!\!\!\!\!
& &\!\!\!\!\!
\dd\wt{\Upsilon}_{@}
(\dd l_{(H(u,\eps)\cdot H(u,0)\inverse)}\dd\iota_\star a(u,0),
\Upsilon(a(u,0))_{\gamma_-(u)}\\
&=&\!\!\!\!\!
\dd\wt{\Upsilon}_{@}
(\dd l_{H(u,\eps)\cdot H(u,0)\inverse}\dd\iota_\star a(u,0),
\dd\wt{\Upsilon}_{(e_\star ,\gamma_-(u))}
(a(u,0),0_{\gamma_-(u)}))\\
&=&\!\!\!\!\!
\dd\wt{\Upsilon}_{@}
(\dd\mu^\star_{(H(u,\eps)\cdot H(u,0)\inverse,e_\star)}
(\dd l_{(H(u,\eps)\cdot H(u,0)\inverse)}\dd\iota_\star a(u,0),
a(u,0)),0_{\gamma_-(u)})\\
&=&\!\!\!\!\!
\dd\wt{\Upsilon}_{@}
(0_{H(u,\eps)\cdot H(u,0)\inverse},0_{\gamma_-(u)})=0
\ee
by equivariance, where we have set $@=(H(u,\eps)\cdot H(u,0)\inverse,\gamma_-(u))$ to
simplify the expressions; therefore
\be
\pa_u X(u,\eps)
&=&
\dd\wt{\Upsilon}_{((H(u,\eps)\cdot H(u,0),\gamma_-(u))}
(\dd r_{H(u,\eps)\cdot H(u,0)\inverse}a(u,\eps),0_{\gamma_-(u)})\\
&=&
\dd\wt{\Upsilon}_{(e_\star,H(u,\eps)\cdot H(u,0)\inverse *\gamma_-(u))}
(a(u,\eps),0_{H(u,\eps)\cdot H(u,0)\inverse *\gamma_-(u)})\\
&=&\Upsilon(a(u,\eps))_{X(u,\eps)}\hs{3}.
\ee
We just have shown that
\be
\dd X_{(u,\eps)}
&=&
\Upsilon(\delta_r^\eps H(u,\eps))_{X(u,\eps)}\cdot d\eps 
+
\Upsilon(\delta_r^u H(u,\eps))_{X(u,\eps)}\cdot du\\ 
&=&
\Upsilon(h(u,\eps))_{X(u,\eps)}
\ee
and this concludes the proof.	
\end{proof}

\chapter{Morphic actions}\label{chapiii}
\spa On the one hand Poisson manifolds behave well under reduction for
actions as general as those of Poisson groupoids: a Poisson bivector
always descends to the quotient by a free and proper compatible
action. On the other hand, Poisson manifolds are not
always integrable to symplectic groupoids; the natural question is
thus:
%
{\em
Are quotients of integrable Poisson manifolds also integrable?\\
}
%
A positive answer was recently given by Fernandes-Ortega-Ratiu \cite{for07} in the case of Lie group actions and Lu
\cite{lu07} gave a construction  of symplectic groupoids for certain Poisson homogeneous spaces. In the case of Poisson group
actions with a complete moment map, Xu described in \cite{xu92} a reduction procedure on a lifted moment map; when the
reduced space of the latter is smooth, it is a symplectic groupoid for the quotient Poisson structure.\\
In fact, lifting processes naturally produce stronger symmetries out of weaker ones, often associated with suitable moment
maps which can be used to compute information on the original action.  Consider for instance the action of a Lie group $G$ on
an ordinary manifold $M$. It was already remarked by Smale \cite{smale} that it is possible to \emph{cotangent lift}  the
action of $G$ on $M$ to a symplectic action of $G$ on $T^*M$; such a lift is  always endowed  with an equivariant moment map
$\jh:T^*M\to\frag^*$ in the sense of Marsden-Weinstein \cite{moment},  morally the dual 
\bgn{equation}\label{cotomo}
\pair{\jh(\theta_q)}{\xi}=\pair{\theta_g}{\sigma(\xi)}
\end{equation}
of the infinitesimal action \cite{mw88}.
\emph{Path prolongation} is another lifting procedure which generates moment maps. Let a Lie group $G$ act on a symplectic
manifold $(M,\ohm)$ by symplectic diffeomorphisms. The space $\Pi(M)$ of paths in $M$ up to homotopy relative to the
endpoints carries a natural symplectic form and a $G$-action. Mikami and Weinstein observed in \cite{mw88} that there always
exist an equivariant moment map $J:\Pi(M)\to\frag^*$ 
$$
\pair{J([\gamma])}{\xi}=\underset{\gamma}\int\iota_{\sigma(\xi)}\ohm
$$
associated with this action. An analogous lifting procedure via Lie algebroid paths 
can be applied to a Poisson action of a Lie group on an
integrable Poisson manifold $P$, with symplectic groupoid $\Lambda$, to obtain a symplectic action on $\Lambda$, 
with equivariant moment map $\Lambda\to\frag^*$ \cite{fs06}; 
this generalization of Mikami and Weinstein's moment map, was used in
\cite{for07} to produce an integration of $P/G$.\\
It is thus quite a general phenomenon that lifting processes tend to enhance symmetries;
this effect is more transparent in the general setting of Poisson geometry and
when more structure is around. Indeed cotangent lifting and path prolongation are
best understood as instances of ``duality'' between Poisson geometric objects and
Lie algebroids-Lie groupoids, and lead to interesting examples of double Lie structures.\\ 
Suppose now a Poisson group $G$ acts on a Poisson manifold $P$. 
The cotangent lifted moment map $\jh:T^*P\to \frag^*$ (\ref{cotomo})
preserves the Poisson structures, in the sense that it is a morphism of Lie
bialgebroids (fig. 2), and it is still equivariant, provided the suitable cotangent lift of
$G$ is correctly identified. In fact $\jh$ is equivariant for the coadjoint action
on the codomain, but, on the domain, the classical cotangent lift needs to be replaced with an action,
which is symplectic in the sense of Mikami and Weinstein \cite{mw88}, of the
symplectic groupoid $\poidd{T^*G}{\frag^*}$. Remarkably, the compatibility of the Poisson group action with 
the Poisson bivector fields is recovered in the dual description in terms of Lie algebroids.
$$
\bcat
\bgn{array}{cccccc}
\xy
*+{}="0",    <-0.7cm,0.7cm>
*+{T^*G}="1", <0.7cm,0.7cm>
*+{\frag^*}="2", <-0.7cm,-0.7cm>
*+{G}="3", <0.7cm,-0.7cm>
*+{\bullet}="4",
\ar  @ <-0.07cm>   @{->} "1";"2"^{}
\ar  @ <0.07cm>    @{->} "1";"2"^{}  
\ar  		   @{->} "1";"3"^{}
\ar                @{->} "2";"4"_{}
\ar  @ <0.07cm>    @{->} "3";"4"^{}
\ar  @ <-0.07cm>   @{->} "3";"4"_{}
\endxy
&\qquad\qquad&
\xy
*+{}="0",    <-0.7cm,0.7cm>
*+{T^*P}="1", <0.7cm,0.7cm>
*+{\frag^*}="3", <-0.7cm,-0.7cm>
*+{P}="2", <0.7cm,-0.7cm>
*+{\bullet}="4",
\ar    @{->} "1";"2"_{} 
\ar  		   @{->} "1";"3"^{\jh}
\ar                @{->} "2";"4"_{}
\ar  @ <0.07cm>    @{->} "3";"4"
\endxy
&\qquad\qquad&
\xy
*+{}="0",    <-1.7cm,0.7cm>
*+{T^*G\acts T^*P}="1", <0.7cm,0.7cm>
*+{T^*P}="2", <-1.7cm,-0.7cm>
*+{G\acts P}="3", <0.7cm,-0.7cm>
*+{P}="4",
\ar  @ <-0.07cm>   @{->} "1";"2"^{}
\ar  @ <0.07cm>    @{->} "1";"2"^{}  
\ar  		   @{->} "1";"3"^{}
\ar                @{->} "2";"4"_{}
\ar  @ <0.07cm>    @{->} "3";"4"^{}
\ar  @ <-0.07cm>   @{->} "3";"4"_{}
\endxy\\
\overset{\quad}Figure\:1&&
Figure\:2&&
\quad Figure\:3\\
\end{array}
\ecat
$$
In fact, $\poidd{T^*G}{\frag^*}$ extends to the cotangent prolongation \tla-groupoid (Fig. 1), and the cotangent 
lifted action induces the \tla-groupoid in figure 3, obtained by Mackenzie in \cite{mkz00b}, which
completely encodes a  Poisson action dually in the category of Lie algebroids and, roughly, lifts it to 
the  category of symplectic manifolds. It was also observed in \cite{mkz00b} that the
zero level reduction of $\jh$ can be used to compute the Lie bialgebroid of the quotient Poisson bivector on $P/G$.

\spa The study of a Poisson action of a Poisson groupoid on a Poisson manifold, as we shall see in Section \ref{macla},
naturally leads to -- and is equivalent to --
a compatible action of the cotangent lifted \tla-groupoid, respectively, under integrability conditions, 
a compatible action of a symplectic double groupoid, on a canonical ``moment morphism'' of Lie algebroids, respectively of Lie groupoids. 
This motivates the study we carry on in this Chapter of the reduction of morphic actions in the categories of Lie
algebroids and Lie groupoids, namely groupoid actions in these categories for which the associated action groupoid is an
object.\\
In Section \ref{macla} we introduce morphic actions of \tla-groupoids and perform the cotangent lift of a Poisson groupoid 
action to a morphic action in the category of Lie algebroids (theorem \ref{ctglift}); by specializing to the
Poisson group case, we obtain an alternative construction for Mackenzie's lift of \cite{mkz00b}\footnote{The chance of 
extending the results of \cite{mkz00b} to Poisson actions of Poisson groupoids and the method used below for this generalization 
emerged from private conversations with Kirill Mackenzie (January 2007).}.\\
In Section \ref{rmalag} we study the reduction of morphic actions of \tla-groupoids and prove the following general result:
\bgn{thm}[\bf\ref{generalred}] Let $(\Ohm,\calG;A,M)$ be an \tla-groupoid acting morphically on a morphism of Lie
algebroids $\jh:B\to A$ over $j:N\to M$. If the action is free and proper (so that $B/\Ohm$ and $N/\calG$ are
smooth manifolds) then there exists a unique Lie algebroid on $B/\Ohm\to N/\calG$ making the quotient
projection $B\to B/\Ohm$ a strong \tla-fibration over $N\to N/\calG$. 
\end{thm}
Thereafter we consider the reduction of the cotangent lifted moment morphism and give a characterization of the Koszul
algebroid associated with a quotient Poisson manifold as a quotient Lie algebroid (proposition \ref{poissonreduction}).\\
In Section
\ref{ima} we consider the reduction of morphic actions of double Lie groupoids and the integrability of morphic actions of
\tla-groupoids. 
The symplectic case is dealt with in Section \ref{iqps}, where we obtain the following
\bgn{thm}[\bf\ref{quoproqui}] Let $(\calS,\calG;\calG^\bullet,M)$ be a symplectic double groupoid acting 
morphically on a morphism of Lie groupoids $\calJ:\Lambda\to\calG^\bullet$ over $j:P\to M$, where
$\poidd{\Lambda}{P}$ is a symplectic groupoid, in such a way that $\Lambda$ is a symplectic $\calS$-space. 
If $\calJ$ is source submersive and the side action is free and proper, then
\medskip\\
$i)$ The reduced kernel groupoid $\poidd{\calJ\inverse(\eps_\bullet(M))/\calG}{P/\calG}$ carries a unique symplectic
form making the projection  $\calJ\inverse(\eps_\bullet(M))\to\calJ\inverse(\eps_\bullet(M))/\calG$  a
Poisson submersion;
\medskip\\
$ii)$ $\poidd{\calJ\inverse(\eps_\bullet(M))/\calG}{P/\calG}$ is a symplectic groupoid for the quotient Poisson
manifold $P/\calG$.
\end{thm}
Last result can be applied to an integrable Poisson $\calG$-space $P$, provided the Poisson groupoid $\calG$ is integrable to
a symplectic double groupoid, in order to produce a symplectic groupoid for the quotient Poisson bivector; this approach is effective,
for instance, in the case of complete Poisson group actions, discussed at the end of this Chapter.
We obtain Xu's reduction of \cite{xu92} and Fernandes-Ortega-Ratiu integration of \cite{for07} as special cases of
our approach to the integration of quotients of complete Poisson group actions, developed in \cite{07b}.\\ 
In the last Section we also use an \tla-path-prolongation of the cotangent lifted action associated with a  
a Poisson $\calG$-space $P$ to derive our main application
\bgn{thm}[\bf\ref{quiproquo}] Let a Poisson groupoid $\gpdm$ act freely and properly on $j:P\to M$. If the action is
Poisson, then $P/\calG$ is integrable to a symplectic groupoid if{f} $\ker{\jh}$ is an integrable Lie algebroid.
\end{thm}
Here $\jh:T^*P\to A^*$ is a lifted moment map canonically associated with the original action, where $A$ is the Lie algebroid of $\calG$ 
and $\ker{\jh}\subset T^*P$ is the vertical subbundle for the infinitesimal action of $\Gamma(A)$.\\
In particular, we give a positive answer to the above question under most natural assumptions.
\bgn{cro}[\bf\ref{happy}] If $P$ is integrable, then so is $P/\calG$.
\end{cro}
\newpage
\subsection*{Definitions and remarks on groupoid actions in a category}\hfill

\spa Let $\gpdm$ be a groupoid acting on a moment map $j:N\to M$; denote with 
$\sigma:\calG\fib{\sor}{j}N\to N$ the action map. The action
 can  be fully described in terms of diagrams
\bgn{equation}\label{compaact}
\barr{c}
\bcat
\xy
*+{}="0",    <0cm,1.2cm>
*+{\calG\fib{\sor}{j}N}="t", <0cm,-1.2cm>
*+{M}="b", <-1.2cm,0cm>
*+{\calG}="l", <+1.2cm,0cm>
*+{N}="r",  
\ar     @ {->} "t";"l"_{\pr_1} 
\ar     @ {->} "t";"r"^{\sigma}  
\ar     @ {->} "l";"b"_{\tar}
\ar     @ {->} "r";"b"^{j}
\endxy
\ecat\\
\\
\hbox{\tsf{compatibility}}\\
\hbox{\tsf{with the moment map}}\\
\earr
\:\:\:
\barr{c}
\bcat
\xy
*+{}="0",    <-0.7cm,1.7cm>
*+{N}="tl",  <0.7cm,1.7cm>
*+{N}="tr",    <-1.2cm,0.5cm>
*+{\Delta_N}="l", <1.2cm,0.5cm>
*+{\calG\fib{\sor}{j}N}="r", <0cm,-0.5cm>	
*+{M\fib{}{j}N}="b",
\ar     @ {->} "tl";"l"_{\Delta_N}
\ar     @ {->} "r";"tr"_{\quad\sigma}  
\ar     @ {->} "l";"b"_{j\times\id_N}  
\ar     @ {->} "b";"r"_{\eps\times\id_N}
\ar     @ {->} "tl";"tr"^{\id_N}
\endxy\\
\\
\hbox{\tsf{unitality}}\\
\earr
\:\:\:
\barr{c}
\bcat
\xy
*+{}="0",    <0cm,1.2cm>
*+{\calG^{(2)}\fib{\sor\comp\pr_1}{j}N
   =
   \calG\fib{\sor}{\sigma}(\calG\fib{\sor}{j}N)}="t", <0cm,-1.2cm>
*+{M}="b", <-1.2cm,0cm>
*+{\calG\fib{\sor}{j}N}="l", <+1.2cm,0cm>
*+{\calG\fib{\sor}{j}N}="r",  
\ar     @ {->} "t";"l"_{\pr_1} 
\ar     @ {->} "t";"r"^{\sigma}  
\ar     @ {->} "l";"b"_{\tar}
\ar     @ {->} "r";"b"^{j}
\endxy
\ecat\hspace{-9.26318pt}.\\
\\
\hbox{\tsf{multiplicativity}}\\
\earr
\end{equation}\\
This leads the abstract notion of a morphic action, namely a
groupoid action in a category.
\bgn{definition}
Let $\bfC$ be a small category with direct and fibered products,
\pmbgpdm  a groupoid object in $\bfC$ and ${\bf j}:\bfN\to\bfM$ an arrow. A
groupoid action of \pmbgpdm on $\bfj$ is a \tsf{morphic action} if the action map
$
\pmbs:\pmbg\fib{\bfs}{\bfj}\bfN\to\bfN
$
is an arrow.
\end{definition}
Morphic actions are precisely those for which the associated action groupoid
stays in the category.
\bgn{proposition}\label{morphiccat} Let $\mathbf{C}$ be a category with fibred products. 
A groupoid action of a groupoid object\: \pmbgpdm on an arrow 
${\bf j}:\bfN\to\bfM$ is morphic if{f} the action groupoid
$\poidd{\pmbg\acts\bfN}{\bfN}$ is a groupoid object.	
\end{proposition}
\bgn{proof} The implication to the left is clear, since the action map is the
target of the action groupoid. For the implication to the right consider that
the source map $\bfs^\acts:\pmbg\acts\bfN\to\bfN$ is the projection to the second
 factor of the
fibered product $\pmbg\fib{\bfs}{\bfj}\bfN$, hence an arrow. Moreover, the unit
section $\pmbe^\acts:\bfN\to\pmbg\acts\bfN$ admits a factorization
$$
\bcat
\xy
*+{}="0",    <-1cm,0.7cm>
*+{\bfN}="1", <1cm,0.7cm>
*+{\pmbg\acts\bfN}="3", <-1cm,-0.7cm>
*+{\Delta_\bfN}="2", <1cm,-0.7cm>
*+{\bfM\fib{}{\bfj}\bfN}="4",
\ar @ {->} "1";"2"^{\Delta}  
\ar @ {->} "1";"3"^{\pmbe^\acts\quad}
\ar @ {->} "2";"4"_{\bfs\times\id_\bfN\quad}
\ar @ {->} "4";"3"_{\pmbe\times\id_\bfN}
\endxy
\ecat
$$
in terms of arrows, therefore it is also an arrow. According to proposition \ref{obietto}, it now
suffices to show that the division map
$
\pmbd^\acts:
\pmbg\acts\bfN\fib{\bfs^\acts}{\bfs^\acts}\pmbg\acts\bfN\to\pmbg\acts\bfN
$,
$$
\qquad
\pmbd^\acts((g_+,n),(g_-,n))=(g_+\pmb{\cdot} g_-\inverse, g_-*n)
=
(\pmbd(g_+,g_-),\pmbs(g_-,n))
\qquad,
$$
is an arrow; this holds true, since there is a factorization
\bgn{equation}\label{exagon}
\bcat
\xy
*+{}="0"
,<-2.5cm,1.9cm>
*+{\pmbg\acts\bfN\fib{\bfs^\acts}{\bfs^\acts}\pmbg\acts\bfN}="1",
<-4cm,0cm>
*+{(\pmbg\fib{\bfs}{\bfs}\pmbg)\fib{\bfs\times\bfs}{\bfj\times\bfj}\Delta_\bfN}="2"
, <-2.5cm,-1.9cm>
*+{(\pmbg\fib{\bfs}{\bfs}\pmbg)\fib{\bfs\comp\pr_2}{\bfj}\bfN}="3"
,<2.5cm,-1.9cm>
*+{(\pmbg\fib{\bfs}{\bfs\comp\pr_1}\Delta_{\pmbg})\fib{\bfs\comp\pr_3}{\bfj}\bfN}="4"
,<4cm,0cm>
*+{(\pmbg\fib{\bfs}{\bfs}\pmbg)\fib{\bfs\comp\pr_2}{\bfs\comp\pr_1}
\pmbg\fib{\bfs}{\bfj}\bfN}="5"
,<2.5cm,1.9cm>
*+{\pmbg\acts\bfN}="6",
\ar @ {->} "1";"2"_{\thicksim}  
\ar @ {->} "2";"3"_{\thicksim}
\ar @ {->} "3";"4"_{\thicksim\quad\:\:}
\ar @ {->} "4";"5"_{\thicksim}
\ar @ {->} "5";"6"_{\pmbd\times\pmbs}
\ar @ {->} "1";"6"^{\qquad\quad\pmbd^\acts}
\endxy
\ecat
\end{equation}
in terms of arrows.
\end{proof}
The same remark as in the introduction to Chapter \ref{chapii} applies: if $\bfC$ does not posses general fibered products,
the statement of last proposition still holds true, up to restricting to a suitable class of groupoid objects for which the
relevant fibered products exist. As we shall see in Sections \ref{macla} and \ref{ima}, this is the case for the categories
of Lie algebroids and Lie groupoids, restricting to \tla-groupoids and double Lie groupoids. 
\vs{1}
\section{Morphic actions in the category of Lie algebroids}\label{macla}
\begin{quotation} 
In this Section we introduce morphic actions in the category of Lie algebroids, namely compatible  actions of
\tla-groupoids on morphisms of Lie algebroids. Furthermore, we produce our main example, the cotangent lift of a
Poisson groupoid action. As it was remarked by He-Liu-Zhong \cite{hlz01}, to a Poisson $\calG$-space $P$ for a
Poisson groupoid $\calG$ one can always associate a canonical morphism of Lie bialgebroids $\jh$,  which is to be
regarded as a moment map for an action of the Lie bialgebroid of $\calG$. We show that $\jh$ is indeed a moment map
for an action of the cotangent prolongation groupoid and such cotangent lifted action determines a morphic action of
the cotangent prolongation \tla-groupoid of $\calG$ on the Koszul algebroid of $P$. As a consequence to a Poisson
$\calG$-space we can associate an action \tla-groupoid, fully encoding the compatibility of the original action, 
in the category of Lie algebroids.
\end{quotation}
\vs{0.5}
\subsection{Morphic actions of Lie algebroids and action
\tla-groupoids}\label{mola}\hfill

\vs{0.1}
\spa Let $(\Ohm,\calG;A,M)$ be an \tla-groupoid and $\sf j$ a morphism of Lie algebroids
$$
\bcat
\sf{\Ohm}:=
\xy
*+{}="0",    <-0.7cm,0.7cm>
*+{\Omega}="1", <0.7cm,0.7cm>
*+{A}="2", <-0.7cm,-0.7cm>
*+{\calG}="3", <0.7cm,-0.7cm>
*+{M}="4",
\ar  @ <-0.07cm>   @{->} "1";"2"^{\Hat{}} 
\ar  @ <0.07cm>    @{->} "1";"2"^{}  
\ar  		   @{->} "1";"3"^{}
\ar                @{->} "2";"4"_{}
\ar  @ <0.07cm>    @{->} "3";"4"^{}
\ar  @ <-0.07cm>   @{->} "3";"4"_{}
\endxy
\qquad\qquad
\sf{j}:=
\xy
*+{}="0",    <-0.7cm,0.7cm>
*+{B}="1", <0.7cm,0.7cm>
*+{A}="3", <-0.7cm,-0.7cm>
*+{N}="2", <0.7cm,-0.7cm>
*+{M}="4",
\ar    @{->} "1";"2"_{} 
\ar  		   @{->} "1";"3"^{\jh}
\ar                @{->} "2";"4"_{j}
\ar  @ <0.07cm>    @{->} "3";"4"
\endxy\qquad.
\ecat
$$
Suppose that 
$\poidd{\Ohm}{A}$ acts on $\jh$ and 
$\poidd{\calG}{M}$ acts on $j$ in such a way that the corresponding action
maps
$$
\msigma:=
\bcat\xy
*+{}="0",    <-1.1cm,0.7cm>
*+{\Ohm\fib{\sh}{\jh}B}="1", <0.7cm,0.7cm>
*+{B}="3", <-1.1cm,-0.7cm>
*+{\calG\fib{\sor}{j}N}="2", <0.7cm,-0.7cm>
*+{M}="4",
\ar    @{->} "1";"2"_{} 
\ar  		   @{->} "1";"3"^{\quad\Hat{\sigma}}
\ar                @{->} "2";"4"_{\quad\sigma}
\ar  @ <0.07cm>    @{->} "3";"4"
\endxy\ecat\qquad\qquad
$$
define a vector bundle map. Note that the vector bundle on the left in the last diagram is always well
defined, since $\sh$  is an \tla-fibration over $\sor$, and it carries a fibered product Lie
algebroid making it a Lie subalgebroid of the direct product Lie algebroid $\Ohm\times
A\rightarrow\calG\times N$, since both $\sh$ and $\jh$ are morphisms of Lie algebroids. Moreover,
it is always  possible to form a diagram
$$
\bcat
\sf{\Ohm\acts B}:=
\xy
*+{}="0",    <-0.9cm,0.7cm>
*+{\Ohm\acts B}="1", <0.7cm,0.7cm>
*+{B}="2", <-0.9cm,-0.7cm>
*+{\calG\acts N}="3", <0.7cm,-0.7cm>
*+{N}="4",
\ar  @ <-0.07cm>   @{->} "1";"2"^{\quad\raise5pt\hbox{$\acts$}\quad} 
\ar  @ <0.07cm>    @{->} "1";"2"^{}  
\ar  		   @{->} "1";"3"^{}
\ar                @{->} "2";"4"_{}
\ar  @ <0.07cm>    @{->} "3";"4"^{}
\ar  @ <-0.07cm>   @{->} "3";"4"_{}
\endxy
\ecat
$$
for the action Lie groupoids $\poidd{\Ohm\acts B}{B}$ and  $\poidd{\calG\acts N}{N}$.
\bgn{proposition}\label{actlagpd} $\sf{\Ohm\acts B}$ is an \tla-groupoid if{f} $\pmb\sigma$ is morphic.
\end{proposition}
\bgn{proof} The implication to the right is clear, since $\sigma$ is the target top map of ${\sf{\Ohm\acts B}}$. Suppose that
$\pmb\sigma$ is morphic. According to proposition (\ref{morphiccat}), to have  $\sf{\Ohm\acts B}$ a groupoid
object, it suffices to show that the fibered product Lie algebroid $\Ohm\acts B\fib{\sh^\acts}{\sh^\acts}\Ohm\acts
B$ exists;  in that case the diagram (\ref{exagon}) commutes in the category of Lie algebroids. By proving that
$\sh^\acts$ is an \tla-fibration over  $\sor^\acts$, we also show that $\sf{\Ohm\acts B}$ is an \tla-groupoid, i.e.
$\sh^\acts$ is fibrewise surjective. 
In order to check fibrewise submersivity of $\sh^\acts$  one can first fix $(g,n)\in\calG\acts N$, $\ohm\in\Ohm_g$ and 
$b\in B_n$ with $\sh(\ohm)=\jh(b)$, then repeat the argument above for the tangent maps.
\end{proof}
We shall call $\sf{\Ohm\acts B}$ the \tsf{action \tla-groupoid} associated with
the morphic action. An \tla-groupoid shall be called \tsf{$($isotropy$)$ free}, respectively \tsf{proper}, if so are its top
and side groupoids. 
\bgn{remark} %
 For any morphic action of an \tla-groupoid the associated action \tla-groupoid is free,  respectively proper, if{f} so
are the top and side groupoid actions.  
\end{remark}
\spa The following corollary, which we get for free from proposition  \ref{actlagpd}, is important in the applications.
\bgn{corollary}\label{strong} Any morphic action $\pmb\sigma:\sf{\Ohm\acts B}\rightarrow\sf{B}$ is a strong 
\tla-fibration. 
\end{corollary}
\bgn{proof} $\pmb\sigma$ is the top target morphism of an \tla-groupoid. 
\end{proof}
%
Let us describe the prototypical example of an action $\mathcal{LA}$-groupoid. 
\bgn{example}\label{fpla}
Let $\gpdm$ be a Lie groupoid, $j:N\rightarrow M$ a smooth map and $\sigma: \calG\fib{\sor}{j}N\rightarrow N$ a
smooth groupoid action. Then the tangent action map  $\dd\sigma:T( \calG\fib{\sor}{j}N)=
T\calG\fib{\dd\sor}{\dd j}TN\to TN$ is a morphic action of the tangent  prolongation groupoid  $\poidd{T\calG}{TM}$ on $\dd
j: TN\rightarrow TM$. The associated action \tla-groupoid
$$
\bcat
\qquad
\xy
*+{}="0",    <-1.1cm,0.7cm>
*+{T\calG\acts TN}="1", <0.9cm,0.7cm>
*+{TN}="2", <-1.1cm,-0.7cm>
*+{\calG\acts  N}="3", <0.9cm,-0.7cm>
*+{N}="4",
\ar  @ <-0.07cm>   @{->} "1";"2"^{} 
\ar  @ <0.07cm>    @{->} "1";"2"^{}  
\ar  		   @{->} "1";"3"^{}
\ar                @{->} "2";"4"_{}
\ar  @ <0.07cm>    @{->} "3";"4"^{}
\ar  @ <-0.07cm>   @{->} "3";"4"_{}
\endxy
\qquad,
\ecat
$$
is also the tangent prolongation \tla-groupoid of the action groupoid $\calG\acts N$. 
\end{example}
\bgn{remark}
Corollary \ref{strong}, implies that the tangent lifted action map is fibrewise surjective; we shall use this fact
in the next Subsection to perform the cotangent lift of a Lie groupoid action. This fact can be also checked
directly: for all $\delta x\in T_{g*p} P$, pick any bisection $\Sigma_g$ through $g$, thus $\delta x=\delta
g*\delta p$, with  $\delta g:=\dd\calR^{\Sigma_g}\dd j(\delta x)$ and $\delta p=\delta g\inverse *\delta x$.
\end{remark}
\spa Next we shall introduce a characteristic morphism associated with action
\tla-groupoids,  the \tsf{moment morphism} 
$\sf J:\sf{\Ohm\acts B}\rightarrow \sf{\Ohm}$,
\bgn{equation}\label{mommmo}
\qquad
\bcat\xy
*++{}="0",    <-1cm,0.7cm>
*++{\Omega\acts B}="1", <0.7cm,0.7cm>
*++{B}="2", <-1cm,-0.7cm>
*++{\calG\acts N}="3", <0.7cm,-0.7cm>
*++{N}="4",     <2.5cm,-0.7cm>
*++{\Omega}="1'", <4.2cm,-0.7cm>
*++{A}="2'", <2.5cm,-2.1cm>
*++{\calG}="3'", <4.2cm,-2.1cm>
*++{M}="4'",  <5.2cm,-0.7cm>
*++{,}="2''"
\ar  @ <-0.07cm>   @{->} "1";"2"^{} 
\ar  @ <0.07cm>    @{->} "1";"2"^{}  
\ar  		   @{->} "1";"3"^{}
\ar                @{->} "2";"4"_{}
\ar  @ <0.07cm>    @{->} "3";"4"^{}
\ar  @ <-0.07cm>   @{->} "3";"4"_{}
\ar  @ <-0.07cm>   @{->} "1'";"2'"^{} 
\ar  @ <0.07cm>    @{->} "1'";"2'"^{}  
\ar  		   @{->} "1'";"3'"^{}
\ar                @{->} "2'";"4'"_{}
\ar  @ <0.07cm>    @{->} "3'";"4'"^{}
\ar  @ <-0.07cm>   @{->} "3'";"4'"_{}
\ar  		   @{->} "1";"1'"^{\:\:\pr_\Ohm}
\ar  		   @{->} "2";"2'"^{\quad\jh}
\ar  		   @{->} "3";"3'"^{\:\:\pr_\calG}
\ar  		   @{->} "4";"4'"^{\quad j}
\endxy\ecat
\end{equation}
where $\pr_\Ohm$ and $\pr_\calG$ denote the projections to the first components.
According  to lemma \ref{kerla}, the regularity of the kernel $\sf{K}^\acts$ of $\sf J$ is
controlled by that of 
\be
\qquad
(\pr_\Ohm,\jh):\Ohm\acts B\rightarrow \Ohm\fib{\sh}{\jh}B
\qquad,
\ee
which is the identity map in this case. Therefore, whenever $\jh$ is an 
\tla-fibration, 
\bgn{equation}\label{cocact}
{\sf
 K^\acts}:=
\bcat\xy
*+{}="0",    <-0.7cm,0.7cm>
*+{\hat{K}}="1", <0.7cm,0.7cm>
*+{ K}="2", <-0.7cm,-0.7cm>
*+{\calG}="3", <0.7cm,-0.7cm>
*+{M}="4",
\ar  @ <-0.07cm>   @{->} "1";"2"^{} 
\ar  @ <0.07cm>    @{->} "1";"2"^{}  
\ar  		   @{->} "1";"3"^{}
\ar                @{->} "2";"4"_{}
\ar  @ <0.07cm>    @{->} "3";"4"^{}
\ar  @ <-0.07cm>   @{->} "3";"4"_{}
\endxy\ecat
\qquad,
\end{equation}
is a sub-\tla-groupoid of $\sf{\Ohm\acts B}$, where $\hat{K}$ and $K$ are the wide
subalgebroids on  $\Pr_K:\ker\,\pr_\Ohm\rightarrow \calG\acts N$ and 
$\pr_K:\ker\,\jh\rightarrow N$ and it is sufficient to have $\jh$ a map of constant
rank for $\sf{K^\acts}$ to be an \tla-groupoid. Actually, $\hat{K}$ is to be identified with the
pullback bundle $\sor_\acts^{\pmb{+}}\ker\,\jh\rightarrow\calG\acts N$ and, as a Lie groupoid,
it is to be thought of as an action groupoid $\calG\acts K$. 
In fact, the action
of $\Ohm$ on $B$ restricts to an action 
$\sigma_K:\calG\fib{\sor}{\pr\comp\jh}K=\calG\fib{\sor}{j\comp\pr_K}K
\rightarrow K$, $(g,k)\mapsto g*_Kk$,
$$
\qquad
g*_Kk:=0^\Ohm_g* k=\sigmah(0^\Ohm_g,k)\qquad,
$$
of $\calG$ on $K$; note that $\sigma_K$ is well defined, since
$\jh(g*_Kk)=\th(0^\Ohm_g)$ vanishes by fibrewise linearity of $\th$. The action map $\sigma_K$
is indeed the restriction of $\sigmah$ to the wide Lie subalgebroid 
$0^\Ohm_\calG\fib{\sh}{\jh}K\subset\Ohm\fib{\sh}{\jh}B$ and, in this sense, the
action of $\calG$ on $K$ is fibrewise linear.
\vs{0.5}
\subsection{The cotangent lift of a Poisson groupoid action}\hfill

\vs{0.1}
\spa If a Poisson groupoid acts on a Poisson manifold, the domain of the action map does not carry any natural
Poisson bivector and the compatibility of the groupoid action with the Poisson structures cannot be formulated in the
requirement that the action map be Poisson. There is however a natural replacement for this condition.  
\bgn{definition}
Let $\gpdm$ be a Poisson groupoid act on a moment map  $j:P\rightarrow M$,  where $P$ is a Poisson manifold.
A \tsf{Poisson action} of $\calG$ on $j$ is a groupoid action $\sigma:\calG\fib{\sor}{j}P\rightarrow P$, such that
\bgn{equation}\label{compact}
\gr{\sigma}\subset \calG\times P\times\ol{P}\qquad\hbox{is coisotropic.} 
\end{equation}
For any Poisson action such as above $j:P\to M$ is called a \tsf{Poisson $\calG$-space}. 
\end{definition}
Any Lie groupoid action is a Poisson action for the zero Poisson structures. When $\calG$ is a Lie
group $G$ and $j=P\to\bullet$ the compatibility conditions amounts to asking the action map 
$\sigma:G\times P\to P$ to be Poisson; a Poisson action of a symplectic groupoid on a symplectic manifold is a
\emph{symplectic groupoid action} in the sense of Mikami and Weinstein \cite{mw88}, since the graph of the action map
must be Lagrangian by dimensional constraints.\\  
\bgn{remark} It follows from the definition and from proposition \ref{holz} that the moment map $j:P\to M$ has to be
anti-Poisson for the Poisson bivector induced by $\calG$ on $M$. Note that this is the case, e.g. for the action of a
Poisson groupoid on itself by left translation: $j\equiv\tar:\calG\to M$ and $\sigma\equiv\mu:\calG^{(2)}\to\calG$.
\end{remark}
To any Lie groupoid action of $\gpdm$ on $j:N\to M$, where $N$ is any manifold, one can associate a Poisson groupoid
action of the cotangent prolongation groupoid $\poidd{T^*\calG}{A^*}$ on the moment map
$\jh:T^*N\to A^*$, defined by
\bgn{equation}\label{momo}
\pair{\jh(\alpha_n)}{a_{j(n)}}
:=
\pair{\alpha_n}{\dd\sigma(a_{j(n)}, 0_n)}
\quad,\quad \alpha_n\in
T^*_nN
\quad,\quad 
a_{j(n)}\in A_{j(n)}\quad,
\end{equation}
i.e. dualizing the infinitesimal action of $A$; note that for all sections $a\in\Gamma(A)$ and $\alpha\in\Ohm^1(N)$,
$$
\pair{\jh\comp\alpha}{a\comp j}
:=
\pair{\alpha}{X^a}
$$
for the infinitesimal action $X^\bullet:\Gamma(A)\to\frax(N)$, or, in terms of fibrewise linear functions, 
$\jh^*F_a=F_{{X^{a}}}$.
\bgn{lemma}\label{op} The moment map $\jh:T^*N\rightarrow A^*$ defined above is 
a Poisson map for the dual Poisson structures; moreover, if the action is locally free
$\jh$ has maximal rank.
\end{lemma}
\bgn{proof}
The pullback map $\jh^*:\cif(A^*)\to\cif(T^*N)$ maps fibrewise linear, respectively fibrewise constant, functions to
fibrewise linear, respectively fibrewise  constant, functions. It suffices to show that
$\jh^*\poib{F}{G}=\poib{\jh^*F}{\jh^* G}$ in the cases %
\medskip\\
 ($i$) $F=F_{a_+}$ and $G=F_{a_-}$, $a_\pm\in\Gamma(A)$,
\medskip\\
 ($ii$) $F=F_{a}$ and $G=\pr^* f$, $a\in\Gamma(A)$, $f\in\cif(M)$,
\medskip\\
 ($iii$) $F=\pr^* f_+$ and $G=\pr^* f_-$, $f_\pm\in\cif(M)$,
\medskip\\
the result follows by the Leibniz rule. ($i$) We have
\be
\poib{\jh^*F_{a_+}}{\jh^*F_{a_-}}
&=&
\poib
{F_{X^{a_+}}}
{F_{X^{a_-}}}
=
F_{
\brak
{X^{a_+}}
{X^{a_-}}
}
=
\jh^*F_{\brak{a_+}{a_-}}\\
&=&
\jh^*\poib{F_{a_+}}{F_{a_-}}\hs{5.5}.
\ee
Since for all $n\in N$,
$$
\quad
X^{a}_{n}(j^* f)=\pair{\dd f_{j(n)}}{\dd j\dd\sigma(a_{j(n)}, 0_n)}=\rho_{j(n)}(a)(f)
\quad,
$$
($ii$) follows:
\be
\poib{\jh^*F_a}{\jh^*\pr^* f}
&=&\poib{F_{{X^{a}}}}{\pr^* j^* f}=\pr^*({X^{a}}(j^*f))=j^*(\rho(a)(f))\\
&=&\jh^*\poib{F_a}{\pr^* f}
\hs{7}.
\ee
Condition ($iii$) holds trivially, since both sides vanish. The second part of the statement is 
a direct consequence of the definition (\ref{momo}) of
$\jh$.
\end{proof}
We are now ready to dualize the tangent lift of a groupoid action.
\bgn{theorem}\label{memento} Let $\gpdm$ be a Lie groupoid and $\sigma:\calG\fib{\sor}{j}N\rightarrow N$ a Lie
groupoid action on $j:N\rightarrow M$. Then $\sigma$ lifts to a Poisson action of $\poidd{T^*\calG}{A^*}$ on 
$\jh:T^*N\to A^*$ with action map  $\sigmah:T^*\calG\fib{\sh}{\jh}T^*N\rightarrow T^*N$, $(\theta_g,\alpha_p)\mapsto
\theta_g\hhh\alpha_p$. The cotangent lifted action is uniquely determined in terms of the tangent lifted action $*$
by the formula
\bgn{equation}\label{ctglift}
\langle\,\theta_g\hhh\alpha_p,\delta g*\delta p\,\rangle
=
\langle\,\theta_g\,,\,\delta g\,\rangle
+
\langle\,\alpha_p\,,\,\delta p\,\rangle
\end{equation}
for all $(\theta_g,\alpha_p)\in T^*\calG\fib{\sh}{\jh}T^*N$ and $(\delta g, \delta n)\in
T\calG\fib{\dd\sor}{\dd j} TN$.
\end{theorem}
\bgn{proof} First of all note that $N^*\gr{\sigma}$ is the graph of a map. Consider that the restriction of the
first projection  $(T^*\calG\times T^*N)\times T^*N\rightarrow T^*\calG\times T^*N$ to $N^*\gr{\sigma}$ takes values 
in $T^*\calG\fib{\sh}{\jh}T^*N$. For all $a_m\in A$ with $m=\sor(g)=j(n)$, the pair
\be
\delta g&:=&\dd l_g(a_m-\dd\eps\rho(a_m))\\
\delta n&:=&\dd\sigma(-a_m,0_m)
\ee
belongs to the domain of the tangent action map,
\be
\dd\sor\,\delta g &=&\dd\sor(a_m-\dd\eps\rho(a_m))\\
&=&-\rho(a_m) \equiv \dd\tar( -a_m)\\
&=&\dd j\,\dd\sigma (-a_m,0_m)\qquad,
\ee
and
\be
\delta g*\delta n&=&[0_g\bullet(a_m-\dd\eps\rho(a_m))]*[(-a_m)*0_n)]\\
&=&0_g*[(a_m-\dd\eps\rho(a_m))\bullet(0_m-a_m)]*0_n\\
&=&0_g*[a_m\bullet 0_m-\dd\eps\rho(a_m)\bullet a_m]*0_n\\
&=&0_g*[a_m-a_m]*0_n\\
&=&0_{g*n}\qquad.
\ee
Therefore, for all $(\theta_g,\alpha_p,\beta_{g*n})\in N^*\gr{\sigma}$, and $a_n\in A$
\be
\pair{\sh(\theta_g)-\jh(\alpha_n)}{a_m}
&=&
\pair{\theta_g}{l_g(a_m-\dd\eps\rho(a_m)}-\pair{\alpha_p}{\dd\sigma(-a_m,0_m)}\\
&=&
\pair{\theta_g}{\delta g}+\pair{\alpha_n}{\delta n}\\
&=&
-\pair{\beta_{g*n}}{\delta g*\delta n}\\
&=&0\qquad\qquad.
\ee
The restriction $\lambda: N^*\gr{\sigma}\to T^*\calG\fib{\sh}{\jh}T^*N$ of the projection above is an isomorphism of
vector bundles: by counting dimensions
$$
\qquad
\rk\,N^*\gr{\sigma}=\dim N +\dim M=
\rk\, T^*\calG\fib{\sh}{\jh}T^*N
\qquad,
$$
moreover
$$
\quad
\ker_{(g,p,g*n)}\lambda=\{(0_g,0_n,\alpha_{g*n})\:\:|\:\:\
\alpha_{g*n}(\delta\,g*\delta n)=0\:\:,\:\:
(\delta\,g,\delta n)\in T(\calG\acts N)\}=\{0\}\quad, 
$$
by fibrewise surjectivity of the tangent lifted action map.
Hence, one can set
$$
\gr{\sigmah}:=(\id_{T^*\calG}\times\id_{T^*N}\times -\id_{T^*N})\,N^*\gr{\sigma}
      	\subset{T^*\calG}\times{T^*N}\times{T^*N}
$$
to define a bundle map $\sigmah$ over $\sigma$ satisfying (\ref{ctglift}); the graph of $\sigmah$ is Lagrangian in
${T^*\calG}\times{T^*N}\times{\ol{T^*N}}$ by construction.  The properties of an action map follow: the following
expressions are to be understood whenever they make sense\medskip\\
$Compatibility\ with\ the\ moment\ map:$ 
\be
\pair{\jh(\theta_g\hhh\alpha_{n})}{a_{j(g*n)}}
&=&
\pair{\theta_g\hhh\alpha_n}{\dd\sigma (a_{j(g*n)},0_{g*n})}\\
&=&
\pair{\theta_g\hhh\alpha_n}{(\dd\,r_g\, a_{j(g*n)})*0_n}\\
&=&
\pair{\theta_g}{\dd\,r_g\, a_{j(g*n)}}+\pair{\alpha_n}{0_n}\\
&=&
\pair{\th(\theta_g)}{a_{j(g*n)}}\hs{2},
\ee
where we regarded $a_{j(n)}\in A_{j(n)}$ as a vector tangent to $\calG$ in the second and third step;
\medskip\\
$Unitality:$ since $\epsh$ is the identification $A^*\simeq N^*M$,
\be
\pair{\epsh(\jh(\alpha_n))\hhh\alpha_n}{\delta n}
&=&
\pair{\epsh(\jh(\alpha_n))\hhh\alpha_n}{\dd\eps(\dd j(\delta n))*\delta n}\\
&=&
\pair{\epsh(\jh(\alpha_n))}{\dd\eps(\dd j(\delta n))}
+
\pair{\alpha_n}{\delta n}\\
&=&
\pair{\alpha_n}{\delta n}\hs{0.5};
\ee
$Multiplicativity:$ 
\be
\pair{\theta_g\hhh(\theta_h*\alpha_n)}{\delta g*\delta h*\delta n}
&=&
\pair{\theta_g}{\delta g}
+
\pair{\theta_h}{\delta h}
+
\pair{\alpha_p}{\delta p}\\
&=&
\pair{\theta_g\hdot \theta_h}{\delta g\bullet \delta h}
+
\pair{\alpha_p}{\delta n}\\
&=&
\pair{(\theta_g\hdot \theta_h)\hhh\alpha_n}{\delta g\bullet \delta h*\delta n}\\
&=&
\pair{(\theta_g\hdot \theta_h)\hhh\alpha_n}{\delta g*\delta h*\delta n}\hs{2},
\ee
then $(\theta_g\hdot \theta_h,\alpha_p ,(\theta_g\hdot \theta_h)*\alpha_p)$ and 
$(\theta_g\hdot \theta_h,\alpha_p ,\theta_g*(\theta_h*\alpha_p))$ are the same
covector, since the tangent lifted action is fibrewise surjective. 
\end{proof}
\bgn{remark} The canonical symplectic form on $T^*\calG$ induces the opposite dual Poisson structure on $A^*$; then
$\jh$ is indeed anti-Poisson as a moment map for the cotangent lifted action.
\end{remark}
Note that the moment map $\jh$ is a bundle map over $j$ and, by construction, the cotangent lifted
action map is linear over the original action map; that is, one can regard $\sigmah$ as a  morphic
action over $\sigma$, for the abelian Lie algebroids on $T^*\calG$, $A^*$ and $T^*P$. When the
cotangent lift is applied to a Poisson groupoid action this fact still holds true.
\bgn{proposition}\label{holz} Let $\gpdm$ be a Poisson groupoid acting on  $j:P\rightarrow M$ with action map 
$\sigma:\calG\fib{\sor}{j}P\rightarrow P$. Let  $\jh:T^*P\rightarrow A^*$
denote the moment map of the cotangent lifted action. Then, if $\sigma$ is a Poisson action
\medskip\\
$i)$ The moment map $\jh$ is a morphism of algebroids $T^*P\rightarrow A^*$;
\medskip\\
$ii)$ The cotangent lifted action $\sigmah$ is morphic over $\sigma$. 
\end{proposition}
($i$) was proved in \cite{hlz01} and together with lemma \ref{op} implies that $\jh$ is a morphism
of Lie bialgebroids (this was also remarked in \cite{hlz01}). We shall give a simpler proof of 
this fact.
%
\bgn{proof}
By construction  $\gr{\sigmah}\subset{T^*\calG}\times{T^*P}\times{T^*P}$ is a Lie subalgebroid whenever
$\gr{\sigma}\subset\calG\times P\times \ol{P}$ is coisotropic. We may identify $\gr{\jh}$ with the subbundle
$\{\epsh(\jh(\alpha)),\alpha,\alpha)\}\subset\gr{\sigmah}$, which is the preimage of the Lie subalgebroid $A^*\times
T^*P$ of $T^*\calG\times T^*P$ under the projection on the first two factors $\gr{\sigmah}\to T^*\calG\times T^*P$.
Thus it follows from corollary \ref{subalsubal} that  $\gr{\jh}$ is canonically isomorphic to a Lie subalgebroid of
$\gr{\sigmah}$, therefore to a Lie subalgebroid of ${T^*\calG}\times{T^*P}\times{T^*P}$. Apply corollary
\ref{caraci}  to conclude that $\jh$ is a morphism of Lie algebroids. It follows that
${T^*\calG}\fib{\sh}{\jh}{T^*P}$ carries a fibred product Lie algebroid making  $\sigmah$ a morphism of Lie
algebroids.
\end{proof}
As a consequence, to every Poisson action of a Poisson groupoid one can associate an action \tla-groupoid
\bgn{equation}\label{lagpdpoissonaction}
\bcat
\xy
*+{}="0",    <-1.1cm,0.8cm>
*+{T^*\calG\acts T^*P}="1", <1.1cm,0.8cm>
*+{T^*P}="2", <-1.1cm,-0.8cm>
*+{\calG\acts  P}="3", <1.1cm,-0.8cm>
*+{P}="4",
\ar  @ <-0.07cm>   @{->} "1";"2"^{} 
\ar  @ <0.07cm>    @{->} "1";"2"^{}  
\ar  		   @{->} "1";"3"^{}
\ar                @{->} "2";"4"_{}
\ar  @ <0.07cm>    @{->} "3";"4"^{}
\ar  @ <-0.07cm>   @{->} "3";"4"_{}
\endxy
\ecat\qquad.
\end{equation}
\bgn{example} In the case of a Poisson action of a Poisson group $G$, formulas
 (\ref{momo}) and
(\ref{ctglift}) allow to compute $\jh$ and $\sigmah$ explicitly:
$$
\qquad
\jh(\alpha_p)={\sigma_p}^*\alpha_p
\qquad\hbox{ and }\qquad
\theta_g\hhh\alpha_p=\sigma_{g\inverse}{}^*\alpha_p\qquad,
$$
for all $\alpha_p\in T^*P$ and $\theta_g\in T^*_gG$ with
$r_g^*\theta_g=\sigma^{p}{}^*\alpha_p$. The associated action \tla-groupoid was
constructed in \cite{mkz00b}.
\end{example}
\bgn{remark} Note that the action of a Poisson groupoid $\calG$ on $j:P\to M$ is Poisson if{f} 
the fibred product $T^*\calG\fib{\sh}{\jh}T^*P$ carries a Lie groupoid over $T^*P$ making (\ref{lagpdpoissonaction})
an \tla-groupoid.
\end{remark}
\vs{1}
\section{Reduction of morphic actions of \tla-groupoids}\label{rmalag}
\begin{quotation} 
This Section is devoted to the study of the reduction of morphic actions of \tla-groupoids. We prove (theorem
\ref{generalred}) that quotients with respect to  suitably free and proper such actions are always Lie algebroids.
Thereafter we discuss the reduction, in the spirit of a categorification of Marsden-Weinstein zero level reduction,
of the moment morphism (of \tla-groupoids) canonically associated with a morphic action (theorem
\ref{marwascat}). In the special case of the cotangent lift of a Poisson groupoid action, the
reduction of the moment morphism produces the Koszul algebroid of the quotient Poisson bivector 
(proposition \ref{poissonreduction}).  
\end{quotation}
%
%
%
%
%
%
The main purpose of this Section is to prove the following general reduction result:\nolinebreak
\bgn{theorem}\label{generalred} Let $(\Ohm,\calG;A,M)$ be an \tla-groupoid acting morphically on a morphism of Lie
algebroids $\jh:B\to A$ over $j:N\to M$. If the action is free and proper (so that $B/\Ohm$ and $N/\calG$ are
smooth manifolds), then there exists a unique Lie algebroid on $B/\Ohm\to N/\calG$ making the quotient
projection $B\to B/\Ohm$ a strong  \tla-fibration over $N\to N/\calG$. 
\end{theorem}
We split the proof in a topological part (Subsection \ref{ss1}), where we push the vector bundle forward over the
quotient, and in an algebraic part (Subsection \ref{ss2}), where we push down the  Lie brackets and anchor. In order
to achieve the second goal we characterize the Lie-Rinehart algebra on the quotient bundle, equivalent to a Lie
algebroid, as a Lie-Rinehart algebra naturally associated with the \tla-groupoid $\sf{\Ohm\acts B}$ encoding the 
action. Namely, sections of the top Lie algebroid of $\sf{\Ohm\acts B}$, which are also functors for the horizontal
groupoids (morphic sections), form a Lie-Rinehart algebra which descends to the quotient by the natural 
equivalence relation given by (categorical) natural transformations. The space of morphic sections modulo equivalence
is isomorphic to the space of projectable sections of $B\to N$ modulo equivalence, therefore, under 
the regularity assumptions on the top and side actions, sections of the quotient inherit a Lie-Rinehart algebra.
Even when the top and side actions are neither free nor proper, but in particular in the non free case,
when the quotient bundle lives in the category of stratified manifolds, one can regard the Lie-Rinehart algebra of 
morphic sections modulo equivalence as a desingularization or a model of a Lie algebroid over the pathological quotient
bundle.\\
Provided only the side action is free and proper, the kernel of the moment morphism behaves however well with respect
to reduction under natural regularity assumptions on the top moment map; in fact the restriction of the top action to  the
kernel \tla-groupoid is essentially a fibrewise linear lift of the side action.\eject
Remarkably, the  moment morphism associated with the cotangent lift of a Poisson
$\calG$-space  always satisfies the regularity requirement and its ``kernel reduction'' procedure yields a quotient Lie
algebroid which is canonically isomorphic to the Koszul algebroid of the quotient Poisson bivector. This fact shall
be used in the last Section to derive applications of the reduction procedure developed here to the integrability of quotient Poisson
manifolds.

\vs{0.5}
\subsection{Free and proper morphic actions: topological reduction}\label{ss1}\hfill

\vs{0.1}
\spa  For any morphic action of an \tla-groupoid $(\Ohm,\calG;A,M)$ on a morphism of Lie algebroids $\jh:B\to A$ over
$j:N\to M$, by linearity of the top action map,  the vector bundle projection $\pr_B:B\to N$ is equivariant under the
actions of $\calG$ and $\Ohm$, i.e. $\pr_B(\ohm_g*b_n)=g*n$ for all  $(\ohm_g,b_n)\in\Ohm\fib{\sh}{\jh}B$. Then,
whenever the quotients $B/\Ohm$ and $N/\calG$ are smooth, $\pr_B$ descends to a smooth map $\ul{\pr}_B:B / \Ohm\to N
/ \calG$.  We shall show that  $\ul{\pr}_B$
is actually a vector bundle under the natural
hypotheses.
\bgn{proposition}\label{topo} Let $(\Ohm,\calG;A,M)$ be an \tla-groupoid acting morphically on a morphism of Lie
algebroids $\jh:B\to A$ over $j:N\to M$. If both the top and side actions are free and proper, 
there exists a unique vector bundle of rank  $\rk B +\rk A -\rk \Ohm$ on the induced map
$\ul{\pr}_B:B/\Ohm\to N/\calG$ making the quotient projection $B\to B/\Ohm$ a strong 
\tvb-fibration over $N\to N/\calG$.
\end{proposition}
According to proposition \ref{actlagpd} to a morphic action such as that in the statement we can
associate a free and proper action \tla-groupoid $(\Ohm\acts B,\calG\acts N;B,N)$. Thus we have to
show that for any free and proper \tla-groupoid $(\Ohm,\calG;A,M)$ the induced map
$\ul{\pr}_A:A/\Ohm\to M/\calG$ carries a vector bundle making the quotient projection 
$A\to A/\Ohm$ a \tvb-fibration, where the quotients are taken for the actions by left translation.\\
We split the proof in three parts:
\medskip\\
\emph{ Step 1.} (Lemma \ref{Step 1.}) We mod out by the action of $\Ohm$ along the
fibres of $A$, which is governed by a wide normal subgroupoid $\Ohm\vip$ of $\poidd{\Ohm}{A}$;
\medskip\\
\emph{ Step 2.} (Lemma \ref{Step 2.}) We show that  $A\to M$ and $\Ohm\to\calG$ descend to vector 
bundles $A/\Ohm\vip\to M$ and $\Ohm/\Ohm\vip\to\calG$;
\medskip\\
\emph{ Step 3.} (Conclusion of the proof) We use the groupoid $\poidd{\Ohm/\Ohm\vip}{A/\Ohm\vip}$
to push $A/\Ohm\vip\to M$ forward to a vector bundle over $M/\calG$ by identifying the fibres along 
the orbits of $\calG$. 
\medskip\\
The quotient vector bundle is then essentially obtained via a reduction 
by stages
$$
\qquad
A/\Ohm=\frac{A/\Ohm\vip}{\Ohm/\Ohm\vip}\qquad;
$$
we shall make this precise in the conclusion of the proof.
\bgn{lemma}\label{Step 1.} Let $(\Ohm,\calG;A,M)$ be an \tla-groupoid
and $\Ohm\vip:=\Pr^{-1}(\eps(M))$ be the restriction of $\Ohm\to\calG$ over $M$. 
Then%
\medskip\\
$i)$ $\poidd{\Ohm\vip}{A}$ is a wide normal Lie subgroupoid of $\Ohm$
\medskip\\
and, if $\poidd{\Ohm}{A}$ is isotropy free
\medskip\\
$ii$) The quotient groupoid $\poidd{\Ohm/\Ohm\vip}{A/\Ohm\vip}$ is an isotropy free Lie groupoid.
\end{lemma}
Note that the quotient is taken by the equivalence relation induced by a normal subgroupoid, \emph{not} by the action
of $\Ohm\vip$ on $\Ohm$ by left translation. On the other hand the quotient $\Ohm/\Ohm\vip$ coincides with the orbit
space of the natural action  $({\ohm}\vip,\ol{\ohm}\vip;\ohm)\mapsto{\ohm}\vip\cdot\ohm\cdot\ol{\ohm}\vip$, of 
$\poidd{\Ohm\vip\times\op\Ohm\vip}{A\times A}$ on $\chih:\Ohm\to A\times A$; it turns out that this action is  free
and proper, whenever $\poidd{\Ohm}{A}$ is isotropy free and proper.
\bgn{proof} ($i$) $\poidd{\Ohm\vip}{A}$ is the kernel groupoid of a strong \tlg-fibration, apply
corollary \ref{pregpd}. ($ii$) The action of $\Ohm$ on $A$ restricts to a free and proper action of $\Ohm\vip$, thus
$A/\Ohm\vip$ is smooth. First we show that $\Ohm/\Ohm\vip$ is also smooth. The action of 
$\Ohm\vip\times\op\Ohm\vip$is free: if 
${\ohm}\vip_+\cdot\ohm\cdot\ol{\ohm}\vip_+={\ohm}\vip_-\cdot\ohm\cdot\ol{\ohm}\vip_-$, 
$\sh({\ohm}\vip_+)=\th(\ohm)=\sh({\ohm}\vip_-)$ and $\th({\ohm}\vip_+)=\th({\ohm}\vip_-)$, 
thus ${\ohm}\vip_+={\ohm}\vip_-$, by injectivity of $\chih$. The same argument applies for immersivity
of the anchor of $(\Ohm\vip\times\op\Ohm\vip)\acts \Ohm$; that is, the action is free.
To see that the action of $\Ohm\vip\times\op\Ohm\vip$ is proper we have to show that for all sequences
$\{({\ohm}\vip_k,\ol{\ohm}\vip_k;\ohm_k)\}\in(\Ohm\vip\times\op\Ohm\vip)\fib{\sh\times\op\sh}{\chih}\Ohm$ such 
that $\{\ohm_k\}$ converges to some $\ohm_\infty$ and $\{{\ohm}\vip_k\cdot\ohm_k\cdot\ol{\ohm}\vip_k\}$ converges to
some $\wt{\ohm}_\infty$, both $\{{\ohm}\vip_k\}$ and $\{\ol{\ohm}\vip_k\}$ have convergent subsequences; since
\be
\chih(\ol{\ohm}\vip_k)&=&(\sh(\ohm_k),\sh({\ohm}\vip_k\cdot\ohm_k\cdot{\ohm}\vip_k))\longrightarrow 
(\sh(\ohm_\infty), \sh(\wt{\ohm}_\infty))\hbox{ and }\\
\chih({\ohm}\vip_k)&=&(\th({\ohm}\vip_k\cdot\ohm_k\cdot{\ohm}\vip_k),\th(\ohm_k))\longrightarrow
(\th(\wt{\ohm}_\infty),\th(\ohm_\infty))\qquad\qquad,
\ee
that is true, by properness of $\chih$. The quotient groupoid by a normal subgroupoid is always well defined; 
the induced source map $\Ohm/\Ohm\vip\to A/\Ohm\vip$ is a surjective submersion since so are the quotient projection
$A\to A/\Ohm\vip$ and the source map $\Ohm\to A$, thus
$\poidd{\Ohm/\Ohm\vip}{A/\Ohm\vip}$ is a Lie groupoid. Next we show that $\poidd{\Ohm/\Ohm\vip}{A/\Ohm\vip}$
is isotropy free. Assume $\ul{\chi}(\ul{\ohm})=\ul{\chi}(\ul{\wt{\ohm}})$, where 
$\ul{\chi}:\Ohm/\Ohm\vip\to A/\Ohm\vip\times A/\Ohm\vip$ is the induced groupoid anchor, and pick any
representatives $\ohm\in\ul{\ohm}$, $\wt{\ohm}\in\ul{\wt{\ohm}}$. Then
$$
\th(\wt{\ohm})={\ohm}\vip_{\tar}*\th(\ohm)=\th({\ohm}\vip_{\tar})
\qquad\qquad
\sh(\wt{\ohm})={\ohm}\vip_{\sor}*\sh(\ohm)=\th({\ohm}\vip_{\sor})
$$
for some ${\ohm}\vip_{\tar}$,${\ohm}\vip_{\sor}\in\Ohm\vip$
and we can form the composition
$x={\ohm}\vip_{\tar}{}\inverse\cdot\wt{\ohm}\cdot{\ohm}\vip_{\sor}\cdot\ohm\inverse$; since 
$\th(x)=\sh({\ohm}\vip_{\tar})=\th(\ohm)=\sh(x)$, $x$ must be a unit, thus 
$\wt{\ohm}\cdot{\ohm}\vip_{\sor}={\ohm}\vip_{\tar}\cdot\ohm$, i.e. $\ul{\wt{\ohm}}=\ul{\ohm}$.
\end{proof}
Note that, if $\sf{\Ohm}$ is an isotropy free \tla-groupoid, $\chih(\Ohm\vip)$ is a vector subbundle of $A\oplus A$,
since, for all $\ohm_g\in\Ohm_g$ $\th(\ohm_g)=\sh(\ohm_g)$ implies that $g$ is an isotropy, therefore a unit. 
Moreover, the graph of the equivalence relation induced by $\Ohm\vip$ on $\Ohm$ is a subbundle of $\Ohm\oplus\Ohm$:
$\Pr(\ul{\ohm}\vip\cdot\ohm\cdot{\ohm}\vip)=\Pr(\ohm)$, for all composable $\ul{\ohm}\vip$, ${\ohm}\vip\in\Ohm\vip$ 
and $\ohm\in\Ohm$. That is, both actions of $\Ohm\vip$ on $A$ and of $\Ohm\vip\times\op\Ohm\vip$ on $\Ohm$ are fibre
preserving.
\bgn{lemma} Let $\pr:E\to M$ be a vector bundle of rank $e$ and $\thicksim$ a regular  {\em \tsf{ linear
equivalence relation}} on $E$, in the sense that $\gr{\thicksim}\subset E\oplus E$ is a vector subbundle. Then
the induced projection $\ul{\pr}:E\mod\to M$ carries a unique (up to bundle isomorphisms) vector bundle making 
the quotient projection  $\pi:E\to E\mod$ a \tvb-fibration.
\end{lemma}
\bgn{proof} Fix any trivializing atlas $\{U_\alpha, \tau^\alpha\}$ for $E$. For all $u\in U_\alpha$, setting
$$
\quad
[e_+]+\lambda\cdot[e_-]=[e_++\lambda\cdot e_-]\qquad,\qquad e_\pm\in\pr\inverse(u)\quad,\quad
\lambda\in\RR
\quad, 
$$
endows $\ul{\pr}\inverse(u)$ with a well defined linear structure. The submersion
$\psi:=\pi\comp\tau^\alpha{}\inverse$ induces a regular equivalence relation on $U_\alpha\times\RR^e$, also denoted
by $\thicksim$, and it is easy to see that  $(U_\alpha\times\RR^e)/\thicksim=U_\alpha\times
\RR^e/\ker_u\,\psi=U_\alpha\times\RR^{e-k}$ for some $k\leq e$.
$$
\bcat
\xy
*+{}="0",    <-5cm,0.7cm>
*+{\pr^{-1}(U_\alpha)}="2", <1cm,0.7cm>
*+{U_\alpha\times\RR^e}="1", <-5cm,-0.7cm>
*+{\ul{\pr}^{-1}(U_\alpha)\equiv\pr^{-1}(U_\alpha)/\thicksim}="4", <1cm,-0.7cm>
*+{U_\alpha\times\RR^n/\thicksim{}\simeq U_\alpha\times\RR^{e-k}}="3",
\ar  @{->} "2";"1"^{\tau^\alpha}  
\ar  		   @{->} "1";"3"^{}
\ar                @{->} "2";"4"_{\pi}
\ar   @{->} "4";"3"_{\ul{\tau}^\alpha}
\endxy
\ecat
$$
One can then endow $\pr: E\mod\to M$ with a vector bundle structure by considering the maximal trivializing atlas
containing $\{U_\alpha,\ul{\tau}^\alpha\}$, where $\ul{\tau}^\alpha:\ul{\pr}\inverse(U_\alpha)\to
U_\alpha\times\RR^{e-k}$, $\ul{\tau}^\alpha([e]_{E\mod}):=[\tau^\alpha(e)]_{(U_\alpha\times\RR^e)/\thicksim}$ is a
well defined fibrewise linear diffeomorphism. Uniqueness is clear.
\end{proof}
The equivalence relations on $A$ and $\Ohm$ induced by $\Ohm\vip$ are both linear in the sense of last lemma, 
due to linearity of the top groupoid multiplication of $\Ohm$;  we shall denote with $\mathrm{q}$ and $\qh$ the
projections $A/\Ohm\vip\to M$ and $\Ohm/\Ohm\vip\to \calG$, respectively. Then we have
\bgn{lemma}\label{Step 2.} For any isotropy free and proper \tla-groupoid $(\Ohm,\calG;A,M)$, the following statements 
hold:
\medskip\\ $i)$ The projection $q:A/\Ohm\vip\to M$ carries a vector bundle making the quotient map
$A\to A/\Ohm\vip$ a \tvb-fibration;
\medskip\\ $ii)$ The projection $\qh:\Ohm/\Ohm\vip\to\calG$ carries a vector bundle making the quotient map 
$\Ohm\to \Ohm/\Ohm\vip$ a \tvb-fibration;
\medskip\\ $iii)$ Let $\ul{\sor}$ and $\ul{\tar}$ denote source and target of
$\poidd{\Ohm/\Ohm\vip}{A/\Ohm\vip}$, then the induced maps 
$\ul{\sor}^!:\Ohm/\Ohm\vip\to \sor\pbk (A/\Ohm\vip)$ and 
$\ul{\tar}^!:\Ohm/\Ohm\vip\to \tar\pbk (A/\Ohm\vip)$ are vector bundle isomorphisms.
\end{lemma}
\bgn{proof} It remains to prove the last statement. Note that $\ul{\sor}^!$ is fibrewise surjective, since so are 
$\sh$ and the quotient projections. For $\ohm_g$, $\wt{\ohm}_g\in\Ohm$ such that 
$\ul{\sor}(\ul{\ohm}_g)=\ul{\sor}(\ul{\wt{\ohm}}_g)$, there is some $\ohm\vip\in\Ohm\vip$ such that
$x=\wt{\ohm}_g\cdot\ohm\vip{\inverse}\cdot\ohm_g\inverse$ is defined; since $\Pr(x)=g\cdot g\inverse$, $x\in\Ohm\vip$
and $\ul{\ohm}_g=\ul{\wt{\ohm}}_g$, that is, $\ul{\sor}^!$ is fibrewise a linear isomorphism. The same reasoning
shows that $\ul{\tar}$ is a bundle isomorphism. 
\end{proof}
We are ready to conclude the proof of  proposition \ref{topo}.
\bgn{proof}[End of proof of proposition \ref{topo}] It is always possible to define a section
$\theta\in\Gamma(\calG,\tar\pbk (A/\Ohm\vip)\oplus\sor\pbk (A/\Ohm\vip)^*)$ by setting
$$
\quad
\langle\theta,\xi^\tar\oplus \ul{a}\rangle
:=
\langle\xi^\tar,\ul{\tar}^!\comp\ul{s}^!{\inverse}\comp\ul{a}\rangle
\quad,\quad 
\xi^\tar\oplus \ul{a}\in\Gamma(\calG,\tar\pbk (A/\Ohm\vip)^*
\oplus 
\sor\pbk (A/\Ohm\vip)^*)\quad.
$$ 
In other words $\{\theta^\sharp_g\}_{g\in\calG}$ is a smooth family of linear isomorphisms 
$$
\theta^\sharp_g:(A/\Ohm\vip)_{\sor(g)}\overset{\thicksim}\longrightarrow(A/\Ohm\vip)_{\tar(g)}
$$ 
that can be characterized by the following property: for all $g\in\calG$ and $[a]_{A/\Ohm\vip}$,
\bgn{equation}\label{descent}
\qquad
\theta^\sharp_g([a]_{A/\Ohm\vip})=[\ohm_g]_{\Ohm/\Ohm\vip}*[a]_{A/\Ohm\vip}
\qquad,
\end{equation}
where $[\ohm_g]_{\Ohm/\Ohm\vip}$ is the unique element in $\ul{\sor}\inverse([a]_{A/\Ohm\vip})$
with $\ul{\Pr}([\ohm_g]_{\Ohm/\Ohm\vip})=g$. From equation (\ref{descent}) is easy to see that
$\theta^\sharp_\bullet$ enjoys the pseudo-group properties
$$
\theta^\sharp_{gh}=\theta^\sharp_g\cdot\theta^\sharp_h
\qquad\hbox{ and }\qquad
\theta^\sharp_{\eps(q)}=\id_{(A/\Ohm\vip)_q} 
$$ 
for all composable $g,h\in\calG$ and $q\in M$.
In the language of \cite{mackbook}, $\theta^\sharp_\bullet$ is a linear action of $\im\chi$
on $A/\Ohm\vup\rightarrow M$; as a consequence \cite{mackbook}.
There exists a unique smooth vector bundle $X\rightarrow M/\calG$, such that
$\mathrm{q}^!X=A/\Ohm\vup$. Note that $\mathrm{q}^!:A/\Ohm\vup\rightarrow X$ is
fibrewise a diffeomorphism and
\be
\dim X
&=&
\dim M/\calG
+
\dim A/\Ohm\vup -\dim M\\
&=&
\dim M -\dim \calG
+
2\dim A - \dim\Ohm\vip \\
&=&
2\dim A -\dim\Ohm\\
&=&
\dim A/\Ohm\qquad\qquad\qquad\qquad\qquad\qquad\qquad\qquad\qquad.
\ee
Concretely, $X$ is the quotient of $A/\Ohm\vip$ for the regular 
equivalence relation $\thicksim_\theta$ induced by $\theta_\bullet$,
$$
[a_+]_{A/\Ohm\vup}\thicksim [a_-]_{A/\Ohm\vup}\quad\hbox{if{f}}\quad
\theta_g([a_+]_{A/\Ohm\vup})=[a_-]_{A/\Ohm\vup}\quad\hbox{for some (unique)
$g\in\calG$}\quad,
$$
which allows identifying the fibres over the same $\calG$-orbit.
Since the quotient projection $A\rightarrow A/\Ohm$ is $\Ohm\vup$-invariant, it
descends to a surjective submersion $A/\Ohm\vip\rightarrow A/\Ohm$; last map
is stable under $\thicksim_\theta$, thus it descend to a surjective submersion 
$X\rightarrow A/\Ohm$. To conclude the proof, it suffices
to show that last map is injective, which is a straightforward check.
\end{proof}\eject
%
%
%
\subsection{The Lie-Rinehart algebras of morphic and pseudoinvariant sections: algebraic reduction}\label{ss2}\hfill

\vs{0.1}
\spa Before endowing the quotient vector bundle of proposition \ref{topo} with a Lie algebroid structure we need an
algebraic digression; in fact we are going to obtain the quotient Lie algebroid via a reduction procedure on
Lie-Rinehart algebras arising from the associated action \tla-groupoid.

\bigskip

\spa Lie-Rinehart algebras provide a  minimal algebraic model for Lie algebroids.
\bgn{definition}\cite{rh63} Let $\KK$ be a commutative ring and $\calA$ a $\KK$-algebra. A \tsf{Lie-Rinehart
algebra} over $\calA$ is given by
\medskip\\
1. A $\KK$-Lie algebra $\mathfrak{L}$
\medskip\\
2. A left $\calA$-module 
$\diamond:\calA\otimes_{\KK}\mathfrak{L}\to\mathfrak{L}$;
\medskip\\
3. A representation $\sf{r}:\mathfrak{L}\to\sf{Der}\calA$, in the Lie algebra of derivations of $\calA$, which is
compatible with the natural $\calA$-module $\cdot:\calA\otimes_{\KK}\sf{Der}\calA\to\sf{Der}\calA$, in the
sense that
$$
\qquad
\sf{r}(a\diamond l)=a\cdot\sf{r}(l)\qquad,\qquad a\in\calA\quad\hbox{and}\quad
l\in\mathfrak{L}\qquad,
$$
and characterizies the defect for the Lie bracket $\brak{}{}$ of 
$\mathfrak{L}$ to be $\calA$-bilinear via the Leibniz rule
$$
\qquad
\brak{l_+}{a\diamond l_-}-a\diamond\brak{l_-}{l_+}=\sf{r}(l_+)(a)\diamond l_-
\qquad,\qquad l_\pm\in\mathfrak{L}\quad\hbox{and}\quad a\in\calA\qquad.
$$ 
\end{definition}
We shall say that the data above constitute an $\calA$-Lie-Rinehart algebra on a $\KK$-module $\mathfrak{L}$, for
short an $\calA$-LR-algebra and call $\sf r$ the \tsf{Lie-Rinehart anchor},  respectively $\calA$ the \tsf{base}
(\tsf{algebra}) of a Lie-Rinehart algebra. Typically, we shall consider examples where $\calA$ is a ring of functions
and call it \tsf{base ring} accordingly.
\bgn{remark} Let $E\to M$ be a vector bundle.  By replacing $\KK$ with $\RR$, $\mathfrak{L}$ with $\Gamma(E)$  and
$\calA$ with $\cif(M)$, one can see that a Lie Algebroid on $E$ is equivalent to a $\cif(M)$-LR-algebra;
in fact, since the LR-anchor takes values in ${\sf{Der}\,\cif(M)}\simeq\frax(M)$, thus it induces
a vector bundle map $\rho:E\to TM$ by setting $\rho(e_q):=\sf r(e)(q)$, for any section $e\in\Gamma(E)$,
$q\in M$ (here $e_q$ denotes the value of $e$ in $q$).
\end{remark}
Let us discuss in some detail a typical example of a Lie-Rinehart algebra.
\bgn{example} \tsf{The Lie-Rinehart algebra of multiplicative vector fields}\\
A vector field $X$ on a Lie groupoid   $\poidd{\calG}{M}$ is \textsf{multiplicative} if 
$$\qquad
X_{gh}=X_g\bullet X_h\qquad,\qquad(g,h)\in\calG^{(2)}
\qquad,
$$
for the cotangent multiplication $\bullet$ on $\poidd{T\calG}{TM}$.  The multiplicativity condition is equivalent
to asking a vector field $X$ to be a morphism of groupoids from $\poidd{\calG}{M}$ to the tangent prolongation
$\poidd{T\calG}{TM}$, over the base map $X^M$ defined by setting
\bgn{equation}\label{multifield}
\qquad\dd\eps X^M_m=X_{\eps{(m)}}\qquad,\qquad m\in M\qquad;
\end{equation}
it is straightforward to see that $X^M$ is then a smooth vector field over $M$, which is $\sor$- and $\tar$-related
to $X$. Multiplicative vector fields are
precisely those, whose local flows are multiplicative \cite{mx98}, in the sense that the equation $\phi_u(gh)=\phi_u(g)\cdot\phi_(h)$ holds 
for all $(g,h)\in\calG^{(2)}$, provided both sides are defined. As a consequence, the space 
$\frax_\mu(\calG)$ of multiplicative vector fields is easily seen to be a Lie subalgebra of $\frax(\calG)$ and the
base map of a bracket of multiplicative vector fields is
$$
[X,Y]^M=[X^M,Y^M]\qquad,\qquad X,Y\in\frax_\mu(\calG)
\qquad.
$$
Moreover, $\frax_\mu(\calG)$ is naturally endowed with a $\calC^\infty(M)^\calG$-(bi)module structure, for the
subring of $\calG$-\tsf{invariant functions}, i.e. those functions which  are constant along the orbits of
$\calG$: $F\in\calC^\infty(M)^\calG$ if{f} $\sor^*F=\tar^*F$. The module is given by
$$
F\diamond X:=(\sor^*F)\cdot X\qquad,\qquad F\in\calC^\infty(M)^\calG\:,\:X\in\frax_\mu(\calG)
$$
The map $\rho:\frax_\mu(\calG)\rightarrow\sf{Der}\,\calC^\infty(M)^\calG$,  $X\mapsto X^M$ is a Lie algebra
homomorphism; note that $\rho$ is takes values in $\sf{Der}\,\calC^\infty(M)^\calG$ indeed,
\be
X^M_{\sor(g)}(F)
&=&\pair{\dd F_{\sor(g)}}{\dd \sor X_{\eps(\sor(g))}}
=X_g(\sor^*F)
=X_g(\tar^*F)
=\pair{\dd F_{\tar(g)}}{\dd \sor X_{\eps(\tar(g))}}\\
&=&X^M_{\tar(g)}(F)\hs{10.25},
\ee
for all $g\in\calG$ and $F\in\calC^\infty(M)^\calG$. It is straightforward to check that the compatibility between
the $\calC^\infty(M)^\calG$-module and the Lie algebra on $\frax_\mu(\calG)$ is encoded by the usual Leibniz rule;
in other words, $\frax_\mu(\calG)$ is a $\calC^\infty(M)^\calG$-Lie-Rinehart algebra over $\RR$.
\end{example}

\spa In a general \tla-groupoid multiplicative sections play a role analogous to that of multiplicative  vector
fields in the tangent prolongation \tla-groupoid. Multiplicative sections have no natural
characterization in terms of flows; nevertheless, one can show that they form a LR-algebra by using the dual
description, in terms of fibrewise linear functions.
\bgn{definition} Let $(\Ohm,A;\calG,M)$ be an \tla-groupoid. We shall call
$\ohm\in\Gamma(\calG,\Ohm)$ a \textsf{morphic section} if 
\bgn{equation}\label{morphiccondition}
\ohm(gh)=\ohm(g)\,\hat{\cdot}\,\ohm(h)\qquad,\qquad (g,h)\in\calG^{(2)}
\qquad,
\end{equation}
equivalently if $\ohm$ a morphism of Lie groupoids over the (uniquely determined) base section
$$
\qquad
a=\sh\comp\ohm\comp\eps=\th\comp\ohm\comp\eps:M\rightarrow A
\qquad.
$$
\end{definition}
Note that the base section of a morphic section is equivariant for the actions of $\calG$ on $M$
and of $\Ohm$ on $A$ by left translation: for all $g\in\sor^{-1}(m)$, $m\in M$.
$$\quad
a(g*m)=a(\tar(g))=\th(\ohm(g))=\ohm(g)*\sh(\ohm(g))=\ohm(g)*a(\sor(g))=\ohm(g)*a(m)
\quad.
$$
\bgn{remark} For any \tla-groupoid $(\Ohm,\calG;A,M)$ the zero section $0^\Ohm$ is a morphic section over the zero
section $0^A$. Morphic sections are precisely the smooth functors $\calG\to \Ohm$ which are left inverses to the
(bundle) projection functor $\Ohm\to \calG$.	
\end{remark}
The following easy lemma gives an effective characterization of morphic sections.
\bgn{lemma}\label{morpheus}
Let $(\Ohm,\calG;A,M)$ be an \tla-groupoid and $\ohm$ be a section of the top Lie algebroid.
Then the following are equivalent:
\medskip\\
$i)$ $\ohm$ is  morphic;
\medskip\\
$ii)$ The fiberwise linear function 
$\widetilde{\ohm}\in\cif(\Ohm^*\times\Ohm^*\times\Ohm^*)$,
$$\qquad
\widetilde{\ohm}(\xi^1,\xi^2,\xi^3)
:=
\ohm(\xi^1)+\ohm(\xi^2)+\ohm(\xi^3)
\qquad,\qquad
\xi^{1},\xi^{2},\xi^{3}\in\Ohm^*\qquad,
$$
vanishes on $\Gamma(\muh)^o$.
\end{lemma}
\begin{proof}  $\ohm$ is morphic if{f}
$\ohm\times \ohm\times\ohm:\Gamma(\mu)\rightarrow\Gamma(\muh)=\Gamma(\muh)^{oo}$,  where
$\Gamma(\muh)^{o}\subset\Ohm^*\times\Ohm^*\times\Ohm^*$  is the annihilator of $\Gamma(\muh)$ and 
$\Gamma(\muh)^{oo}=(\Gamma(\muh)^o)^o$.
\end{proof}
From now on we shall denote with $\mgamma$ the space of morphic sections of an
\tla-groupoid $(\Ohm,\calG;A,M)$.
It is sometimes convenient to regard a morphic section $\ohm$ with base section
$a$ as a \tsf{morphic pair} $(a,\ohm)$; this point of view
allows computing the base section of a morphic section in the dual picture using
the following tautological lemma, which can be proved in the same way as lemma
\ref{morpheus}.
\bgn{lemma}\label{mopa}
Let $(\Ohm,\calG;A,M)$ be an \tla-groupoid,and consider sections 
$\ohm\in\mgamma$ and $a\in\Gamma(M,A)$. Then, %
\medskip\\
$i)$ $(a,\ohm)$ is a morphic pair;
\medskip\\
$ii)$ The fiberwise linear function 
$E(a,\ohm)\in\cif(A^*\times\Ohm^*)$,
$$\qquad
E(a,\ohm)(\alpha,\xi)
=
a(\alpha)+\ohm(\xi)\qquad,\qquad
\alpha\in A^*\quad,\xi\in\Ohm^*\qquad,
$$
vanishes on $\Gamma(\epsh)^o$.
\end{lemma}
Morphic sections $\mgamma\subset\Gamma(\calG,\Ohm)$ form an  $\RR$-linear subspace\fn{Note that the top
multiplication $\muh_{(g,h)}:\Ohm^{(2)}_{(g,h)}\to \Ohm_{gh}$ is fibrewise linear for the linear structure on 
$\Ohm^{(2)}_{(g,h)}$ by $\Ohm_g\times \Ohm_h$. It is \emph{not} linear in the separate components.}
naturally endowed with a $\cif(M)^\calG$-(bi)module:
$$
\qquad
F\diamond\ohm=(\sor^*F)\cdot\ohm\qquad,\qquad F\in\calC^\infty(M)^\calG\:,
\:\ohm\in\mgamma
\qquad;
$$
if $a$ is the base section of $\ohm$, $F\cdot a$ is the base section of 
$F\diamond\ohm$. If for all $g\in\calG$, $(\Ohm,\calG;A,M)$ admits a morphic section $(\ohm,a)$ with $a(\sor(g))\neq
0\neq a(\tar(g))$, it is easy to see, using the defining condition (\ref{morphiccondition}), that 
$\cif(M)^\calG\subset\cif(\calG)$
is the largest subring for which  the $\cif(\calG)$-module on 
$\Gamma(\calG,\Ohm)$ restricts to $\mgamma$. 
For any morphic pair $(\ohm,a)$ 
$$\qquad
\dd\sor\comp\rhoh\comp\ohm=\rho\comp a\comp \sor
\qquad\hbox{ and }\qquad 
\dd\tar\comp\rhoh\comp\ohm=\rho\comp a\comp \tar\qquad,
$$
thus the mapping $\Tilde{\rho}:\mgamma\rightarrow\frax(M)$, $\ohm\longmapsto \rho\comp a$
takes values in $\sf{Der}\,\cif(M)^{\calG}$: for all $f\in\cif(M)$,
$
\rho(a)_{\tar(g)}(f)
=
[\dd\tar\comp(\rhoh(\ohm))]_{\tar(g)}(f)
=
(\rhoh(\ohm))_g(\tar^*f)
$,
and, analogously, $\rho(a)_{\sor(g)}(f)=(\rhoh(\ohm))_g(\sor^*f)$; thus, if $f$ is a $\calG$-invariant function,
$\rho(a)_{\tar(g)}(f)=\rho(a)_{\sor(g)}(f)$ and $\rho(a)(f)$ is also $\calG$-invariant.
\bgn{theorem}\label{lrm}
The space of morphic sections $\mgamma$ of an \tla-groupoid $(\Ohm,\calG;$ $A,M)$ is a
$\cif(M)^\calG$-Lie-Rinehart algebra for the Lie bracket induced from
$\Gamma(\calG,\Ohm)$ and the anchor 
$\Tilde{\rho}:\mgamma\rightarrow\sf{Der}\,\cif(M)^{\calG}$
defined above. In particular, for any  morphic pairs $(\ohm_\pm,a_\pm)$, 
$([\ohm_+,\ohm_-], [a_+,a_-])$ is also a morphic pair.
\end{theorem}
\begin{proof} If the second statement holds $\mgamma$ is an $\RR$-Lie subalgebra and a 
$\cif(M)^\calG$-submodule of $\Gamma(\calG,\Ohm)$; moreover $\Tilde{\rho}$ is a morphism of
$\RR$-Lie algebras. The top multiplication
$\muh:\Ohm\fib{\sh}{\th}\Ohm\rightarrow\Ohm$ is a
morphism of Lie algebroids, hence $\Gamma(\muh)\subset\threetimes{\Ohm}$
is a Lie subalgebroid and $\Gamma(\muh)^o\subset\threetimes{\Ohm^*}$ a
coisotropic submanifold. Since  $\Tilde{\ohm}_\pm\in\calI_{\Gamma(\muh)^o}$,
$\{\Tilde{\ohm}_+,\Tilde{\ohm}_-\}\in\calI_{\Gamma(\muh)^o}$; moreover, for all $\xi^{1}$, $\xi^{2}$, 
$\xi^{3}\in\Ohm^*$,
\be
\{\widetilde{\ohm}_+,\widetilde{\ohm}_-\}(\xi_1,\xi_2,\xi_3)
&=&
\sum_{i=1}^{3}\{{\ohm}_+,{\ohm}_-\}(\xi_i)\quad=\quad
\sum_{i=1}^{3}\langle[{\ohm}_+,{\ohm}_-],\xi_i\rangle\\
&=&
\widetilde{[{\ohm}_+,{\ohm}_-]}(\xi_1,\xi_2,\xi_3)	
\qquad,
\ee
then $[{\ohm}_+,{\ohm}_-]$ is morphic. Let us check the second statement. Since $(a_\pm,\ohm_\pm)$ are morphic
pairs $\{E(a_+,\ohm_-),E(a_+,\ohm_-)\}$ vanishes on $\Gamma(\epsh)^o$; being 
$\Gamma(\epsh)$ coisotropic and,
for all $(\alpha,\xi)\in A^*\times\Ohm^*$,
\be
\{E(a_+,\ohm_+),E(a_-,\ohm_-)\}(\alpha,\xi)
&=&
\{a_+,a_-\}(\alpha)
+\{{\ohm}_+,{\ohm}_-\}(\xi)\\
&=&
\langle[a_+,a_-],\alpha\rangle+
\langle[{\ohm}_+,{\ohm}_-],\xi\rangle\\
&=&
E([a_+,a_-],[{\ohm}_+,{\ohm}_-])(\alpha,\xi)\qquad\:,
\ee
$E([a_+,a_-],[{\ohm}_+,{\ohm}_-])\in\calI_{\Gamma(\epsh)^o}$, 
i.e. $[a_+,a_-]$ is the base section of $[{\ohm}_+,{\ohm}_-]$. 
The Leibniz rule follows: 
\be
[\ohm_+,f\diamond\ohm_-]_g
&=&
[\ohm_+,\sor^*f\cdot \ohm_-]_g
=
(\sor^*f\cdot [\ohm_+,\ohm_-])_g
+
[(\rhoh(\ohm_+))]_g(\sor^*f)\cdot\ohm_-\\
&=&(f\diamond[\ohm_+,f\diamond\ohm_-])_g
+  \sor^*\Tilde{\rho}(\ohm_+)(g)\ohm_-(g)\\
&=&(f\diamond[\ohm_+,\ohm_-]+ \Tilde{\rho}(\ohm_+)\diamond\ohm_-)_g
\ee
holds, for all $f\in\cif(M)^\calG$ and $g\in\calG$.
\end{proof}
Apart from multiplicative vector fields, other examples of morphic sections have already appeared in literature.
\bgn{example} Morphic sections of the cotangent prolongation \tla-groupoid of a Lie groupoid are the 
multiplicative 1-forms considered in \cite{mx98}.
\end{example}

Morphic sections are biinvariant in a sense we are about to specify. A section $\ohm\in\Gamma(\calG,\Ohm)$ is
\tsf{right pseudoinvariant}, respectively 
\tsf{left pseudoinvariant}, if, for all $(g,h)\in\calG^{(2)}$
\bgn{equation}\label{rl}
\:\ohm(gh)=\ohm(g)\,\hat{\cdot}\,R^\ohm(h)
\quad,\quad\hbox{respectively}\quad,\quad
\ohm(gh)=L^\ohm(g)\,\hat{\cdot}\,\ohm(h)\qquad,
\end{equation}
for some sections $R^\ohm,L^\ohm\in\Gamma(\calG,\Ohm)$. The formulas (\ref{rl}) determine $R^\ohm$ and $L^\ohm$
uniquely:
$$
R^\ohm(g)=\ohm(\eps(\tar(g)))^{-1}\hdot\ohm(g)\quad,\quad\hbox{respectively}\quad,
\quad
L^\ohm(g)=\ohm(g)\hdot\ohm(\eps(\sor(g)))^{-1}
$$
and it is easy to show that $R^\ohm$ and  $L^\ohm$ are both morphic: for example
$$
\:
R^\ohm(gh)=\ohm(\eps(\tar(h)))^{-1}\hdot\ohm(gh)
=
\ohm(\eps(\tar(h)))^{-1}
\hdot
\ohm(g)\,\hat{\cdot}\,R^\ohm(h)
=R^\ohm(g)\hdot R^\ohm(h)
\quad.
$$
The base section $r^\ohm$, respectively $l^\ohm$, of the
morphic section associated with a right, respectively left, pseudoinvariant section is
$$
\quad
r^\ohm(m)=\sh(\ohm(\eps(m)))\quad,
\quad\hbox{respectively}\quad,
\qquad l^\ohm(m)=\th(\ohm(\eps(m)))\quad,
$$
$m\in M$.
If $\ohm$ is a pseudoinvariant section with
associated morphic section $X$, we shall say that $(\ohm,X)$ is an 
\tsf{invariant pair}; we shall also say that $X$ is the \tsf{morphic component} of
$\ohm$.
\bgn{example} Right, respectively left, invariant vector fields on a Lie groupoid are pseudoinvariant sections of
the tangent prolongation \tla-groupoid  with zero morphic component. Morphic sections of an \tla-groupoid are
right and left pseudoinvariant sections, coinciding with their invariant components.
\end{example}
In the same way as morphic sections, right and left pseudoinvariant sections admit a dual geometric characterization.
\bgn{lemma}\label{righorlef}
Let $(\Ohm,\calG;A,M)$ be an \tla-groupoid and $\ohm,R\in\Gamma(\calG,\Ohm)$.
Then the following are equivalent:
\medskip\\
$(i)$ $(\ohm,R)$ is a right invariant pair;
\medskip\\
$(ii)$ The fiberwise linear function 
$\rinv{\ohm}\in\cif_{\sf{lin}}(\Ohm^*\times\Ohm^*\times\Ohm^*)$
$$
\qquad
\rinv{\ohm}(\xi^1,\xi^2,\xi^3):=\ohm(\xi^1)+R^\ohm(\xi^2)+\ohm(\xi^3)
\qquad
$$
vanishes on $\Gamma(\muh)^o$.
\medskip\\
The analogous statements hold for left invariant functions.
\end{lemma}
\begin{proof} Adapt the proof of theorem \ref{lrm}.
\end{proof}
We shall denote with $\gammar{\calG,\Ohm}$ and $\gammal{\calG,\Ohm}$ the $\RR$-linear spaces of right and left
invariant sections of an \tla-groupoid $(\Ohm,\calG;A,M)$; the natural $\RR$-linear structure of $\Gamma(\calG,\Ohm)$
restricts to the invariant sections, due to linearity of the top multiplication $\muh$. A $\cif(M)^\calG$-module can
be defined both on  $\gammar{\calG,\Ohm}$, respectively $\gammal{\calG,\Ohm}$, precisely in the same way  as for
morphic sections and the anchor
$$\qquad
\rinv{\rho}:\gammar{\calG,\Ohm}\longrightarrow\sf{Der}\,\cif(M)^\calG\quad,\quad
\ohm\longmapsto \rho(r^\ohm)
\qquad,$$
respectively
$$
\qquad
\linv{\rho}:\gammal{\calG,\Ohm}\longrightarrow\sf{Der}\,\cif(M)^\calG\quad,\quad
\ohm\longmapsto \rho(l^\ohm)\qquad,
$$
is compatible with the multiplication by $\calG$-invariant functions; the images of $\rinv{\rho}$ and $\linv{\rho}$
actually lay in $\sf{Der}\,\cif(M)^\calG$, since, for any right, respectively left, invariant section $\ohm$,
$r^\ohm$, respectively $l^\ohm$  is the base section of a morphic section.\\ 
Note that $\gammar{\calG,\Ohm}$ and $\gammal{\calG,\Ohm}$ are linearly isomorphic over $\cif(M)^\calG$, the
isomorphism being given (either way)  by composition with the top inversion map of $(\Ohm,\calG;A,M)$. More
precisely,  under the mapping
$$\qquad\ol{\,\cdot\,}:
\gammar{\calG,\Ohm}\longrightarrow\gammal{\calG,\Ohm}\qquad,
\qquad\ohm\longmapsto\ol{\ohm}:=\iotah\comp\ohm\comp\iota
\qquad,
$$
we have
\be
\ol{\ohm}(gh)&=&\ohm(h^{-1}g^{-1})^{-1}
=
(\ohm(h^{-1})\hdot R^\ohm(g^{-1}))^{-1}
=
(R^\ohm(g^{-1}))^{-1}\hdot(\ohm(h^{-1}))^{-1}\\
&=&
(R^\ohm(g^{-1}))^{-1}\hdot\ol{\ohm}(h)
=
R^\ohm(g)\hdot\ol{\ohm}(h)\qquad,
\ee
that is, $\ol{\ohm}$ is left pseudoinvariant indeed and $L^{\ol{\ohm}}=\ol{R^\ohm}=R^\ohm$, since morphic sections are
$\ol{\,\cdot\,}$-stable. Then $\ol{\,\cdot\,}$  is compatible with the anchors, in the sense that
$\linv{\rho}\comp\ol{\cdot}=\rinv{\rho}$:
$$\qquad
\rinv{\rho}(\ohm)=\widetilde{\rho}(R^\ohm)=\widetilde{\rho}(L^{\ol{\ohm}})
=
\linv{\rho}(\ol{\ohm})\qquad,\qquad \ohm\in\rgamma\qquad.
$$
Thinking in terms of linear functions on $\Ohm^*$, 
$$\qquad
\ol{\ohm}=(\iotah^{\sf{t}})^*\ohm\qquad\qquad\hbox{ and }\qquad\qquad
L^{\ol{\ohm}}=(\iotah^{\sf{t}})^*R^\ohm=R^\ohm\qquad
$$
for the Poisson automorphism $\iotah^{\sf{t}}:\Ohm^*\rightarrow\Ohm^*$. The
isomorphism inverse to $\ol{\,\cdot\,}$ enjoys the same properties as 
$\ol{\,\cdot\,}$ and we shall denote it also by the same symbol.
\bgn{proposition}\label{lrinv} For any \tla-groupoid $(\Ohm,\calG;A,M)$
the space of right, respectively left, invariant sections $\gammar{\calG,\Ohm}$, respectively
$\gammal{\calG,\Ohm}$, is a $\cif(M)^\calG$-Lie-Rinehart algebra over $\RR$
for the Lie bracket induced from
$\Gamma(\calG,\Ohm)$ and the anchor 
$\rinv{\rho}:\gammar{\calG,\Ohm}\rightarrow\sf{Der}\,\cif(M)^{\calG}$, respectively
$\linv{\rho}:\gammal{\calG,\Ohm}\rightarrow\sf{Der}\,\cif(M)^{\calG}$
defined above. Moreover, 
$\ol{\,\cdot\,}:\gammal{\calG,\Ohm}\rightarrow\gammar{\calG,\Ohm}$ is
an isomorphism of  $\cif(M)^\calG$-Lie-Rinehart algebras.
\end{proposition}
\bgn{corollary} For any \tla-groupoid $(\Ohm,\calG;A,M)$,
$$
\:
[\gammar{\calG,\Ohm},\mgamma]\subset\gammar{\calG,\Ohm}
\quad\hbox{and}\quad
[\gammal{\calG,\Ohm},\mgamma]\subset\gammal{\calG,\Ohm}\:.
$$
\end{corollary}
\begin{proof} 
The proof goes on the same lines of that of theorem \ref{lrm}. The following 
facts hold\medskip\\
($i$) For any  $\ohm_\pm\in\gammar{\calG,\Ohm}$ with associated morphic
	sections $R^\ohm_\pm\in\mgamma$, $[\ohm_+,\ohm_-]$ 
	is right invariant  with associated morphic section
	section $[R^{\ohm_+},R^{\ohm_-}]$;\medskip\\
($ii$) The anchor map $\rinv{\rho}$ is a morphism of Lie algebras;\medskip\\
As in theorem \ref{lrm}, ($i$) can be proved using lemma \ref{righorlef},
($ii$) follows and the check of the Leibniz rule amounts to a straightforward computation.
The similar statements for $\gammal{\calG,\Ohm}$ can be 
shown using the top groupoid
inversion map, at the same time proving the isomorphism 
$\gammal{\calG,\Ohm}\simeq\gammar{\calG,\Ohm}$. 
For any $\ohm_\pm\in\lgamma$ and $\xi\in\Ohm^*\quad$
\be
\langle[\ohm_+,\ohm_-],\xi\rangle
&=&
\{\ohm_+,\ohm_-\}(\xi)
=
\{(\iotah^{\sf{t}})^*\ol{\ohm}_+,(\iotah^{\sf{t}})^*\ol{\ohm}_-\}(\xi)
=
(\iotah^{\sf{t}})^*\{\ol{\ohm}_+,\ol{\ohm}_-\}(\xi)\\
&=&
\langle\ol{[\ol{\ohm}_+,\ol{\ohm}_-]},\xi\rangle
\ee
and
$$\qquad
L^{[\ohm_+,\ohm_-]}
=
L^{\ol{[\ol{\ohm}_+,\ol{\ohm}_-]}}
=
R^{[\ol{\ohm}_+,\ol{\ohm}_-]}
=
[R^{\ol{\ohm}_+},R^{\ol{\ohm}_-}]
=
[L^{{\ohm}_+},L^{{\ohm}_-}]
\qquad
$$
as it follows from the remarks above.
\end{proof}

\bigskip

\spa Morphic sections of an \tla-groupoid $(\Ohm,\calG;A,M)$ are smooth functors $\calG\rightarrow\Ohm$, then the
obvious notion of equivalence on $\mgamma$ is described by natural transformations. However, we shall say that  two morphic
sections $\ohm_\pm\in\mgamma$ are \textsf{equivalent} if there exists a smooth natural transformation $\eta$ from
$\ohm_-$ to $\ohm_+$, which is moreover compatible with the projection functor $\Pr:\Ohm\to\calG$, i.e. a smooth map $\eta:
M\rightarrow\Ohm$, $m\rightarrow\eta_m$, such that $\Pr\comp\eta=\eps$ and
\bgn{equation}\label{natural}
\qquad
\ohm_+(g)\hdot\eta_{\sor(g)}
=
\eta_{\tar(g)}\hdot\ohm_-(g)
\qquad,\qquad g\in\calG\qquad.
\end{equation}
An equivalence $\eta$ from $\ohm_-$ to $\ohm_+$ shall be denoted also by  $\eta:\ohm_-\Rightarrow\ohm_+$. Note
that, for any equivalence of morphic sections as above, the base sections $a_\pm\in\Gamma(M,A)$ are related by the
formula 
$$
\qquad
a_+(m)=\eta_{m}*a_-(m)
\qquad,\qquad m\in M\qquad,
$$
for the action of $\Ohm$ on $A$ by left translation.\\ 
The extra compatibility condition with the projection functor makes an equivalence of morphic sections a section
of the restriction of $\Ohm$ to $M$ and, technically, it is required to make the groupoid
$\poidd{\sf{Nat}(\calG,\Ohm)}{\mgamma}$ on the set of equivalences of morphic sections a groupoid object in the
category of $\cif(M)^\calG$-modules. For any pair of equivalences $\eta^{1,2}:\ohm^{1,2}_-\Rightarrow\ohm^{1,2}_+$
of morphic sections and $F\in\cif(M)^\calG$, setting
$$
\eta_m:=\eta^1_m +_{\eps(m)}F(m)\cdot\eta^2_m
$$
yields an equivalence 
$
\eta=(\eta^1 +F\diamond\eta^2):
\ohm^1_- + F\diamond\ohm^2_-
\Rightarrow
\ohm^1_++ F\diamond\ohm^2_+
$,
 where the sum is taken in the fibre of $\Ohm$ over $\eps(m)$.  The algebraic properties of a
$\cif(M)^\calG$-module can be easily proved using the commutation relation (\ref{natural}).\\ 
There are interesting cases in which the extra compatibility condition is fulfilled by all natural
transformations:
\bgn{remark}\label{equo} Let $\eta:\ohm_-\Rightarrow\ohm_+$ be a natural transformation of morphic sections. 
One can compute, for all $g\in\calG$,
$$
\tar(\Pr\,\eta_{\sor(g)})=\pr\,\th(\eta_{\sor(g)})=\pr\,\sh(\ohm_+(g))=\sor(g)
$$
and similarly $\sor(\Pr\,\eta_{\tar(g)})=\tar(g)$; thus for $g=\eps(q)$, 
$\sor(\Pr\,\eta_{q})=q=\tar(\Pr\,\eta_{q})$. That is, whenever all the isotropy groups of $\calG$ are trivial, the
compatibility with the projection functor is automatically satisfied.  
\end{remark} 
\bgn{theorem}
For any \tla-groupoid $(\Ohm,\calG;A,M)$, the  $\cif(M)^\calG$-Lie-Rinehart algebra on  $\mgamma$ descends to the
quotient  $\ul{\mgamma}:=\mgamma\mod$.
\end{theorem}
\bgn{proof}
The $\cif(M)^\calG$-module descends to the quotient since $\sf{Nat}(\calG,\Ohm)$ is a $\cif(M)^\calG$-module. The
anchors of  equivalent morphic sections coincide as derivations of $\cif(M)^\calG$:  for any equivalence
$\eta:\ohm_-\Rightarrow\ohm_+$ of morphic sections and $F\in\cif(M)$
$$
\widetilde{\rho}_m(\ohm_+)(f)
=
\rho_m(\sh_{\eps(m)}\ohm_+)(f)
=
\langle\dd\,f,(\rho\comp\th)_{\eps(m)}\eta_m\rangle
=
\langle\dd\,\tar^*f,\rhoh_{\eps(m)}\eta_m\rangle
$$
and, similarly, 
$
\widetilde{\rho}_m(\ohm_-)(f)
=\langle\dd\,\sor^*f,\rhoh_{\eps(m)}\eta_m\rangle
$; 
then, $\widetilde{\rho}_m(\ohm_+)(f)=\widetilde{\rho}_m(\ohm_-)(f)$, whenever $f$ is $\calG$-invariant. Then the
anchor descends to a morphism of  $\cif(M)^\calG$-modules $\ul{\mgamma}\to{\sf{Der}}\,{\cif(M)}^\calG$; if the Lie bracket descends, the anchor is then
automatically a morphism of Lie algebras and the Leibniz rule holds. Let 
$\eta^{1,2}:\ohm^{1,2}_-\Rightarrow\ohm^{1,2}_+$ be a pair of equivalences of morphic sections. Then, setting 
$$\qquad
(\ohm^{1,2}_+\star\eta^{1,2})(g)
:=
\ohm^{1,2}_+(g)\hdot\eta^{1,2}_{\eps(\sor(g))}
=
\eta^{1,2}_{\eps(\tar(g))}\hdot\ohm^{1,2}_-(g)\qquad,\qquad
g\in\calG\qquad,$$ defines
biinvariant sections with right and left morphic components
$$\qquad
L^{\ohm^{1,2}_+\star\eta^{1,2}}=\ohm^{1,2}_+
\qquad\hbox{ and }\qquad R^{\ohm^{1,2}_+\star\eta^{1,2}}=\ohm^{1,2}_-
\qquad.
$$
According to the proof of theorem \ref{lrinv},  $[\ohm^1_+\star\eta^1,\ohm^2_+\star\eta^2]$ is a
biinvariant section  with morphic components
$$\qquad
L^{[\ohm^1_+\star\eta^1,\ohm^2_+\star\eta^2]}
=[\ohm^1_+,\ohm^2_+]
\qquad\hbox{ and }\qquad R^{[\ohm^1_+\star\eta^1,\ohm^2_+\star\eta^2]}
=[\ohm^1_-,\ohm^2_-]
\qquad,
$$
as a consequence
$
[\ohm^1_+\star\eta^1,\ohm^2_+\star\eta^2]_g
=
[\ohm^1_+,\ohm^2_+]_g
\hdot
[\ohm^1_+\star\eta^1,\ohm^2_+\star\eta^2]_{\eps(\sor(g))}
$,
by left invariance, on the other hand,
$
[\ohm^1_+\star\eta^1,\ohm^2_+\star\eta^2]_g
=[\ohm^1_+\star\eta^1,\ohm^2_+\star\eta^2]_{\eps(\tar(g))}
\hdot
[\ohm^1_-,\ohm^2_-]_g
$,
by right  invariance. Note that $\eta^{12}:=[\ohm^1_+\star\eta^1,\ohm^2_+\star\eta^2]\comp\eps$  is, by construction
compatible with the projection functor, therefore we have found an  equivalence of morphic sections
$\eta^{12}:[\ohm^1_-,\ohm^2_-]\Rightarrow[\ohm^1_+,\ohm^2_+]$.
\end{proof}

\bigskip

\spa Next we shall see that, for any \tla-groupoid $(\Ohm,\calG;A,M)$, $\mgamma\mod$ is isomorphic to the
Lie-Rinehart algebra of $\Gamma(M/\calG, A/\Ohm)$ for the groupoid actions by left translation, whenever the
quotient Lie algebroid exists. Note that the construction of $\mgamma\mod$ is independent of any assumption on
$(\Ohm,\calG;A,M)$, beside the \tla-groupoid structure; hence, even when the quotient Lie algebroid on the
topological bundle $A/\Ohm\rightarrow A/\calG$ does not exist, $\mgamma\mod$ is to be regarded as a
model or a desingularization thereof.

\spa From now on we shall assume that $(\Ohm,\calG;A,M)$ is an isotropy free and proper \tla-groupoid, so that the
quotient vector bundle $A/\Ohm\to M/\calG$ exists. A section  $a\in\Gamma(M,A)$ is
\tsf{projectable} if there exists a section $\ul{a}\in\Gamma(M/\calG,A/\Ohm)$, such that
$$
\bcat
\xy
*+{}="0",    <-0.7cm,0.7cm>
*+{A}="1", <0.7cm,0.7cm>
*+{A/\Ohm}="2", <-0.7cm,-0.7cm>
*+{M}="3", <0.7cm,-0.7cm>
*+{M/\calG}="4",
\ar                @{->} "1";"2"^{} 
\ar  		   @{->} "3";"1"^{a}
\ar                @{->} "4";"2"_{\ul{a}}
\ar                @{->} "3";"4"_{}
\endxy
\ecat
$$
commutes for the quotient projections on the horizontal edges. Clearly, a section is projectable if{f} it is
equivariant for the groupoid actions by left translation, that is, for all $m\in M$ and $g\in\sor^{-1}(m)$, there
exists an $\ohm_g^a\in\Ohm_g$, such that
\bgn{equation}\label{project}
\qquad a(g*m)=\ohm^a_g*a(m)\qquad.
\end{equation}
Since $\poidd{\Ohm}{A}$ is free, $\ohm_g^a$ is uniquely determined by $g$;  since it is proper, the factorization
$
\ohm^a=\chih\inverse\comp(a\times a)\comp \chi
$
defines a smooth section $\ohm^a\in\Gamma(\calG,\Ohm)$. It turns out that $\ohm^a\in\mgamma$; recall that the base
section of a morphic section is always equivariant. Due to fibrewise linearity of the top groupoid action,
projectable sections form an $\RR$-linear space $\pgamma$; it is easy to see that $\pgamma$ is also naturally a
$\cif(M)^\calG$-module.
\bgn{proposition}\label{moods}
For any free and proper \tla-groupoid $(\Ohm,\calG;A,M)$, 
\medskip\\
$i$) The projection $\mgamma\rightarrow\pgamma$ is an isomorphism of $\cif(M)^\calG$-modules;
\medskip\\
$ii$) $\pgamma\subset\Gamma(M,A)$ is a Lie-Rinehart subalgebra.
\end{proposition}  
Note that $\pgamma$ is a $\cif(M)^\calG$-LR algebra, thus a LR-subalgebra 
with change of base.
\bgn{proof} %
($i$) We only have to prove that, for any $a\in\pgamma$, the section $\ohm^a$ defined by (\ref{project}) is morphic:
for any $(g,h)\in\calG^{(2)}$, applying (\ref{project}) and associativity of the  groupoid actions,
\be
\ohm^a(gh)*a(m)&=&a(gh*m)=a(g*h*m)=\ohm^a(g)*\ohm^a(h)*a(m)\\
              &=&\ohm^a(g)\hdot\ohm^a(h)*a(m)\hs{5};
\ee
thus $\ohm^a(gh)*a(m)=\ohm^a(g)\hdot\ohm^a(h)$, since  $\poidd{\Ohm}{A}$ is isotropy free. ($ii$) Straightforward
consequence of ($i$) and theorem \ref{lrm}.
\end{proof}

\spa The natural notion of equivalence for projectable sections is  $a_+\thicksim a_-$, $a_\pm\in\pgamma$, if
$\ul{a}_+=\ul{a}_-$; this is equivalent to requiring that, for all $m\in M$, there exists an $\eta_m\in\Ohm$ such
that
\bgn{equation}\label{natta}
\qquad
a_+(m)=\eta_m*a_-(m)
\qquad,
\end{equation}
by equivariance. Once again, since $\poidd{\Ohm}{A}$ is free and proper, equation (\ref{natta}) determines a smooth
map $\eta:M\rightarrow \Ohm$, actually a natural transformation $\eta:\ohm^{a}_-\Rightarrow\ohm^{a}_+$. To see this,
consider that,  for all $m\in M$ and $g\in\sor^{-1}(m)$,
$$
\quad
a_+(g*m)=\ohm^{a_+}(g)*a_+(m)=\ohm^{a_+}(g)*\eta_{\sor(g)}*a_-(m)=(\ohm^{a_+}(g)\hdot\eta_{\sor(g)})*a_+(m)
\quad
$$
and
$$
\quad
a_+(g*m)=\eta_{\tar(g)}*a_-(g*m)=\eta_{\tar(g)}*\ohm^{a_-}(g)*a_-(m)=(\eta_{\tar(g)}\hdot\ohm^{a_-}(g))*a_-(m)
\quad.
$$
It follows then from remark \ref{equo} that $\eta$ is an equivalence of morphic sections, since $\poidd{\calG}{M}$ is
free.

\medskip

We are ready to conclude the proof of theorem \ref{generalred}.
%
%
\bgn{proof}[Proof of theorem \ref{generalred}] %
 Note that the base sections of equivalent morphic sections are equivalent in the sense of  (\ref{natta}). Then the
isomorphism $\hat{\Gamma}(\calG\acts N,\Ohm\acts B)\to\Gamma^\downarrow (N,B)$  descends to an isomorphism of
$\cif(N)^\calG$-modules 
$$
\qquad
\ul{\hat{\Gamma}(\calG\acts N,\Ohm\acts B)}\overset{\thicksim}\longrightarrow\Gamma^\downarrow(N,B)\mod\equiv\Gamma(N/\calG,B/\Ohm)\qquad,
$$ 
endowing the quotient vector bundle $B/\Ohm\to N/\calG$ with a Lie algebroid structure. For all 
$\beta\in \Gamma(N/\calG,B/\Ohm)$ and $b\in\Gamma^\downarrow(N,B)$ with $\qh\comp b=\beta\comp\mathrm{q}$, the
equality 
$$
\rho_{B/\Ohm}\comp\beta=\rho_B\comp b
$$
holds in the space of derivations on $\cif(N)^\calG$, and similarly, for all 
$\beta'\!\in \Gamma(N/\calG,B/\Ohm)$ and $b'\in\Gamma^\downarrow(N,B)$ with $\qh\comp b'=\beta'\comp\mathrm{q}$,
by definition,
$$
\qquad
\brak{\beta}{\beta'}_{B/\Ohm}([n]_{N/\calG})=[\brak{b}{b'}_B(n)]_{B/\Ohm}\qquad,\qquad n\in N\qquad,
$$
that is, the Lie bracket $\brak{b}{b'}$ is $(\qh,\mathrm{q})$-related to the Lie bracket $\brak{\beta}{\beta'}$; 
the quotient projection is then by construction a fibrewise surjective, base submersive and surjective morphism of 
Lie  algebroids hence a strong \tla-fibration. Uniqueness is clear.
\end{proof}
\bgn{example}\label{tangentq} %
Consider the tangent prolongation \tla-groupoid,  and the actions of $T\calG$ on $TM$, respectively $\calG$ on $M$ by
left translation. It easy to see that the reduction of the (tangent lift of) the action of $\calG$ on $M$ by left
translation yields  $TM/T\calG\simeq T(M/\calG)$, as a vector bundle. Since the quotient projection is submersive
for any $\ul{\delta m}_{\ul{m}}\in(TM/T\calG)_{\ul{m}}$, one can pick a representative $\delta m\in T_mM$, for all
$m\in\ul{m}$; define  $TM/T\calG\rightarrow T(M/\calG)$,  $\ul{\delta m}_{\ul{m}}\mapsto [\delta m_m]$. Such a map
does not depend on the choice of representatives since, for all $\delta m_\pm\in T_m M$, with $\delta m_+\thicksim
\delta m_-$, $\delta m_+=\eta_m*\delta m_-$, for some $\eta\in T_{\eps(m)}\calG$; for any class $C^1$ path
$\gamma^\pm$, such that $\dot{\gamma}^\pm(0)=\eta_m$, set $\gamma^+:=\tar\comp\gamma^\pm$ and
$\gamma^-:=\sor\comp\gamma^\pm$. Thus $\dot{\gamma}^{+,-}(0)=\delta m_{+,-}$,  $[\gamma^{+,-}]$ is a $C^1$ path with
tangent vector $[\delta m_{+,-}]$ at $0$ and by definition $[\gamma^+]=[\gamma^-]$.  The mapping 
$TM/T\calG\rightarrow
T(M/\calG)$ is fibrewise linear and injective, since, if  $[\delta m]=0$, then $\delta m$ is tangent to the orbits of
$\calG$, i.e. $\delta m=\dd\,\tar\,\delta g=\delta g*0_m$ for some $T_{\eps(m)}^\sor\calG$, and $\ul{\delta m}=0$.
Then $TM/T\calG\simeq T(M/\calG)$, since both bundles have the same rank. Moreover, consider that the quotient
projection $M\to M/\calG$ is invariant for the action of $\calG$ on $M$, then so is the tangent map $TM\to
T(M/\calG)$ under the action of $T\calG$; the induced map  $TM/T\calG\to T(M/\calG)$ is precisely the map defined
above, which is then a smooth isomorphism of vector bundles. Up to this identification, the fiberwise linear map on
the quotient associated with the Lie-Rinehart anchor is the identity and the induced Lie bracket, namely that of
multiplicative vector fields, is the canonical bracket on $\frax(M/\calG)$.
\end{example}
\bgn{example}
From last example and the remarks in example \ref{fpla}, one can see that, for
any free and proper action of a Lie group $G$ on a manifold $M$, $T(M/G)$ is to
be canonically identified with $TM/TG$ for the tangent lifted action.
\end{example}
\vs{0.5}
\subsection{Reduction of the moment morphism: quotient Lie algebroids and Poisson structures}\hfill

\vs{0.1}
\spa The kernel of the moment morphism (\ref{mommmo}) associated with a morphic action is well behaved under
reduction under mild regularity assumptions; in particular, in the case of the cotangent lift of a Poisson 
groupoid action, these assumptions are met whenever the action is free and proper. 
\bgn{theorem}\label{marwascat} Let  $(\Ohm,\calG;A,M)$ be an \tla-groupoid acting morphically on morphism of
Lie algebroids $\jh:B\rightarrow A$ over $j: N\rightarrow M$ of maximal rank. Then
\medskip\\ 
$i)$ The kernel of the moment morphism is an \tla-groupoid of the form
$$
\qquad
\bcat
\xy
*+{}="0",    <-0.7cm,0.7cm>
*+{\calG\acts K}="1", <0.7cm,0.7cm>
*+{ K}="2", <-0.7cm,-0.7cm>
*+{\calG}="3", <0.7cm,-0.7cm>
*+{M}="4",
\ar  @ <-0.07cm>   @{->} "1";"2"^{} 
\ar  @ <0.07cm>    @{->} "1";"2"^{}  
\ar  		   @{->} "1";"3"^{}
\ar                @{->} "2";"4"_{}
\ar  @ <0.07cm>    @{->} "3";"4"^{}
\ar  @ <-0.07cm>   @{->} "3";"4"_{}
\endxy
\ecat
\qquad,
$$
where $K=\ker\jh$;
\medskip\\ 
$ii)$ If $\calG$ acts freely and properly on $N$, the action on $K$ is also free and proper and there exists a unique
Lie algebroid of rank 
\bgn{equation}\label{rankquo}
\rank K/\calG=\rank K=\rank B-\rank A
\end{equation}
on $K/\calG\to N/\calG$ making the quotient projection $K\to K/\calG$  an  \tla-fibration.
%
%
\end{theorem}
\bgn{proof} %
($i$) was already remarked in Subsection \ref{mola}. ($ii$) For all $(g_\pm,k)\in\calG\acts K$, by fibrewise
linearity of the top action, we have
$$
g_+*k=g_-*k\quad\Leftrightarrow\quad g\inverse_+g_-*\pr (k)=\pr (k)\quad\Rightarrow g_+=g_-
$$ 
since the action of $\calG$ on $N$ is free. For all sequences $\{(g_n,k_n)\}\subset\calG\acts K$, such that $k_n$
converges to some $k_\infty$ and $g_n*k_n$ converges to some $\wt{k}_\infty$, we have
$$
\qquad
g_n*\pr(k_n)=\pr(g_n*k_n)\longrightarrow \pr(\wt{k}_\infty)\qquad\hbox{ and }\qquad \pr(k_n)\longrightarrow
\pr(k_\infty)
\qquad, 
$$
thus $g_n$ has a convergent subsequence, due to properness of the action of $\calG$ on $N$; ($ii$) follows
specializing theorem \ref{generalred}.
\end{proof}
The quotient Lie algebroid $K/\calG\rightarrow N/\calG$ is, in a sense, the push forward of $K$ under the quotient
projection $\mathrm{q}:N\to N/\calG$. In fact the action $\sigma_K$, in the language of \cite{mackbook}, 
defines a linear action of $\mathrm{q}$ on $K$ and $K/\calG$ is the unique vector bundle $\ul{K}$
over $N/\calG$ such that $K=\mathrm{q}^{\sf{+}}\ul{K}$. Last proposition shows that such a vector bundle is canonically
endowed with a compatible Lie algebroid structure.\\
Note that the restriction of a morphic action to the kernel \tla-groupoid can be free and proper 
even when the top action is not.
We shall present below an example of this phenomenon arising from Poisson groupoid actions.

\bigskip

\spa Let us briefly review the basic facts about quotient Poisson structures.
Consider a Poisson action $\sigma: G\times P\rightarrow P$ of a 
Poisson group $G$.
If the action is free and proper $P/G$, is a smooth manifold that can be easily
endowed with a Poisson structure for which the quotient projection is a Poisson
submersion. For any $F\in\cif(G\times P)$, denote with $F^G_p\in\cif(G)$ and
$F^P_g\in\cif(P)$, $(g,p)\in G\times P$, 
$$
\qquad
(F^G_p)(g):=F(g,p)=:(F^P_g)(p)
\qquad,
$$
the restrictions to the $P$- and $G$-direction(s).
The Poisson bracket $\{\,,\,\}_{G\times P}$ of $G\times P$ can be expressed as
$$
\quad
\{F,H\}_{G\times P}(g,p)=\{F^G_p,H^G_p\}_G(g) + \{F^P_g,H^P_g\}_P(p)
\qquad,\qquad(g,p)\in G\times P\quad,
$$
in terms of the Poisson brackets $\{\,,\,\}_G$ of $G$ and 
$\{\,,\,\}_P$ of $P$. Smooth functions on $P/G$ are to be identified with
$G$-invariant functions on $P$; for any $f\in\cif(P/G)$ and 
$(g,p)\in G\times P$, $f(g*p)=f(p)$; then $(\sigma^*f)^G_p$ is constant on $G$ 
for all
$p$ and $(\sigma^*f)^P_g=\pr_B^*f$ for all $g$.
Since $\sigma$ is Poisson, for all $f,h\in\cif(P)^G$  and $(g,p)\in G\times P$,
\be
\{f,h\}_P(g*p)&=&\{(\sigma^*f)^P_g,(\sigma^*h)^P_g\}_P(p)
+
\{(\sigma^*f)^G_p,(\sigma^*h)^G_p\}_G(g)\\
&=&
\{f,h\}_P(p)
\quad,
\ee
i.e. $\{\,,\,\}_P$ restricts to a biderivation of $\cif(P)^G$. Thus, upon
identifying $\cif(P/G)$ with $\cif(P)^G$, $\{\,,\,\}_P$ yields a Poisson bracket
$\{\,,\,\}_{P/G}$ on $P/G$; note that Jacobi identity for $\{\,,\,\}_P$
implies the same property for $\{\,,\,\}_{P/G}$ and the quotient map is Poisson
by construction.\\
It is well known that quotient structures also exist for free and proper Poisson groupoid 
actions.
\bgn{theorem}\label{quangto}
Let a Poisson groupoid $\poidd{\calG}{M}$ act on  $P\rightarrow M$ freely and properly (so that the quotient
manifold $P/\calG$ is smooth). If the action is Poisson, there exists a unique Poisson structure on $P/\calG$, 
such that the quotient projection is a Poisson submersion.
\end{theorem}
%
%
\bgn{proof}
For any $f\in\cif(M)$, define $F\in\cif(\calG\times P\times\ol{P})$ setting
$F(g,p,q):=f(p)-f(q)$; then $f\in\cif(P)^\calG$ if{f} $F\in\calI_{\gr{\sigma}}$.
For all $f_\pm\in\cif(P)^\calG$,
$$
\qquad
\poib{F_+}{F_-}_{\calG\times P\times\ol{P}}(g,p,q)=\poib{f_+}{f_-}_P(p)
-
\poib{f_+}{f_-}_P(q) 
\qquad,
$$
since $\gr{\sigma}$ is coisotropic both sides of last equation vanish
identically on the graph of the action map, i.e.
$\poib{f_+}{f_-}\in\cif(P)^\calG$. 
\end{proof}
Next we shall describe quotient Poisson structures within the framework of reduction of
\tla-groupoids. Consider the moment morphism
$$
\qquad
\bcat\xy
*++{}="0",    <-1.7cm,0.7cm>
*++{T^*\calG\acts T^*P}="1", <0.7cm,0.7cm>
*++{T^*P}="2", <-1.7cm,-0.7cm>
*++{\calG\acts P}="3", <0.7cm,-0.7cm>
*++{P}="4",     <2cm,-0.7cm>
*++{T^*\calG}="1'", <4.2cm,-0.7cm>
*++{A^*}="2'", <2cm,-2.1cm>
*++{\calG}="3'", <4.2cm,-2.1cm>
*++{M}="4'",  <5.2cm,-0.7cm>
*++{}="2''"
\ar  @ <-0.07cm>   @{->} "1";"2"^{} 
\ar  @ <0.07cm>    @{->} "1";"2"^{}  
\ar  		   @{->} "1";"3"^{}
\ar                @{->} "2";"4"_{}
\ar  @ <0.07cm>    @{->} "3";"4"^{}
\ar  @ <-0.07cm>   @{->} "3";"4"_{}
\ar  @ <-0.07cm>   @{->} "1'";"2'"^{} 
\ar  @ <0.07cm>    @{->} "1'";"2'"^{}  
\ar  		   @{->} "1'";"3'"^{}
\ar                @{->} "2'";"4'"_{}
\ar  @ <0.07cm>    @{->} "3'";"4'"^{}
\ar  @ <-0.07cm>   @{->} "3'";"4'"_{}
\ar  		   @{->} "1";"1'"^{}
\ar  		   @{->} "2";"2'"^{\quad\jh}
\ar  		   @{->} "3";"3'"^{}
\ar  		   @{->} "4";"4'"^{\quad j}
\endxy\ecat
$$
associated with the cotangent lift of a Poisson $\calG$-space $j:P\to M$.
\bgn{proposition}\label{poissonreduction}
Let $P\rightarrow M$ a Poisson $\calG$-space. If $\calG$ acts freely and properly,  the reduced kernel  $K/\calG$
of the associated action \tla-groupoid  is a Lie algebroid canonically isomorphic to the  Koszul algebroid on
$T^*(P/\calG)$.
\end{proposition}
Note that for all $(g,\kappa)\in\calG\acts K_{(g,p)}$ and
 $(\delta g,\delta p)\in T\calG\acts TP_{(g,p)}$
\be
\pair{g*_K\kappa}{\delta g*\delta p}&=&\pair{0_g}{\delta g}+\pair{\kappa}{\delta p}\\
&=&
\pair{\kappa}{\delta p}
\qquad\qquad,
\ee
i.e the canonical pairing of $T^*P$ with $TP$ restricts to a $(\calG,T\calG)$-invariant pairing of  $K$ with $TP$,
which might be degenerate. However, passing to the quotients, we obtain a nondegenerate pairing $\ppair{}{}$,
$$
\qquad
\ppair{[\kappa]_{K/\calG}}{[\delta p]_{T(P/\calG)}}:=\pair{\kappa}{\delta p}\qquad,\qquad
(\kappa,\delta p)\in K\oplus TP
\qquad,
$$
of $K/\calG$ with $T(P/\calG)\equiv TP/T\calG$. The corresponding sharp map  $K/\calG\to T^*(P/\calG)$ is clearly
injective and according to formula (\ref{rankquo}) 
$$
\qquad
\rank K/\calG=\dim P-(\dim\calG-\dim M)=\dim P/\calG
\qquad,
$$ 
thus it is an isomorphism. It is now easy to prove proposition \ref{poissonreduction}.
\bgn{proof} For all $f\in\cif(P/\calG)\equiv\cif (P)^\calG$ and $a\in\Gamma(A)$,
$$
\qquad
\pair{\jh(\dd f)}{a}=\pair{\dd f}{X^a}=0
\qquad,
$$
for the infinitesimal action $X^\bullet$.
Thus $\dd f\in\Gamma(K)\subset\Ohm^1(P)$ represents the class
$\dd_{P/\calG}f\in\Ohm^1(P/\calG)\equiv\Gamma^\downarrow(P,K)\mod$ and for all 
$f_\pm\in\cif (P)^\calG$, $\poib{f_+}{f_-}_{P/\calG}$ is represented by $\poib{f_+}{f_-}$.
We have
\be
\rho_{K/\calG}\dd_{P/\calG} f&=&[\rho_{K}\dd f]_{T(P/\calG)}=[\pi^\sharp\dd f]_{T(P/\calG)}\\
&=&\pi_{P/\calG}^\sharp\dd_{P/\calG} f
\ee  
and for all $p\in P$
\be
\brak{\dd_{P/\calG}f_+}{\dd_{P/\calG}f_-}_{K/\calG}([p]_{P/\calG})
&=&
[\brak{\dd f_+}{\dd f_-}_{K}(p)]_{K/\calG}\\
&=&
[\dd \poib{f_+}{f_-}_p]_{K/\calG}\\
&=&
\dd_{P/\calG} (\poib{f_+}{f_-}_{P/\calG})_{[p]_{P/\calG}}\qquad.
\ee
Then, up to the identification of $K/\calG$ with $T^*(P/\calG)$ provided by $\ppair{}{}$, $K/\calG$
coincides with the Koszul algebroid. 
\end{proof}
\vs{1}
\section{Integrability of morphic actions}\label{ima}
\begin{quotation} 
We discuss in this Section morphic actions in the category of  Lie groupoids, introduced in \cite{bm92}, and develop an integrated version of the 
reduction procedures studied in the last Section. To obtain a quotient Lie groupoid from a morphic action of a double Lie
groupoid, beside the natural requirements for the quotients to be smooth, one has to further assume that the double
source map of the double Lie groupoid which is acting be surjective (theorem \ref{mainred}). Nevertheless a kernel
reduction procedure on moment morphisms (of double Lie groupoids) is effective under natural assumptions (proposition
\ref{kkk}). On the one hand, a free and proper morphic action of a double Lie groupoid always differentiates to a
morphic action of the associated \tla-group-oid and we further show that the quotient Lie groupoid of the original
action is a Lie groupoid integrating the quotient Lie algebroid of the induced morphic action of \tla-groupoid
(proposition \ref{oop}). On the other hand,  under a suitable completeness  condition in terms of the \tla-homotopy
lifting conditions of Chapter \ref{chapii}, a morphic action of an integrable \tla-groupoid $\sf{\Ohm}$ on an
integrable Lie algebroid  can be integrated to a morphic action in  the category of  Lie groupoids 
(theorem \ref{poo}).
\end{quotation}

\vs{0.1}
\spa Consider a double Lie groupoid $\pmb\calD$ and a morphism of Lie groupoids $\pmb\calJ$
$$
\bcat
\sf{D}:=
\xy
*+{}="0",    <-0.7cm,0.7cm>
*+{\calD}="1", <0.7cm,0.7cm>
*+{\calV}="2", <-0.7cm,-0.7cm>
*+{\calH}="3", <0.7cm,-0.7cm>
*+{M}="4",
\ar  @ <0.07cm>   @ {->} "1";"2"^{} 
\ar  @ <-0.07cm>  @ {->} "1";"2"_{}  
\ar  @ <-0.07cm>  @ {->} "1";"3"_{}
\ar  @ <0.07cm>   @ {->} "1";"3"^{}
\ar  @ <0.07cm>   @ {->} "2";"4"^{}
\ar  @ <-0.07cm>  @ {->} "2";"4"_{}
\ar  @ <-0.07cm>  @ {->} "3";"4"_{}
\ar  @ <0.07cm>   @ {->} "3";"4"^{}
\endxy
\qquad\qquad
\pmb{\calJ}:=
\xy
*+{}="0",    <-0.7cm,0.7cm>
*+{\calG}="1", <0.7cm,0.7cm>
*+{\calV}="3", <-0.7cm,-0.7cm>
*+{N}="2", <0.7cm,-0.7cm>
*+{M}="4",
\ar  @ <0.07cm>   @ {->} "1";"2"^{} 
\ar  @ <-0.07cm>  @ {->} "1";"2"_{}  
\ar  @ <-0.07cm>  @ {->} "3";"4"_{}
\ar  @ <0.07cm>   @ {->} "3";"4"^{}
\ar  		  @ {->} "1";"3"^{\calJ}
\ar  		  @ {->} "2";"4"_{j}
\endxy\qquad.
\ecat
$$
It is easy to see that, if actions of $\poidd{\calD}{\calV}$ on $\calJ$ and of $\calH$ on $j$ are given, the diagram
\bgn{equation}\label{actdoubgp}
\bcat
\xy
*+{}="0",    <-0.9cm,0.7cm>
*+{\calD\acts\calG}="1", <0.9cm,0.7cm>
*+{\calG}="2", <-0.9cm,-0.7cm>
*+{\calH\acts N}="3", <0.9cm,-0.7cm>
*+{N}="4",
\ar  @ <0.07cm>    @ {->} "1";"2"
\ar  @ <-0.07cm>   @ {->} "1";"2"
\ar  @ <-0.07cm>   @ {->} "1";"3"_{}
\ar  @ <0.07cm>    @ {->} "1";"3"^{}
\ar  @ <0.07cm>    @ {->} "2";"4"^{}
\ar  @ <-0.07cm>   @ {->} "2";"4"_{}
\ar  @ <-0.07cm>   @ {->} "3";"4"
\ar  @ <0.07cm>    @ {->} "3";"4"
\endxy
\ecat
\end{equation}
is a double Lie groupoid if{f} the action maps, respectively $\wt{\sigma}$ and $\sigma$, form a morphism of Lie
groupoids
$$
\pmb{\sigma}:=\bcat
\xy
*+{}="0",    <-1.5cm,0.7cm>
*+{\calD\fib{\sth}{\calJ}\calG}="1", <0.9cm,0.7cm>
*+{\calG}="3", <-1.5cm,-0.7cm>
*+{\calH\fib{\ssh}{j}N}="2", <0.9cm,-0.7cm>
*+{N}="4",
\ar  @ <0.07cm>   @ {->} "1";"2"^{} 
\ar  @ <-0.07cm>  @ {->} "1";"2"_{}  
\ar  @ <-0.07cm>  @ {->} "3";"4"_{}
\ar  @ <0.07cm>   @ {->} "3";"4"^{}
\ar  		  @ {->} "1";"3"^{}
\ar  		  @ {->} "2";"4"_{}
\endxy
\ecat\qquad.
$$
In particular the top horizontal source map $\nsth^\acts$ of (\ref{actdoubgp}) is clearly source submersive since
$$
\dd\ntth^\acts:T^{\stv^\acts}\calD\acts\calG=T^{\stv}\calD\fib{\dd\nsth}{\dd\calJ}T^\sor\calG\longrightarrow
T^\sor\calG
$$
is the restriction of the first projection and $\nsth$ is an \tlg-fibration. Then a statement analogous to
proposition \ref{actlagpd} holds true also in the case of a morphic actions of double Lie groupoids.
\bgn{proposition}\label{sef} With the above notations diagram ($\ref{actdoubgp}$) is a double Lie groupoid if{f}
$\pmb{\sigma}$ is morphic in the category of Lie groupoids. In that case $\pmb{\sigma}$ is an \tlg-fibration.
\end{proposition}

\spa For any morphic action $\pmb\sigma$ such as described above, we shall say that  $\pmb{\calD\acts\calG}$ is the
associated \tsf{action double groupoid}. If both  top and side actions are free and proper and provided a 
suitable regularity condition on the double source map of $\sf{D}$ is met, the top reduced space
carries a natural Lie groupoid over the side reduced space.
\bgn{theorem}\label{mainred}
Let $(\calD,\calH;\calV,N)$ a double Lie groupoid  with surjective double source map act morphically on a 
morphism of Lie groupoids $\calJ:\calG\to \calV$ over $j:N\to M$.  If the action is free and  proper 
(so that $\calG/\calD$ and $N/\calH$ are smooth manifolds), 
then there exists a unique Lie groupoid $\poidd{\calG/\calD}{N/\calH}$ making the quotient projection
$\calG\to\calG/\calD$ a strong \tlg-fibration  over $N\to N/\calH$. 
\end{theorem}
\bgn{proof} Thanks to proposition \ref{sef} it is sufficient to prove the statement for the  morphic action of
$\calD$ on its side vertical groupoid by left translation when $\poidd{\calD}{\calV}$ and $\poidd{\calH}{N}$ are free
and proper groupoids; the general case follows considering the action groupoid associated with the morphic action.
Note that the double source map $\mathbb{S}^\acts:\calD\acts\calG\to (\calH\acts N)\fib{\ssh^\acts}{\calJ}\calG$ is 
also a surjective, since it admits the factorization
$$
\bcat
\xy
*+{}="0",    <-3.5cm,1cm>
*+{\calD\fib{\sth}{\calJ}\calG}="1", <1.5cm,1cm>
*+{(\calH\acts N)\fib{\pr_2}{\calJ}\calG}="3", <-3.5cm,-1cm>
*+{(\calH\fib{\ssh}{\ssv}\calV)\fib{\pr_2}{\calJ}\calG}="2", <1.5cm,-1cm>
*+{\calH\fib{\ssh}{\ssv\comp\pr_1}\gr{\calJ}}="4",
\ar  @ {->} "1";"2"^{\:\simeq}_{\mathbb{S}\times\id_\calG\:}  
\ar  @ {->} "4";"3"_{\simeq}
\ar  @ {->} "1";"3"^{\mathbb{S}^\acts}
\ar  @ {->} "2";"4"_{\:\simeq}
\endxy
\ecat\qquad.
$$
For all $(d,v)\in\calD\fib{\sth}{}\calV$,
\be
\ssv(d*v)
&=&
\ssv(\ntth(d))
=
\tsh(\nstv(d))
=
\nstv(d)*\ssh(\nstv(d))
=
\nstv(d)*\ssv(\nsth(d))\\
&=&
\nstv(d)*\ssv(v)\hs{3},
\ee
that is, the source map of $\calV$ descends to the quotient; by similar arguments one can show that target and
inversion also descend. Note that for all $v^\pm\in\calV$, with $\ssv(v^+)=h*\tsv(v^-)$, for some $h\in\calH$, 
there always exist composable representatives: since there exist an element $d\in\calD$ such that 
$\nstv(d)=h\inverse$ and $\nsth(d)=v$, thanks to the regularity  condition on the double
source map, $d*v^+\equiv\ntth(d)$ is composable with $v^-$. Then one can define a multiplication
on the graph $(\calG/\calD, M/\calH)$ by picking composable representatives $w^\pm\in\ul{v^\pm}$ and setting, 
 $\ul{v^+}\cdot\ul{v^-}:=\ul{w^+\cdot_v w^-}$
where $\ul{v}$ denotes the class of 
$v\in\calV$.
%
Moreover, for all pairs 
$w^\pm_{1,2}\in\ul{v^\pm}$ of composable representatives such that $w^\pm_{2}=d^\pm*w^\pm_{1}$, one has
$$
\qquad
\ssv(w^+_2)=\nstv(d^+)*\ssv(w^+_1)\qquad\hbox{ and }\qquad\tsv(w^-_2)=\nttv(d^-)*\tsv(w^-_1)
\qquad,
$$
therefore $d^+$ and $d^-$ are vertically composable elements, since $\calH$ acts freely on $N$, and
\be
\ul{w^+_{2}}\cdot\ul{w^+_{2}}
&=&\ul{\ntth(d^+)\cdot_v\ntth(d^-)}=\ul{\ntth(d^+\cdot\vup d^-)}=\ul{(d^+\cdot\vup d^-)*(w^+_1\cdot_v w^-_1)}\\
&=&
\ul{w^-_{1}}\cdot\ul{w^-_{1}}\hs{4};
\ee
that is, the multiplication on $\calG/\calD$ does not depend on the choice of composable representatives.
It is straightforward to check that the induced source map is submersive. By the snake lemma one can see that
the quotient projection is an \tlg-fibration if{f} the top vertical source map of $\calD$ is orbitwise submersive,
i.e. if $\dd\nstv: T_v\calO^{\calD}_v\to T_{\ssv(v)}\calO^{\calH}_{\ssv(v)}$ is onto; this can be 
checked easily using submersivity of the double source map, since
$\calO^{\calD}_v$ and $\calO^{\calH}_{\ssv(v)}$ are locally diffeomorphic respectively to $\nsth\inverse(v)$ and
$\ssh\inverse(\ssv(v))$.
\end{proof}

\spa By functoriality, a morphic action of a double Lie groupoid differentiates to a a morphic
action on the vertically induced \tla-groupoid; in particular, all the fibred products relevant to the
definition of the induced morphic action exist, since \tlg-fibrations differentiate to \tla-fibrations. Namely,
when a morphic action as above is given, the action morphism differentiates to a morphism of Lie algebroids 
$$
\qquad
\bcat\xy
*+{}="0",    <-1.9cm,0.7cm>
*+{A\vup(\calD)\fib{\sh}{\jh}A(\calG)}="1", <0.9cm,0.7cm>
*+{A(\calG)}="3", <-1.9cm,-0.7cm>
*+{\calH\fib{\sor}{j}N}="2", <0.9cm,-0.7cm>
*+{N}="4",
\ar    @{->} "1";"2"_{} 
\ar  		   @{->} "1";"3"^{}
\ar                @{->} "2";"4"_{}
\ar  @ <0.07cm>    @{->} "3";"4"
\endxy\ecat\qquad,
$$
where we identify the Lie algebroid $A\vup(\calD\acts \calG)$ of $\poidd{\calD\acts \calG}{\calH\acts N}$ with 
the fibred product $A\vup(\calD)\fib{\sh}{\jh}A(\calG)$ for the morphisms 
$\sh\colon A\vup(\calD)\to A(\calV)$ and $\jh\colon A(\calG)\to A(\calV)$, differentiating the top horizontal source map
 and  the moment map of the top action. It is immediate to check that,
the associated  action double Lie groupoid differentiates to the action \tla-groupoid for the induced morphic
action of Lie algebroids. 
\bgn{proposition}\label{oop} 
Let $(\calD,\calH;\calV,N)$ be a double Lie groupoid acting morphically on a morphism of Lie algebroids
$\calJ:\calG\to \calV$ over $j:N\to M$ and assume that the quotient Lie groupoid $\poidd{\calG/\calD}{N/\calH}$
exists and makes the quotient projection a strong \tlg-fibration. Then the quotient Lie algebroid
$A(\calG)/A\vup(\calD)\to N/\calH$ for the induced morphic  action of Lie algebroids exists and is canonically
isomorphic to the Lie algebroid of  $\poidd{\calG/\calD}{N/\calH}$.
\end{proposition}
\bgn{proof}
As in the proof of theorem \ref{mainred} it is sufficient to consider the case of horizontally free and proper double
Lie groupoids for the action on the side vertical groupoid by left translation. Since the tangent prolongation
groupoid $\poidd{T\calD}{T\calV}$ is also  free and proper, $A\vup(\calD)\subset T\calD$ and $A(\calV)\subset T\calV$
are embedded as normal bundles and the induced groupoid anchor $\ol{\chi}: A\vup(\calD)\to A(\calV)\times A(\calV)$
is the restriction of the top horizontal tangent anchor $\dd\chi\hup:T\calD\to T\calV\times T\calV$,
the Lie groupoid $\poidd{A\vup(\calD)}{A(\calV)}$ is a free and proper and the quotient Lie algebroid
$A(\calV)/A\vup(\calD)\to M/\calH$ exists. The quotient projection $\calV\to \calV/\calD$ is a strong \tlg-fibration
over $M\to M/\calH$, thus it differentiates to a strong \tla-fibration $\wp:A(\calV)\to A(\calV/\calD)$.
For all $q\in M$, define 
\bgn{equation}\label{sdrr}
\psi_{q}:(A(\calV)/A\vup(\calD))_{\ul{q}}\to A(\calV/\calD)_{\ul{q}}
\end{equation}
as $\psi_{q}([a_q])=\wp_q(a_q)$, by picking any $a_q\in A(\calV)_q$ such that $a_q\in[a_q]$; $\psi_{q}$ is well
defined since the commuting diagram
$$
\bcat
\xy
*+{}="0",    <0cm,1.2cm>
*+{\calD\fib{\sth}{}\calV}="t", <0cm,-1.2cm>
*+{\calV/\calD}="b", <-1.2cm,0cm>
*+{\calV}="l", <+1.2cm,0cm>
*+{\calV}="r",  
\ar     @ {->} "t";"l"_{\pr_2} 
\ar     @ {->} "t";"r"^{\sigma}  
\ar     @ {->} "l";"b"_{}
\ar     @ {->} "r";"b"^{}
\endxy
\ecat
\qquad\hbox{ differentiates to }
\qquad
\bcat
\xy
*+{}="0",    <0cm,1.2cm>
*+{A\vup(\calD)\fib{\sh}{}A(\calV)}="t", <0cm,-1.2cm>
*+{A(\calV/\calD)}="b", <-1.2cm,0cm>
*+{A(\calV)}="l", <+1.2cm,0cm>
*+{A(\calV)}="r",  
\ar     @ {->} "t";"l"_{\pr_2} 
\ar     @ {->} "t";"r"^{\sigmah}  
\ar     @ {->} "l";"b"_{\wp}
\ar     @ {->} "r";"b"^{\wp}
\endxy
\ecat
$$
Note that for all elements $h*q$ in the $\calH$-orbit through $q$,
$$
\psi_{h*q}([a_{h*q}])=\wp(a_{h*q})=\wp(\ohm_h*a_{q})=\psi_{q}([a_{q}])=\psi_{q}([a_{h*q}])
$$
for the unique $\ohm_h\in A\vup(\calD)_h$, such that $\ol{\chi}(\ohm_h)=((a_{h*q}),a_q)$. Then $\psi_{q}$ does not 
depend on the choice of $q$ and it is surjective, since $\wp$ is a strong \tla-fibration; by counting dimensions, 
one can see that it is actually a linear isomorphism, inducing a bundle isomorphism $A(\calV/\calD)\to
A(\calV)/A\vup(\calD)$ over the identity of $M/\calH$; up to this identification, the Lie algebroid on
$A(\calV/\calD)$ coincides with that on  $A(\calV)/A\vup(\calD)$ by uniqueness (theorem \ref{generalred}).
\end{proof}

\spa Next we shall consider the integrability of morphic actions of \tla-groupoids. Under natural conditions for
an \tla-groupoid $(\Ohm,\calG;A,M)$ to have a source 1-connected integration $\sf\Gamma$, all of its morphic
actions  on morphisms of integrable Lie algebroids also integrate to morphic actions of $\sf\Gamma$.
\bgn{theorem}\label{poo} Let $\sf\Ohm:=(\Ohm,\calG;A,M)$ be an \tla-groupoid and $\jh:B\to A$ a morphism of Lie
algebroids  over $j:N\to M$. Assume that $(\Ohm,\calG;A,M)$ has integrable top Lie algebroid and $B$ is also
integrable with source 1-connected Lie groupoid $\calB$. Then, if the top source map of $(\Ohm,\calG;A,M)$
satisfies the \tla-homotopy lifting conditions of definition \ref{liftingprop},  any morphic action of $\sf\Ohm$ on $\jh$
integrates to a morphic action of the vertically  source 1-connected double Lie groupoid
$\sf\Gamma:=(\Gamma;\calG,\calA,M)$ on the integration  $\calJ:\calB\to \calA$.   
\end{theorem}
\bgn{proof} Note that, under the assumptions, the top source map $\nsth:\Gamma\to \calA$ of $\sf\gamma$ is a
\tlg-fibration and the top source map $\sh$ of $\Ohm$ and $\jh$ are strongly transversal; therefore the fibered
products Lie groupoids
$$
\Gamma\fib{\sth}{\calJ}\calB
\qquad\hbox{ and }
\qquad
\Gamma\fib{\sth}{\calJ}(\Gamma\fib{\sth}{\calJ}\calB)
\simeq
\Gamma\hup^{\!\!\!\!\hbox{\tiny{$(2)$}}}\fib{\sth\comp\pr_1}{\calJ}\calB
$$
are well defined and source 1-connected for the integration $\wt{\sigma}:\Gamma\fib{\sth}{\calJ}\calB\to\calB$ of
the top action map $\sigmah:\Ohm\fib{\sh}{\jh}B\to B$. The compatibility diagrams \ref{compaact} for
$\wt{\sigma}$ to be an action map compatible with $\calJ$ commute then by functoriality. 
\end{proof}
It follows from the last theorem and proposition \ref{oop} above that the reduction $\poidd{\calB/\Gamma}{N/\calG}$, 
provided it exists, is an integration of the reduction $B/\Ohm\to N/\calG$.

\bigskip

\spa As in the case of morphic actions of \tla-groupoids there is a moment morphism associated with each
morphic action of a double Lie groupoid whose kernel is well behaved under reduction.\\
According to lemma \ref{kerla} the \tsf{moment morphism} 
$$
\bcat
\xy
*++{}="0",    <-1.6cm,0.7cm>
*++{\calD\acts\calG}="1", <0.7cm,0.7cm>
*++{\calG}="2", <-1.6cm,-0.7cm>
*++{\calH\acts N}="3", <0.7cm,-0.7cm>
*++{N}="4",     <1.9cm,-0.7cm>
*++{\calD}="1'", <4.2cm,-0.7cm>
*++{\calV}="2'", <1.9cm,-2.1cm>
*++{\calH}="3'", <4.2cm,-2.1cm>
*++{M}="4'",  <5.2cm,-0.7cm>
\ar  @ <-0.07cm>   @{->} "1";"2"^{} 
\ar  @ <0.07cm>    @{->} "1";"2"^{}  
\ar  @ <0.07cm>   @{->} "1";"3"_{}
\ar  @ <-0.07cm>  @{->} "1";"3"_{}
\ar  @ <0.07cm>   @{->} "2";"4"_{}
\ar  @ <-0.07cm>  @{->} "2";"4"_{}
\ar  @ <0.07cm>    @{->} "3";"4"^{}
\ar  @ <-0.07cm>   @{->} "3";"4"_{}
\ar  @ <-0.07cm>   @{->} "1'";"2'"^{} 
\ar  @ <0.07cm>    @{->} "1'";"2'"^{}  
\ar  @ <0.07cm>   @{->} "1'";"3'"_{}
\ar  @ <-0.07cm>  @{->} "1'";"3'"_{}
\ar  @ <0.07cm>   @{->} "2'";"4'"_{}
\ar  @ <-0.07cm>  @{->} "2'";"4'"_{}
\ar  @ <0.07cm>    @{->} "3'";"4'"^{}
\ar  @ <-0.07cm>   @{->} "3'";"4'"_{}
\ar  		   @{->} "1";"1'"^{\pr_\calD}
\ar  		   @{->} "2";"2'"^{\calJ}
\ar  		   @{->} "3";"3'"^{}
\ar  		   @{->} "4";"4'"^{j}
\endxy
\ecat
$$
has a kernel double Lie groupoid if{f} $\calJ$ is source submersive, typically if it is an \tlg-fibration. In
that case it is easy to identify the kernel with an action double Lie groupoid
$$
\bcat
\xy
*+{}="0",    <-1cm,0cm>
*+{\calH\acts\calK}="1", <1.4cm,0cm>
*+{\calK}="2", <-1cm,-1.4cm>
*+{\calH\acts N}="3", <1.4cm,-1.4cm>
*+{M}="4",
\ar  @ <0.07cm>   @{->} "1";"2"^{} 
\ar  @ <-0.07cm>  @{->} "1";"2"_{}  
\ar  @ <-0.07cm>  @{->} "1";"3"_{}
\ar  @ <0.07cm>   @{->} "1";"3"^{}
\ar  @ <0.07cm>   @{->} "2";"4"^{}
\ar  @ <-0.07cm>  @{->} "2";"4"_{}
\ar  @ <-0.07cm>  @{->} "3";"4"_{}
\ar  @ <0.07cm>   @{->} "3";"4"^{}
\endxy
\ecat
$$
for the restriction $\sigma_{\calK}$ of the top action,
$$
\quad
h*_{\calK}\kappa\equiv\sigma_{\calK}(h,\kappa):=\wt{\sigma}(\etv(h),\kappa)
\quad,\quad \esv(\ssh(h))=\calJ(\kappa)
\quad,\quad (h,\kappa)\in\calH\times\calK\quad,
$$
to the kernel groupoid $\calK:=\ker\calJ$. Remarkably, not only is $\sigma_\calK$ a morphic action of the
double Lie groupoid trivially associated to $\calH$ over
$\sigma$ along $\calJ_\calK:=\ssv\comp J|_\calK=\tsv\comp J|_\calK$, but it is also a compatible action in
the sense of definition \ref{compaction}. As a consequence the kernel double Lie groupoid always has a quotient, provided
the side action is free and proper.
\bgn{proposition}\label{kkk} Let $(\calD,\calH;\calV,N)$ be a double Lie groupoid acting morphically on a
morphism of Lie algebroids $\calJ:\calG\to \calV$ over $j:N\to M$.  If $\calJ$ is source submersive 
(so that the kernel double Lie groupoid exists), and the side action is free and proper, then:
\medskip\\
$i)$ The induced action of $\calH$ on $\calK=\ker\calJ$ is a compatible groupoid action in the sense
of definition \ref{compaction};
\medskip\\
$ii)$ The quotient $\calK/\calH$ carries a unique Lie groupoid structure over $N/\calH$ making the
quotient projection a strong \tlg-fibration.
\medskip\\
$iii)$ The source fibres of $\calK/\calH$ have the same homotopy type as those of $\calK$
\end{proposition}
\bgn{remark}
Note that no surjectivity condition on the double source map is necessary. Moreover, one can see from the proof
below that \emph{the second and third statements remain true when only a compatible groupoid action, in the sense of
definition \ref{compaction}, of $\calH$
on $\calK$ is given}. 
\end{remark}
\bgn{remark}
It follows from proposition \ref{oop} that the Lie algebroid of $\calK/\calH$ is canonically
isomorphic to the reduction $\ker\jh/\calH$ associated with the moment morphism of the induced action of
\tla-groupoids.
\end{remark}
\bgn{proof} ($i$) Clearly the induced action is compatible with source, target and inversion of $\calK$.
Moreover, for all composable $\kappa_\pm\in\calK$,
$$
\qquad
\calJ_\calK(\kappa_+)=j(\sor_\calK (\kappa_+))=j(\tar_\calK (\kappa_-))=\calJ_\calK(\kappa_-)
\qquad,
$$
thus
\be
h*_{\calK}(\kappa_+\cdot\kappa_-)
&=&
(\esv(h)\cdot_V\esv(h))*(\kappa_+\cdot\kappa_-)\\
&=&
(\esv(h)*\kappa_+)\cdot_V(\esv(h)*\kappa_-)\\
&=&
(h*_{\calK}\kappa_+)\cdot(h*_{\calK}\kappa_-)\hs{3}.
\ee
($ii$) Since the side action is free and proper, $N\to N/\calH$ is a principal $\calH$ bundle; apply lemma 
\ref{pgbred}.
Alternatively the statement can be proved using a nonlinear version of the proof of lemma \ref{Step 2.}. For
all $n\in N$ and $h\in\ssh\inverse(j(h))$, set 
$$
\qquad
\theta_h:\sor_\calK\inverse(n)\to \sor_\calK\inverse(h*n)\qquad,\qquad \kappa\mapsto h*_\calK\kappa 
\qquad,
$$
this yields a smooth family of diffeomorphisms enjoying the usual pseudo-group property. It is easy to see  that the
obvious equivalence relation $\thicksim_\theta$ induced by $\theta_\bullet$ is regular (properness of  $\sigma_\calK$
is implied by properness of the base action map). Being $\sigma_\calK$ an action map,  one can check that actually
$\gr{\thicksim_\theta}\subset \calK\times\calK$ is a Lie subgroupoid of the pair groupoid, as well as the graph 
$\gr{\thicksim}\subset N\times N$ of the equivalence relation associated with the action on $N$ and
$\poidd{\gr{\thicksim_\theta}}{\gr{\thicksim}}$ a Lie subgroupoid of the direct product $\twice{\calK}$. That is, in
the language \cite{mackbook}, $(\gr{\thicksim_\theta},\gr{\thicksim})$  defines a congruence on $\poidd{\calK}{N}$,
which therefore provides the descent data to push the Lie groupoid  forward on $N/\calH$ along the orbits of the side action
\cite{mackbook}\footnote{A congruence in \cite{mackbook} is further required to fulfill certain
surjectivity requirements which are not met here; the extra condition is only needed to make the quotient projection
a \emph{strong} \tlg-fibration, i.e. a fibration of Lie groupoids in the language of Mackenzie.}.
($iii$) Note that for all $n\in N$ and $\kappa_\pm\in\sor\inverse_\calK(n)$, $h*_\calK\kappa_-=\kappa_+$ implies that
$(h,n)\in\calH\acts N$ is an isotropy, therefore a unit; that is, $\calH$ acts transversally to the source fibres and
the quotient projection $\calK\to\calK/\calH$ is sourcewise a diffeomorphism.
\end{proof}
In complete analogy with the reduction of the moment morphism of a morphic action in the category of Lie
algebroids, the reduced kernel groupoid is really a ``pushforward'' of the kernel groupoid, obtained by
identifying the source fibres along the $\calG$-orbits on $P$.
\vs{1}
\section{Integration of quotient Poisson structures}\label{iqps}
\begin{quotation} 
We derive here two approaches to the integration of quotient Poisson structures for compatible Poisson groupoid actions. The
first consist in the integration of the action \tla-groupoid associated with the cotangent lifted action and in the kernel
reduction of the corresponding moment morphism (theorem \ref{quoproqui}). This approach is not always effective but, when it
is, it has the advantage of providing an explicit description of the symplectic form on an integration of the quotient Poisson
bivector and control on the connectivity of its source fibres (corollary \ref{erty}). The second approach allows us to 
obtain necessary and sufficient conditions for the integrability of quotient Poisson bivector fields 
(theorem \ref{quiproquo}), but no explicit description of the symplectic form on the integrations.
The result is obtained roughly by a  prolongation of the cotangent lifted action to a compatible action on the space of Lie
algebroid homotopies for the kernel of the cotangent lifted moment morphism.
We finally consider a class of examples, i.e. the case of complete Poisson group actions, where both approaches apply.   
\end{quotation}
\spa Let us consider the reduction of the moment morphism $\pmb{\calJ}$ associated with a compatible morphic
action of a symplectic double groupoid $\pmb{\calS}$
$$
\bcat
\pmb{\calS}=
\xy
*+{}="0",    <-0.7cm,0.7cm>
*+{\calS}="1", <0.7cm,0.7cm>
*+{\calG^\bullet}="2", <-0.7cm,-0.7cm>
*+{\calG}="3", <0.7cm,-0.7cm>
*+{M}="4",
\ar  @ <0.07cm>   @{->} "1";"2"^{} 
\ar  @ <-0.07cm>  @{->} "1";"2"_{}  
\ar  @ <-0.07cm>  @{->} "1";"3"_{}
\ar  @ <0.07cm>   @{->} "1";"3"^{}
\ar  @ <0.07cm>   @{->} "2";"4"^{}
\ar  @ <-0.07cm>  @{->} "2";"4"_{}
\ar  @ <-0.07cm>  @{->} "3";"4"_{}
\ar  @ <0.07cm>   @{->} "3";"4"^{}
\endxy
\qquad\qquad
\pmb{\calJ}=
\xy
*+{}="0",    <-0.7cm,0.7cm>
*+{\Lambda}="1", <0.7cm,0.7cm>
*+{\calG^\bullet}="3", <-0.7cm,-0.7cm>
*+{P}="2", <0.7cm,-0.7cm>
*+{M}="4",
\ar  @ <0.07cm>   @ {->} "1";"2"^{} 
\ar  @ <-0.07cm>  @ {->} "1";"2"_{}  
\ar  @ <-0.07cm>  @ {->} "3";"4"_{}
\ar  @ <0.07cm>   @ {->} "3";"4"^{}
\ar  		  @ {->} "1";"3"^{\calJ}
\ar  		  @ {->} "2";"4"_{j}
\endxy\qquad,
\ecat
$$
where $\poidd{\Lambda}{P}$ is a symplectic groupoid, $\calJ$ an anti-Poisson map 
(hence a morphism of Poisson groupoids $\Lambda\to\ol{\calG^\bullet}$) 
and both the top and side actions are Poisson; in particular, since both $\calS$ and $\Lambda$ are symplectic, the
graph of the top action is a Lagrangian subgroupoid of $\calS\times\Lambda\times\ol{\Lambda}$ for the vertical groupoid
of $\calS$. Let us denote with $\Ohm$ the symplectic form of $\calS$ and with $\ohm$ that of $\Lambda$, the last
requirement is equivalent to the multiplicativity condition
\bgn{equation}\label{multisigma}
\wt{\sigma}^*\ohm=\pr^*_{\calS}\Ohm+\pr^*_{\Lambda}\ohm
\end{equation}
on $\calS\acts\Lambda$. Explicitly, (\ref{multisigma}) reads
$$
\qquad
\ohm_{s*\lambda}(\delta s_+*\delta \lambda_+,\delta s_+*\delta \lambda_+)
=
\Ohm_{s}(\delta s_+,\delta s_+)
+
\ohm_{\lambda}(\delta \lambda_+,\delta \lambda_+)
\qquad,
$$
where we have used the symbol $*$ for the tangent lift of $\wt{\sigma}$, for all composable 
$\delta s_\pm\in T\calS$ and $\delta \lambda_\pm\in T\Lambda$. In this setting the kernel double groupoid of
the moment morphism takes the form
$$
\qquad
\bcat
\xy
*+{}="0",    <-1.9cm,0cm>
*+{\calG\acts\calJ\inverse(\eps_\bullet(M))}="1", <1.4cm,0cm>
*+{\calJ\inverse(\eps_\bullet(M))}="2", <-1.9cm,-1.4cm>
*+{\calG\acts P}="3", <1.4cm,-1.4cm>
*+{M}="4",
\ar  @ <0.07cm>   @{->} "1";"2"^{} 
\ar  @ <-0.07cm>  @{->} "1";"2"_{}  
\ar  @ <-0.07cm>  @{->} "1";"3"_{}
\ar  @ <0.07cm>   @{->} "1";"3"^{}
\ar  @ <0.07cm>   @{->} "2";"4"^{}
\ar  @ <-0.07cm>  @{->} "2";"4"_{}
\ar  @ <-0.07cm>  @{->} "3";"4"_{}
\ar  @ <0.07cm>   @{->} "3";"4"^{}
\endxy
\ecat
\qquad,
$$
where $\eps_\bullet:M\to \calG^\bullet$ is the unit section (here we assume that $\calJ$ is regular enough).
Under the hypothesis of theorem \ref{kkk} for the underlying morphic action the reduction of the \tla-groupoid
above produces a Lie groupoid $\poidd{\calJ\inverse(\eps_\bullet(M))/\calG}{P/\calG}$; 
the reduction procedure is compatible with the side Poisson action in the sense of the following
\bgn{theorem}\label{quoproqui} Let $(\calS,\calG;\calG^\bullet,M)$ be a symplectic double groupoid acting 
morphically on a morphism of Lie groupoids $\calJ:\Lambda\to\calG^\bullet$ over $j:P\to M$, where
$\poidd{\Lambda}{P}$ is a symplectic groupoid, in such a way that $\Lambda$ is a symplectic $\calS$-space. 
If $\calJ$ is source submersive and the side action is free and proper, then
\medskip\\
$i)$ The reduced kernel groupoid $\poidd{\calJ\inverse(\eps_\bullet(M))/\calG}{P/\calG}$ carries a unique symplectic
form making the quotient projection  $\pr:\calJ\inverse(\eps_\bullet(M))\to\calJ\inverse(\eps_\bullet(M))/\calG$  a
Poisson submersion;
\medskip\\
$ii)$ $\poidd{\calJ\inverse(\eps_\bullet(M))/\calG}{P/\calG}$ is a symplectic groupoid for the quotient Poisson
manifold $P/\calG$.
\end{theorem}
\bgn{proof} The quotient $\calJ\inverse(\eps_\bullet(M))/\calG$ is smooth and carries a Lie groupoid on $P/\calG$
thanks to proposition \ref{kkk}.  Set $\calK:=\calJ\inverse(\eps_\bullet(M))$.  ($i$) On the one hand, 
by coisotropicity of $\eps_\bullet(M)\subset\calG^\bullet$, $\calK\subset\calS$ is also coisotropic, since $\calJ$ is
anti-Poisson, and evaluating the multiplicativity
condition $\ref{multisigma}$ on $\calG\acts\calK$ yields
$$
\qquad
\ohm_{g*\kappa}(\delta g_+*\delta \kappa_+,\delta g_+*\delta \kappa_+)
=
\ohm_{\kappa}(\delta \kappa_+,\delta \kappa_+)
\qquad,
$$
being $\calG\subset\calS$ Lagrangian, for all composable 
$\delta g_\pm\in T\calG$ and $\delta \kappa_\pm\in T\calK$.
Setting
\bgn{equation}\label{thy}
\qquad
\ul{\ohm}_{[\kappa]}([\delta\kappa_+],[\delta\kappa_-])
:=
\ohm_{{\kappa}}({\delta\kappa_+},{\delta\kappa_-})
\qquad
\end{equation}
for any representatives  $\delta\kappa_\pm\in[\delta\kappa_\pm]\in T_{[k]}\calK/\calG$  over the same
$\kappa\in[\kappa]$, yields a 2-form $\ul{\ohm}$ on $\calK/\calG$. The left hand side of (\ref{thy}) does not depend
on the choice of $\delta\kappa_\pm$, provided they are tangent to $\calK$ at the same $\kappa$ (in fact there are
unique such representatives) and,  by changing the representative of $[\kappa]$, we have
$$
\qquad
\ohm_{g*\kappa}(\kappa_+',\delta \kappa_+')
=
\ohm_{\kappa}(\delta g\inverse_+*\delta \kappa_+',\delta g\inverse_-*\delta \kappa_+')
\qquad,
$$
for all $\delta g_\pm$ such that $\dd\eps_\bullet(\dd\tar(\delta
g_\pm))=\dd\calJ(\dd\kappa_\pm')$. Then $\ul{\ohm}$ is well defined and multiplicative;
by construction $\pr^*\ul{\ohm}=\iota^*\ohm$,  thus it is closed. We claim that the
characteristic distribution $\Delta$ of ${\calK}$ spans $T_\kappa\calO$ to the
$\calG$-orbit $\calO$ through $\kappa$ at each $\kappa\in\calK$; it follows that
$\ul{\ohm}$ is nondegenerate. To see this, note that all vectors $\delta o\in
T_\kappa\calO$ are those of the form $a*0_k$ with $a\in T^\sor_{\eps(j(k))}\calG$,
thus
$$
\qquad
\ohm_{\kappa}(\delta o,\delta \kappa)
=
\ohm_{\kappa}(a*0_\kappa,0_{\eps(j(\kappa))}*\delta \kappa)=0
\qquad,\quad\hbox{ for all }\quad \delta \kappa\in T_\kappa\calK,
$$
that is, $T_\kappa\calO$ is contained in the symplectic orthogonal 
$T_\kappa^\ohm\calK=\Delta_\kappa$; the two
spaces coincide since
\be
\dim\calO&=&\dim\calG-\dim M
\quad=\quad
\dim\Lambda - (\dim M +\dim\Lambda -\dim\calG^\bullet)\\
&=&
\rank T^\ohm\calK.
\ee
($ii$) Since the characteristic leaves of $\calK$ are the connected components of the
$\calG$-orbits, $\cif(\calK)^\calG\subset\cif(\calK)^\Delta$; thus
all extensions $F_\pm\in\cif(\Lambda)$ of $f_\pm\in\cif(\calK)^\calG$ are in the 
normalizer of $\calI_{\calK)}$ and the restriction of the Hamiltonian vector fields
$X^{F_\pm}$ are tangent to $\cif(\calK)^\calG$. Let $X^F_\pm=\dd\iota Y^{F_\pm}$ for
some (in general non smooth) vector fields $Y^{F_\pm}$ on $\calK$; we have
$$
\qquad
(\dd\pr^{\sf{t}}_{\pr(\kappa)}\comp\ohm^\sharp_{\pr(\kappa)}\comp\dd\pr_{k})Y^{F_\pm}
=
(\dd\iota^{\sf{t}}_{\iota(\kappa)}\comp\ohm^\sharp_k)
X^{F_\pm}_{\iota(\kappa)}=\dd_{P}f_\pm
=
\dd\pr^{\sf{t}}_{\pr(\kappa)}\dd_{P/\calG}f_\pm
\qquad,
$$
i.e. the Hamiltonian vector fields $\ul{X}^{f_\pm}$ of $f_\pm\in\cif(\calK/\calG)$ are
given by $\ul{X}^{f_\pm}_{\pr(\kappa)}=\dd\pr Y^{F_\pm}_\kappa$ and the Poisson bracket
$\poib{}{}_{\calK/\calG}$ associated with $\ul{\ohm}$ can be computed using extensions:
\be
\pr^*\poib{f_+}{f_-}_{\calK/\calG}
&:=&
\ul{\ohm}(\ul{X}^{f_-},\ul{X}^{f_+})\comp\pr
=
\pr^*\ul{\ohm}(Y^{F_-},Y^{F_+})
=
\ohm(X^{F_-},X^{F_+})\comp\iota\\
&=:&\iota^*\poib{F_+}{F_-}_{\Lambda}\hs{3}.
\ee
Let now $\poib{}{}'$ be the Poisson bracket induced by the symplectic groupoid of ($i$) on $P/\calG$ and
$u_\pm\in\cif(P)^\calG$, since
\be
\poib{u_+}{u_-}'([p]_{P/\calG})
&:=&
\poib{\sor_{\calK/\calG}^*u_+}{\sor_{\calK/\calG}^*u_-}_{\calK/\calG}(\eps_{\calK/\calG}([p]_{P/\calG}))\\
&=&
\poib{\sor_\Lambda^*u_+}{\sor_\Lambda^*u_-}_{\Lambda}(\eps_\Lambda(p)))\\
&=&
\poib{u_+}{u_-}_{P}(p)\\
&=:&
\poib{u_+}{u_-}_{P/\calG}([p]_{P/\calG})\hs{3},
\ee
for all $p\in P$, ($ii$) follows by uniqueness (theorem \ref{quangto}).
\end{proof}
It follows from proposition \ref{kkk} and the proof above that the source fibres of 
$\calJ\inverse(\eps_\bullet(M))/\calG$ have the same homotopy type as those of $\calJ\inverse(\eps_\bullet(M))$. 
Therefore  we obtain a condition for  $\calJ\inverse(\eps_\bullet(M))$ to be source 1-connected in terms of the
infinitesimal data only, i.e. a condition for integration to commute with reduction, as an 
application of corollary \ref{lupo}.
\bgn{corollary}\label{erty}
Assume that $\Lambda$ and $\calG^\bullet\equiv\calG^\star$ are source 1-connected. Then the source connected
component of $\calJ\inverse(\eps_\star(M))/\calG$ is the source 1-connected integration of $P/\calG$ if{f} the loop
groups $\KK_\bullet(\jh)$ are trivial.
\end{corollary}
The condition of lemma \ref{erty} was recently considered in \cite{for07}, along the lines of \cite{fs06}, in the special case
of Poisson  actions of Lie groups and in \cite{07b}, in the case of Poisson actions of Poisson groups.
\bgn{remark} Assume that a Poisson groupoid $\calG$ induces a top source map satisfying the \tla-homotopy lifting conditions 
of definition \ref{liftingprop} on the cotangent prolongation \tla-groupoid, and therefore is integrable to a symplectic
double groupoid.  Then last result, together with propositions \ref{poissonreduction}, \ref{poo} implies that, for any integrable
Poisson  $\calG$-space $P$, the quotient Poisson bivector on $P/\calG$ is also integrable, since the integrated action map
defines a symplectic action on the integration $\calJ$ of the cotangent lifted moment map $\jh:T^*P\to A^*$ (this follows
reasoning in the same way as in the proof of theorem \ref{double}). In fact an integration -- generally not source
(1-)connected -- is given by the quotient $\calJ\inverse(\eps_\star(M))/\calG$.
\end{remark} 

\spa By a suitable path-lifting procedure of the cotangent lifted action of $\poidd{T^*\calG}{A^*}$ we can prove
integrability of quotient Poisson manifolds arising from Poisson $\calG$-spaces, 
independently of the existence of a double (or any integration ) of $\calG$. 
\bgn{theorem}\label{quiproquo} Let a Poisson groupoid $\gpdm$ act freely and properly on $j:P\to M$. If the action is
Poisson, then $P/\calG$ is integrable to a symplectic groupoid if{f} $\ker{\jh}$ is an integrable Lie algebroid.
\end{theorem}
In particular, we give a positive answer to the above question under most natural assumptions: since a Lie subalgebroid of an 
integrable Lie algebroid is also integrable (theorem \ref{1}), we have:
\bgn{corollary}\label{happy} If $P$ is integrable, then so is $P/\calG$.
\end{corollary}
%
%
We show in two steps that integrability of $\ker{\jh}$ is a sufficient condition:
\medskip\\{\em Step 1.} Through the cotangent lift of the $\calG$-action on $j$, we obtain a compatible groupoid
action of $\calG$ on the Weinstein groupoid $\calW(K)$ of the kernel Lie algebroid $K:=\ker\jh$  of the moment
map $\jh:T^*P\to A^*$
\medskip\\{\em Step 2.} We identify the quotient $\calW(K)/\calG$ with the symplectic groupoid of $P/\calG$
\medskip\\ 
Finally we explain why the integrability condition is also necessary.

\spa Let us fix some notations. For $i=1,2$, denote with $\Delta_i$ the diagonal morphism 
$TI^{\times i}\to TI^{\times i}\times TI^{\times i}$ and regard $0_g\in T^*\calG$ as the constant morphism of
Lie algebroids $TI^{\times i}\to T^*\calG$, denoted $g_{(i)}$.
\bgn{proof}[Proof of theorem \ref{quiproquo}]
({\em Step 1.}) Note that, for all morphisms of Lie algebroids $h_i:TI^{\times i}\to K$ over 
$\gamma_i:I^{\times i}\to\calG$,
$$
\qquad
\dd j\comp\dd\gamma_i=\dd j\comp\pish\comp h_i=\rho_{A^*}\comp\jh\comp h_i=0
\qquad,
$$
where $\pi$ is the Poisson bivector of $P$,
thus the image of the base map $\gamma_i$ is contained in some $j$-fibre. If $j\comp\gamma_i\equiv \sor(g)$,
$g\in\calG$, setting
$$
g*h_i:=\sigmah\comp(g_{(i)}\times h_i)\comp\Delta_i
$$
yields a morphism of Lie algebroids $TI^{\times i}\to K$, since 
$\jh\comp (g*h_i)=\th\comp g_{(i)}\equiv 0_{\tar(g)}$. The moment map $J:\calW(K)\to M$ for the lifted action
is induced by post composition with $\jh$ of representatives of classes of $K$-paths, equivalently
$$
\qquad
J([\kappa]):=j(\pr_{K}\kappa(0))\qquad,
$$ 
$[\kappa]\in{\calW(K)}$; the action map 
$\sigma_{\calW(K)}:\calG\fib{\sor}{J}\calW(K)\to\calW(K)$
is given by
\bgn{eqnarray}\label{wellpos}
\nn\qquad
\sigma_{\calW(K)}(g,[\kappa\,])&:=\:&[g_{(1)}*_K\kappa\,]_{\calW(K)}\\
&\:=:&g\circledast[\kappa\,]
\qquad.
\end{eqnarray}
For all representatives $\kappa_\pm$ of the same class and $K$-homotopy $h$ from $\kappa_-$ to $\kappa_+$,
$g_{(2)}*h$ is a $K$-homotopy from $g_{(1)}*\kappa_-$ to $g_{(1)}*\kappa_+$, thus $\sigma_{\calW(K)}$ is well defined by
(\ref{wellpos}); that it is an action map compatible with $J$, follows straightforwardly from the cotangent
lifted action being morphic. Whenever $g\circledast[\kappa\,]$ is defined, we have
\be
\sor_{\calW(K)}(g\circledast[\kappa\,])&=&\pr_{K}(g*\kappa(0))=g*\pr_K(\kappa(0))\\
&=&g*\sor_{\calW(K)}([\kappa\,])
\ee
and similarly
$\tar_{\calW(K)}(g\circledast[\kappa\,])=g*\tar_{\calW(K)}([\kappa\,])$. For all composable classes 
$[\kappa_\pm]\in\calW(K)$ we may choose compactly supported smooth representatives to compute
$g\circledast([\kappa_+]\cdot [\kappa_+])$:
$$
\quad
g\circledast([\kappa_+]\cdot [\kappa_-])
=
\left[
{
\bgn{array}{ll}
2\cdot g*\kappa_-(2u)& 0\leq u\leq 1/2\\
2\cdot g*\kappa_+(2u-1)& 1/2\leq u\leq 1
\end{array}
}
\right]
=
(g\circledast[\kappa_+])\cdot(g\circledast[\kappa_-])
\quad,
$$
that is, the lifted action is compatible with the concatenation of $K$-paths.\\
({\em Step 2.})  Under the integrability assumptions $\calW(K)$ is the source 1-connected integration of $K$; thus it
is possible to form an action \tla-groupoid
$$
\bcat
\xy
*+{}="0",    <-1cm,0cm>
*+{\calG\acts\calW(K)}="1", <1.4cm,0cm>
*+{\calW(K)}="2", <-1cm,-1.4cm>
*+{\calG\acts P}="3", <1.4cm,-1.4cm>
*+{P}="4",
\ar  @ <0.07cm>   @{->} "1";"2"^{} 
\ar  @ <-0.07cm>  @{->} "1";"2"_{}  
\ar  @ <-0.07cm>  @{->} "1";"3"_{}
\ar  @ <0.07cm>   @{->} "1";"3"^{}
\ar  @ <0.07cm>   @{->} "2";"4"^{}
\ar  @ <-0.07cm>  @{->} "2";"4"_{}
\ar  @ <-0.07cm>  @{->} "3";"4"_{}
\ar  @ <0.07cm>   @{->} "3";"4"^{}
\endxy
\ecat
$$ 
and, according to theorem \ref{mainred} (with $\calD=\calG\acts\calW(K)$, $\calJ=\id_\calW(K)$, for the action by left
translation), to push  $\poidd{\calW(K)}{P}$ forward to a source 1-connected Lie groupoid 
$\calW(K)/\calG\equiv\calW(K)/(\calG\acts\calW(K))$ over $P/\calG$.  Since the Lie algebroid of $\calW(K)/\calG$ is 
$K/\calG$, thanks to proposition \ref{oop}, which in turn is isomorphic to $T^*(P/\calG)$ (proposition
\ref{poissonreduction}), $\calW(K)/\calG$ is a \emph{Lie} groupoid   integrating the Lie algebroid of the
quotient Poisson structure on $P/\calG$. It follows by the uniqueness of quotient Poisson bivector fields 
(theorem \ref{inty}) that the quotient Poisson manifold is integrable to a symplectic groupoid.

Conversely, assume that $P/\calG$ is integrable to some symplectic groupoid $\Lambda$. The projection
$\pr:P\to P/\calG$ is a surjective submersion, then the pullback algebroid
$\pr\daga\, T^*(P/\calG)$ always exists and $K$ can be embedded (as a Lie algebroid) in it 
by factoring the projection $K\to T^* (P/\calG)$ along the identity of $P$ (proposition \ref{indio}). Note that the Lie algebroid on  
$\pr\daga\, T^*(P/\calG)$ coincides with the fibred product 
$TP\fib{\dd\pr}{\ul{\pi}^\sharp}T^*(P/\calG)$.
The pullback Lie groupoid $\poidd{\pr\daga\Lambda}{P}$, namely the fibred product $(P\times
P)\fib{\pr\times\pr}{\chi_\Lambda}\Lambda$, exists thanks to the regularity of $\pr$ and 
integrates $TP\fib{\dd\pr}{\ul{\pi}^\sharp}T^*(P/\calG)$. Thus $K$ is integrable by 2nd Lie's
theorem.
\end{proof}
%
%
%
%
\bgn{example} Consider a Poisson groupoid $\gpdm$ and the Poisson action of $\calG$ on itself 
by left translation. The quotient Poisson manifold $\calG/\calG\simeq M$ exists and the 
quotient projection is given by the source map. Therefore
\emph{for any Poisson groupoid $\gpdm$ with integrable Poisson structure, the Poisson
structure induced on $M$ is integrable}. Note that $T^*M\subset T^*\calG$ 
is a Lie subalgebroid, being the core Lie algebroid of the cotangent prolongation \tla-groupoid; thus its integrability follows straightforwardly from
that of $\calG$. 
\end{example}
%
%
%

\spa Our proof of the integrability of quotient Poisson bivector fields, even though constructive, does not
produce explicitly a \emph{symplectic groupoid} for the quotient Poisson structure. In fact we obtain a \emph{Lie
groupoid}  integrating the corresponding Koszul algebroid, which is in general source 1-connected (the projection
$\calW(K)\to\calW(K)/\calG$ is sourcewise a diffeomorphism), thus it carries a compatible symplectic form, thanks to
Mackenzie and Xu's theorem \ref{inty}. Nevertheless we have no explicit characterization of the symplectic form. On the
other hand the integration of quotient Poisson bivector fields ``via symplectic  double groupoids'' (theorem
\ref{quoproqui}) produces an integrating Lie groupoid \emph{and} a canonical symplectic form.\\ 
We conclude by discussing a class of examples where the latter method is effective.
%
%
%
%
%
\subsection{The case of complete Poisson group actions}\hfill

\vs{0.1}
\spa The functorial approach of Chapter \ref{chapii} to the integration of \tla-groupoids can be applied to a wide
class of \tla-groupoids, namely action \tla-groupoids associated with the cotangent lift of a Poisson action of a
Poisson group $(G,\Pi)$ (including all actions  of Lie groups by Poisson diffeomorphisms in the case $\Pi=0$).
\bgn{proposition}\label{mainapp}\cite{07b} Let $G$ be a complete Poisson group and $P$ an integrable
Poisson manifold with source 1-connected symplectic groupoid $\Lambda$.  If $P$ is a
Poisson $G$-space, then
\medskip\\
$i)$ $\calJ:\Lambda\to G^\star$ is a symplectic $\calS$-space for the top
horizontal groupoid $\poidd{\calS}{G^\star}$ of the vertically source 1-connected
double of $G$;
\medskip\\
$ii)$ The action map 
$
\wt{\sigma}:\calS\fib{\sth}{\calJ}\Lambda\to\Lambda
$ 
is a morphism of Lie groupoids over the action map $\sigma :G\times P\to P$;
\medskip\\
$iii)$ The action double groupoid of the integrated action
\bgn{equation}\label{lgsga}
\bcat\xy
*+{}="0",    <-0.7cm,0.7cm>
*+{\calS\acts \Lambda}="1", <0.7cm,0.7cm>
*+{\Lambda}="2", <-0.7cm,-0.7cm>
*+{G\acts P}="3", <0.7cm,-0.7cm>
*+{P}="4",
\ar  @ <0.07cm>   @{->} "1";"2"^{} 
\ar  @ <-0.07cm>  @{->} "1";"2"_{}  
\ar  @ <-0.07cm>  @{->} "1";"3"_{}
\ar  @ <0.07cm>   @{->} "1";"3"^{}
\ar  @ <0.07cm>   @{->} "2";"4"^{}
\ar  @ <-0.07cm>  @{->} "2";"4"_{}
\ar  @ <-0.07cm>  @{->} "3";"4"_{}
\ar  @ <0.07cm>   @{->} "3";"4"^{}
\endxy\ecat
\end{equation}
is the vertically source 1-connected double Lie groupoid integrating the $\lagpd$
{\em (\ref{lagpdpoissonaction})} associated with the Poisson action.
\end{proposition}
\bgn{proof} It follows from corollary \ref{bipo} that $\calJ$ is anti-Poisson, since it integrates a morphism of Lie
bialgebroids $(T^*P,TP)\to (A^*, A)$; theorem \ref{poo} and lemma \ref{lift} imply that the cotangent lift of the
given action integrates to a morphic action of the vertically 1-connected double of $G$. ($i$) It remains to show
that the graph $\gr{\wt{\sigma}}$ of the integrated action map is coisotropic, by the same reasoning as in the proof
of theorem \ref{double} one can show that it is actually Lagrangian;
%
%
($ii$) holds by construction and ($iii$) is now obvious. 
\end{proof}
We are therefore in the condition to apply theorem \ref{quoproqui} and corollary \ref{erty} to obtain an
integrating symplectic groupoid for quotients of Poisson group actions
\bgn{corollary}\label{raspa} Under the hypotheses of proposition \ref{mainapp}, if $G$ acts freely and properly on $P$, then 
\medskip\\
i$)$ $\poidd{\jmod}{P/G}$ is a symplectic groupoid for the quotient Poisson structure;
\medskip\\
i$)$ The source connected component of $\poidd{\jmod}{P/G}$ is source 1-connected if{f} the
groups $\KK_\bullet(\jh)$ are trivial.
\end{corollary}
\bgn{remark}
Theorem \ref{raspa} generalizes a result by Xu (\cite{xu92}, theorem
4.2), regarding Poisson actions with a complete moment map. A
\tsf{moment map} \cite{luth} for a Poisson $G$-space $(P,\pi)$ is a
Poisson map  $j:P\to G^\star$ such that $\sigma(x)=\pish
j^*\linv{x}$, for the infinitesimal action
$\sigma(\cdot):\frag\to\frax(P)$ and the left invariant 1-form
$\linv{x}$ on $G^\star$ associated with  $x\in\frag\simeq\frag^{**}$;
such a Poisson map is called complete, when the Hamiltonian vector
field of $j^*f$ is complete for all compactly supported 
$f\in\cif(G^\star)$. If $G$ is 1-connected, so that the right
dressing action of $G$ on $G^\star$ is globally defined, $j$ is
always equivariant. When the action admits a complete moment map 
$j$, defefine $J:\Lambda\to G^\star$ and 
$J(\lambda)=j(\tar(\lambda))\cdot j(\sor(\lambda))\inverse$ (this space carries a natural action of $G$).
Assuming that $J\inverse(e_\star)/G$ is smooth, the construction of \cite{xu92}, produces a groupoid
$\poidd{J\inverse(e_\star)/G}{P/G}$, when $G$ is
complete and 1-connected.  Note that $J$ is by construction a
morphism of Lie groupoids and one can check \cite{xu95} that it
differentiates to $\jh$, therefore $J$ coincides with our $\calJ$; it
is easy to see that the $G$ action on $J\inverse(e_\star)$ of
\cite{xu92} is the same as that induced by the
$\poidd{\calS}{G^\star}$-action on $\Lambda$   ($\calS\simeq
G^\star{\ltimes\hs{-0.3}\rtimes} G$, under the assumptions).
Specializing last theorem to the case considered in \cite{xu92}, shows
that Xu's quotient is always smooth. 
\end{remark}
\bgn{remark}
A Lie group $G$ is trivially ($\Pi=0$) a complete Poisson group, with
the abelian group $\ol{\frag^*}$ as a dual Poisson group; in this
case a Poisson group action is an action by Poisson diffeomorphisms
and our approach  reproduces
the ``symplectization functor'' treatment of Fernandes \cite{fs06} and
Fernandes-Ortega-Ratiu \cite{for07} from the viewpoint of double
structures. 
The construction given in \cite{for07} (proposition 4.6) of a
symplectic groupoid for  the quotient Poisson manifold is precisely
the construction of theorem  \ref{quoproqui} in the special case $\Pi=0$.
\end{remark}

\bibliographystyle{habbrv}
\bibliography{biblione}

\def\cprime{$'$} \def\cprime{$'$} \def\cprime{$'$}
\begin{thebibliography}{10}

\bibitem{amrbook}
R.~Abraham, J.~E. Marsden, and T.~Ratiu.
\newblock {\em Manifolds, tensor analysis, and applications}, volume~75 of {\em
  Applied Mathematical Sciences}.
\newblock Springer-Verlag, New York, second edition, 1988.

\bibitem{aspp}
I.~Androulidakis and S.~G.
\newblock The holonomy groupoid of a singular foliation.
\newblock Preprint arXiv.org:math/0612370, 2007.

\bibitem{balan}
R.~Balan.
\newblock A note about integrability of distributions with singularities.
\newblock {\em Boll. Un. Mat. Ital. A (7)}, 8(3):335--344, 1994.

\bibitem{bm94}
J.~V. Beltr{\'a}n and J.~Monterde.
\newblock Poisson-{N}ijenhuis structures and the {V}inogradov bracket.
\newblock {\em Ann. Global Anal. Geom.}, 12(1):65--78, 1994.

\bibitem{bv88}
K.~H. Bhaskara and K.~Viswanath.
\newblock Calculus on {P}oisson manifolds.
\newblock {\em Bull. London Math. Soc.}, 20(1):68--72, 1988.

\bibitem{bh78}
R.~Brown and P.~J. Higgins.
\newblock On the connection between the second relative homotopy groups of some
  related spaces.
\newblock {\em Proc. London Math. Soc. (3)}, 36(2):193--212, 1978.

\bibitem{bm92}
R.~Brown and K.~C.~H. Mackenzie.
\newblock Determination of a double {L}ie groupoid by its core diagram.
\newblock {\em J. Pure Appl. Algebra}, 80(3):237--272, 1992.

\bibitem{bs76}
R.~Brown and C.~B. Spencer.
\newblock Double groupoids and crossed modules.
\newblock {\em Cahiers Topologie G\'eom. Diff\'erentielle}, 17(4):343--362,
  1976.

\bibitem{cds01}
A.~Cannas~da Silva.
\newblock {\em Lectures on symplectic geometry}, volume 1764 of {\em Lecture
  Notes in Mathematics}.
\newblock Springer-Verlag, Berlin, 2001.

\bibitem{cw99}
A.~Cannas~da Silva and A.~Weinstein.
\newblock {\em Geometric models for noncommutative algebras}, volume~10 of {\em
  Berkeley Mathematics Lecture Notes}.
\newblock American Mathematical Society, Providence, RI, 1999.

\bibitem{cat04}
A.~S. Cattaneo.
\newblock On the integration of {P}oisson manifolds, {L}ie algebroids, and
  coisotropic submanifolds.
\newblock {\em Lett. Math. Phys.}, 67(1):33--48, 2004.

\bibitem{cafe01}
A.~S. Cattaneo and G.~Felder.
\newblock Poisson sigma models and symplectic groupoids.
\newblock In {\em Quantization of singular symplectic quotients}, volume 198 of
  {\em Progr. Math.}, pages 61--93. Birkh\"auser, Basel, 2001.

\bibitem{chpr}
V.~Chari and A.~Pressley.
\newblock {\em A guide to quantum groups}.
\newblock Cambridge University Press, Cambridge, 1995.
\newblock Corrected reprint of the 1994 original.

\bibitem{conlon}
L.~Conlon.
\newblock {\em Differentiable manifolds}.
\newblock Birkh\"auser Advanced Texts: Basler Lehrb\"ucher. [Birkh\"auser
  Advanced Texts: Basel Textbooks]. Birkh\"auser Boston Inc., Boston, MA,
  second edition, 2001.

\bibitem{cdw87}
A.~Coste, P.~Dazord, and A.~Weinstein.
\newblock Groupo\"\i des symplectiques.
\newblock In {\em Publications du D\'epartement de Math\'ematiques. Nouvelle
  S\'erie. A, Vol.\ 2}, volume~87 of {\em Publ. D\'ep. Math. Nouvelle S\'er.
  A}, pages i--ii, 1--62. Univ. Claude-Bernard, Lyon, 1987.

\bibitem{crfs03}
M.~Crainic and R.~L. Fernandes.
\newblock Integrability of {L}ie brackets.
\newblock {\em Ann. of Math. (2)}, 157(2):575--620, 2003.

\bibitem{crfs04}
M.~Crainic and R.~L. Fernandes.
\newblock Integrability of {P}oisson brackets.
\newblock {\em J. Differential Geom.}, 66(1):71--137, 2004.

\bibitem{dfd93}
V.~G. Drinfel{\cprime}d.
\newblock On {P}oisson homogeneous spaces of {P}oisson-{L}ie groups.
\newblock {\em Teoret. Mat. Fiz.}, 95(2):226--227, 1993.

\bibitem{ev98}
P.~Etingof and A.~Varchenko.
\newblock Geometry and classification of solutions of the classical dynamical
  {Y}ang-{B}axter equation.
\newblock {\em Comm. Math. Phys.}, 192(1):77--120, 1998.

\bibitem{fs02}
R.~L. Fernandes.
\newblock Lie algebroids, holonomy and characteristic classes.
\newblock {\em Adv. Math.}, 170\quad(1):119--179, 2002.

\bibitem{fs06}
R.~L. Fernandes.
\newblock The symplectization functor.
\newblock Preprint arXiv.org:math/0610542, 2006.

\bibitem{for07}
R.~L. Fernandes, J.-P. Ortega, and T.~S. Ratiu.
\newblock The momentum map in {P}oisson geometry.
\newblock Preprint arXiv.org:0705.0562, 2007.

\bibitem{hlz01}
L.-g. He, Z.-J. Liu, and D.-S. Zhong.
\newblock Poisson actions and {L}ie bialgebroid morphisms.
\newblock In {\em Quantization, Poisson brackets and beyond (Manchester,
  2001)}, volume 315 of {\em Contemp. Math.}, pages 235--244. Amer. Math. Soc.,
  Providence, RI, 2002.

\bibitem{hm90a}
P.~J. Higgins and K.~Mackenzie.
\newblock Algebraic constructions in the category of {L}ie algebroids.
\newblock {\em J. Algebra}, 129(1):194--230, 1990.

\bibitem{hm90b}
P.~J. Higgins and K.~C.~H. Mackenzie.
\newblock Fibrations and quotients of differentiable groupoids.
\newblock {\em J. London Math. Soc. (2)}, 42(1):101--110, 1990.

\bibitem{ks86}
M.~V. Karas{\"e}v.
\newblock Analogues of objects of the theory of {L}ie groups for nonlinear
  {P}oisson brackets.
\newblock {\em Izv. Akad. Nauk SSSR Ser. Mat.}, 50(3):508--538, 638, 1986.

\bibitem{kv75}
A.~A. Kirillov.
\newblock Lie algebra structures that have the property of being local.
\newblock {\em Funkcional. Anal. i Prilo\v zen.}, 9(2):75--76, 1975.

\bibitem{ksb95}
Y.~Kosmann-Schwarzbach.
\newblock Exact {G}erstenhaber algebras and {L}ie bialgebroids.
\newblock {\em Acta Appl. Math.}, 41(1-3):153--165, 1995.
\newblock Geometric and algebraic structures in differential equations.

\bibitem{ksb96}
Y.~Kosmann-Schwarzbach.
\newblock The {L}ie bialgebroid of a {P}oisson-{N}ijenhuis manifold.
\newblock {\em Lett. Math. Phys.}, 38(4):421--428, 1996.

\bibitem{km90}
Y.~Kosmann-Schwarzbach and F.~Magri.
\newblock Poisson-{N}ijenhuis structures.
\newblock {\em Ann. Inst. H. Poincar\'e Phys. Th\'eor.}, 53(1):35--81, 1990.

\bibitem{lang}
S.~Lang.
\newblock {\em Fundamentals of differential geometry}, volume 191 of {\em
  Graduate Texts in Mathematics}.
\newblock Springer-Verlag, New York, 1999.

\bibitem{lp05}
L.-C. Li and S.~Parmentier.
\newblock On dynamical {P}oisson groupoids. {I}.
\newblock {\em Mem. Amer. Math. Soc.}, 174(824):vi+72, 2005.

\bibitem{lnw77}
A.~Lichnerowicz.
\newblock Les vari\'et\'es de {P}oisson et leurs alg\`ebres de {L}ie
  associ\'ees.
\newblock {\em J. Differential Geometry}, 12(2):253--300, 1977.

\bibitem{Lie}
S.~Lie.
\newblock Theorie der {T}ransformationsgruppen {I}.
\newblock {\em Math. Ann.}, 16(4):441--528, 1880.

\bibitem{lx96}
Z.-J. Liu and P.~Xu.
\newblock Exact {L}ie bialgebroids and {P}oisson groupoids.
\newblock {\em Geom. Funct. Anal.}, 6(1):138--145, 1996.

\bibitem{luth}
J.-H. Lu.
\newblock {\em Multiplicative and Affine {P}oisson Structures on {L}ie Groups}.
\newblock {PhD} thesis, University of California Berkeley, Department of
  Mathematics, 1990.

\bibitem{lu97}
J.-H. Lu.
\newblock Poisson homogeneous spaces and {L}ie algebroids associated to
  {P}oisson actions.
\newblock {\em Duke Math. J.}, 86(2):261--304, 1997.

\bibitem{lu07}
J.-H. Lu.
\newblock A note on poisson homogeneous spaces.
\newblock arXiv.org:0706.1337, 2007.

\bibitem{lw89}
J.-H. Lu and A.~Weinstein.
\newblock Groupo\"\i des symplectiques doubles des groupes de {L}ie-{P}oisson.
\newblock {\em C. R. Acad. Sci. Paris S\'er. I Math.}, 309(18):951--954, 1989.

\bibitem{mkz92}
K.~C.~H. Mackenzie.
\newblock Double {L}ie algebroids and second-order geometry. {I}.
\newblock {\em Adv. Math.}, 94(2):180--239, 1992.

\bibitem{mkz98}
K.~C.~H. Mackenzie.
\newblock Drinfel\cprime d doubles and {E}hresmann doubles for {L}ie algebroids
  and {L}ie bialgebroids.
\newblock {\em Electron. Res. Announc. Amer. Math. Soc.}, 4:74--87
  (electronic), 1998.

\bibitem{mkz99}
K.~C.~H. Mackenzie.
\newblock On symplectic double groupoids and the duality of {P}oisson
  groupoids.
\newblock {\em Internat. J. Math.}, 10(4):435--456, 1999.

\bibitem{mkz00a}
K.~C.~H. Mackenzie.
\newblock Double {L}ie algebroids and second-order geometry. {II}.
\newblock {\em Adv. Math.}, 154(1):46--75, 2000.

\bibitem{mkz00b}
K.~C.~H. Mackenzie.
\newblock A unified approach to {P}oisson reduction.
\newblock {\em Lett. Math. Phys.}, 53(3):215--232, 2000.

\bibitem{mackbook}
K.~C.~H. Mackenzie.
\newblock {\em General theory of {L}ie groupoids and {L}ie algebroids}, volume
  213 of {\em London Mathematical Society Lecture Note Series}.
\newblock Cambridge University Press, Cambridge, 2005.

\bibitem{mx94}
K.~C.~H. Mackenzie and P.~Xu.
\newblock Lie bialgebroids and {P}oisson groupoids.
\newblock {\em Duke Math. J.}, 73(2):415--452, 1994.

\bibitem{mx98}
K.~C.~H. Mackenzie and P.~Xu.
\newblock Classical lifting processes and multiplicative vector fields.
\newblock {\em Quart. J. Math. Oxford Ser. (2)}, 49(193):59--85, 1998.

\bibitem{mx00}
K.~C.~H. Mackenzie and P.~Xu.
\newblock Integration of {L}ie bialgebroids.
\newblock {\em Topology}, 39(3):445--467, 2000.

\bibitem{moment}
J.~Marsden and A.~Weinstein.
\newblock Reduction of symplectic manifolds with symmetry.
\newblock {\em Rep. Mathematical Phys.}, 5(1):121--130, 1974.

\bibitem{mw88}
K.~Mikami and A.~Weinstein.
\newblock Moments and reduction for symplectic groupoids.
\newblock {\em Publ. Res. Inst. Math. Sci.}, 24(1):121--140, 1988.

\bibitem{mm02}
I.~Moerdijk and J.~Mr\v{c}un.
\newblock On integrability of infinitesimal actions.
\newblock {\em Amer. J. Math.}, 124(3):567--593, 2002.

\bibitem{nij55}
A.~Nijenhuis.
\newblock Jacobi-type identities for bilinear differential concomitants of
  certain tensor fields. {I}, {II}.
\newblock {\em Nederl. Akad. Wetensch. Proc. Ser. A. {\bf 58} Indag. Math.},
  17:390--397, 398--403, 1955.

\bibitem{prd66}
J.~Pradines.
\newblock Th\'eorie de {L}ie pour les groupo\"\i des diff\'erentiables.
  {R}elations entre propri\'et\'es locales et globales.
\newblock {\em C. R. Acad. Sci. Paris S\'er. A-B}, 263:A907--A910, 1966.

\bibitem{prd67}
J.~Pradines.
\newblock Th\'eorie de {L}ie pour les groupo\"\i des diff\'erentiables.
  {C}alcul diff\'erenetiel dans la cat\'egorie des groupo\"\i des
  infinit\'esimaux.
\newblock {\em C. R. Acad. Sci. Paris S\'er. A-B}, 264:A245--A248, 1967.

\bibitem{prd67b}
J.~Pradines.
\newblock Th\'eorie de {L}ie pour les groupo\"\i des diff\'erentiables.
  {C}alcul diff\'erenetiel dans la cat\'egorie des groupo\"\i des
  infinit\'esimaux.
\newblock {\em C. R. Acad. Sci. Paris S\'er. A-B}, 264:A245--A248, 1967.

\bibitem{prd68a}
J.~Pradines.
\newblock G\'eom\'etrie diff\'erentielle au-dessus d'un groupo\"\i de.
\newblock {\em C. R. Acad. Sci. Paris S\'er. A-B}, 266:A1194--A1196, 1968.

\bibitem{prd68b}
J.~Pradines.
\newblock Troisi\`eme th\'eor\`eme de {L}ie les groupo\"\i des
  diff\'erentiables.
\newblock {\em C. R. Acad. Sci. Paris S\'er. A-B}, 267:A21--A23, 1968.

\bibitem{p84}
J.~Pradines.
\newblock Feuilletages: holonomie et graphes locaux.
\newblock {\em C. R. Acad. Sci. Paris S\'er. I Math.}, 298(13):297--300, 1984.

\bibitem{rh63}
G.~S. Rinehart.
\newblock Differential forms on general commutative algebras.
\newblock {\em Trans. Amer. Math. Soc.}, 108:195--222, 1963.

\bibitem{sch40}
J.~A. Schouten.
\newblock Ueber {D}ifferentialkomitanten zweier kontravarianter {G}r\"ossen.
\newblock {\em Nederl. Akad. Wetensch., Proc.}, 43:449--452, 1940.

\bibitem{sch53}
J.~A. Schouten.
\newblock On the differential operators of first order in tensor calculus.
\newblock In {\em Convegno Internazionale di Geometria Differenziale, Italia,
  1953}, pages 1--7. Edizioni Cremonese, Roma, 1954.

\bibitem{smale}
S.~Smale.
\newblock Topology and mechanics. {I}.
\newblock {\em Invent. Math.}, 10:305--331, 1970.

\bibitem{st74}
P.~Stefan.
\newblock Accessible sets, orbits, and foliations with singularities.
\newblock {\em Proc. London Math. Soc. (3)}, 29:699--713, 1974.

\bibitem{07b}
L.~Stefanini.
\newblock Integrability and reduction of poisson group actions.
\newblock {\em Submitted}, Preprint arXiv.org:0710.5753, 2007.

\bibitem{07a}
L.~Stefanini.
\newblock On the integration of {LA}-groupoids and duality for poisson
  groupoids.
\newblock {\em To appear in Travaux Math{\'e}matiqu{\'e}s}, Preprint
  arXiv.org:math/0701231, 2007.

\bibitem{su73}
H.~J. Sussmann.
\newblock Orbits of families of vector fields and integrability of
  distributions.
\newblock {\em Trans. Amer. Math. Soc.}, 180:171--188, 1973.

\bibitem{vm94}
I.~Vaisman.
\newblock {\em Lectures on the geometry of {P}oisson manifolds}, volume 118 of
  {\em Progress in Mathematics}.
\newblock Birkh\"auser Verlag, Basel, 1994.

\bibitem{sev05}
P.~\v{S}evera.
\newblock Some title containing the words ``homotopy'' and ``symplectic'', e.g.
  this one.
\newblock In {\em Travaux math\'ematiques. Fasc. XVI}, Trav. Math., XVI, pages
  121--137. Univ. Luxemb., Luxembourg, 2005. Based on a talk at "Poisson 2000",
  CIRM Marseille, June 2000.

\bibitem{weibel}
C.~A. Weibel.
\newblock {\em An introduction to homological algebra}, volume~38 of {\em
  Cambridge Studies in Advanced Mathematics}.
\newblock Cambridge University Press, Cambridge, 1994.

\bibitem{ws83}
A.~Weinstein.
\newblock The local structure of {P}oisson manifolds.
\newblock {\em J. Differential Geom.}, 18(3):523--557, 1983.

\bibitem{ws87}
A.~Weinstein.
\newblock Symplectic groupoids and {P}oisson manifolds.
\newblock {\em Bull. Amer. Math. Soc. (N.S.)}, 16(1):101--104, 1987.

\bibitem{ws88}
A.~Weinstein.
\newblock Coisotropic calculus and {P}oisson groupoids.
\newblock {\em J. Math. Soc. Japan}, 40(4):705--727, 1988.

\bibitem{ws98}
A.~Weinstein.
\newblock Poisson geometry.
\newblock {\em Differential Geom. Appl.}, 9(1-2):213--238, 1998.
\newblock Symplectic geometry.

\bibitem{xu92}
P.~Xu.
\newblock Symplectic groupoids of reduced {P}oisson spaces.
\newblock {\em C. R. Acad. Sci. Paris S\'er. I Math.}, 314(6):457--461, 1992.

\bibitem{xu95}
P.~Xu.
\newblock On {P}oisson groupoids.
\newblock {\em Internat. J. Math.}, 6(1):101--124, 1995.

\bibitem{xuonvaisman}
P.~Xu.
\newblock Review of `` {L}ectures on the geometry of {P}oisson manifolds'', by
  {I}. {V}aisman.
\newblock {\em Bull. Amer. Math. Soc. (N.S.)}, 33(2):255--261, 1996.

\bibitem{zr90}
S.~Zakrzewski.
\newblock Quantum and classical pseudogroups. {I}. {U}nion pseudogroups and
  their quantization.
\newblock {\em Comm. Math. Phys.}, 134(2):347--370, 1990.

\end{thebibliography}
\newpage
\thispagestyle{empty}
\newpage
\selectlanguage{german}
\newpage 
\chapter*{Curriculum Vitae}

\vspace{\stretch{1}}

\noi {\bf Stefanini Luca} \\
${}$ \\
Geboren am 9. May 1978 in Bergamo  (Italien) \\

\vspace{\stretch{2}}

\noi {\bf Ausbildung}
${}$ \\
\begin{itemize}
\item[] Juli 1997: \emph{Diploma di Maturit\`a Scientifica}, Liceo Lorenzo Mascheroni, Bergamo (Italien)  
\item[] 1997-2002: Studium der Physik, Universit\`a di Pavia, Pavia (Italien)
\item[] September 2002: \emph{Laurea in Fisica Teorica}
\item[] Tesi di Laurea: 
	``Tecniche Algebriche nelle Teorie di Campo Topologiche''
\item[] Leitung: Chiar.mo Prof. Dr. Mauro Carfora (Universit\`a di Pavia)
\end{itemize}

\vspace{\stretch{3}}

\noi {\bf Dissertation}
${}$ \\
\begin{itemize}
\item[] Titel: ``On Morphic Actions and Integrability of \tla-Groupoids''
\item[] Leitung: Prof. Dr. Alberto S. Cattaneo
\end{itemize}

\vspace{\stretch{4}}

\noi {\bf Gegenw\"artige Stellung} \\
${}$ \\
	Seit M\"arz 2003 Doktorand und Assistent am 
	Institut f\"ur Mathematik der Universit\"at Z\"urich. 

\vspace{\stretch{5}}

\quad
\end{document}